\documentclass[reqno]{siamonline1116}

\usepackage{amsmath}
\usepackage{amssymb}
\usepackage{lipsum}
\usepackage{braket,amsfonts}
\usepackage{graphicx}
\usepackage{algorithmic}
\usepackage{multirow}
\usepackage{bold-extra}
\usepackage{fancybox,fancyhdr}

\usepackage{setspace,color}
\usepackage[caption=false]{subfig}

\numberwithin{theorem}{section}
\numberwithin{equation}{section}

\usepackage[mathscr]{eucal}

\newtheorem{remark}[theorem]{Remark}

\newtheorem{example}[theorem]{Example}

\def\R{{\mathbb R}}

\def\N{{\mathbb N}}

\def\Z{{\mathbb Z}}

\def\BB{{\mathbb B}}
\def\OO{{\mathbb O}}

\def\sign{\mathrm{sign}}

\def\pv{\mathrm{pv}}
\def\rd{\mathrm{d}}
\def\Om{\Omega}
\def\om{\omega}

\def\f{\frac}
\def\p{\partial}

\def\na{\nabla}
\def\la{\langle}
\def\ra{\rangle}

\def\bi{{\mathbf i}}

\def\bk{{\boldsymbol k}}

\def\bn{\boldsymbol{n}}

\def\x{\boldsymbol{x}}
\def\y{{\boldsymbol y}}

\def\B{\mathbf{B}}

\def\bH{{\mathbf H}}

\def\M{{\mathbf M}}

\def\mD{{\mathcal D}}
\def\mE{{\mathcal E}}

\def\mH{{\mathcal H}}
\def\mI{{\mathcal I}}

\def\mT{{\mathcal T}}

\def\mX{{\mathcal X}}

\def\msF{{\mathscr F}}

\def\msL{{\mathscr L}}

\def\msX{{\mathscr X}}

\def\aal{\boldsymbol{\alpha}}

\def\ssi{\boldsymbol{\sigma}}

\def\nnu{\boldsymbol{\nu}}

\def\bi{\begin{itemize}} \def\ei{\end{itemize}}
\def\be{\begin{eqnarray*}}
\def\ee{\end{eqnarray*}}

\def\0{{\mathbf 0}}

\def\pv{\mathrm{pv}}

\newcommand{\beq}{\begin{equation}}
\newcommand{\eeq}{\end{equation}}

\def\xxi{{\boldsymbol{\xi}}}

\def\Ps{{\boldsymbol \Psi}}
\def\wt{\widetilde}
\def\wh{\widehat}

\newcommand{\eps}{\varepsilon}
\def\la{\langle}
\def\ra{\rangle}

\def\XXint#1#2#3{{\setbox0=\hbox{$#1{#2#3}{\int}$ }
\vcenter{\hbox{$#2#3$ }}\kern-.55\wd0}}

\newcommand{\TheTitle}{Whole Brain Susceptibility Mapping Using Harmonic Incompatibility Removal}
\newcommand{\TheAuthors}{Chenglong Bao, Jae Kyu Choi, and Bin Dong}
\newcommand{\argmin}{\operatornamewithlimits{argmin}}

\headers{Whole Brain Susceptibility Mapping Using HIRE}{\TheAuthors}

\title{{\TheTitle}\thanks{Submitted to the editors DATE. \funding{The research of the third author is supported by NSFC grant 11831002.}}}

\author{
  Chenglong Bao\thanks{Yau Mathematical Sciences Center, Tsinghua University, Beijing, 100084 China, (\email{clbao@math.tsinghua.edu.cn}).}
  \and
  Jae Kyu Choi\thanks{Corresponding Author. School of Mathematical Sciences, Tongji University, Shanghai, 200092 China, (\email{jaycjk@tongji.edu.cn}).}
  \and
  Bin Dong\thanks{Beijing International Center for Mathematical Research and Laboratory for Biomedical Image Analysis Beijing Institute of Big Data Research, Peking University, Beijing, 100871 China, (\email{dongbin@math.pku.edu.cn}).}
}

\usepackage{amsopn}

\begin{document}
\maketitle

\begin{abstract} Quantitative susceptibility mapping (QSM) aims to visualize the three dimensional susceptibility distribution by solving the field-to-source inverse problem using the phase data in magnetic resonance signal. However, the inverse problem is ill-posed since the Fourier transform of integral kernel has zeroes in the frequency domain. Although numerous regularization based models have been proposed to overcome this problem, the incompatibility in the field data has not received enough attention, which leads to deterioration of the recovery. In this paper, we show that the data acquisition process of QSM inherently generates a harmonic incompatibility in the measured local field. Based on such discovery, we propose a novel regularization based susceptibility reconstruction model with an additional sparsity based regularization term on the harmonic incompatibility. Numerical experiments show that the proposed method achieves better performance than the existing approaches.
\end{abstract}



\section{Introduction}\label{Introduction}

Quantitative susceptibility mapping (QSM) \cite{L.deRochefort2010} is a novel imaging technique that visualizes the magnetic susceptibility distribution from the measured field data associated with magnetization $\M=(M_1,M_2,M_3)$ induced in the body by an MR scanner. The magnetic susceptibility $\chi$ is an intrinsic property of the material which relates $\M$ and the magnetic field $\bH=(H_1,H_2,H_3)$ through $\M=\chi\bH$ \cite{J.K.Seo2014}. As physiological and/or pathological processes alter tissues' magnetic susceptibilities, QSM has been widely applied in biomedical image analysis \cite{J.K.Seo2014}. Applications include demyelination, inflammation, and iron overload in multiple sclerosis \cite{W.Chen2014}, neurodegeneration and iron overload in Alzheimer's disease \cite{J.Acosta-Cabronero2013}, Huntington's disease \cite{J.M.G.vanBergen2016}, changes in metabolic oxygen consumption \cite{E.M.Haacke2010}, hemorrhage including microhemorrhage and blood degradation \cite{J.Klohs2011}, bone mineralization \cite{A.V.Dimov2018}, drug delivery using magnetic nanocarriers \cite{T.Liu2010}.

\begin{figure}[tp!]
\centering
\includegraphics[width=15.6cm]{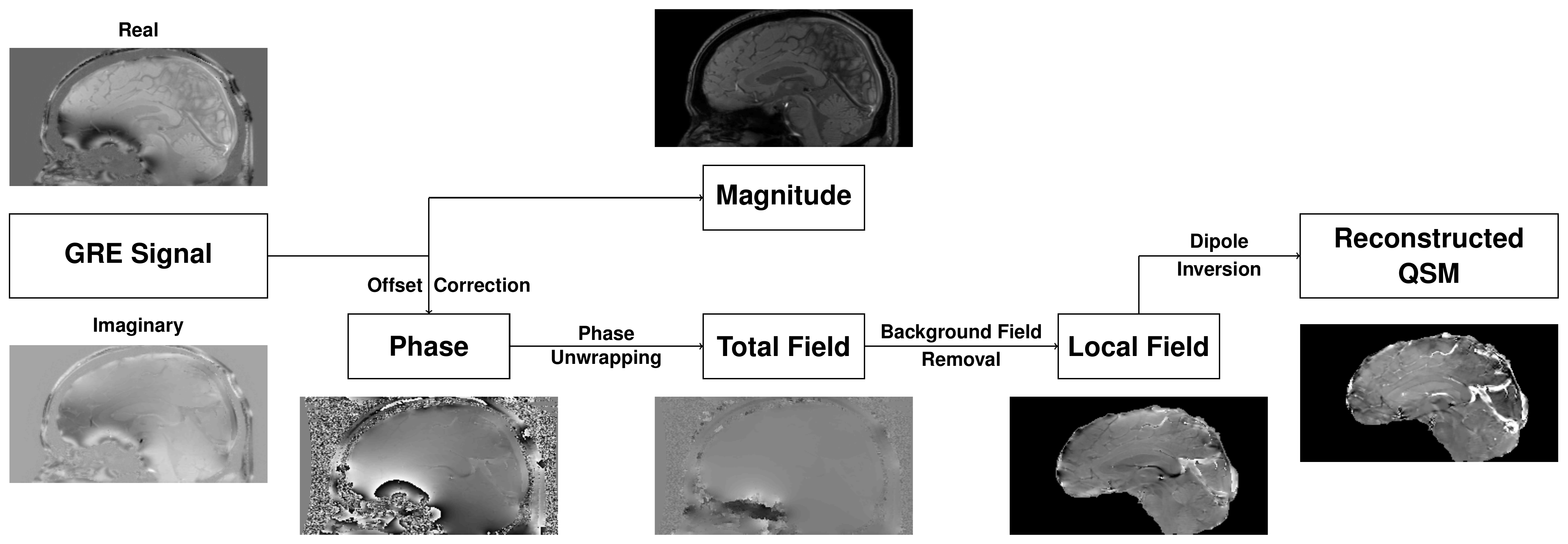}\vspace{-0.20cm}
\caption{Schematic diagram of QSM reconstruction process.}\label{SchematicProcedure}
\end{figure}

QSM uses the phase data of a complex gradient echo (GRE) signal as the phase linearly increases with respect to the field perturbation induced by the magnetic susceptibility distribution in an MR scanner \cite{Y.Wang2015}. More concretely, assume that an object is placed in an MR scanner with the main static magnetic field $\B_0=(0,0,B_0)$ where $B_0$ is a positive constant. Then, for any $\x\in\R^3$, the observed complex GRE signal $I(\x,TE)$ at an echo time $TE\mathrm{sec}$ is modeled as
\begin{align}\label{ComplexMR}
I(\x,TE)=m(\x)\exp\left\{-i\left(b(\x)\omega_0B_0TE+\theta_0(\x)\right)\right\},
\end{align}
where $\omega_0=42.577 \mathrm{MHz/T}$ is the proton gyromagnetic ratio, $b$ is the total field induced by the susceptibility distribution in an MR scanner, and $\theta_0$ is the coil sensitivity dependent phase offset. The magnitude image $m(\x)$ in \eqref{ComplexMR} is proportional to the proton density \cite{Y.Wang2015}, and the phase $\theta(\x)$ in $I(\x,TE)$ is written as
\begin{align}\label{phasemodel}
\theta(\x)=b(\x)\omega_0B_0TE+\theta_0(\x).
\end{align}
Based on the observations $\theta(\x)$, QSM aims at visualizing the susceptibility distribution $\chi(\x)$ in the region of interests (ROI) $\Omega$ which occupies the water and brain tissues. Note that the ROI $\Omega$ can be readily determined by $I(\x,TE)$ (and thus by $m(\x)$) as $m(\x)=|I(\x,TE)|\approx0$ whenever $\x\notin\Omega$ \cite{W.Li2014,F.Schweser2016,Y.Wang2015}. The standard QSM consists of the following four steps: offset correction, phase unwrapping, background field removal and dipole inversion (see \cref{SchematicProcedure} for the overview of the process). The first three steps extract the local field $b_l$ that is contained in the total field $b$: the offset correction removes/corrects $\theta_0(\x)$ from $\theta(\x)$ to obtain $b(\x)\omega_0B_0TE$ (the offset corrected phase) lying in $(-\pi,\pi]$; the phase unwrapping removes the artificial jumps in the offset corrected phase when estimating the total field $b$; the background field removal eliminates the field induced by the susceptibility outside $\Omega$ such as skulls and nasal cavity. Interested readers may refer to \cite{E.M.Haacke2015,F.Schweser2016,Y.Wang2015} and references therein for more details.

Given the local field $b_l$, the dipole inversion recovers the susceptibility distribution $\chi$ in $\Om$ by solving the following convolution relation \cite{T.Liu2011,T.Liu2011a,T.Liu2009}:
\begin{align}\label{QSM_Deconv}
b_l(\x)=\pv\int_{\Om}d(\x-\y)\chi(\y)\rd\y,
\end{align}
where $\pv$ denotes the principal value \cite{E.M.Stein2011} of the singular integral with the kernel $d$:
\begin{align*}
d(\x)=\f{2x_3^2-x_1^2-x_2^2}{4\pi|\x|^5}.
\end{align*}
In the frequency domain, \cref{QSM_Deconv} reads
\begin{align}\label{QSM_Fou}
\msF(b_l)(\xxi)=\mD(\xxi)\msF(\chi)(\xxi)=\left(\f{1}{3}-\f{\xi_3^2}{|\xxi|^2}\right)\msF(\chi)(\xxi)
\end{align}
where $\mD=\msF(d)$ is the Fourier transform of $d$ and $\mD(\0)=0$  by the definition of $\pv$ \cite{L.deRochefort2010,E.M.Haacke2015}. From \eqref{QSM_Fou}, it is easy to see that recovering the susceptibility distribution $\chi$ is ill-posed as $\mD=0$ on the critical manifold $\Gamma_0=\left\{\xxi\in\R^3:\xi_1^2+\xi_2^2-2\xi_3^2=0\right\}$. This ill-posedness leads to the streaking artifacts unless the data $b_l$ satisfies a proper compatibility condition \cite{J.K.Choi2014}.

\subsection{Existing QSM Reconstruction Methods}\label{Existing}

In the literature, various QSM reconstruction methods have been explored to deal with the ill-posed nature of the inverse problem \cref{QSM_Fou}. Early attempts mainly focus on the direct methods based on the modification of \cref{QSM_Fou} near $\Gamma_0$ \cite{Y.Kee2017}. One benchmark method, called the truncated K-space division (TKD) \cite{K.Shmueli2009}, finds the approximate solution to \cref{QSM_Fou} via:
\begin{align}\label{TKD}
\chi_{\hbar}=\msF^{-1}(\mX_{\hbar}),~~~\text{where}~~\mX_{\hbar}(\xxi)=\f{\mathrm{sign}(\mD(\xxi))}{\max\left\{|\mD(\xxi)|,\hbar\right\}}\msF(b_l)(\xxi)
\end{align}
with a threshold level $\hbar>0$. Another method recovers $\chi$ via solving the following Tikhonov regularization \cite{B.Kressler2010}:
\begin{align}\label{Tikhonov}
\min_{\chi}~\f{1}{2}\left\|A\chi-b_l\right\|_2^2+\eps\left\|\chi\right\|_2^2
\end{align}
where $\eps>0$ and $A$ denotes the forward operator that is obtained by discretizing the kernel $\mD$. Recently, some other direct methods are proposed, e.g. the iterative susceptibility weighted imaging and susceptibility mapping \cite{J.Tang2013}, the analytic continuation \cite{Natterer2016} and so on. Even though these direct methods are simple to implement, they can introduce additional artifacts due to the modification of $1/\mD$ near $\Gamma_0$ in the frequency domain \cite{J.K.Choi2014,Y.Kee2017,B.Palacios2017}.

In recent years, the regularization based methods have been proposed and show the superior performance over the direct method \cite{Y.Kee2017,S.Wang2013}. Mathematically, it is formulated as solving the minimization problem:
\begin{align}\label{Regularization}
\min_{\chi}~F(b_l|\chi)+R(\chi),
\end{align}
where $F(b_l|\chi)$ denotes the data fidelity term and $R(\chi)$ is the regularization term which mostly promotes the sparse approximation of $\chi$ under some linear transformation such as total variation and wavelet frames. According to the choices of $F(b_l|\chi)$, the regularization based methods can be classified into the \emph{integral approaches} and the \emph{differential approaches} \cite{Y.Kee2017}. The most widely used integral approaches are based on the convolution relation \cref{QSM_Deconv}. For example, $F(b_l|\chi)=\f{1}{2}\left\|A\chi-b_l\right\|_2^2$ when the data is corrupted by a white Gaussian noise. Even though the integral approach is capable of suppressing streaking artifacts, it is empirically reported in \cite{Y.Kee2017} that the reconstructed image can contain the shadow artifacts in the region of piecewise constant susceptibility. The differential approaches are based on the following partial differential equation (PDE)
\begin{align}\label{QSM_PDE}
-\Delta b_l(\x)=P(D)\chi(\x)=\left(-\f{1}{3}\Delta+\f{\p^2}{\p x_3^2}\right)\chi(\x)~~~~~~~\x\in\Om
\end{align}
which is derived from the Maxwell's equation \cite{E.M.Haacke1999,J.K.Seo2014}. In this case, one typical fidelity term is $F(b_l|\chi)=\f{1}{2}\left\|P(D)\chi+\Delta b_l\right\|_2^2$ by considering $-\Delta b_l$ as a measurement. Compared with the integral approach, the differential approach is able to restore susceptibility image with less shadow artifacts. However, the noise in the data can be amplified by $-\Delta$, which leads to the streaking artifacts \cite{Y.Wang2015}. In \cite{Y.Kee2017}, the differential approach is implemented by incorporating the spherical mean value (SMV) filter $S_r$ with a radius $r>0$ \cite{W.Li2014} into the integral approach:
\begin{align}\label{IntegralSMV}
\min_{\chi}~\f{1}{2}\left\|S_r\left(A\chi-b_l\right)\right\|_2^2+R(\chi).
\end{align}
Since the implementation of $S_r$ causes the erosion of $\Om$ according to the choice of $r$, the loss of anatomical information near $\p\Om$ is inevitable at the cost of the shadow artifact removal \cite{Y.Kee2017}.

\subsection{Motivations and Contributions of Our Approach}\label{OurApproach}

Even though the equations \cref{QSM_Deconv,QSM_PDE} are known to be equivalent \cite{J.K.Choi2014,Y.Kee2017,B.Palacios2017}, it is observed that the local field $b_l$ defined as \cref{QSM_Deconv} is a particular solution of the PDE \cref{QSM_PDE}. Whenever the data acquisition is based on the PDE \cref{QSM_PDE}, the measured local field data will be written as the superposition of $b_l$ in \cref{QSM_Deconv} and the ambiguity of $-\Delta$, which will be referred as the \emph{harmonic incompatibility}. Therefore, there is a need to identify/remove the harmonic incompatibility from the measured local field data for better reconstruction results as it is smooth, analytic and satisfies the mean value property in an open set \cite{Evans2010}, which are different from the noise properties.

It is noted that the background field removal aims at obtaining the local field $b_l$ via solving a Poisson equation with certain boundary condition as the background field is harmonic in $\Om$ \cite{W.Li2014,F.Schweser2010,Y.Wang2015,D.Zhou2014}. In this case, the measured local field $b_l$ is represented by the Green's function associated with the boundary condition \cite{J.K.Choi2014}. Thus, it is inevitable that $b_l$ contains the incompatibility associated with the imposed boundary condition. In this paper, we investigate the incompatibility of the local filed data in QSM and establish that this incompatibility consists of two harmonic functions inside and outside $\Om$ respectively, and its (distributional) Laplacian defines a surface measure on $\p\Om$ (see \cref{Th1} for details and \cref{IllustrateTh1,IllustrateTh1Axial,IllustrateTh1Brain,IllustrateTh1BrainAxial} for illustrations). Therefore, we can establish a new forward model in QSM by taking this harmonic incompatibility into account.



Based on this discovery, we impose a constraint on harmonic incompatibility term in susceptibility reconstruction model. Since our theoretical results suggest that the incompatibility is harmonic except on $\p\Om$, one straightforward approach is to penalize its (discrete) Laplacian on points $\x\notin\p\Om$. However, it is in general difficult to explicitly model this harmonic incompatibility and/or to directly impose its property into the susceptibility reconstruction model due to the complicated geometries of human brains and the limited spatial resolution in real MRI data. Instead, we impose the sparse regularization of the incompatibility as the support of its Laplacian is small compared to the size of image. Combing it with traditional regularization on the susceptibility image, we propose a novel regularization based QSM model by imposing additional constraints on the incompatibility term. Within the new model, we can suppress the incompatibility other than the noise, achieving the whole brain imaging with less artifacts together with the regularization of susceptibility image. Experiments on both brain phantom and vivo MR data consistently show the advantages of the proposed HIRE model which achieves the state-of-the-art performance. Besides, our experiments suggest that tight frame regularization of the susceptibility image can avoid the constant offset \cite{Y.Kee2017} and lead to efficient computation.

\subsection{Organization of Paper}\label{Organization}

In \cref{HIREModel}, we introduce our HIRE model for whole brain susceptibility imaging. More precisely, we first briefly review the biophysics forward model of QSM in \cref{BasicForward}, and characterize the harmonic incompatibility in the local field data in \cref{LocalFieldAcquisition}. Based on the characterization, we introduce the proposed HIRE model in \cref{ProposedHIREModel}, followed by an alternating minimization algorithm in \cref{AlternatingMinimizationAlgorithm}. In \cref{Experiments}, we present experimental results for both brain phantom and in vivo MR data, and the concluding remarks are given in \cref{Conclusion}.

\section{Harmonic Incompatibility Removal (HIRE) Model for Whole Brain Imaging}\label{HIREModel}

\subsection{Preliminaries on Biophysics of QSM}\label{BasicForward}

In an MRI scanner with the main static magnetic field $\B_0=(0,0,B_0)$ where $B_0$ is a positive constant, objects gain a magnetization $\M(\x)$. This magnetization generates a macroscopic field $\B(\x)$ satisfying the following magnetostatic Maxwell's equation \cite{E.M.Haacke1999,J.K.Seo2014}
\begin{align}\label{MagnetostaticMaxwell}
\begin{split}
\na\cdot\B&=0\\
\na\times\B&=\mu_0\na\times\M,
\end{split}
\end{align}
where $\mu_0=8.854\times 10^{-12} \mathrm{F/m}$ is the vacuum permittivity. Since the MRI signal is generated by the microscopic field $\B_{\ell}(\x)$ experienced by the spins of water protons \cite{Y.Kee2017}, we use the following Lorenz sphere correction model \cite{E.M.Haacke1999}:
\begin{align}\label{Lorenz}
\B_{\ell}(\x)=\B(\x)-\f{2}{3}\mu_0\M(\x)
\end{align}
to relate $\B(\x)$ and $\B_{\ell}(\x)$.

Note that since $\M(\x)$ is generated by $\B_0$ field, we have $\M(\x)=(0,0,M(\x))$. Moreover, since we consider the linear magnetic materials with $|\chi|\ll1$, $\chi$ can be approximated as
\begin{align}\label{Approx1}
\chi(\x)\approx\f{\mu_0}{B_0}\M(\x).
\end{align}
Finally, we introduce the total field $b(\x)$ as
\begin{align}\label{Approx2}
b(\x)=\f{B_{\ell3}(\x)-B_0}{B_0}
\end{align}
where $B_{\ell3}(\x)$ denotes the third component of $\B_{\ell}(\x)$.

Combining \cref{MagnetostaticMaxwell,Lorenz,Approx1,Approx2} and taking the third component into account only, we obtain the following relation between $\chi$ and $b$ in the frequency domain:
\begin{align}\label{QSM_Fou2}
|\xxi|^2\msF(b)(\xxi)=\left(\f{1}{3}|\xxi|^2-\xi_3^2\right)\msF(\chi)(\xxi),
\end{align}
which gives
\begin{align}\label{QSM_PDE2}
-\Delta b=P(D)\chi:=\left(-\f{1}{3}\Delta+\f{\p^2}{\p x_3^2}\right)\chi.
\end{align}
Then for a given susceptibility distribution $\chi$ (in $\R^3$), the general solution $b$ which is bounded everywhere in $\R^3$ is expressed as
\begin{align}\label{GeneralSolution}
b(\x)=\int_{\R^3}\Phi(\x-\y)\left(-\f{1}{3}\Delta_{\y}+\f{\p^2}{\p y_3^2}\right)\chi(\y)\rd\y+b_0
\end{align}
where $b_0$ is some constant, and $\Phi(\x)=1/\left(4\pi|\x|\right)$.

In MRI, the phase of a complex GRE MR signal is linear with respect to the total field $b$ in \cref{GeneralSolution} \cite{Y.Wang2015}, and the constant $b_0$ is determined by the coil sensitivity of an MR scanner as the coil sensitivity dependent phase offset is in general assumed to be a constant \cite{E.M.Haacke2015,F.Schweser2016}. However, since we can remove it during the phase estimation from the multi echo GRE signal \cite{L.deRochefort2008}, we assume that $b_0=0$ and
\begin{align}\label{Forward1}
b(\x)=\int_{\R^3}\Phi(\x-\y)\left(-\f{1}{3}\Delta_{\y}+\f{\p^2}{\p y_3^2}\right)\chi(\y)\rd\y
\end{align}
in the rest of this paper. Note that $b$ defined as above is induced by the susceptibility distribution \emph{in the entire space,} which is different from $b_l$ in \cref{QSM_Deconv}.

\begin{remark}\label{RK211}
Since \cite[Proposition A.1.]{J.K.Choi2014} has discussed the equivalence between \cref{Forward1} and the following representation in the literature
\begin{align}\label{Forward_Convention}
b(\x)=\pv\int_{\R^3}d(\x-\y)\chi(\y)\rd\y,
\end{align}
we shall use \cref{Forward1} in the rest of this paper. Note that \cref{Forward1} avoids the singularity of the kernel $d(\x-\y)$ in \cref{Forward_Convention} as $\Phi(\x-\y)$ is locally integrable near $\x=\y$.
\end{remark}

\subsection{Characterization of Harmonic Incompatibility in Local Field Data}\label{LocalFieldAcquisition}

In QSM, the total field $b(\x)$ is obtained from the phase data of a complex GRE MR signal \cite{F.Schweser2016,Y.Wang2015}. In fact, if the information of $b$ is available over the entire space, then we can directly solve inverse problem from the knowledge of $b$ without the background field removal step. However, since the GRE signal is not available outside $\Om$, the information of $b$ is available only inside $\Om$. Moreover, even if $\chi$ is compactly supported, the support of $b$ may not necessarily coincide with that of $\chi$, which inevitably leads to the information loss outside $\Om$ \cite{F.Schweser2016,Y.Wang2015}.

Since the total field $b$ depends on the susceptibility distribution \emph{throughout the entire space} \cite{F.Schweser2016}, it consists of the background field induced from the susceptibility outside $\Om$, which is of no interest, and the local field $b_l$ by the susceptibility in $\Om$ which we aim to visualize. Since the substantial susceptibility sources are usually located outside $\Om$ which makes the background field dominant in $b$ compared to the local field $b_l$, we need to remove the background field from the (incomplete) total field prior to the dipole inversion \cite{F.Schweser2016,Y.Wang2015}.

In the literature, given that the background field is harmonic in $\Om$ \cite{Y.Wang2015,D.Zhou2014}, the background field removal methods take the form of the following Poisson's equation in \cite{D.Zhou2014}:
\begin{align}\label{LBV}
\left\{\begin{array}{rcll}
-\Delta b_l\hspace{-0.65em}&=\hspace{-0.65em}&-\Delta b&\text{in}~\Om \vspace{0.125em}\\
b_l\hspace{-0.65em}&=\hspace{-0.65em}&0&\text{on}~\p\Om.
\end{array}\right.
\end{align}
Under this setting, we present \cref{Th1} which characterizes the relation between \cref{LBV} and the PDE \cref{QSM_PDE2}, and the measured local field obtained by solving \cref{LBV} contains an incompatibility which consists of two harmonic functions both inside and outside $\Om$ due to the imposed boundary condition.

\begin{theorem}\label{Th1} Let $\Om\subseteq\R^3$ be an open and bounded set with $C^1$ boundary $\p\Om$. Let $b$ satisfy \cref{Forward1} for a given $\chi$ compactly supported in $\R^3$, and let $b_l:\overline{\Om}\to\R$ be obtained from \cref{LBV}. If we extend $b_l$ into $\R^3$ by assigning $b_l(\x)=0$ for $\x\notin\Om$, then we have the followings:
\end{theorem}
{\it\begin{enumerate}
\item There exists $v(x)$ such that
\begin{align}\label{NewModeling}
b_l(\x)=\int_{\Om}\Phi(\x-\y)\left(-\f{1}{3}\Delta_{\y}+\f{\p^2}{\p y_3^2}\right)\chi(\y)\rd\y+v(\x)
\end{align}
for $\x\in\R^3$, and $v(\x)$ satisfies
\begin{align}\label{DistLapIncompatibility}
\int_{\R^3}v(\x)\left(-\Delta\varphi\right)(\x)\rd\x=\int_{\p\Om}\left[\f{\p v_{\mathrm i}}{\p\bn}(\x)-\f{\p v_{\mathrm e}}{\p\bn}(\x)\right]\varphi(\x)\rd\ssi(\x)
\end{align}
for $\varphi\in C_0^{\infty}(\R^3)$, where $v_{\mathrm i}$ and $v_{\mathrm e}$ denote the restriction of $v$ in $\overline{\Om}$ and $\R^3\setminus\Om$ respectively, and $\bn$ denotes the outward unit normal vector of $\p\Om$.
\item Moreover, we have
\begin{align}\label{NonVanishingSurfMeas}
\f{\p v_{\mathrm i}}{\p\bn}-\f{\p v_{\mathrm e}}{\p\bn}\neq0~~~~\text{almost everywhere on}~\p\Om
\end{align}
whenever $P(D)\chi\neq0$ in $\Om$. Hence, $-\Delta v=0$ in $\R^3\setminus\p\Om$, and $-\Delta v\neq0$ on $\p\Om$ in this case.
\end{enumerate}}


\begin{proof} Since $-\Delta b=P(D)\chi$, the governing equation in \cref{LBV} becomes
\begin{align*}
\left\{\begin{array}{rcll}
-\Delta b_l\hspace{-0.65em}&=\hspace{-0.65em}&P(D)\chi&\text{in}~\Om \vspace{0.125em}\\
b_l\hspace{-0.65em}&=\hspace{-0.65em}&0&\text{on}~\p\Om.
\end{array}\right.
\end{align*}
Let $G(\x,\y)$ denote the Green's function in $\Om$:
\begin{align*}
G(\x,\y)=\Phi(\y-\x)-H(\x,\y)
\end{align*}
where for each $\x\in\Om$, the corrector function $H(\x,\y)$ satisfies
\begin{align*}
\left\{\begin{array}{rclll}
-\Delta_{\y} H(\x,\y)\hspace{-0.65em}&=\hspace{-0.65em}&0&\text{if}&\y\in\Om \vspace{0.125em}\\
H(\x,\y)\hspace{-0.65em}&=\hspace{-0.65em}&\Phi(\y-\x)&\text{if}&\y\in\p\Om.
\end{array}\right.
\end{align*}
Note that since $G(\x,\y)=G(\y,\x)$ for $\x,\y\in\Om$ and $\Phi(\y-\x)=\Phi(\x-\y)$, we have $H(\x,\y)=H(\y,\x)$ for $\x,\y\in\Om$. Consequently, we have
\begin{align}\label{CorrectorHarmonicSym}
-\Delta_{\x}H(\x,\y)=-\Delta_{\x}H(\y,\x)=0~~~~~\x\in\Om.
\end{align}
Then the solution to \cref{LBV} is represented as
\begin{align*}
b_l(\x)=\int_{\Om}G(\x,\y)P(D_{\y})\chi(\y)\rd\y=\int_{\Om}\Phi(\x-\y)P(D_{\y})\chi(\y)\rd\y+\mH(\x)
\end{align*}
where we used the fact that $\Phi(\y-\x)=\Phi(\x-\y)$, and $\mH(\x)$ is defined as
\begin{align*}\label{HarmonicInterior}
\mH(\x)=-\int_{\Om}H(\x,\y)P(D_{\y})\chi(\y)\rd\y
\end{align*}
for $\x\in\Om$. Then we can see that $\mH(\x)$ satisfies
\begin{align}\label{H_HarmonicEq}
\left\{\begin{array}{rcll}
-\Delta\mH\hspace{-0.65em}&=\hspace{-0.65em}&0&\text{in}~\Om \vspace{0.125em}\\
\mH\hspace{-0.65em}&=\hspace{-0.65em}&-\wt{b}_l&\text{on}~\p\Om
\end{array}\right.
\end{align}
where the first equation of \cref{H_HarmonicEq} comes from \cref{CorrectorHarmonicSym}. Here, $\wt{b}_l$ is induced by the information of $\chi$ only in $\Om$:
\begin{align}\label{TrueLocalField}
\wt{b}_l(\x)=\int_{\Om}\Phi(\x-\y)P(D_{\y})\chi(\y)\rd\y=\int_{\R^3}\Phi(\x-\y)1_{\Om}(\y)P(D_{\y})\chi(\y)\rd\y
\end{align}
with $1_{\Om}$ being the characteristic function of $\Om$.



Based on the fact that $b_l(\x)=0$ for $\x\in\R^3\setminus\overline{\Om}$, we define
\begin{align*}
v(\x)=\left\{\begin{array}{cccll}
v_{\mathrm i}(\x)\hspace{-0.65em}&=\hspace{-0.65em}&\mH(\x)&\text{if}&\x\in\Om \vspace{0.125em}\\
v_{\mathrm e}(\x)\hspace{-0.65em}&=\hspace{-0.65em}&-\wt{b}_l(\x)&\text{if}&\x\notin\Om.
\end{array}\right.
\end{align*}
Hence, we obtain \cref{NewModeling}, and we can further see that $v_{\mathrm i}$ and $v_{\mathrm e}$ satisfy
\begin{align}
-\Delta v_{\mathrm i}&=0~~~~~~\text{in }  \Om \label{v_iharmonic}\\
-\Delta v_{\mathrm e}&=0~~~~~~\text{in }  \R^3\setminus\overline{\Om} \label{v_eharmonic}\\
v_{\mathrm i}=v_{\mathrm e}&=-\wt{b}_l~~~\text{on } \p\Om, \label{vboundaryconti}
\end{align}
respectively, where \cref{v_eharmonic} comes from \cref{TrueLocalField}; $-\Delta\wt{b}_l=1_{\Om}P(D)\chi$, i.e. $-\Delta\wt{b}_l=P(D)\chi$ in $\Om$, and $-\Delta\wt{b}_l=0$ in $\R^3\setminus\overline{\Om}$.

To prove \cref{DistLapIncompatibility}, let $\varphi\in C_0^{\infty}(\R^3)$, and we consider
\begin{align*}
\int_{\R^3}v(\x)\left(-\Delta\varphi\right)(\x)\rd\x&=\int_{\Om}v_{\mathrm i}(\x)\left(-\Delta\varphi\right)(\x)\rd\x+\int_{\R^3\setminus\overline{\Om}}v_{\mathrm e}(\x)\left(-\Delta\varphi\right)(\x)\rd\x\\
&=I_1+I_2
\end{align*}
Using \cref{v_iharmonic,vboundaryconti} and the Green's identity (e.g. \cite{Evans2010}), we have
\begin{align*}
I_1&=\int_{\Om}\left[\varphi(\x)\left(\Delta v_{\mathrm i}\right)(\x)-v_{\mathrm i}(\x)\left(\Delta\varphi\right)(\x)\right]\rd\x\\
&=\int_{\p\Om}\left[\varphi(\x)\f{\p v_{\mathrm i}}{\p\bn}(\x)-v_{\mathrm i}(\x)\f{\p\varphi}{\p\bn}(\x)\right]\rd\ssi(\x)=\int_{\p\Om}\left[\varphi(\x)\f{\p v_{\mathrm i}}{\p\bn}(\x)+\wt{b}_l(\x)\f{\p\varphi}{\p\bn}(\x)\right]\rd\ssi(\x).
\end{align*}
Similarly using \cref{v_eharmonic,vboundaryconti}, we also have
\begin{align*}
I_2&=\int_{\R^3\setminus\overline{\Om}}\left[\varphi(\x)\left(\Delta v_{\mathrm e}\right)(\x)-v_{\mathrm e}(\x)\left(\Delta\varphi\right)(\x)\right]\rd\x\\
&=\int_{\p\Om}\left[-\varphi(\x)\f{\p v_{\mathrm e}}{\p\bn}(\x)+v_{\mathrm e}(\x)\f{\p\varphi}{\p\bn}(\x)\right]\rd\ssi(\x)=-\int_{\p\Om}\left[\varphi(\x)\f{\p v_{\mathrm e}}{\p\bn}(\x)+\wt{b}_l(\x)\f{\p\varphi}{\p\bn}(\x)\right]\rd\ssi(\x)
\end{align*}
where the second equality comes from the fact that we need to compute the inward normal derivatives on $\p\Om$. Hence, combining these two equalities, we obtain \cref{DistLapIncompatibility}.

To prove 2, we assume on the contrary that there exists $\x\in\p\Om$ such that
\begin{align*}
\f{\p v_{\mathrm i}}{\p\bn}=\f{\p v_{\mathrm e}}{\p\bn}
\end{align*}
for some open and connected set $U\subseteq\p\Om$ such that $\x\in U$ and $\ssi(U)>0$. Choose $r>0$ such that $B(\x,r)\cap\p\Om$ is contained in $U$, where $B(\x,r)$ denotes an open ball centered at $\x$ with radius $r$. Then since $-\Delta v=0$ in $B(\x,r)$ and $v=v_{\mathrm e}=-\wt{b}_l$ in $B(\x,r)\cap\left(\R^3\setminus\overline{\Om}\right)$, it follows that $v=-\wt{b}_l$ in $B(\x,r)$ by the analyticity of $v$ in $B(\x,r)$. Since this means that $v=v_{\mathrm i}=-\wt{b}_l$ in $B(\x,r)\cap\Om$, together with the fact that $v_{\mathrm i}$ is harmonic in $\Om$, we have $v_{\mathrm i}=-\wt{b}_l$ in $\Om$ by the analyticity of $v_{\mathrm i}$ in $\Om$. Since $v=-\wt{b}_l$ on $\p\Om$, we have $v=-\wt{b}_l$ in $\R^3$. Hence, $-\Delta\wt{b}_l=0$ in $\R^3$, and thus, $P(D)\chi=0$ in $\Om$, which is a contradiction.
\end{proof}



\begin{remark}\label{RK221} From the proof of \cref{Th1},  the incompatibility $v$ in \cref{NewModeling} is from the boundary condition of \cref{LBV} which is not related to the regularity of $\chi$. More precisely, $v_{\mathrm i}$ inside $\Om$ is generated by the information of the unknown true local field $\wt{b}_l$ on $\p\Om$ so that the boundary condition of \cref{LBV} is matched. In addition, it is obvious that $v_{\mathrm e}$ outside $\Om$ is due to the information loss outside $\Om$.
\end{remark}

\begin{remark}\label{RK222} Notice that $P(D)$ is a ``wave type'' differential operator (by considering $x_3$ as the time variable). Indeed, the proof of \cref{NonVanishingSurfMeas} tells us that if $P(D)\chi=0$ in $\Om$, such $\chi$ has a wave type structure in $\Om$ regardless of its regularity, whereas the susceptibility of human brain does not have such a wave type structure \cite{J.K.Choi2014}. Hence in QSM, it follows that $-\Delta v$ defined as \cref{DistLapIncompatibility} is a nonvanishing surface measure on $\p\Om$, i.e. $-\Delta v=0$ in $\R^3\setminus\p\Om$, but $-\Delta v\neq0$ on $\p\Om$.
\end{remark}

\begin{figure}[tp!]
\centering
\hspace{-0.1cm}\subfloat[ROI $\Om$]{\label{LoganMask}\includegraphics[width=3.60cm]{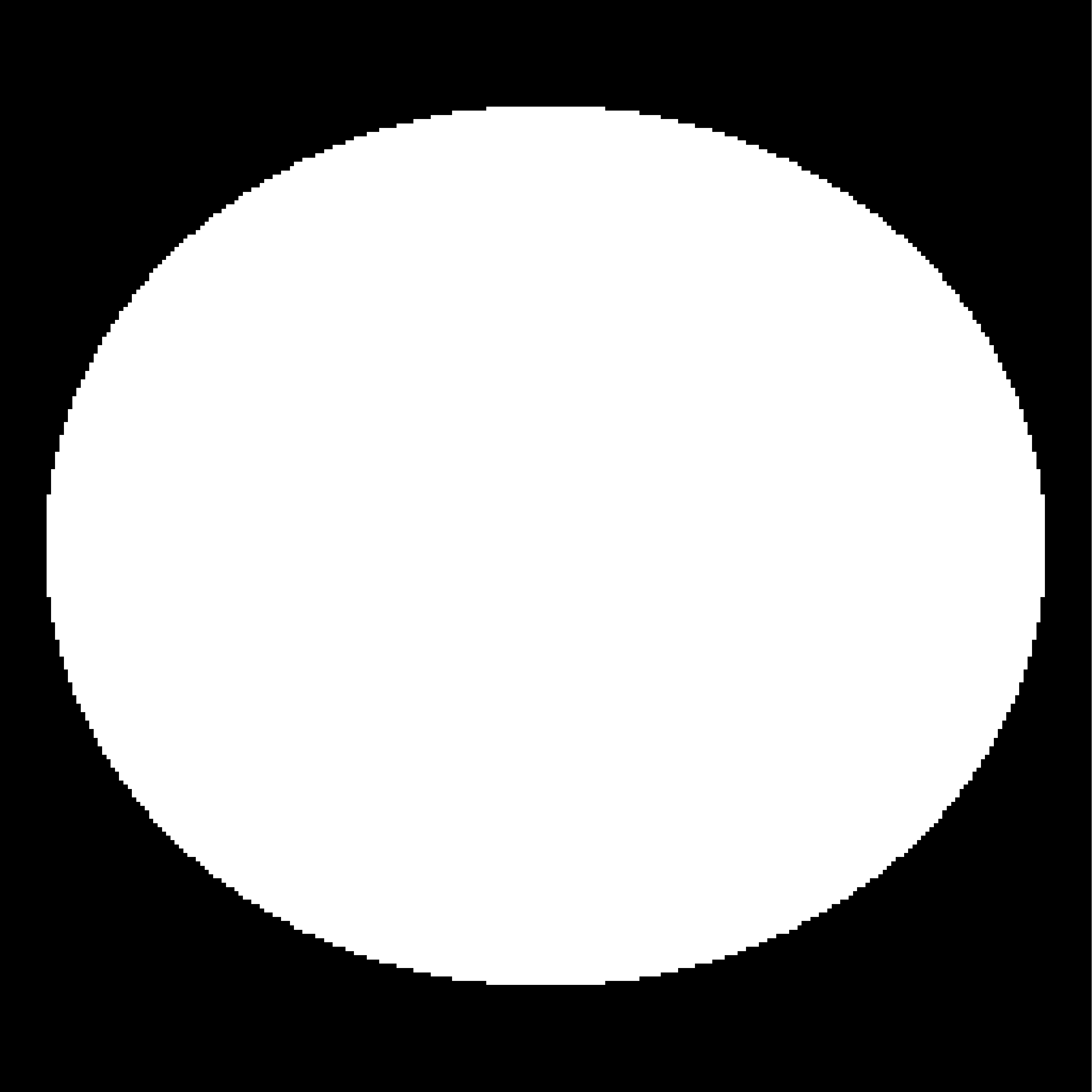}}\hspace{0.005cm}
\subfloat[$\chi$ in $\Om$]{\label{LoganQSM}\includegraphics[width=3.60cm]{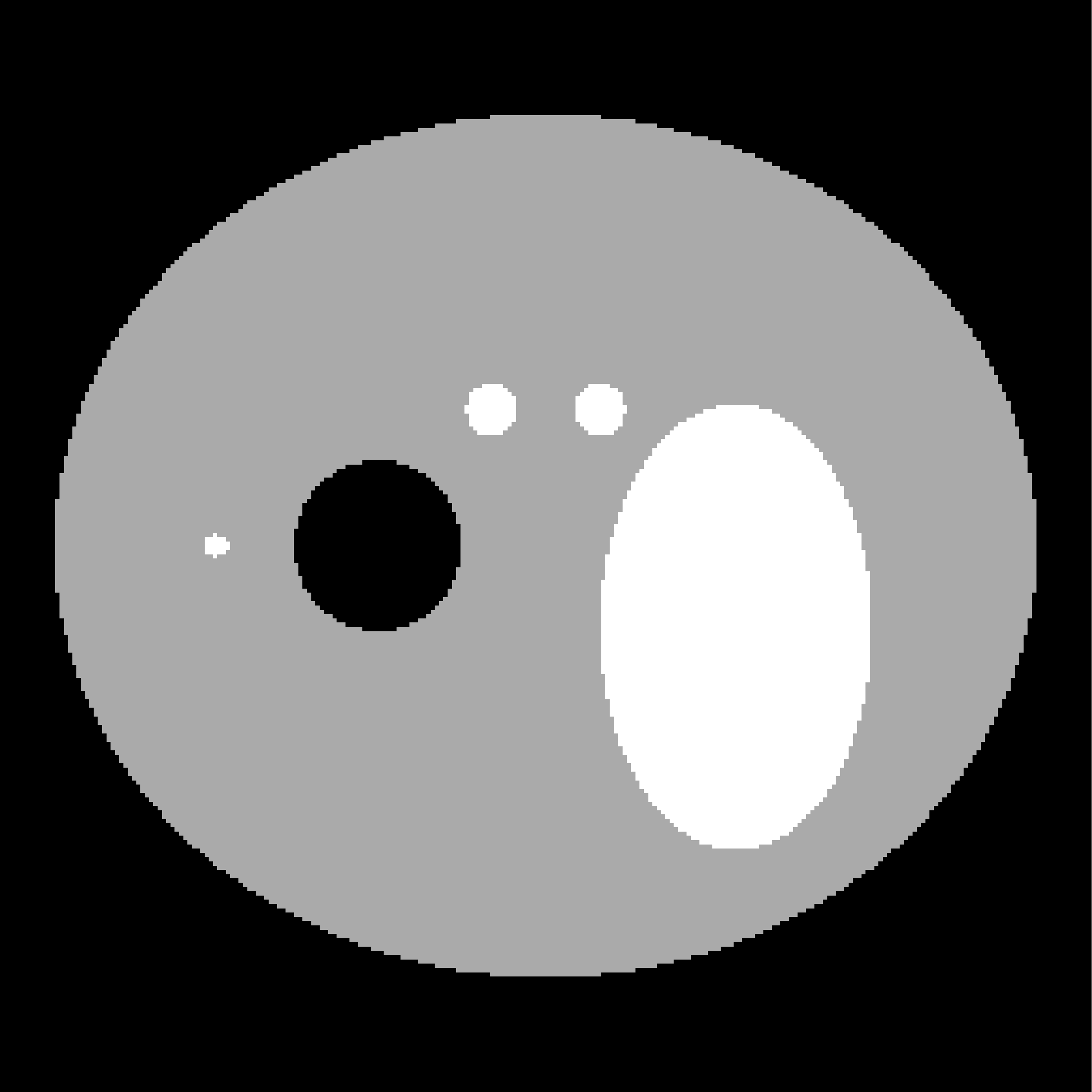}}\hspace{0.005cm}
\subfloat[$\chi$ in $\R^3\setminus\Om$]{\label{LoganSusBG}\includegraphics[width=3.60cm]{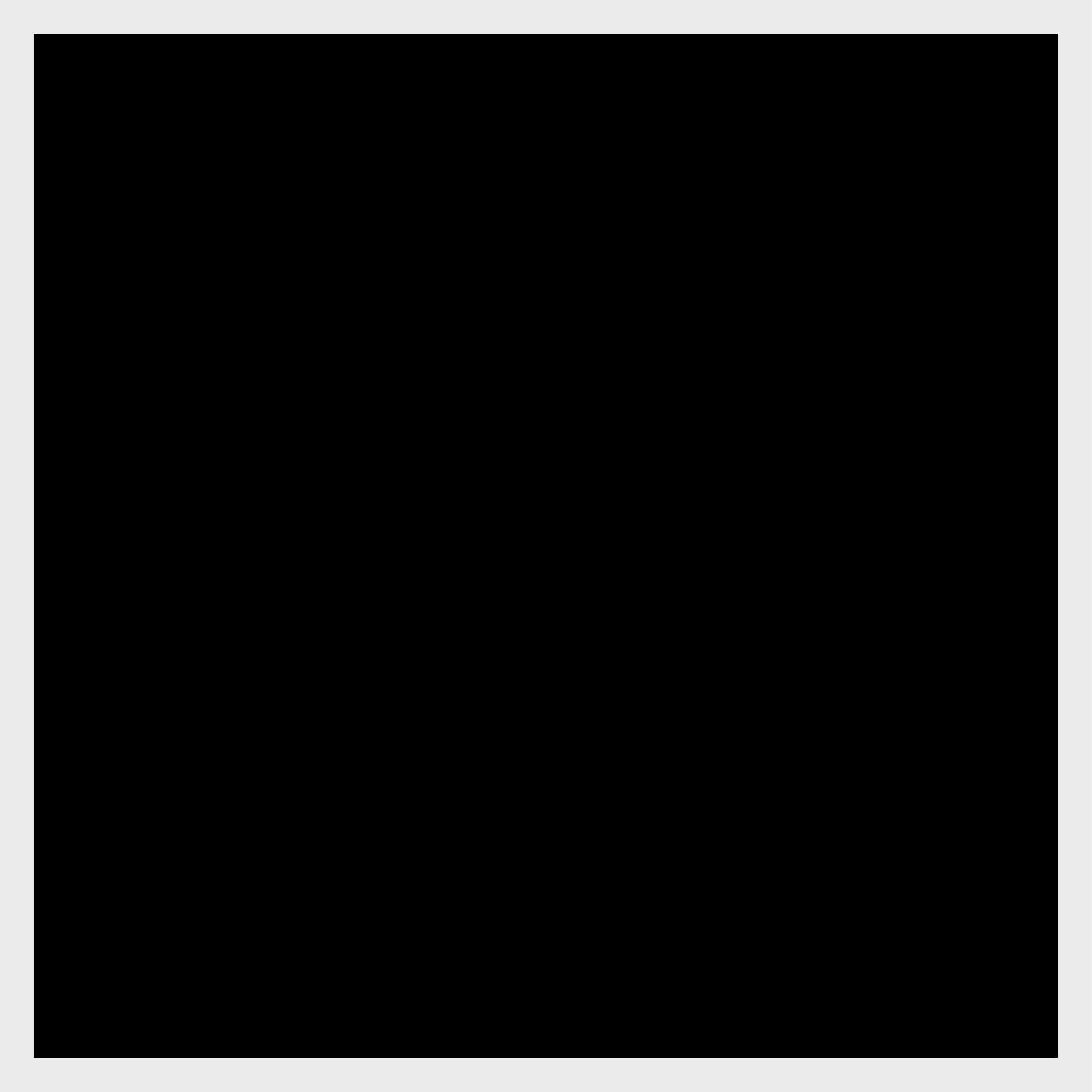}}\hspace{0.005cm}
\subfloat[Simulated total field $b$]{\label{LoganTotal}\includegraphics[width=3.60cm]{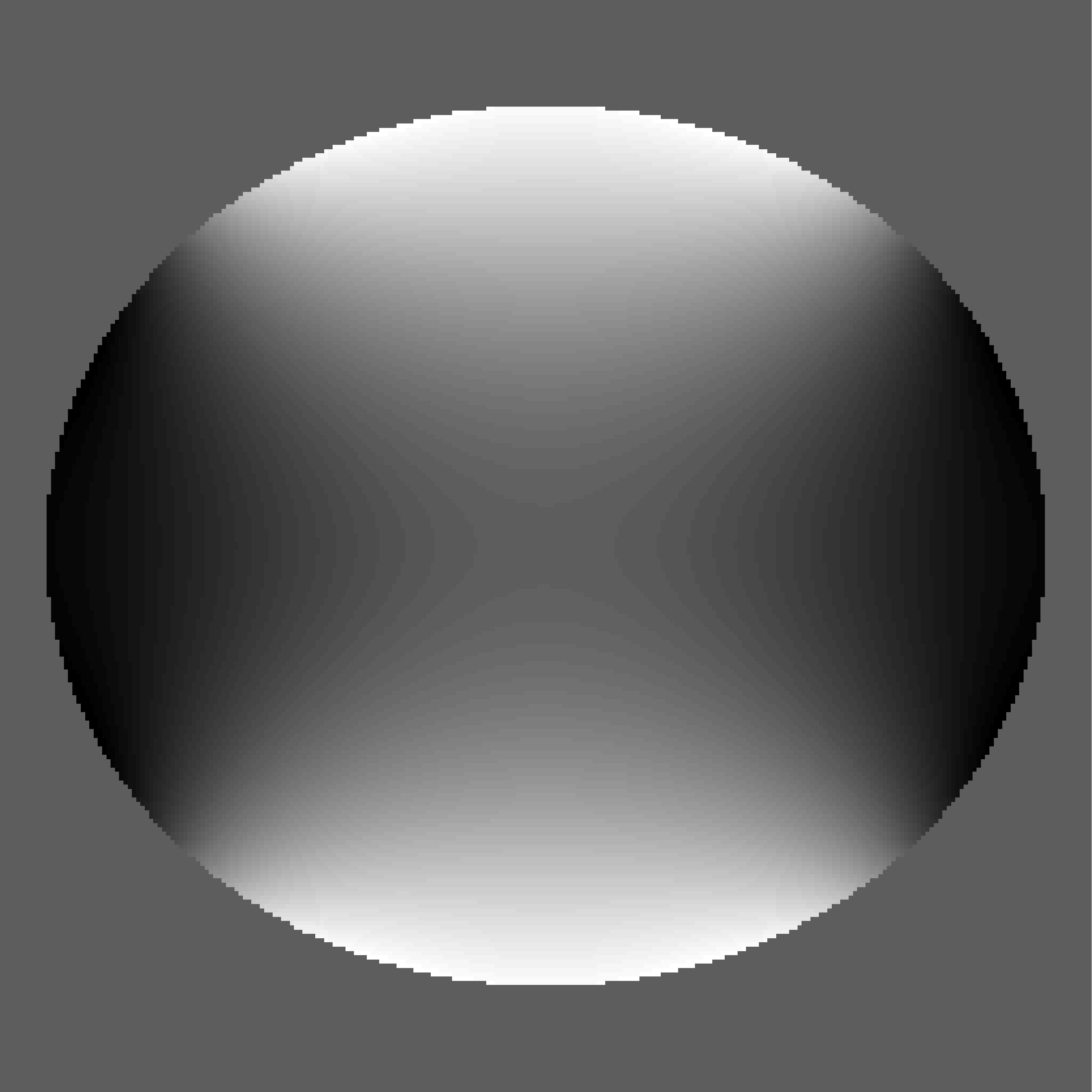}}\vspace{-0.20cm}\\
\subfloat[True local field $\wt{b}_l$]{\label{LoganTrueLocal}\includegraphics[width=3.60cm]{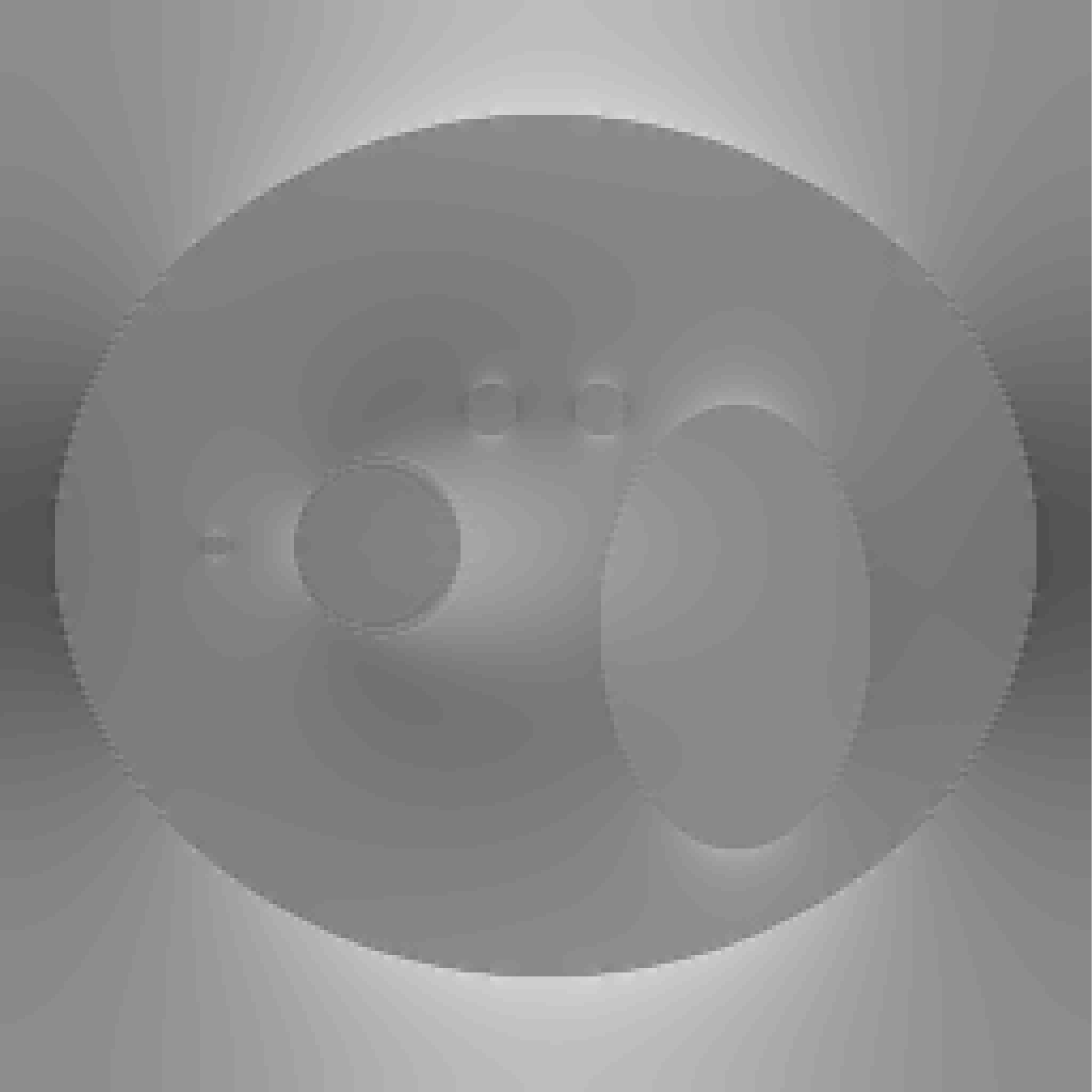}}\hspace{0.005cm}
\subfloat[Measured local field $b_l$]{\label{LoganMeasLocal}\includegraphics[width=3.60cm]{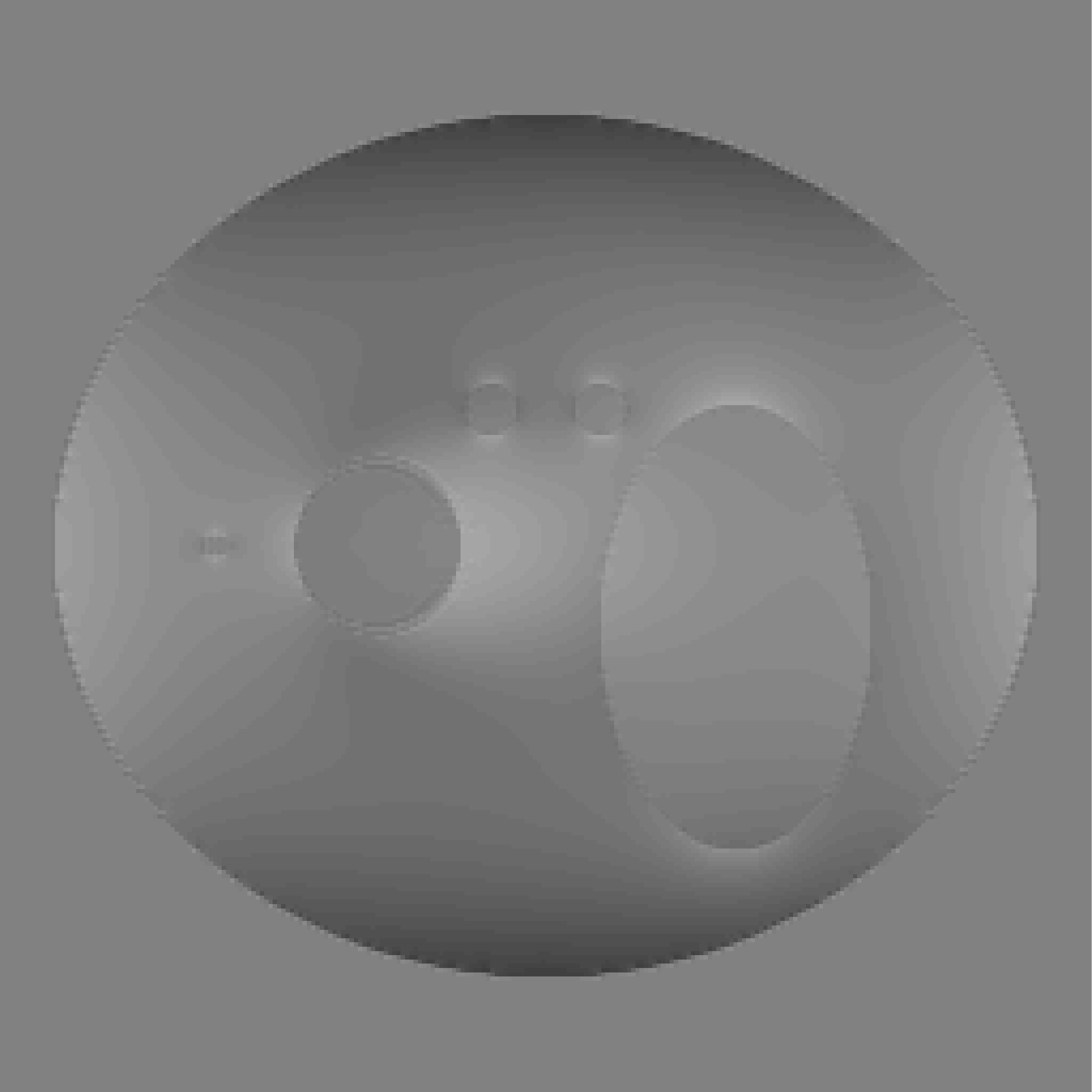}}\vspace{-0.20cm}
\subfloat[$v=b_l-\wt{b}_l$]{\label{LoganIncomp}\includegraphics[width=3.60cm]{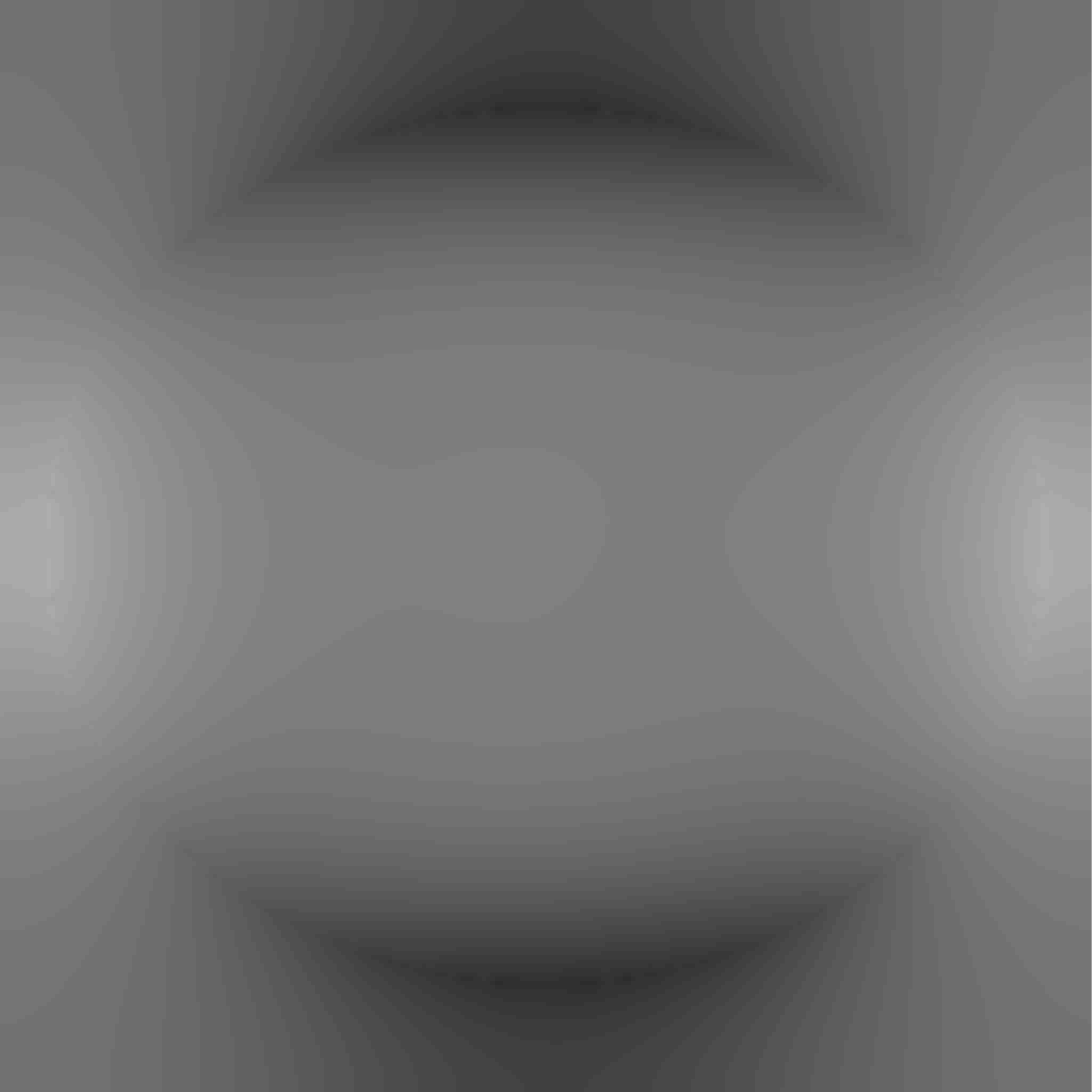}}\hspace{0.005cm}
\subfloat[$\left|-\Delta v\right|$]{\label{LoganAbsLv}\includegraphics[width=3.60cm]{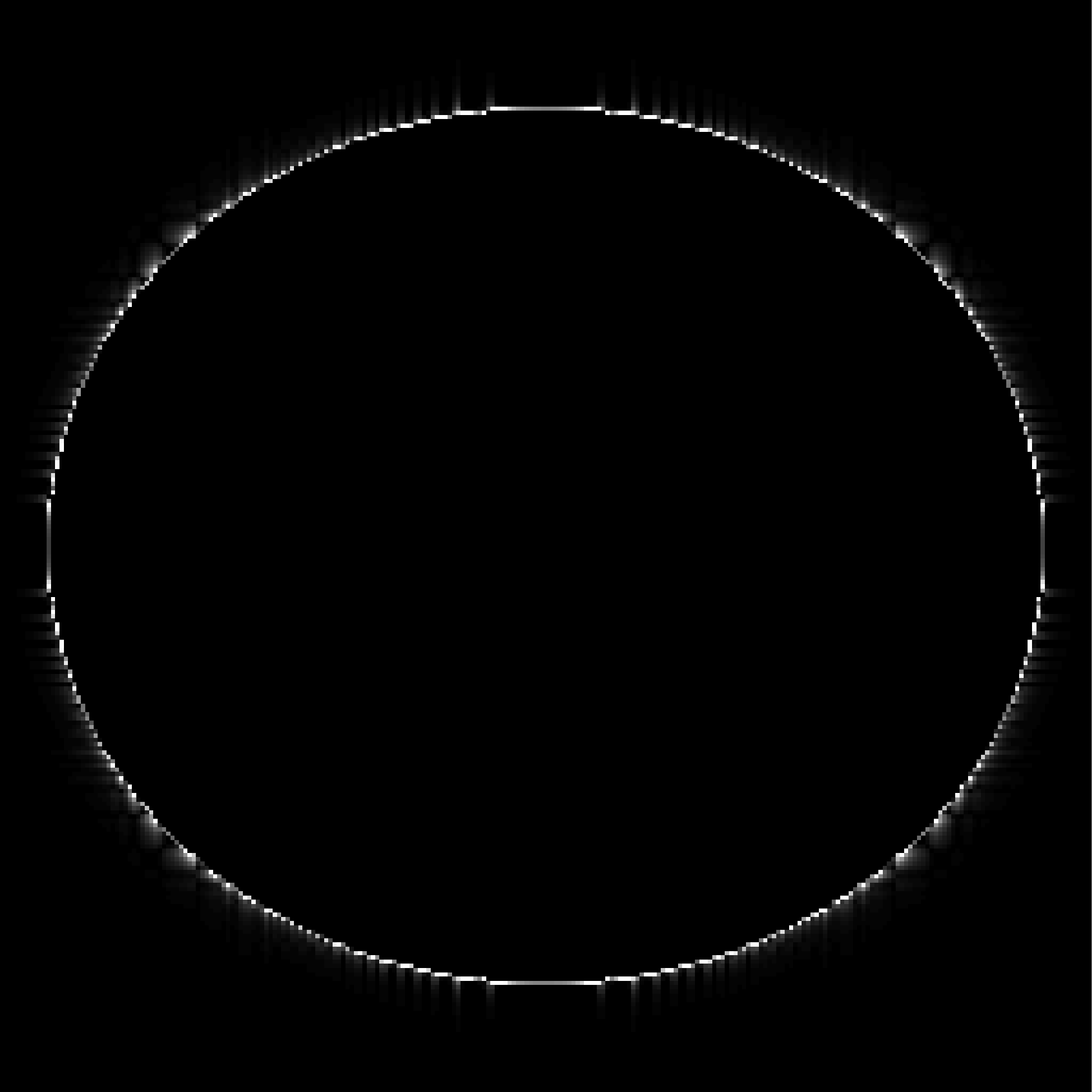}}\vspace{-0.20cm}
\caption{Sagittal slice images of $256\times256\times256$ Shepp-Logan phantom with $1\times1\times1\mathrm{mm}^3$. $\chi$ in $\R^3\setminus\Om$ is displayed in the window level $[0,550]$, $\wt{b}_l$, $b_l$ and $v$ in the window level $[-0.025,0.025]$, and $\left|-\Delta v\right|$ in the window level $[0,0.001]$.}\label{IllustrateTh1}
\end{figure}
\begin{figure}[tp!]
\centering
\hspace{-0.1cm}\subfloat[ROI $\Om$]{\label{LoganMaskAx}\includegraphics[width=3.60cm]{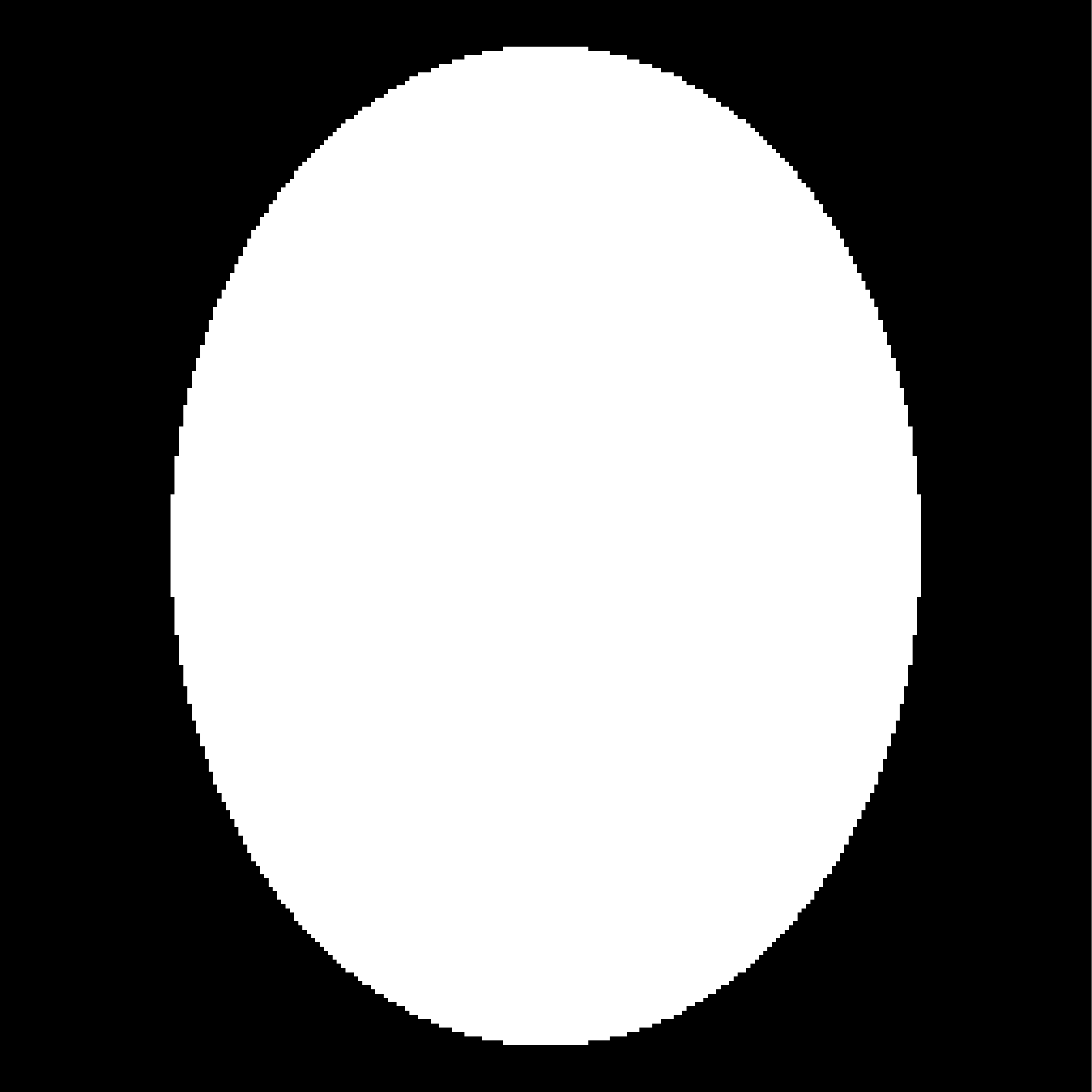}}\hspace{0.005cm}
\subfloat[$\chi$ in $\Om$]{\label{LoganQSMAx}\includegraphics[width=3.60cm]{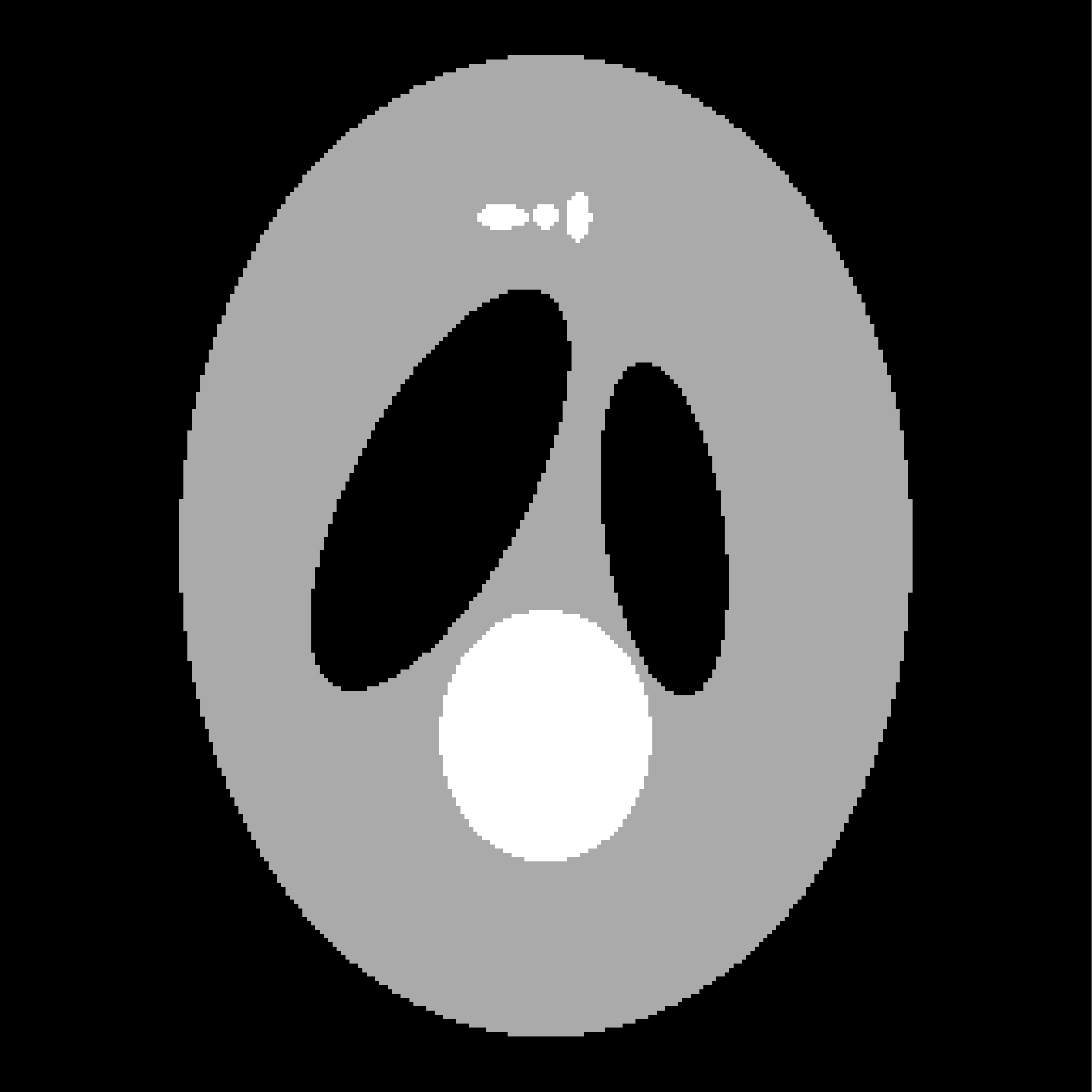}}\hspace{0.005cm}
\subfloat[$\chi$ in $\R^3\setminus\Om$]{\label{LoganSusBGAx}\includegraphics[width=3.60cm]{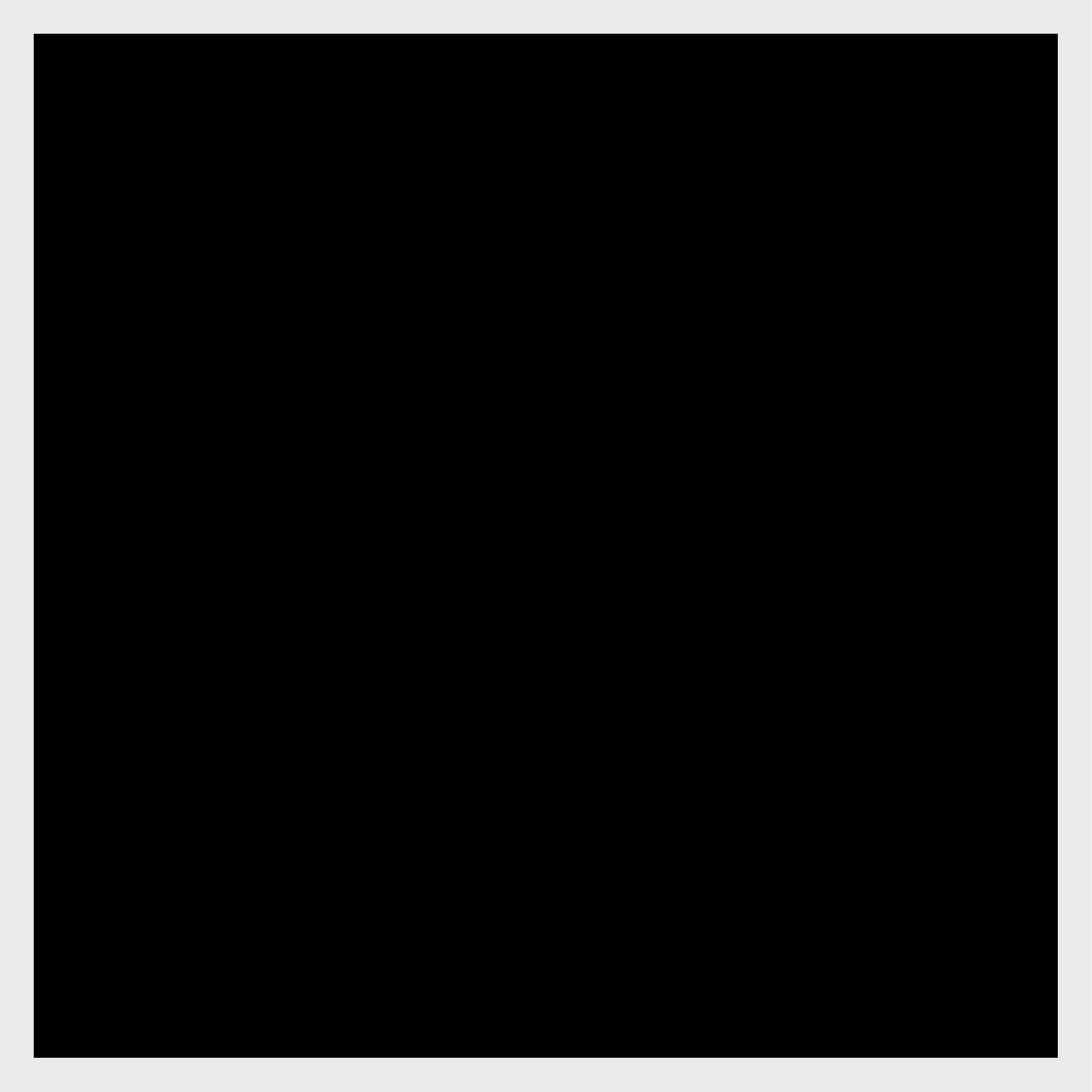}}\hspace{0.005cm}
\subfloat[Simulated total field $b$]{\label{LoganTotalAx}\includegraphics[width=3.60cm]{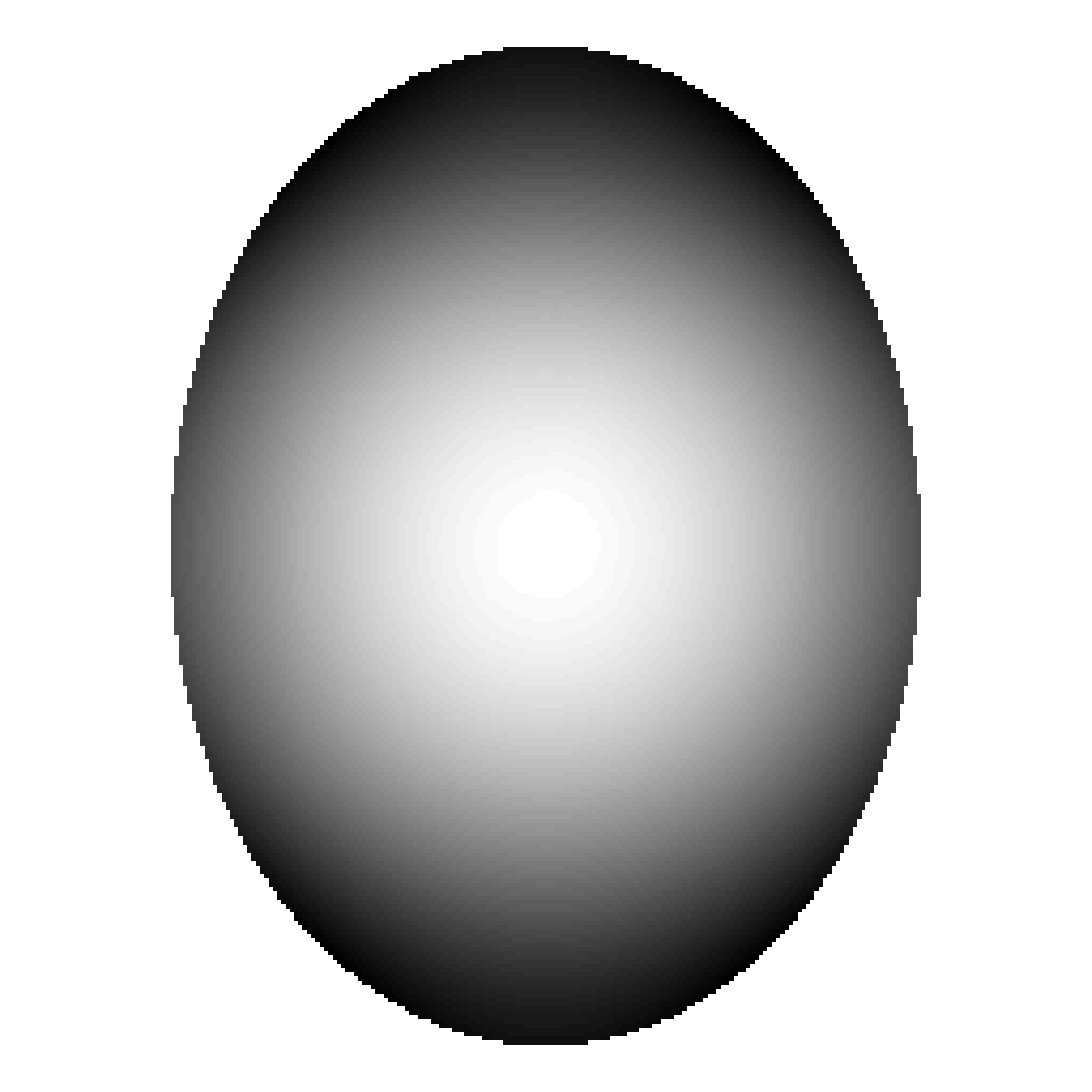}}\vspace{-0.20cm}\\
\subfloat[True local field $\wt{b}_l$]{\label{LoganTrueLocalAx}\includegraphics[width=3.60cm]{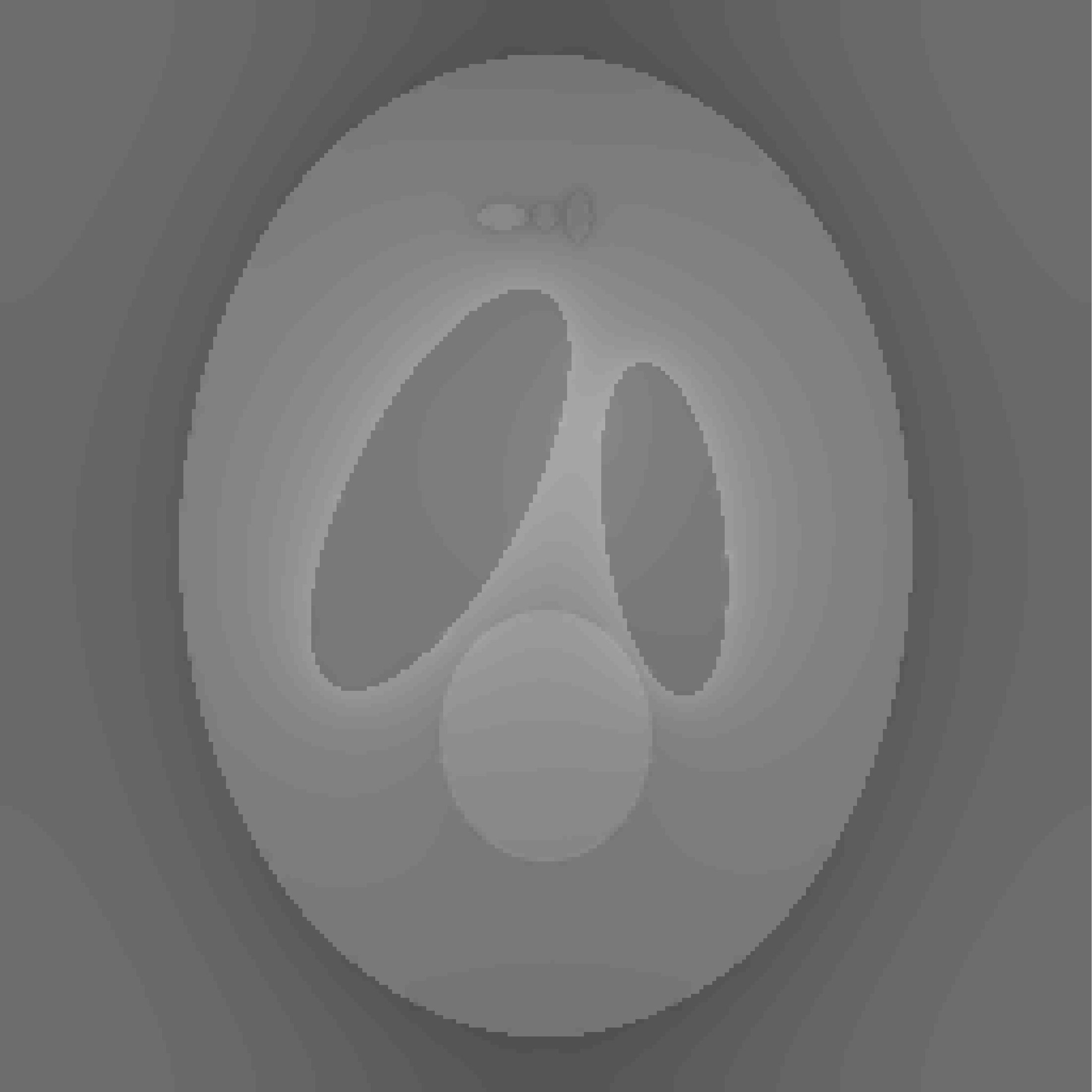}}\hspace{0.005cm}
\subfloat[Measured local field $b_l$]{\label{LoganMeasLocalAx}\includegraphics[width=3.60cm]{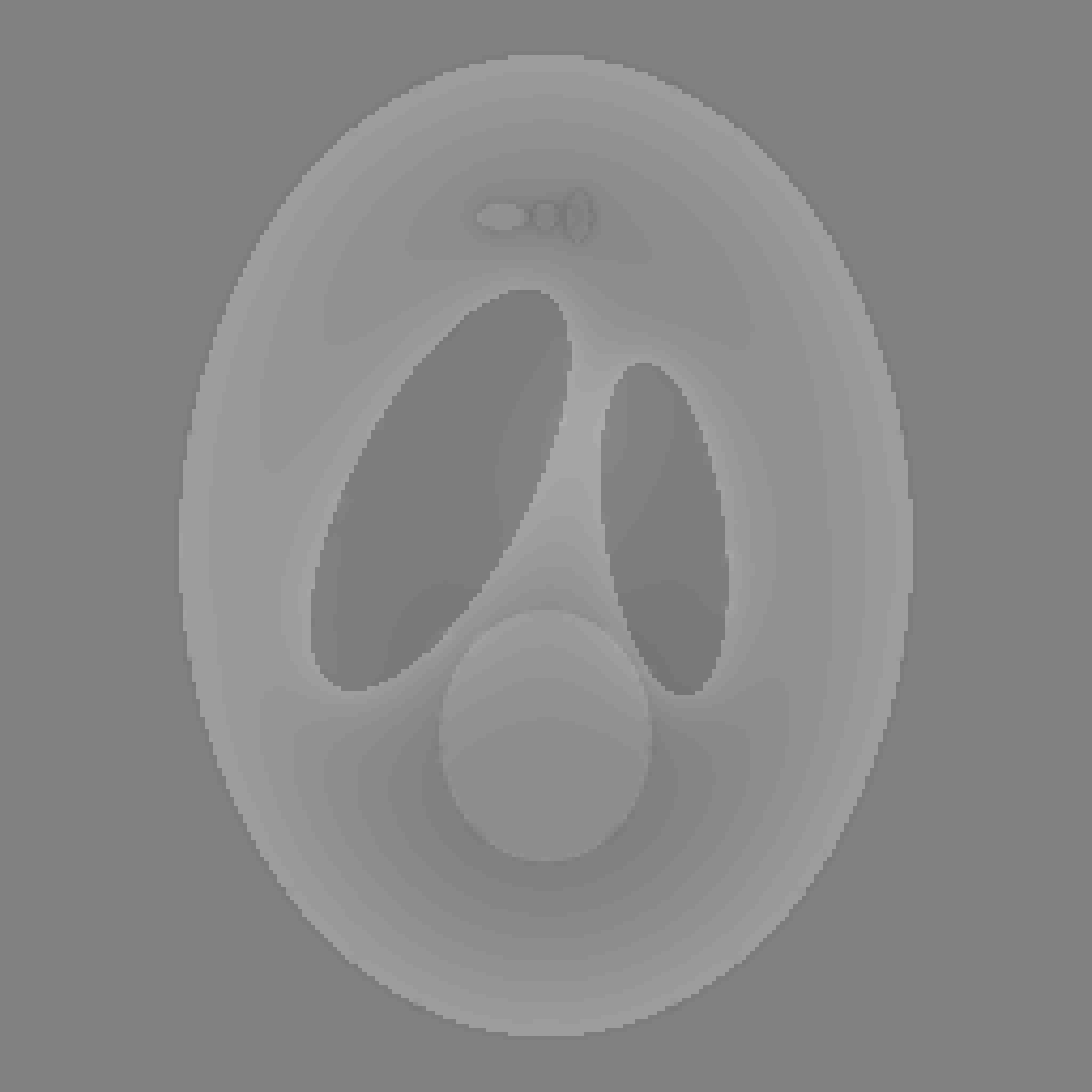}}\hspace{0.005cm}
\subfloat[$v=b_l-\wt{b}_l$]{\label{LoganIncompAx}\includegraphics[width=3.60cm]{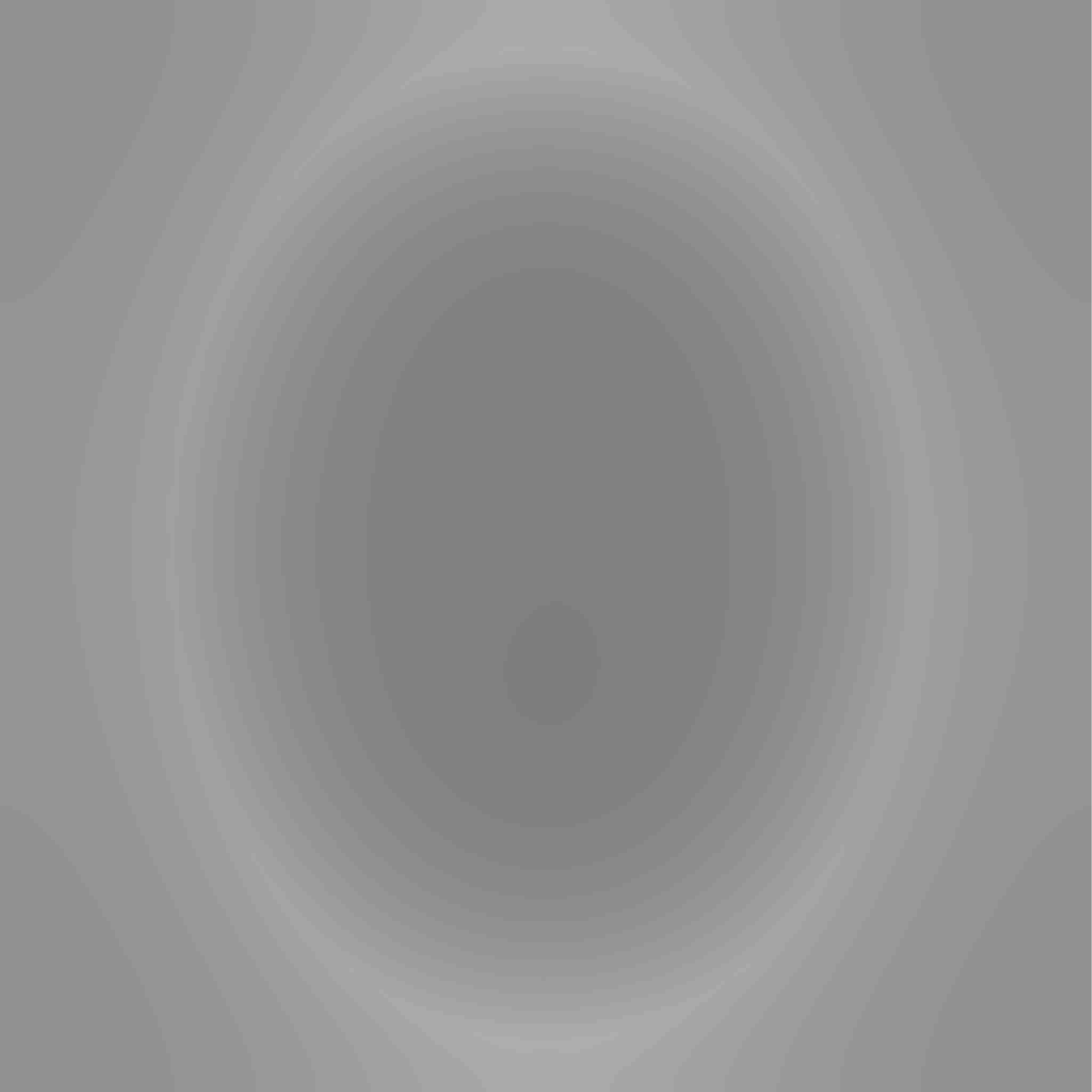}}\hspace{0.005cm}
\subfloat[$\left|-\Delta v\right|$]{\label{LoganAbsLvAx}\includegraphics[width=3.60cm]{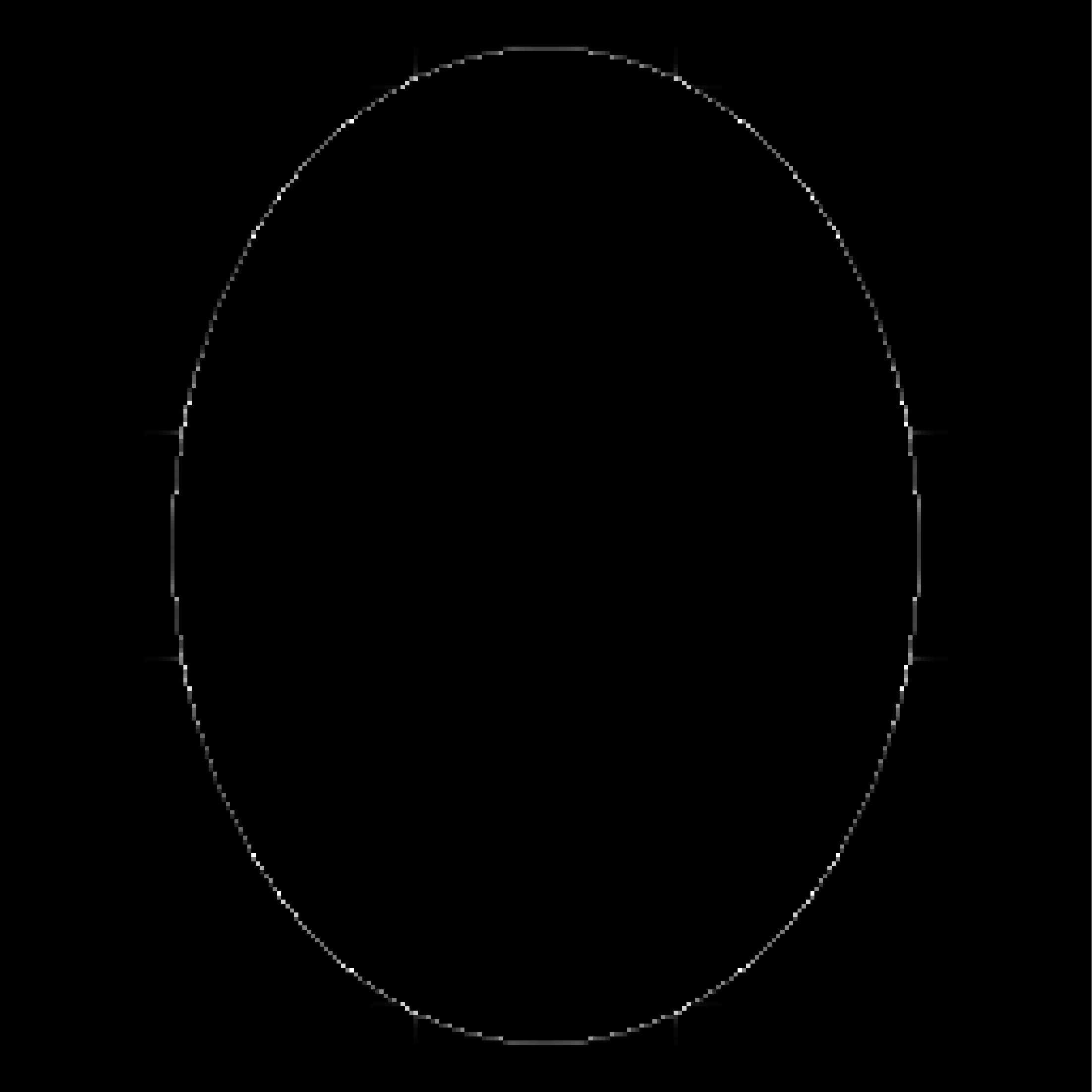}}\vspace{-0.20cm}
\caption{Axial slice images of \cref{IllustrateTh1}. The images of $\chi$ in $\R^3\setminus\Om$, $\wt{b}_l$, $b_l$, $v$, and $\left|-\Delta v\right|$ are displayed in the same window level as \cref{IllustrateTh1}.}\label{IllustrateTh1Axial}
\end{figure}
\begin{figure}[tp!]
\centering
\hspace{-0.1cm}\subfloat[ROI $\Om$]{\label{BrainMask}\includegraphics[width=3.60cm]{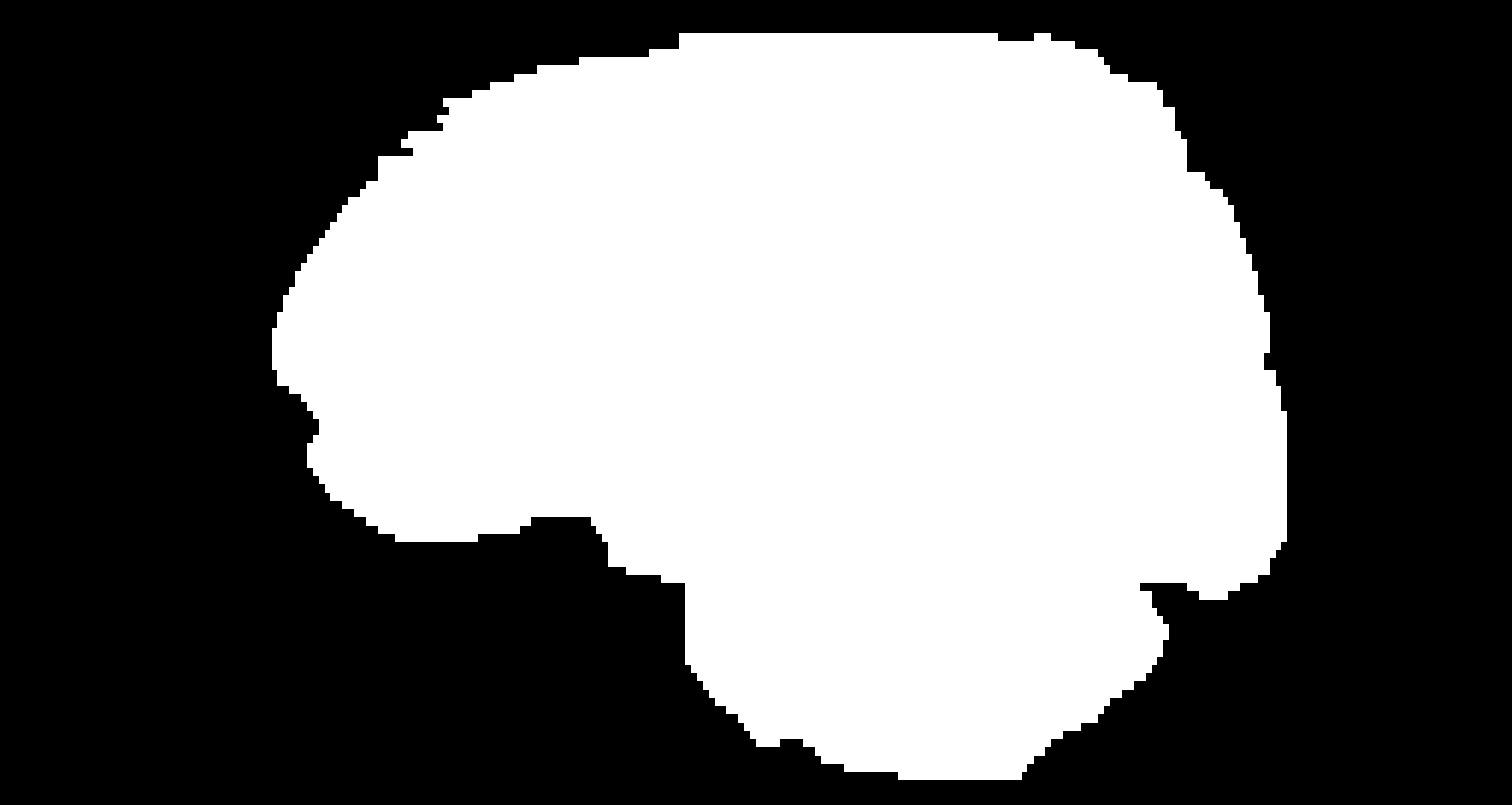}}\hspace{0.005cm}
\subfloat[$\chi$ in $\Om$]{\label{BrainQSM}\includegraphics[width=3.60cm]{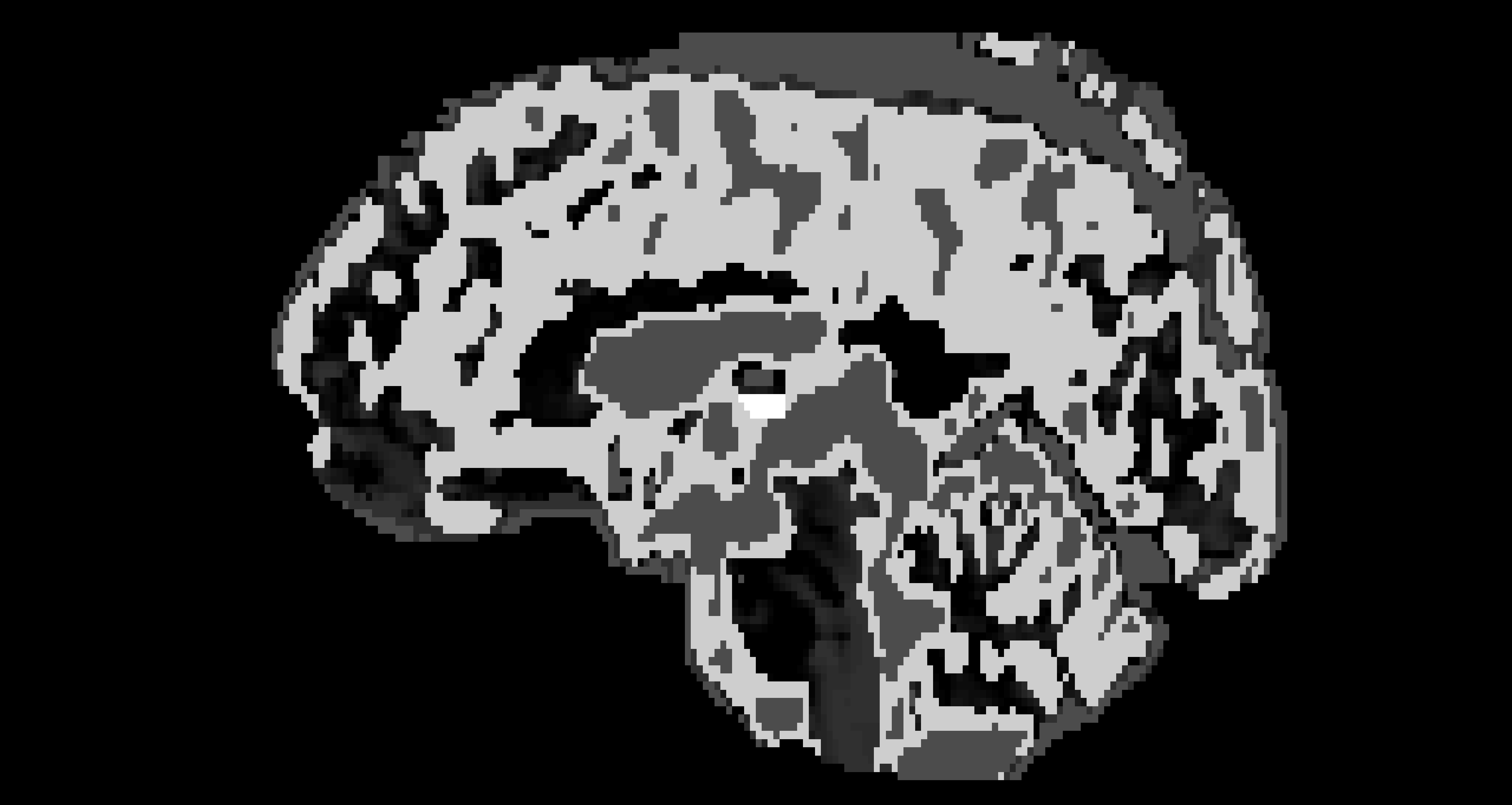}}\hspace{0.005cm}
\subfloat[$\chi$ in $\R^3\setminus\Om$]{\label{BrainSusBG}\includegraphics[width=3.60cm]{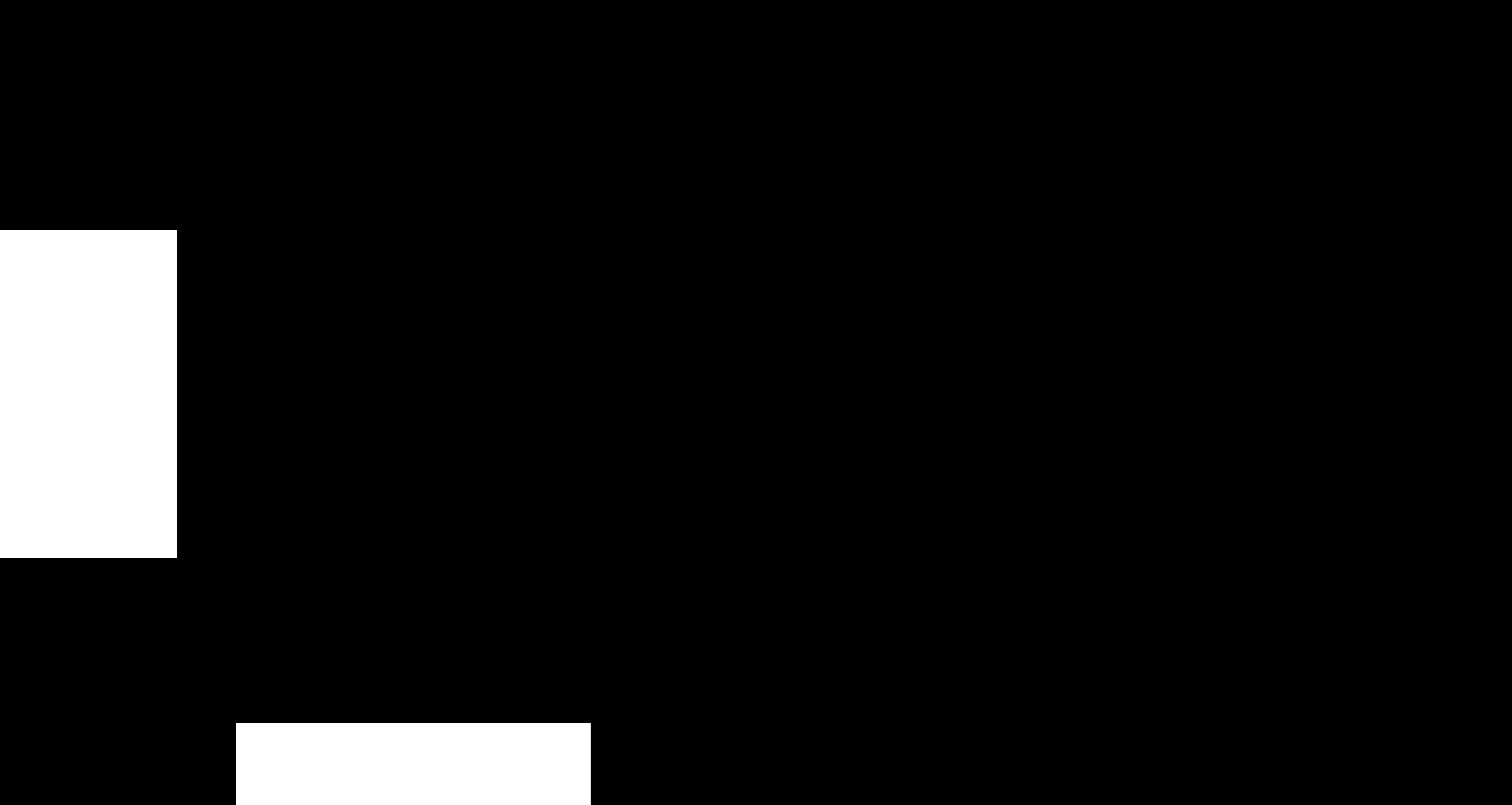}}\hspace{0.005cm}
\subfloat[Simulated total field $b$]{\label{BrainTotal}\includegraphics[width=3.60cm]{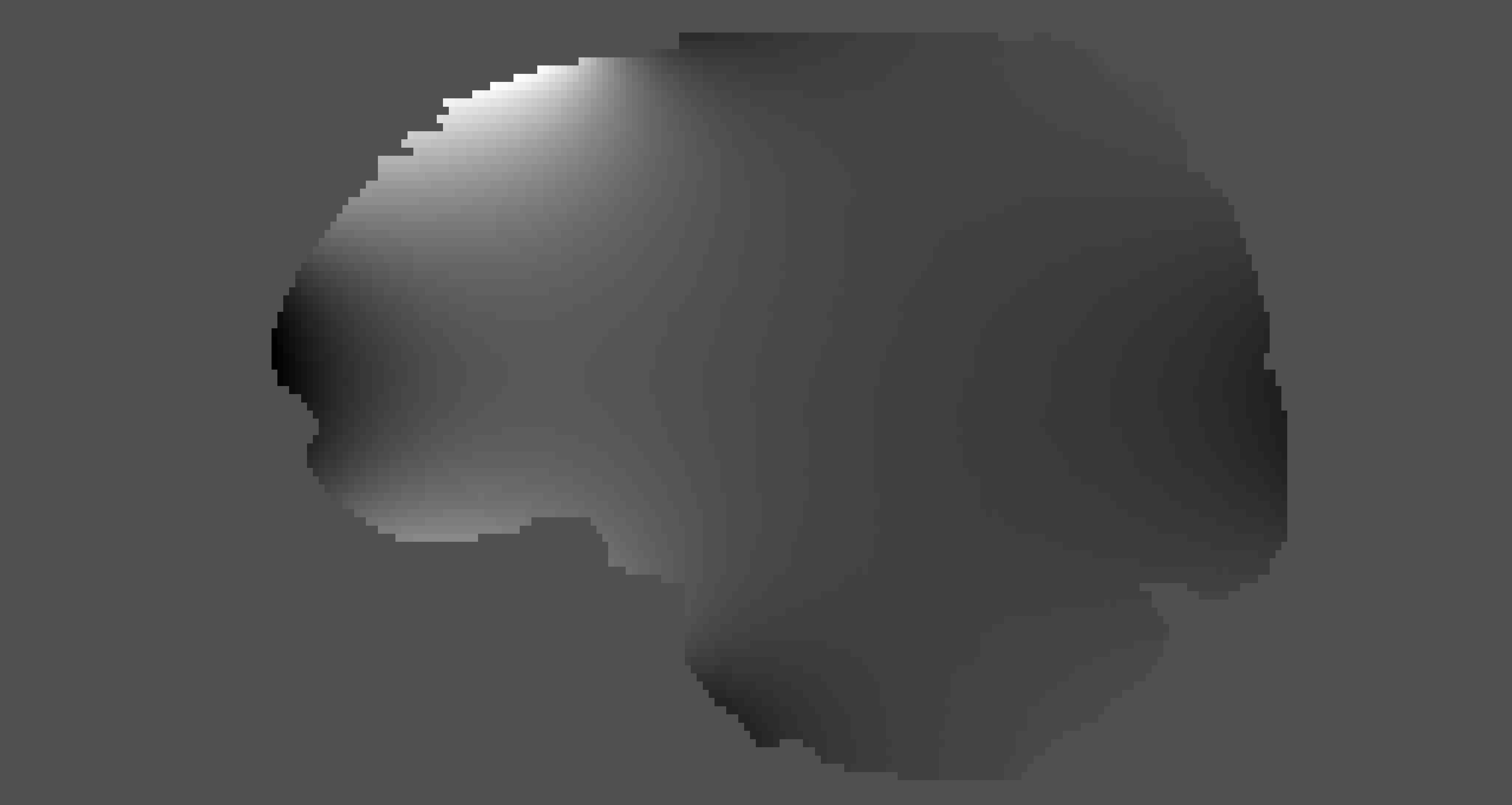}}\vspace{-0.20cm}\\
\subfloat[True local field $\wt{b}_l$]{\label{BrainTrueLocal}\includegraphics[width=3.60cm]{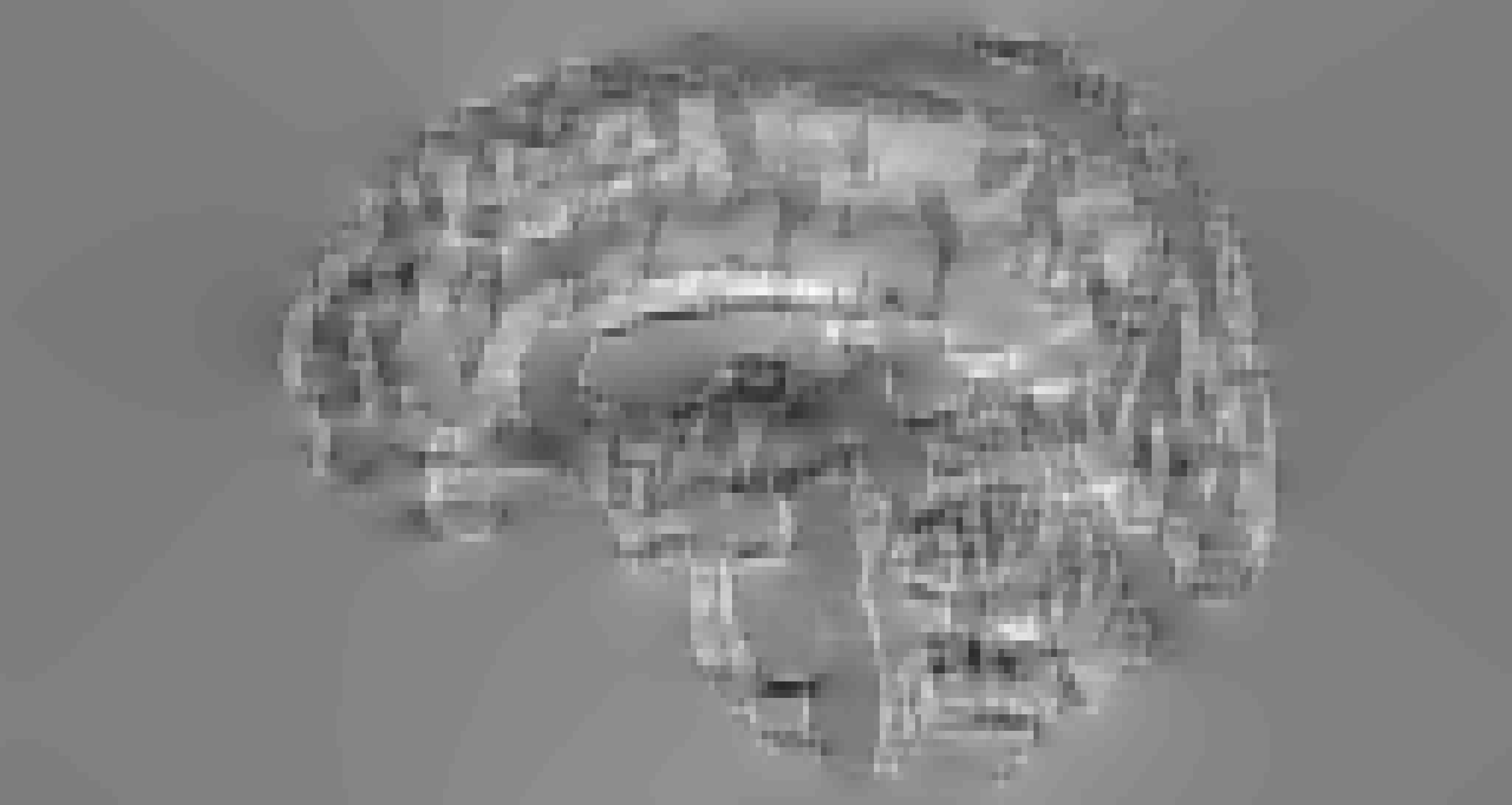}}\hspace{0.005cm}
\subfloat[Measured local field $b_l$]{\label{BrainMeasLocal}\includegraphics[width=3.60cm]{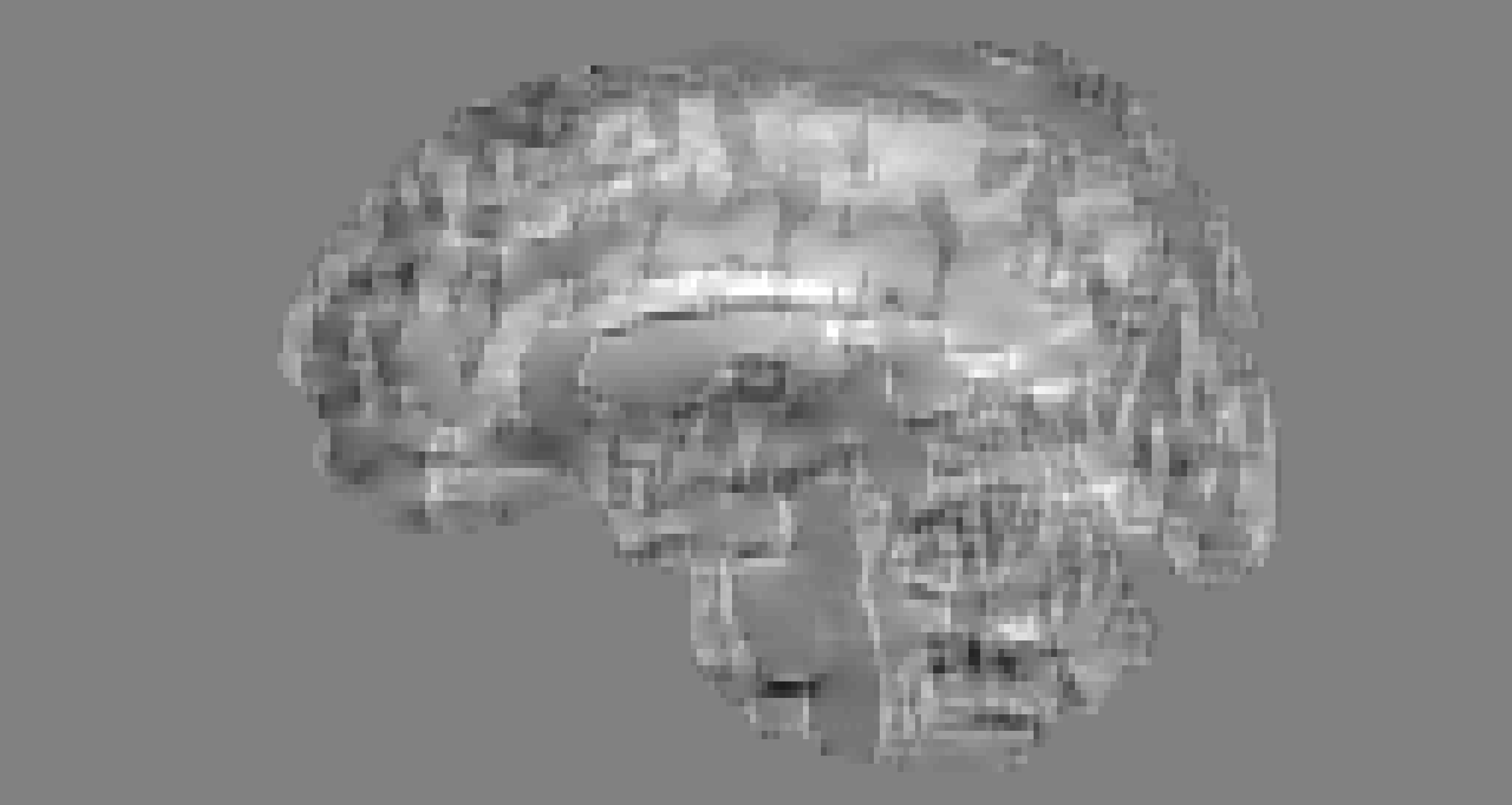}}\hspace{0.005cm}
\subfloat[$v=b_l-\wt{b}_l$]{\label{BrainIncomp}\includegraphics[width=3.60cm]{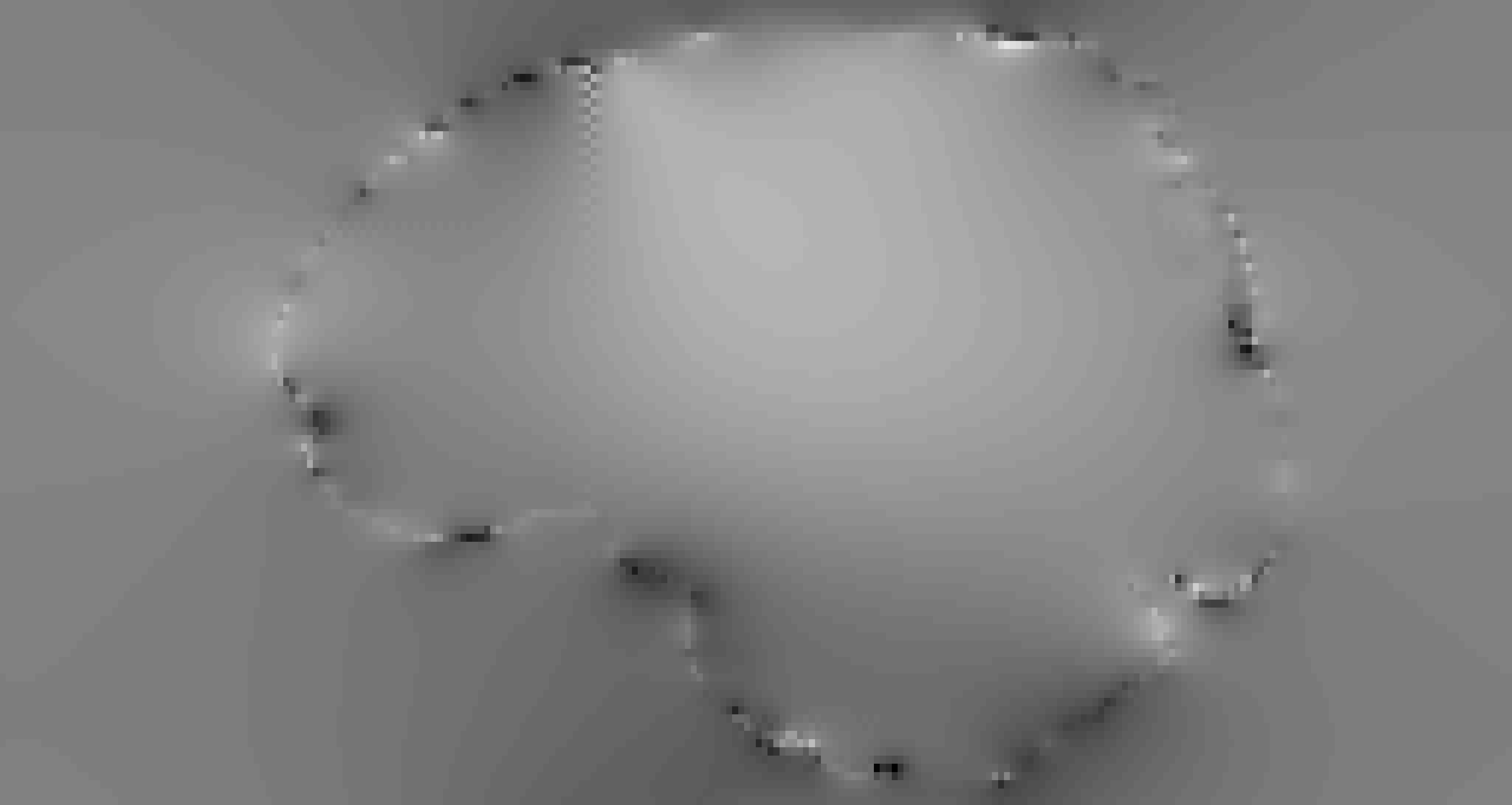}}\hspace{0.005cm}
\subfloat[$\left|-\Delta v\right|$]{\label{BrainAbsLv}\includegraphics[width=3.60cm]{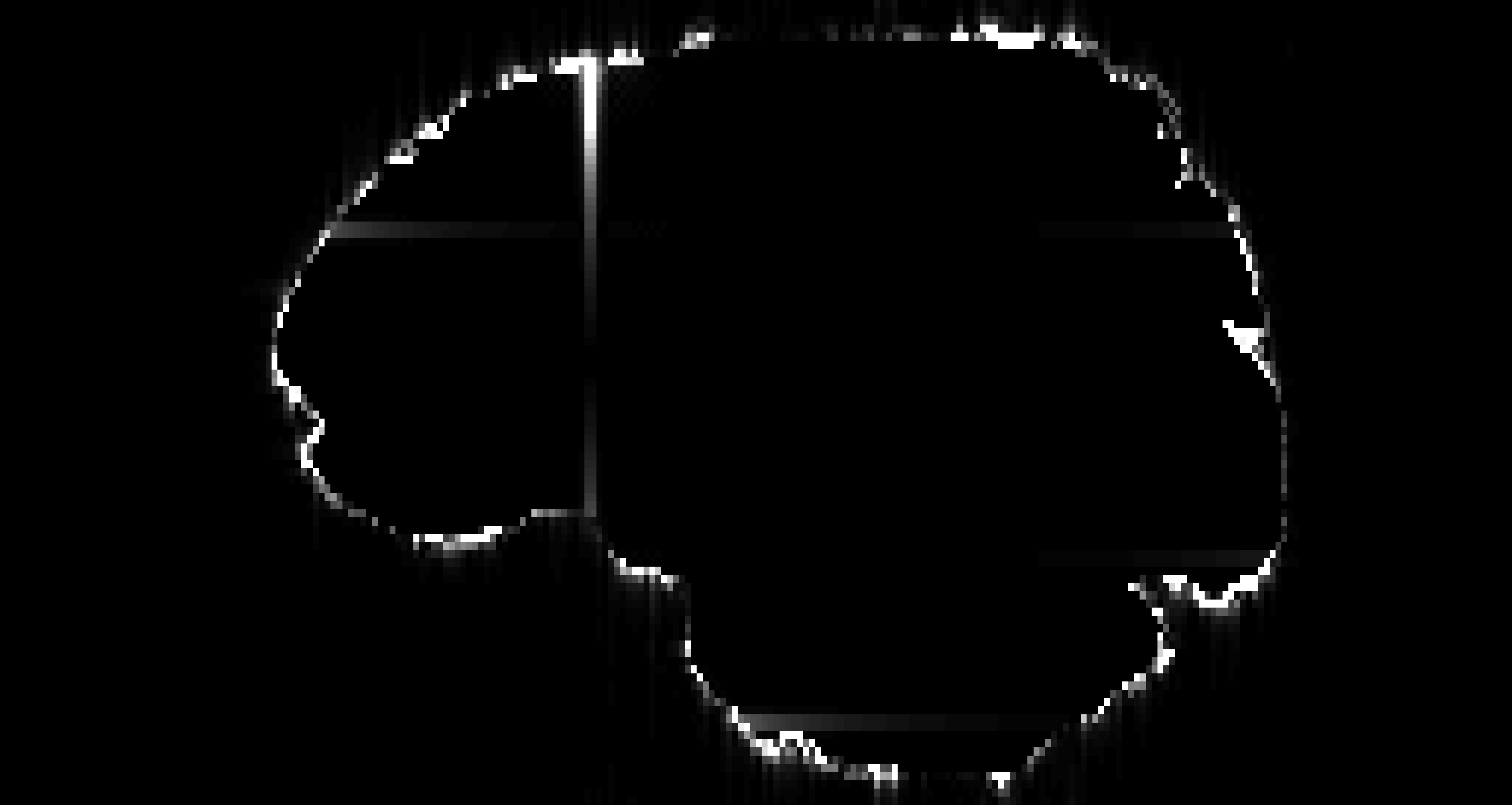}}\vspace{-0.20cm}
\caption{Sagittal slice images of $256\times256\times98$ brain phantom images with $0.9375\times0.9375\times1.5\mathrm{mm}^3$. Image of $\chi$ in $\R^3\setminus\Om$ is displayed in the window level $[0,500]$, $\wt{b}_l$ and $b_l$ in the window level $[-0.05,0.05]$, $v$ in the window level $[-0.025,0.025]$, and $\left|-\Delta v\right|$ in the window level $[0,0.01]$ respectively.}\label{IllustrateTh1Brain}
\end{figure}
\begin{figure}[tp!]
\centering
\hspace{-0.1cm}\subfloat[ROI $\Om$]{\label{BrainMaskAx}\includegraphics[width=3.60cm]{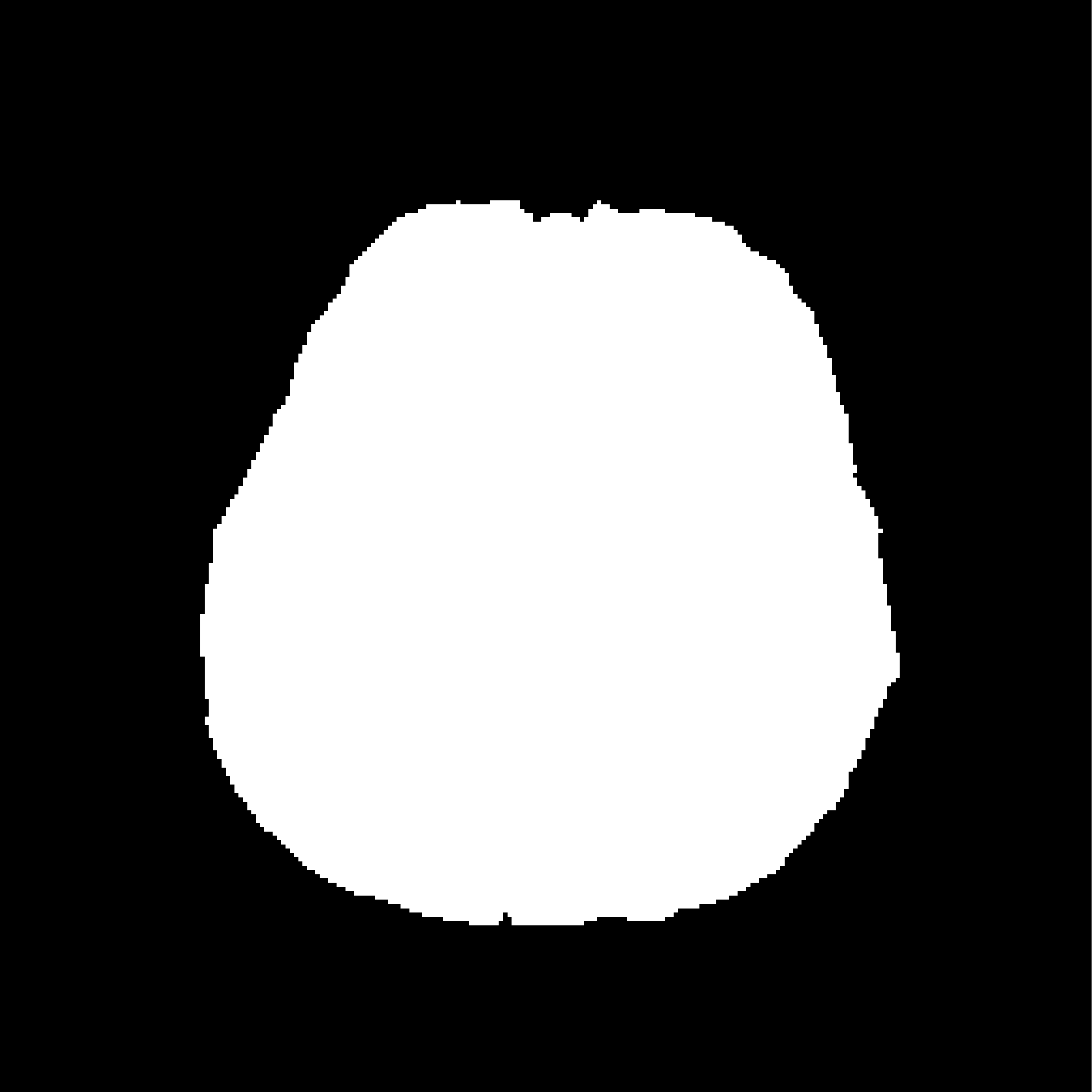}}\hspace{0.005cm}
\subfloat[$\chi$ in $\Om$]{\label{BrainQSMAx}\includegraphics[width=3.60cm]{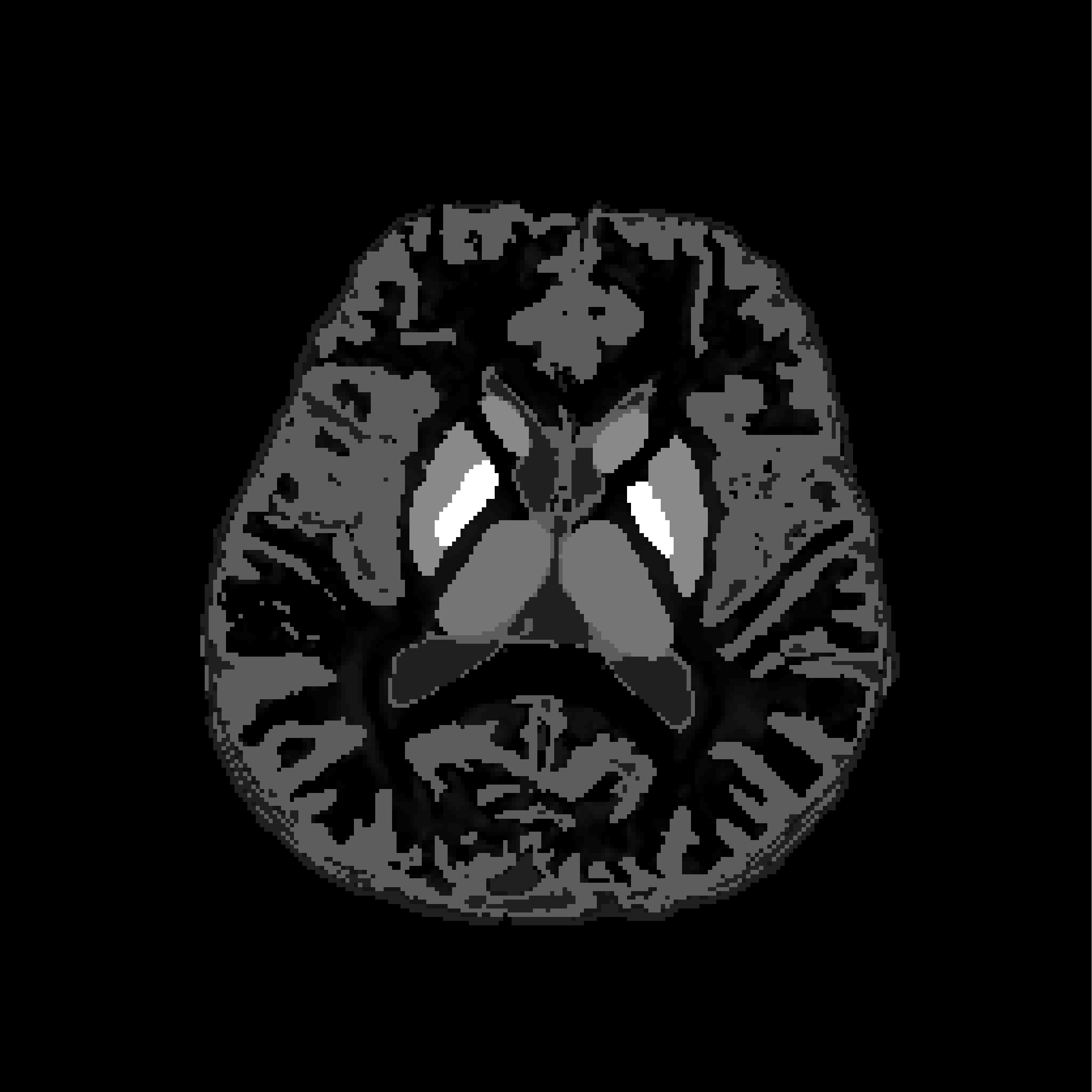}}\hspace{0.005cm}
\subfloat[$\chi$ in $\R^3\setminus\Om$]{\label{BrainSusBGAx}\includegraphics[width=3.60cm]{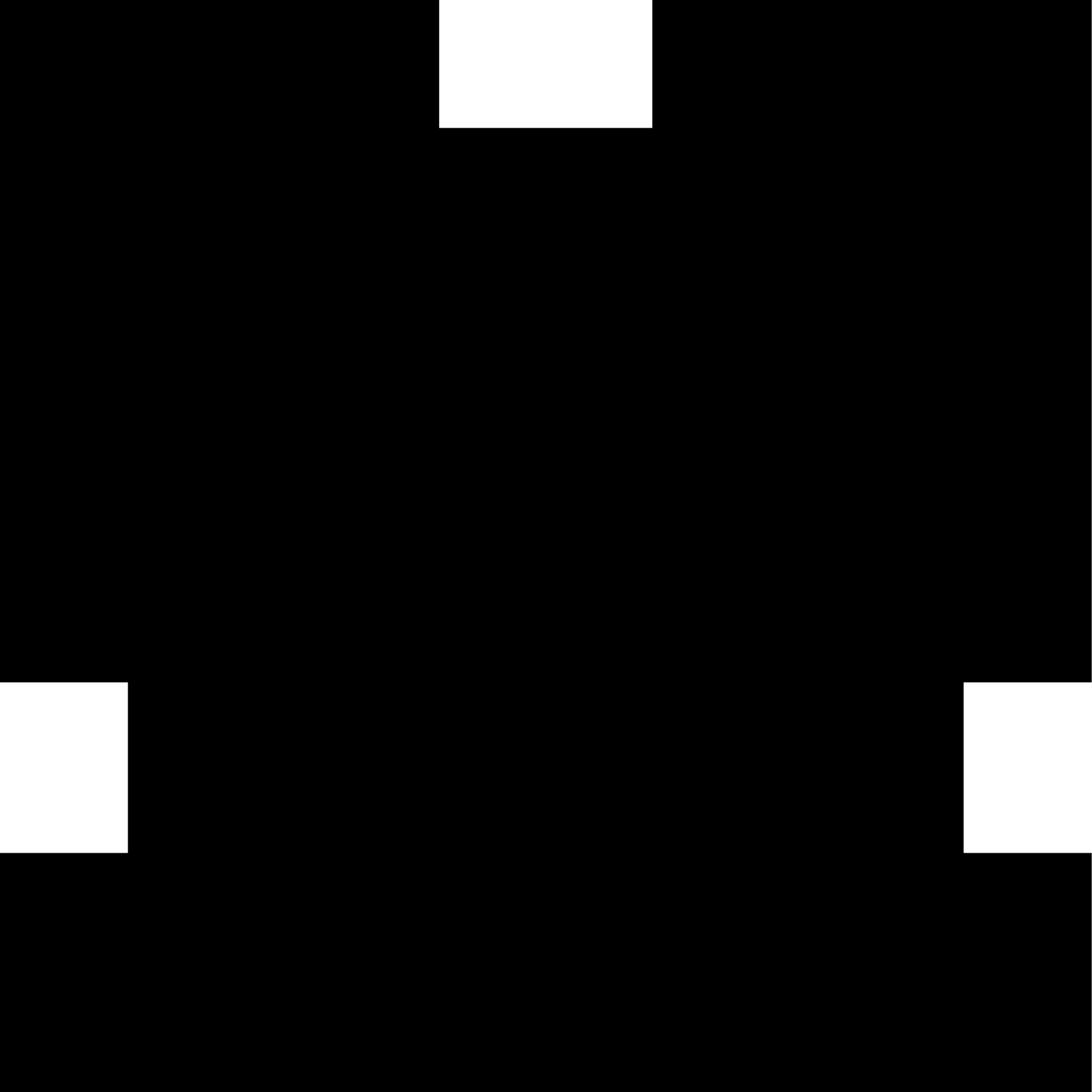}}\hspace{0.005cm}
\subfloat[Simulated total field $b$]{\label{BrainTotalAx}\includegraphics[width=3.60cm]{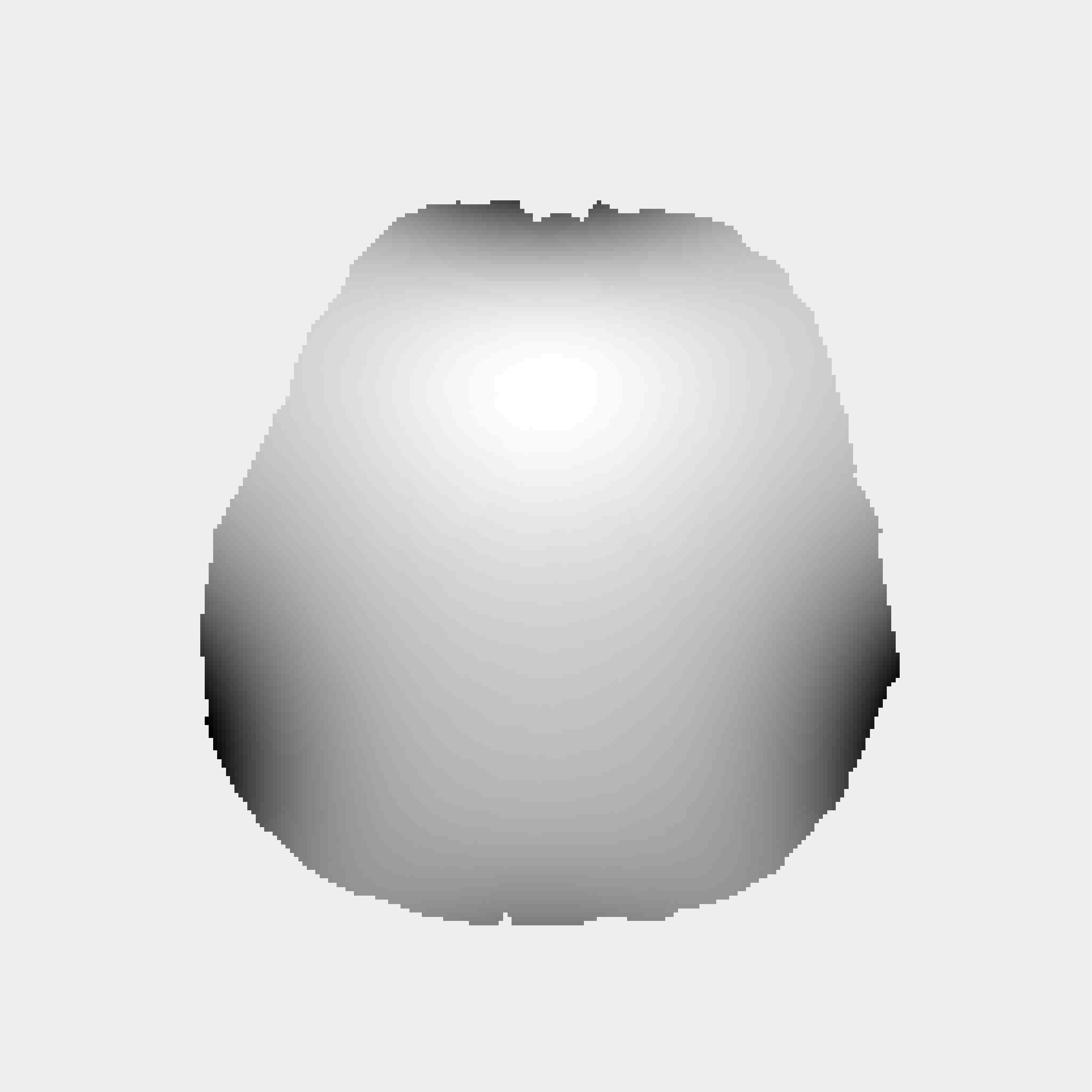}}\vspace{-0.20cm}\\
\subfloat[True local field $\wt{b}_l$]{\label{BrainTrueLocalAx}\includegraphics[width=3.60cm]{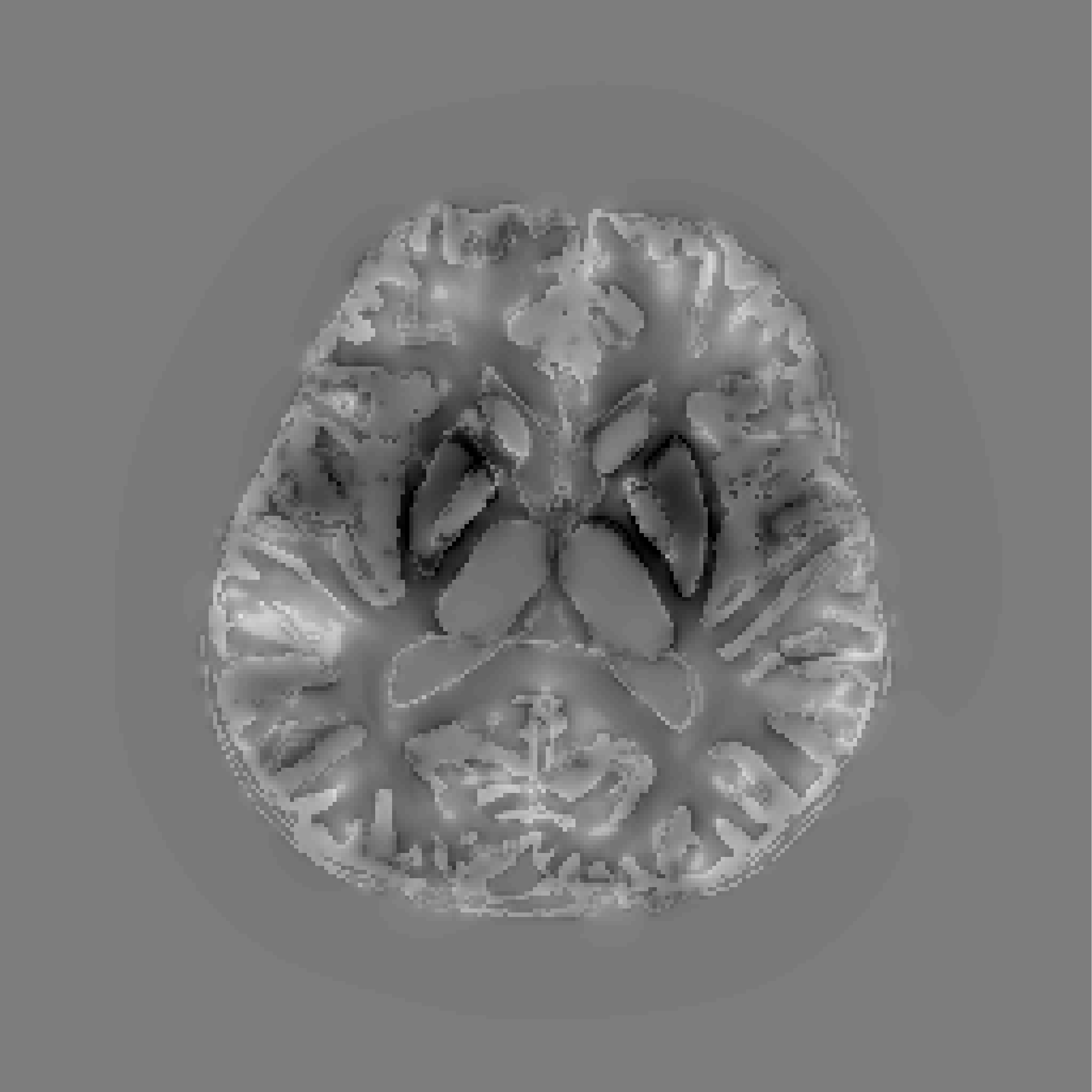}}\hspace{0.005cm}
\subfloat[Measured local field $b_l$]{\label{BrainMeasLocalAx}\includegraphics[width=3.60cm]{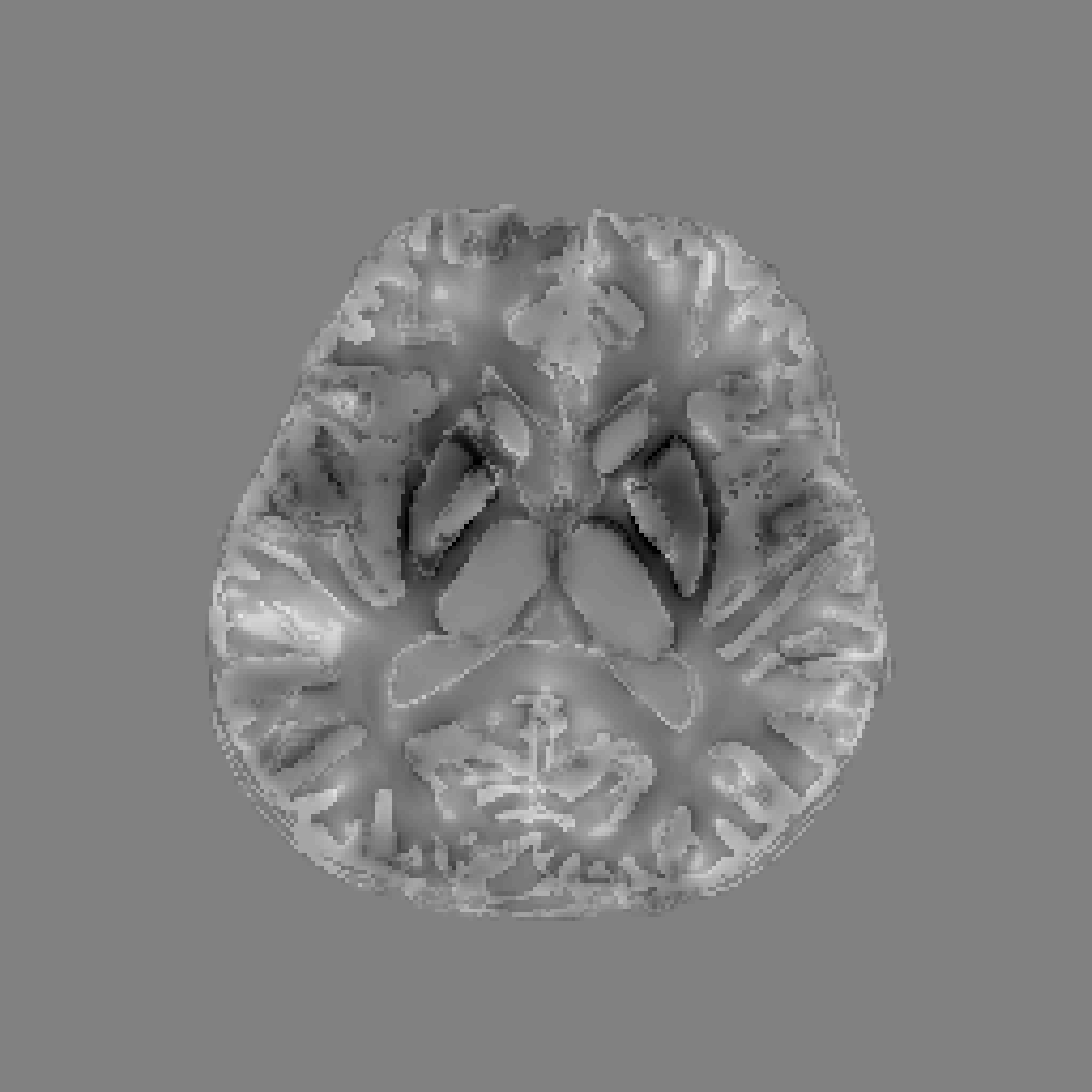}}\vspace{-0.20cm}
\subfloat[$v=b_l-\wt{b}_l$]{\label{BrainIncompAx}\includegraphics[width=3.60cm]{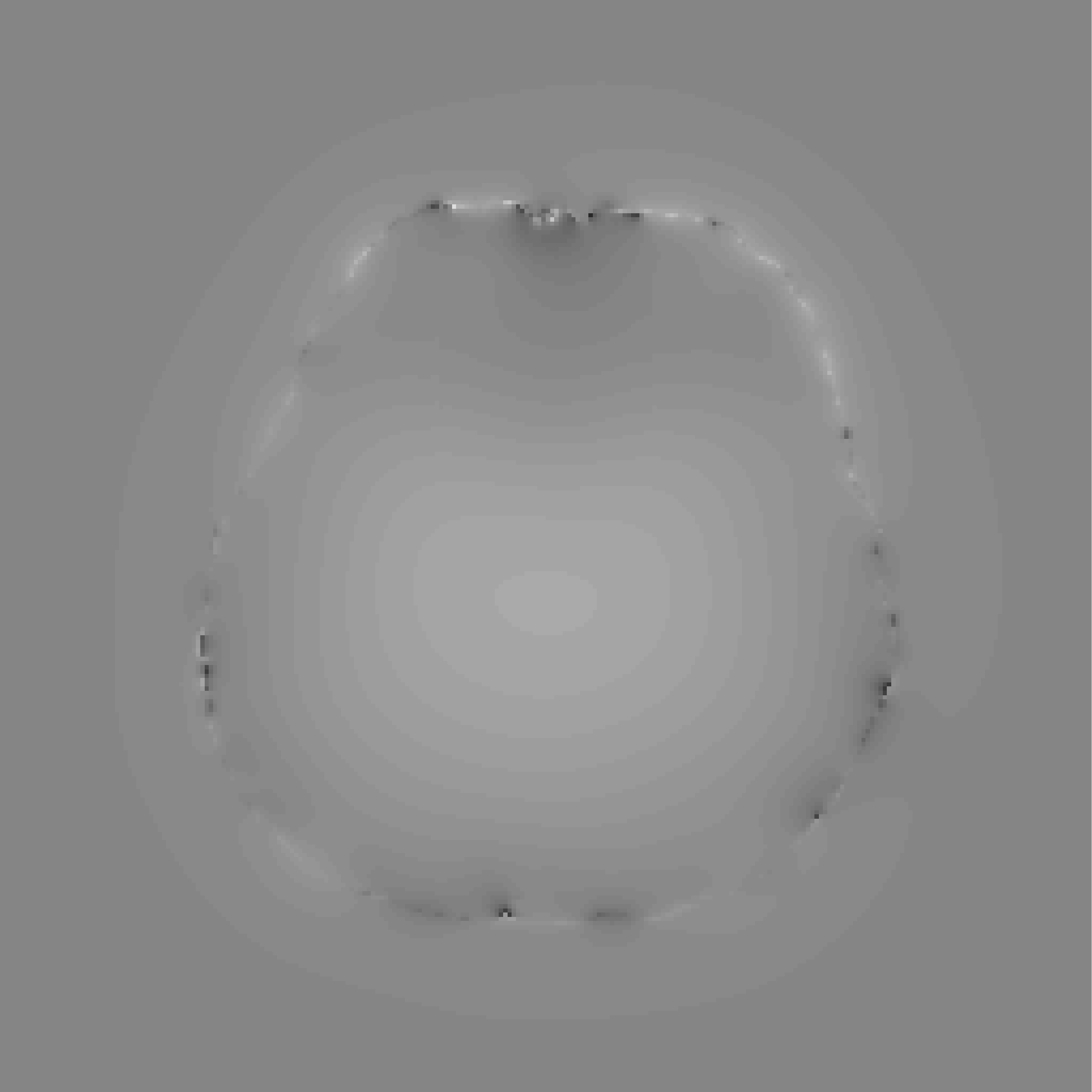}}\hspace{0.005cm}
\subfloat[$\left|-\Delta v\right|$]{\label{BrainAbsLvAx}\includegraphics[width=3.60cm]{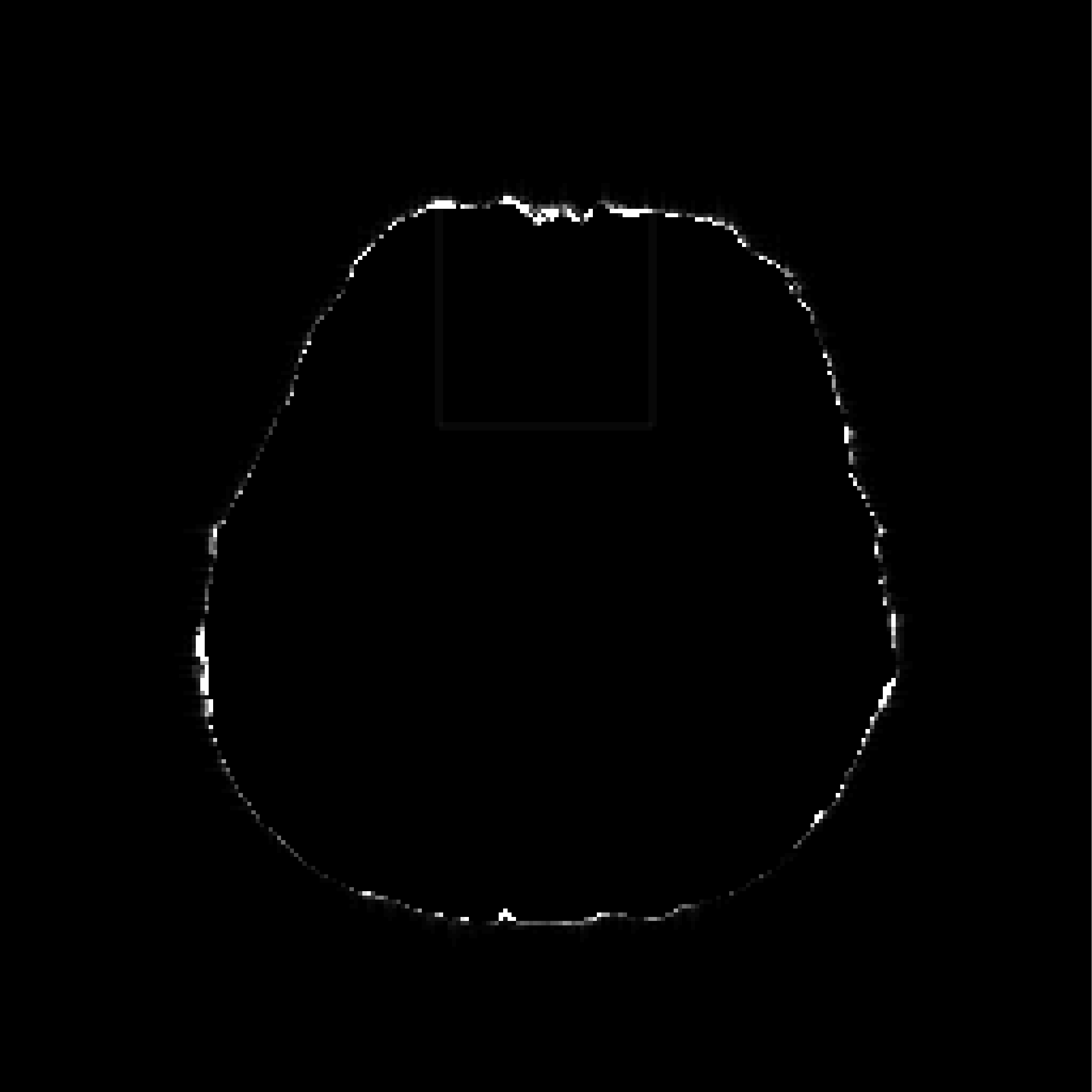}}\vspace{-0.20cm}
\caption{Axial slice images of \cref{IllustrateTh1Brain}. The images of $\chi$ in $\R^3\setminus\Om$, $\wt{b}_l$, $b_l$, $v$, and $\left|-\Delta v\right|$ are displayed in the same window level as \cref{IllustrateTh1Brain}.}\label{IllustrateTh1BrainAxial}
\end{figure}

We present \cref{IllustrateTh1,IllustrateTh1Axial,IllustrateTh1Brain,IllustrateTh1BrainAxial} to illustrate \cref{Th1} by using the Shepp-Logan phantom (\cref{IllustrateTh1,IllustrateTh1Axial}) and the brain phantom (\cref{IllustrateTh1Brain,IllustrateTh1BrainAxial}). Using the limited total field $b$ in \cref{LoganTotal,LoganTotalAx,BrainTotal,BrainTotalAx} which are derived from \cref{Forward1} by placing strong susceptibilities outside $\Om$, we solve \cref{LBV} using multigrid (MG) based the finite difference method \cite{D.Zhou2014} to obtain the measured local field $b_l$ in \cref{LoganMeasLocal,LoganMeasLocalAx,BrainMeasLocal,BrainMeasLocalAx} which are used for the susceptibility reconstruction. We also display the true local field $\wt{b}_l$ obtained from \cref{TrueLocalField} in \cref{LoganTrueLocal,LoganTrueLocalAx,BrainTrueLocal,BrainTrueLocalAx} for the comparison with the measured $b_l$. Finally, $v=b_l-\wt{b}_l$ and $|-\Delta v|$ are displayed in \cref{LoganIncomp,LoganIncompAx,BrainIncomp,BrainIncompAx} and \cref{LoganAbsLv,LoganAbsLvAx,BrainAbsLv,BrainAbsLvAx} respectively for better illustrations. Compared to the Shepp-Logan phantom, the brain phantom shows the artifacts as shown in \cref{BrainAbsLv}. There are two possible reasons of the artifacts. Firstly, since the boundary of human brain is more complicated than the Shepp-Logan phantom, the erroneous boundary values may have affected the background field removal in the case of brain phantom, as pointed out in \cite{D.Zhou2014}. Secondly, unlike the Shepp-Logan phantom with isotropic spatial resolution $(1\times1\times1\mathrm{mm}^3)$, the brain phantom has an anisotropic spatial resolution of $0.9375\times0.9375\times1.5\mathrm{mm}^3$. As pointed out in \cite{M.W.Gee2009}, the MG method has a limitation that errors in certain directions ($x_3$ direction in our case) are not smoothed by standard relaxation and as a consequence, it is inappropriate to coarsen in these directions, which may lead to artifacts in \cref{BrainAbsLv} along the $x_3$ direction. Since the real spatial resolution of phase data is not necessarily isotropic, an efficient and effective numerical solver of \cref{LBV} need to be investigated in the future, which is beyond the scope of this paper at this point.


\subsection{Proposed HIRE Susceptibility Reconstruction Model}\label{ProposedHIREModel}

We begin with introducing some notation. Let $\OO=\left\{0,\cdots,N_1-1\right\}\times\left\{0,\cdots,N_2-1\right\}\times\left\{0,\cdots,N_3-1\right\}$ denote the set of indices of $N_1\times N_2\times N_3$ grids, and let $\Om\subseteq\OO$ denote the set of indices corresponding to the ROI. Denote $\p\Om$ to be the indices of the boundary of ROI $\Om$. Finally, the space of real valued functions defined on $\OO$ is denoted as $\mI_3\simeq\R^{N_1\times N_2\times N_3}$.

Let $b_l\in\mI_3$ be the (noisy) measured local field data obtained from \cref{LBV}, which satisfies $b_l=0$ in $\OO\setminus\Om$. From the viewpoint of \cref{Th1}, we can model it as
\begin{align*}
b_l=A\chi+v+\eta
\end{align*}
where $A=\msF^{-1}\mD\msF$ denotes the discretization of the forward operator in \cref{NewModeling}. Here, $\chi\in\mI_3$ is the unknown true susceptibility image supported in $\Om$, $v\in\mI_3$ is the incompatibility arising from solving \cref{LBV}, and $\eta$ is the additive noise.



We observe that in the discrete setting, \cref{DistLapIncompatibility} in \cref{Th1} can be understood as
\begin{align}\label{DiscreteHarmonicIncompatibility}
\msL v=0\quad\text{in}\quad\OO\setminus\p\Om\quad\text{ and }\quad \msL v\neq0\quad\text{on}\quad\p\Om
\end{align}
with the discrete Laplacian $\msL$, as $\msL v$ is supported on $\p\Om$. However, it is in general difficult to directly impose \cref{DiscreteHarmonicIncompatibility} into the susceptibility reconstruction model (e.g. penalizing $\|(\msL v)_{\OO\setminus\p\Om}\|_2^2$) for the following reasons: 1) the estimation of $\Om$ always contains error due to the complicated geometry of the human brain; 2) the real MRI data may not exactly satisfy \cref{DiscreteHarmonicIncompatibility} due to its spatial resolution \cite{E.M.Haacke1999}; 3) the discretization can introduce the error on the boundary of $\Om$. However, it is a fact that the support of $\msL v$ is small compared to $|\OO|=N_1N_2N_3$, i.e. $\msL v$ is sparse. Consequently, we penalize the $\|\msL v\|_1$ for the incompatibility term $v$. Although the $\|\msL v\|_1$ does not necessarily satisfy the harmonic constraints on $\OO\setminus\p\Om$, it is a good relaxation approach when considering the error source of the forward model in QSM. In addition, motivated by the successful results on the wavelet frame based image restoration (e.g. \cite{J.F.Cai2010,J.F.Cai2012}), we assume the sparse approximation of $\chi$ under a given wavelet transformation $W$, and propose our HIRE model as follows:
\begin{align}\label{WaveletFrameHIRE}
\min_{\chi,v\in\mI_3}~\f{1}{2}\left\|A\chi+v-b_l\right\|_{\Sigma}^2+\lambda\left\|\msL v\right\|_1+\left\|\gamma\cdot W\chi\right\|_{1,2}
\end{align}
where $\left\|\cdot\right\|_{\Sigma}^2=\left\la\Sigma\cdot,\cdot\right\ra$ with the SNR weight $\Sigma$ which is estimated from the MRI \cite{T.Liu2011a,C.Milovic2018}. Here, $\left\|\gamma\cdot W\chi\right\|_{1,2}$ is the isotropic $\ell_1$ norm of the wavelet frame coefficients \cite{J.F.Cai2012} defined as
\begin{align}\label{IsoTropicWaveletFrame}
\left\|\gamma\cdot W\chi\right\|_{1,2}:=\sum_{\bk\in\OO}\sum_{l=0}^{L-1}\gamma_l[\bk]\left(\sum_{\aal\in\BB}\left|\left(W_{l,\aal}\chi\right)[\bk]\right|^2\right)^{1/2}.
\end{align}
(See \cref{WaveletFrame} for the brief introduction on the wavelet frames.)

There are many variational regularizations for the susceptibility image including total variation (TV) \cite{L.I.Rudin1992}, total generalized variation (TGV) \cite{K.Bredies2010,F.Yanez2013}, and weighted TV for morphological consistency \cite{Y.Kee2018}. However, since $\mD(\0)=0$, the $\chi$ subproblem is a rank deficient system matrix, when using the alternating direction method of multipliers (ADMM) methods \cite{J.Eckstein1992} to solve the regularization model. As a consequence, we may need additional prior information such as the zero susceptibility value in the cerebrospinal fluid region \cite{Z.Liu2017} for the stable reconstruction. In contrast, by using the tight frame regularization, the system matrix of $\chi$ subproblem has a full column rank, which can lead to the computational efficiency over the existing variational methods.

\begin{remark}\label{RK231} For better understanding of our HIRE model, we temporarily assume that $\Sigma=I$ and consider
\begin{align*}
\min_{\chi,v\in\mI_3}~\f{1}{2}\left\|A\chi+v-b_l\right\|_2^2+\lambda\left\|\msL v\right\|_1+\left\|\gamma\cdot W\chi\right\|_{1,2}.
\end{align*}
	If $v\equiv0$, our model reduces to the integral approach model:
	\begin{align*}
	\min_{\chi\in\mI_3}~\f{1}{2}\left\|A\chi-b_l\right\|_2^2+\left\|\gamma\cdot W\chi\right\|_{1,2}.
	\end{align*}
	In addition, if we fix $v=b_l-A\chi$, our model reduces to the $\ell_1$ fidelity version of the following differential approach model:
	\begin{align*}
	\min_{\chi\in\mI_3}~\f{1}{2}\left\|\msL b_l-\msL A\chi\right\|_2^2+\left\|\gamma\cdot W\chi\right\|_{1,2}
	\end{align*}
	as $\msL A\chi=\msL b_l$ discretizes the PDE \cref{QSM_PDE2} in the sense of \cite[Proposition A.1.]{J.K.Choi2014}.	
\end{remark}


From \cref{RK231}, we can see that our model considers the incompatibility $v$ and noise separately, thereby providing a more precise forward model for QSM. This is because $b_l$ is obtained from the Poisson's equation \cref{LBV} and it inevitably contains the harmonic incompatibility related to the imposed boundary condition, as described in \cref{Th1}. Even though more rigorous theoretical analysis is needed, we can somehow explain the effect of harmonic incompatibility in this manner; since the standard arguments on the harmonic functions (e.g. \cite{Evans2010}) tell us that $v$ is smooth and satisfies the mean value property except on $\p\Om$, it has slow variations on this region. As a consequence, it mostly affects the low frequency components in $b_l$ compared to the noise which mainly affects the high frequency components. Together with the fact that the critical manifold $\Gamma_0$ forms a conic manifold in the frequency domain, the harmonic incompatibility $v$ in $b_l$ mainly leads to the loss of $\msF(\chi)$ in low frequency components.

As empirically observed in \cite{Y.Kee2017}, the incompatibility in low frequency components of $b_l$ leads to the shadow artifacts in the reconstructed image, while that in high frequency components leads to the streaking artifacts. Therefore, the simultaneous consideration on the incompatibilities in both components is crucial for better susceptibility imaging. The integral approach does not take the harmonic incompatibility in $b_l$ into account, which may not be capable of suppressing the incompatibility in low frequency components of $b_l$, and leads to the shadow artifacts in the reconstructed images. The differential approach can be viewed as a preconditioned integral approach since the harmonic incompatibility in $b_l$ has been removed in advance. However, the noise in $b_l$ can be amplified by $\msL$ at the cost of harmonic incompatibility removal, and this leads to the streaking artifacts propagating from the noise in final image \cite{Y.Wang2015}. In contrast, the HIRE model takes the form of integral approach which explicitly considers the incompatibility $v$ other than the noise by incorporating its sparsity under $\msL$. By doing so, we expect that the HIRE model can suppress both the noise (cause of streaking artifacts) and the harmonic incompatibility (cause of shadow artifacts), so that we can achieve the whole brain imaging with less artifacts.

We would like to mention that the formulation of HIRE model is not limited to \cref{WaveletFrameHIRE}. In fact, we can use the nonlinear fidelity term $F(b_l|\chi,v)=\f{1}{2}\left\|e^{i(A\chi+v)\om_0B_0TE}-e^{ib_l\om_0B_0TE}\right\|_{\Sigma}^2$ to further compensate the errors in phase unwrapping, which will be more coincident with the GRE signal model \cite{Y.Kee2017,T.Liu2013}. However, we will not discuss the details on such nonlinear variants as this is beyond the scope of this paper. We will focus on \cref{WaveletFrameHIRE} throughout this paper.


\subsection{Numerical Algorithm}\label{AlternatingMinimizationAlgorithm}

In the literature, there are numerous algorithms which can solve the proposed HIRE model \cref{WaveletFrameHIRE}. In this paper, we adopt the split Bregman algorithm given in \cite{C.Milovic2018} in the framework of ADMM \cite{J.Eckstein1992} as we can convert \cref{WaveletFrameHIRE} into several subproblems which can be solved efficiently. More precisely, let $d=W\chi$, $e=\msL v$, $f=A\chi$, and $g=v$. Then \cref{WaveletFrameHIRE} is reformulated as follows:
\begin{align*}
&\min_{\chi,v,d,e,f,g}~\f{1}{2}\left\|f+g-b_l\right\|_{\Sigma}^2+\lambda\left\|e\right\|_1+\left\|\gamma\cdot d\right\|_{1,2}\\
&\text{subject to}~d=W\chi,~e=\msL v,~f=A\chi,~\text{and}~g=v.
\end{align*}
Under this reformulation, we summarize the split Bregman algorithm for \cref{WaveletFrameHIRE} in \cref{Algorithm1}.

\begin{algorithm}[tp!]
\caption{Split Bregman Algorithm for \eqref{WaveletFrameHIRE}}\label{Algorithm1}
\begin{algorithmic}
\STATE{\textbf{Initialization:} $\chi^0$, $v^0$, $d^0$, $e^0$, $f^0$, $g^0$, $p^0$, $q^0$, $r^0$, $s^0$}
\FOR{$k=0$, $1$, $2$, $\cdots$}
\STATE{Update $\chi$ and $v$:
\begin{align}
\chi^{k+1}&=\argmin_{\chi}~\f{\beta}{2}\|A\chi-f^k+r^k\|_2^2+\f{\beta}{2}\|W\chi-d^k+p^k\|_2^2 \label{chisubproblem}\\
v^{k+1}&=\argmin_{v}~\f{\beta}{2}\|v-g^k+s^k\|_2^2+\f{\beta}{2}\|\msL v-e^k+q^k\|_2^2 \label{vsubproblem}
\end{align}
Update $d$, $e$, $f$, and $g$:
\begin{align}
d^{k+1}&=\argmin_d~\left\|\gamma\cdot d\right\|_{1,2}+\f{\beta}{2}\|d-W\chi^{k+1}-p^k\|_2^2 \label{dsubproblem}\\
e^{k+1}&=\argmin_e~\lambda\left\|e\right\|_1+\f{\beta}{2}\|e-\msL v^{k+1}-q^k\|_2^2 \label{esubproblem}\\
f^{k+1}&=\argmin_f~\f{1}{2}\|f+g^k-b_l\|_{\Sigma}^2+\f{\beta}{2}\|f-A\chi^{k+1}-r^k\|_2^2 \label{fsubproblem}\\
g^{k+1}&=\argmin_g~\f{1}{2}\|g+f^{k+1}-b_l\|_{\Sigma}^2+\f{\beta}{2}\|g-v^{k+1}-s^k\|_2^2 \label{gsubproblem}
\end{align}
Update $p$, $q$, $r$, and $s$:
\begin{align}
p^{k+1}&=p^k+W\chi^{k+1}-d^{k+1} \label{psubproblem}\\
q^{k+1}&=q^k+\msL v^{k+1}-e^{k+1} \label{qsubproblem}\\
r^{k+1}&=r^k+A\chi^{k+1}-f^{k+1} \label{rsubproblem}\\
s^{k+1}&=s^k+v^{k+1}-g^{k+1} \label{ssubproblem}
\end{align}}
\ENDFOR
\end{algorithmic}
\end{algorithm}

It is easy to see that each subproblem has a closed form solution. The solutions to \cref{chisubproblem,vsubproblem} can be written as
\begin{align}
\chi^{k+1}&=\big(A^TA+I\big)^{-1}\big[A^T(f^k-r^k)+W^T(d^k-p^k)\big] \label{chisubexplicit}\\
v^{k+1}&=\big(I+\msL^T\msL\big)^{-1}\big[g^k-s^k+\msL^T(e^k-q^k)\big]. \label{vsubexplicit}
\end{align}
Since we use the periodic boundary conditions, both \cref{chisubexplicit,vsubexplicit} can be easily solved by using the fast Fourier transform. In addition, the solutions to \cref{dsubproblem,esubproblem} are expressed in terms of the soft thresholding:
\begin{align}
d^{k+1}&=\mT_{\gamma/\beta}\big(W\chi^{k+1}+p^k\big) \label{dsubexplicit}\\
e^{k+1}&=\max\big(|\msL v^{k+1}+q^k|-\lambda/\beta,0\big)\sign\big(\msL v^{k+1}+q^k\big). \label{esubexplicit}
\end{align}
Here, $\mT_{\gamma}$ is the isotropic soft thresholding in \cite{J.F.Cai2012}: given $d$ defined as
\begin{align*}
d=\left\{d_{l,\aal}:(l,\aal)\in\left(\left\{0,\cdots,L-1\right\}\times\BB\right)\cup\left\{(L-1,\0)\right\}\right\}
\end{align*}
and $\gamma=\left\{\gamma_l:l=0,1,\cdots,L-1\right\}$ with $\gamma_l\geq0$, $\mT_{\gamma}\left(d\right)$ is defined as
\begin{align*}
\left(\mT_{\gamma}\left(d\right)\right)_{l,\aal}[\bk]=\left\{\begin{array}{ll}
d_{l,\aal}[\bk],&(l,\aal)=(L-1,\0)\vspace{0.125em}\\
\displaystyle{\max\left(R_l[\bk]-\gamma_l[\bk],0\right)\f{d_{l,\aal}[\bk]}{R_l[\bk]}},&(l,\aal)\in\left\{0,\cdots,L-1\right\}\times\BB
\end{array}\right.
\end{align*}
where $R_l[\bk]=\left(\sum_{\aal\in\BB}\left|d_{l,\aal}[\bk]\right|^2\right)^{1/2}$ for $\bk\in\OO$. Finally, the solutions to \cref{fsubproblem,gsubproblem} are expressed as
\begin{align}
f^{k+1}&=\big(\Sigma+\beta I\big)^{-1}\big[\Sigma(b_l-g^k)+\beta(A\chi^{k+1}+r^k)\big] \label{fsubexplicit}\\
g^{k+1}&=\big(\Sigma+\beta I\big)^{-1}\big[\Sigma(b_l-f^{k+1})+\beta(v^{k+1}+s^k)\big] \label{gsubexplicit}
\end{align}
where $\Sigma+\beta I$ is simply a diagonal matrix and thus, no matrix inversion is needed.

Note that, since our model \cref{WaveletFrameHIRE} is convex, it can be verified that \cref{Algorithm1} converges to the minimizer of \cref{WaveletFrameHIRE} by following the framework of \cite[Theorem 3.2.]{J.F.Cai2009/10}, whenever it has the unique global minimizer.

\section{Experimental Results}\label{Experiments}

In this section, we present some experimental results on brain phantom in \cite{C.Wisnieff2013} and in vivo MR data in \cite{Y.Wang2015}, both of which are available on Cornell MRI Research Lab webpage\footnote{\href{}{http://www.weill.cornell.edu/mri/pages/qsm.html}}, to compare the wavelet frame HIRE model \cref{WaveletFrameHIRE} (Frame-HIRE) with several existing approaches. In this paper, we choose to compare with the TKD method \cref{TKD} in \cite{K.Shmueli2009}, the Tikhonov regularization \cref{Tikhonov} in \cite{B.Kressler2010}, the wavelet frame integral approach (Frame-Int)
\begin{align}\label{WaveletFrameIntegral}
\min_{\chi\in\mI_3}~\f{1}{2}\left\|A\chi-b_l\right\|_{\Sigma}^2+\left\|\gamma\cdot W\chi\right\|_{1,2}
\end{align}
and the wavelet frame differential approach (Frame-Diff)
\begin{align}\label{WaveletFrameDifferential}
\min_{\chi\in\mI_3}~\f{1}{2}\left\|\msL A\chi-\msL b_l\right\|_{\Sigma}^2+\left\|\gamma\cdot W\chi\right\|_{1,2}
\end{align}
where the SNR weight for \cref{WaveletFrameDifferential} is estimated by the method described in \cite{Y.Kee2017}. Moreover, in order to highlight the main focus of this paper-to propose a two system regularization model by identifying a harmonic incompatibility in the measured local field data, we also test the models by replacing $\left\|\gamma\cdot W\chi\right\|_{1,2}$ in \cref{WaveletFrameHIRE,WaveletFrameIntegral,WaveletFrameDifferential} into the following TGV term:
\begin{align}\label{TGVTerm}
\mathrm{TGV}_{\alpha_1,\alpha_0}^2(\chi)=\alpha_1\left\|\na\chi-p\right\|_1+\alpha_0\left\|\mE p\right\|_1,~~~\text{where}~~\mE=\f{1}{2}\left(\na+\na^T\right),
\end{align}
which will be denoted as TGV-HIRE, TGV-Int and TGV-Diff respectively. All experiments are implemented on MATLAB $\mathrm{R}2015\mathrm{a}$ running on a platform with $16\mathrm{GB}$ RAM and Intel(R) Xeon(R) CPU $\mathrm{E}5$-$2609$ $0$ at $2.40\mathrm{GHz}$ with $4$ cores.

In \cref{WaveletFrameHIRE,WaveletFrameIntegral,WaveletFrameDifferential}, we choose $W$ to be the tensor product Haar framelet transform with $1$ level of decomposition to avoid the memory storage problem. Note, however, that the decomposition level and the choice of $W$ will do affect the restoration results. In addition, we use the standard difference for the TGV term, and the standard centered difference for $\msL$ in the HIRE approaches. The stopping criterion for \cref{Algorithm1} is
\begin{align*}
\f{\|\chi^{k+1}-\chi^k\|_2}{\|\chi^{k+1}\|_2}\leq5\times10^{-3},
\end{align*}
and \cref{WaveletFrameIntegral,WaveletFrameDifferential} as well as the TGV models are solved using the split Bregman algorithm presented in \cite{J.F.Cai2009/10,W.Guo2014} with the same stopping criterion as above. For the parameters, we choose $\gamma$ in \cref{IsoTropicWaveletFrame} as $\gamma=\big\{\nu2^{-l}:l=0,\cdots,L-1\big\}$ with $\nu>0$ according to \cite{J.F.Cai2012}. Empirically, we observe that $\alpha_0=2\alpha_1$ for \cref{TGVTerm}, $\lambda=5\nu$ for the Frame-HIRE, and $\lambda=8\alpha_1$ for the TGV-HIRE are good choices. Parameters $\nu$ and $\alpha_1$ vary case by case, and are chosen manually to promote an optimal balance between indices and visual qualities; even though the parameters have few effects on the indices, the reconstructed images contain more artifacts as the parameters become smaller. Finally, we compute the root mean square error (RMSE), the structural similarity index map (SSIM) \cite{Z.Wang2004}, and the computation time of the brain phantom experimental results for the quantitative comparison of each reconstruction model.

\subsection{Experiments on Brain Phantom}\label{BrainPhantom}

\begin{figure}[tp!]
\centering
\hspace{-0.1cm}\subfloat[True $\chi$]{\label{PhantomQSM}\includegraphics[width=3.60cm]{PhantomQSM.pdf}}\hspace{0.005cm}
\subfloat[Magnitude]{\label{PhantomMag}\includegraphics[width=3.60cm]{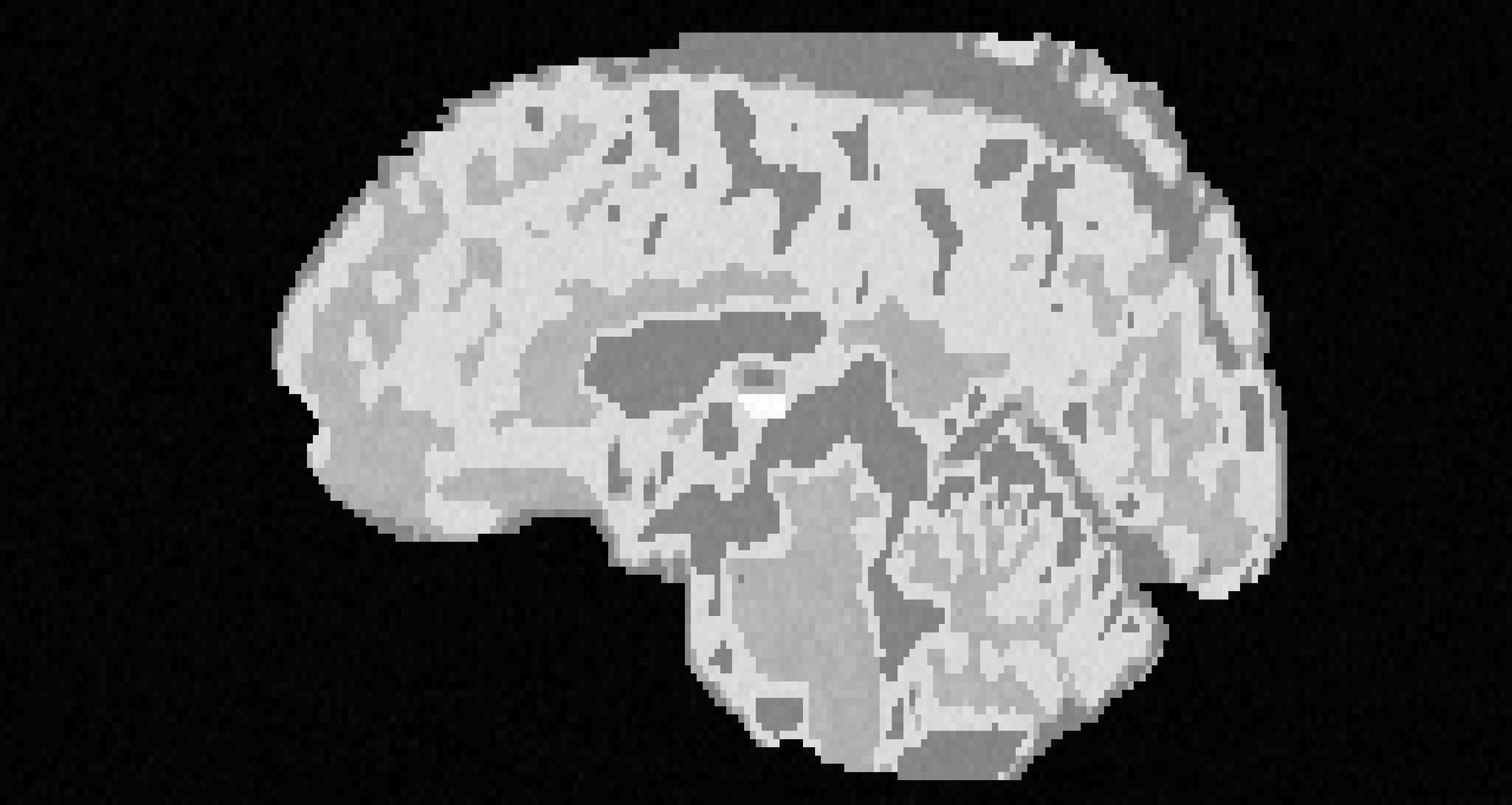}}\hspace{0.005cm}
\subfloat[ROI]{\label{PhantomMask}\includegraphics[width=3.60cm]{PhantomMask.pdf}}\vspace{-0.20cm}\\
\subfloat[Phase]{\label{PhantomPhase}\includegraphics[width=3.60cm]{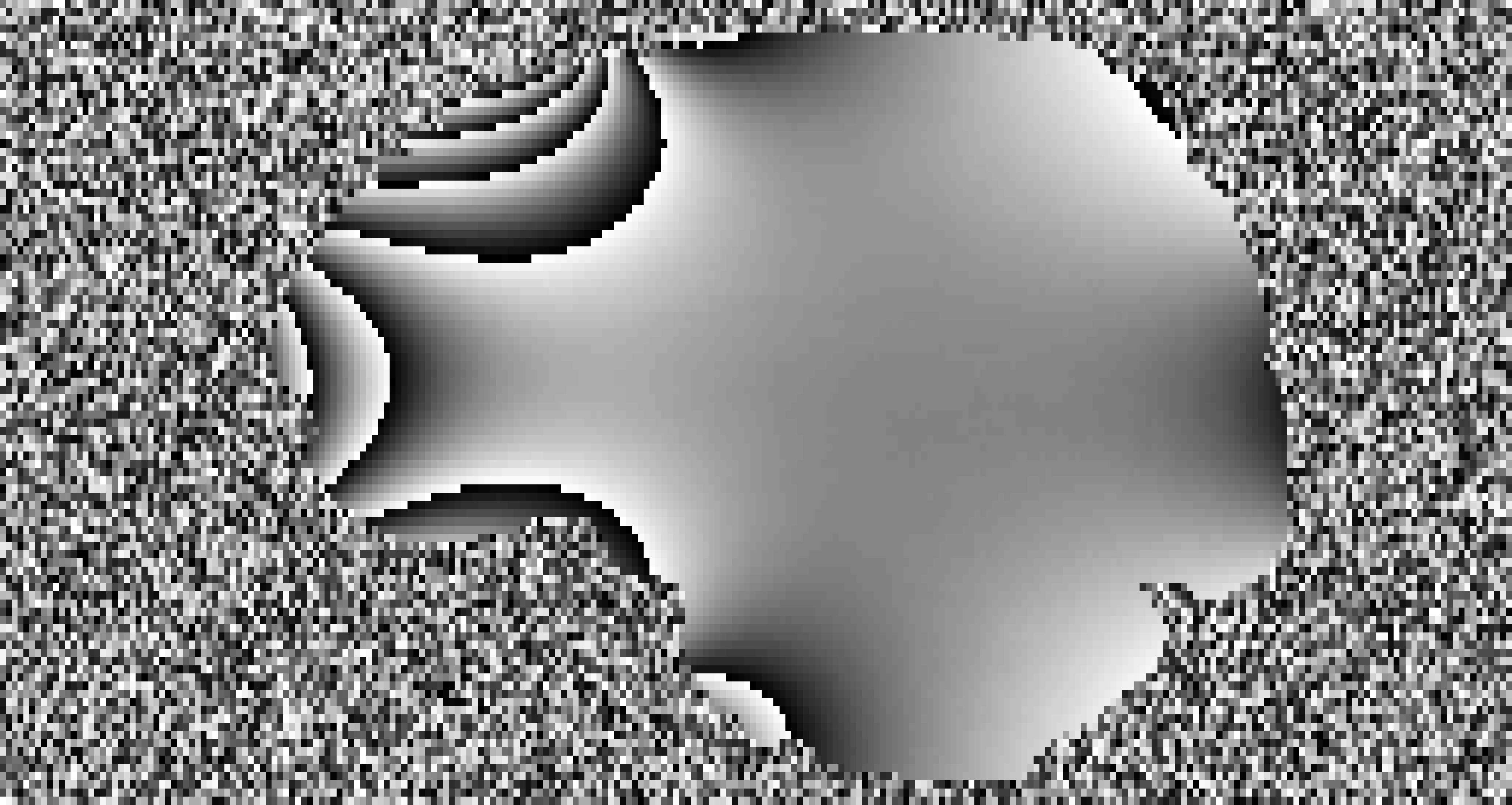}}\hspace{0.005cm}
\subfloat[Total field]{\label{PhantomTotalField}\includegraphics[width=3.60cm]{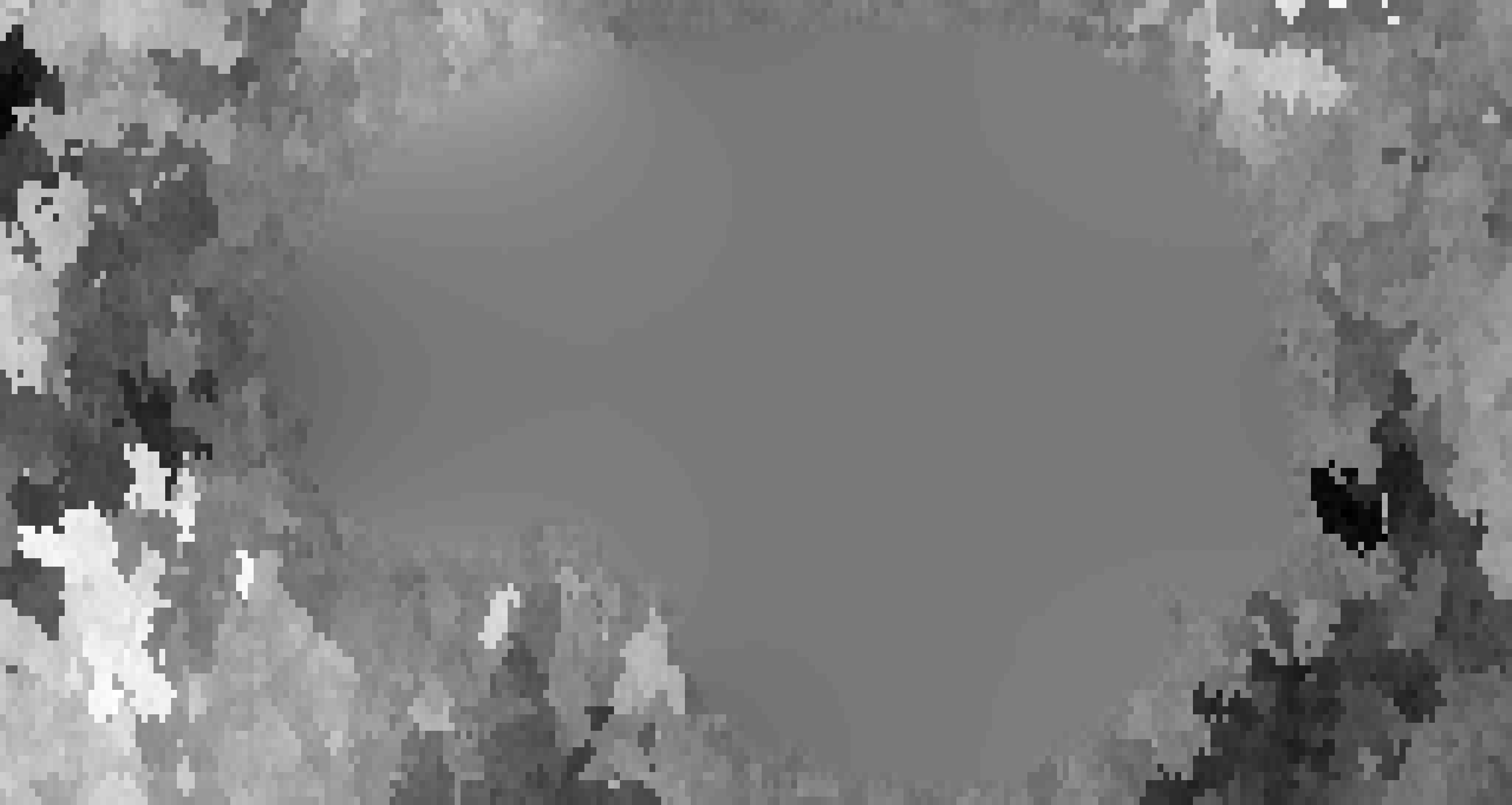}}\hspace{0.005cm}
\subfloat[Local field]{\label{PhantomLocalField}\includegraphics[width=3.60cm]{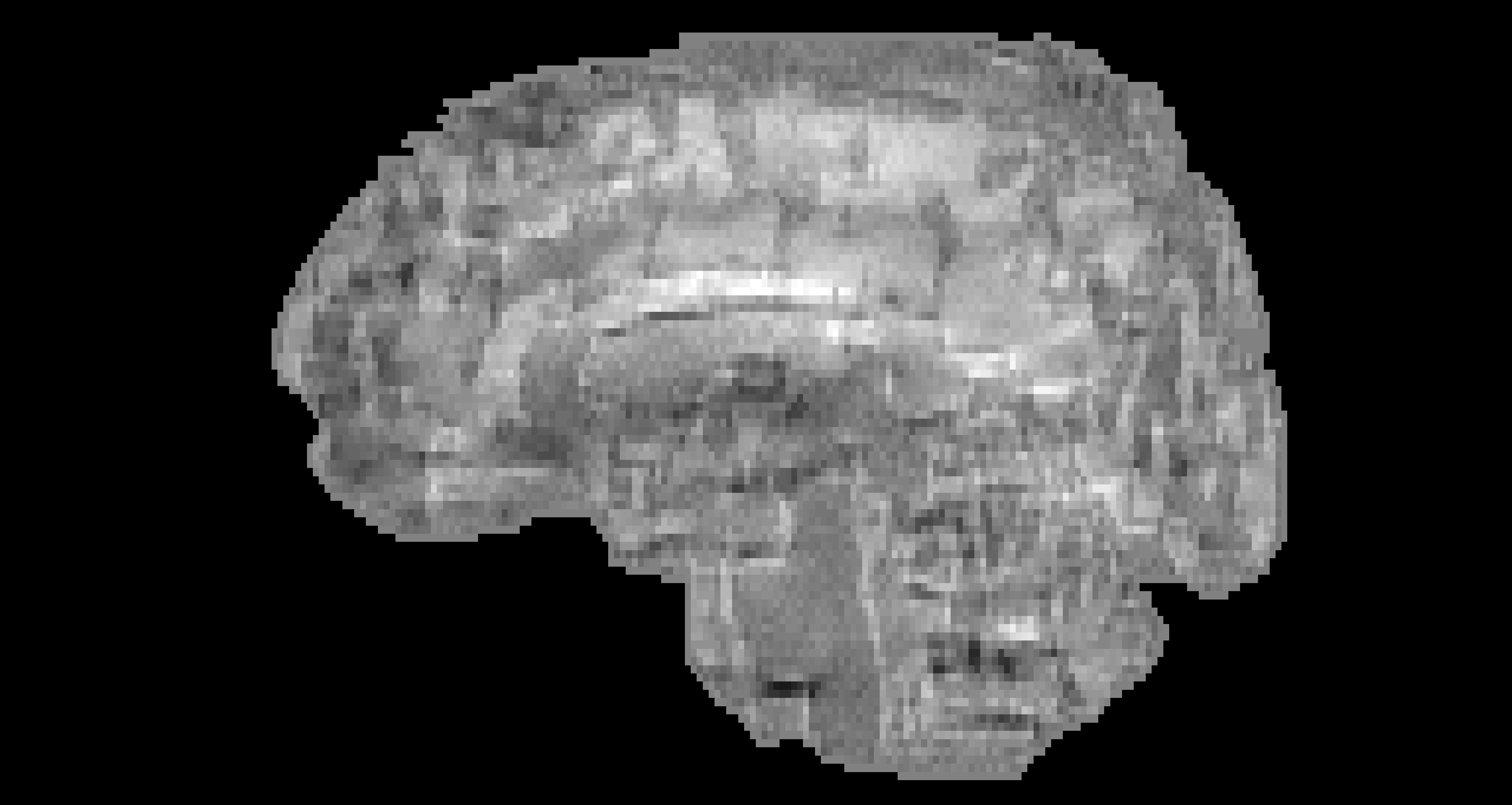}}\vspace{-0.20cm}
\caption{Sagittal slice images of synthesized data sets for the brain phantom experiments.}\label{BrainPhantomDataSet}
\end{figure}

\begin{figure}[tp!]
\centering
\hspace{-0.1cm}\subfloat[True $\chi$]{\label{PhantomQSMAx}\includegraphics[width=3.60cm]{PhantomQSMAxial.pdf}}\hspace{0.005cm}
\subfloat[Magnitude]{\label{PhantomMagAx}\includegraphics[width=3.60cm]{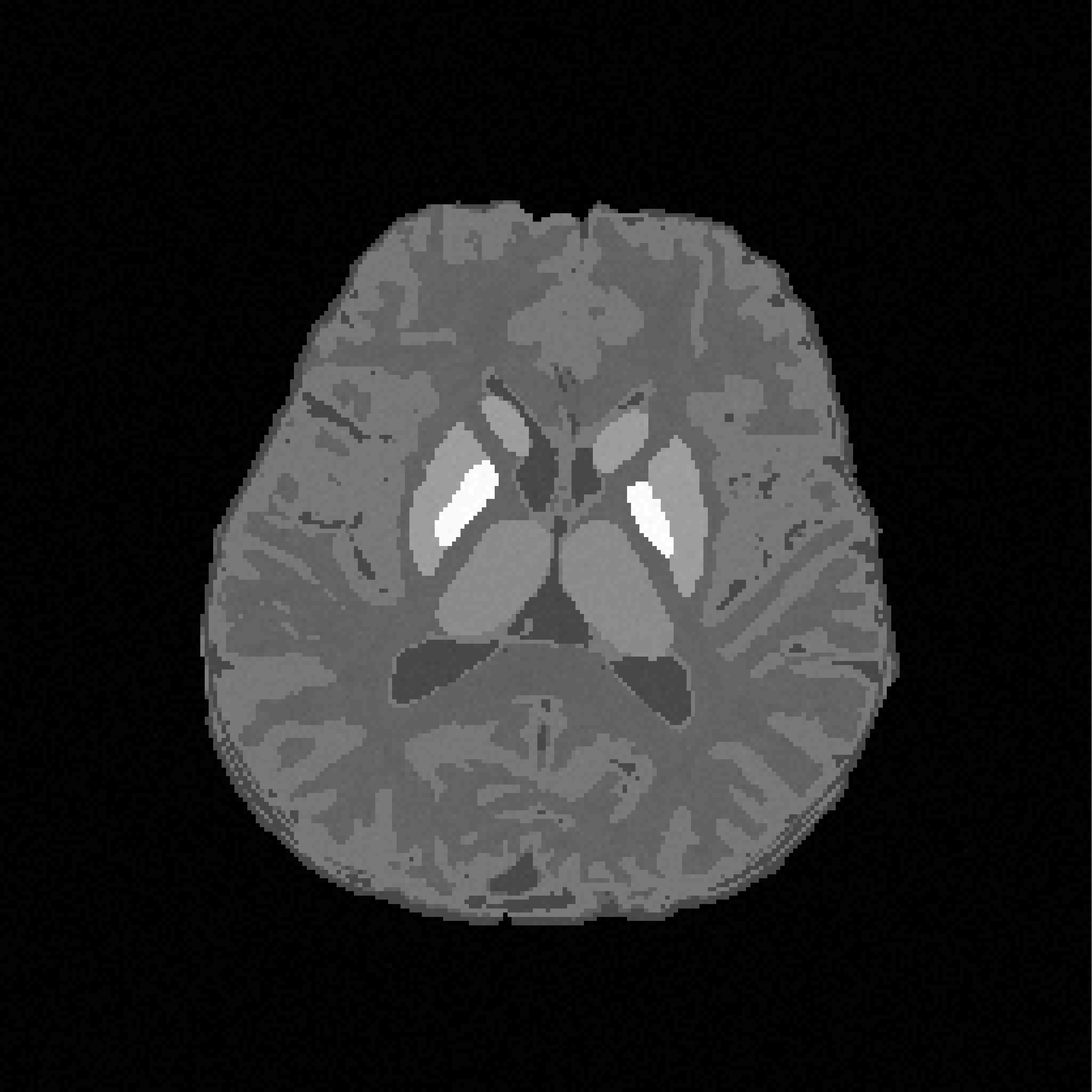}}\hspace{0.005cm}
\subfloat[ROI]{\label{PhantomMaskAx}\includegraphics[width=3.60cm]{PhantomMaskAxial.pdf}}\vspace{-0.20cm}\\
\subfloat[Phase]{\label{PhantomPhaseAx}\includegraphics[width=3.60cm]{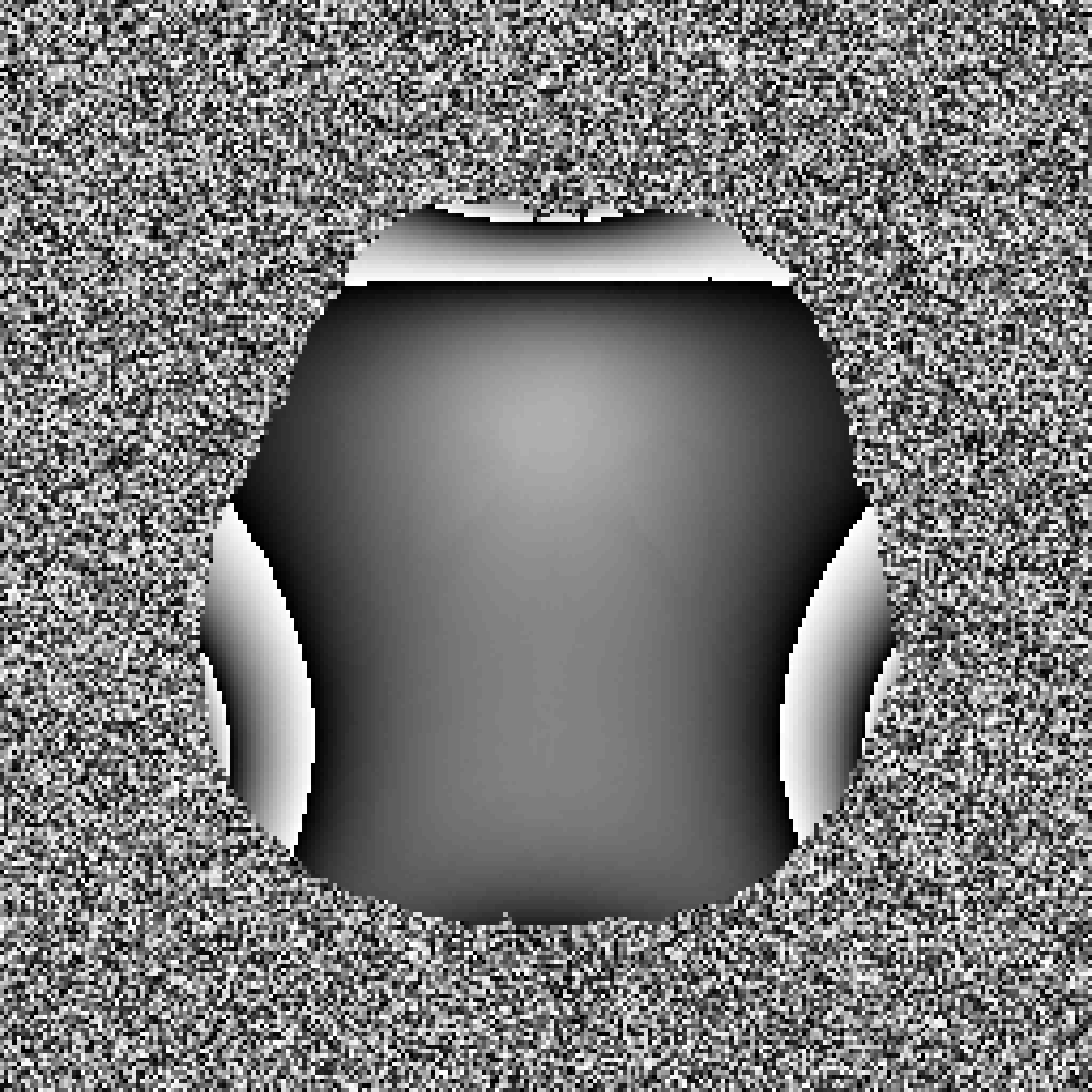}}\hspace{0.005cm}
\subfloat[Total field]{\label{PhantomTotalFieldAx}\includegraphics[width=3.60cm]{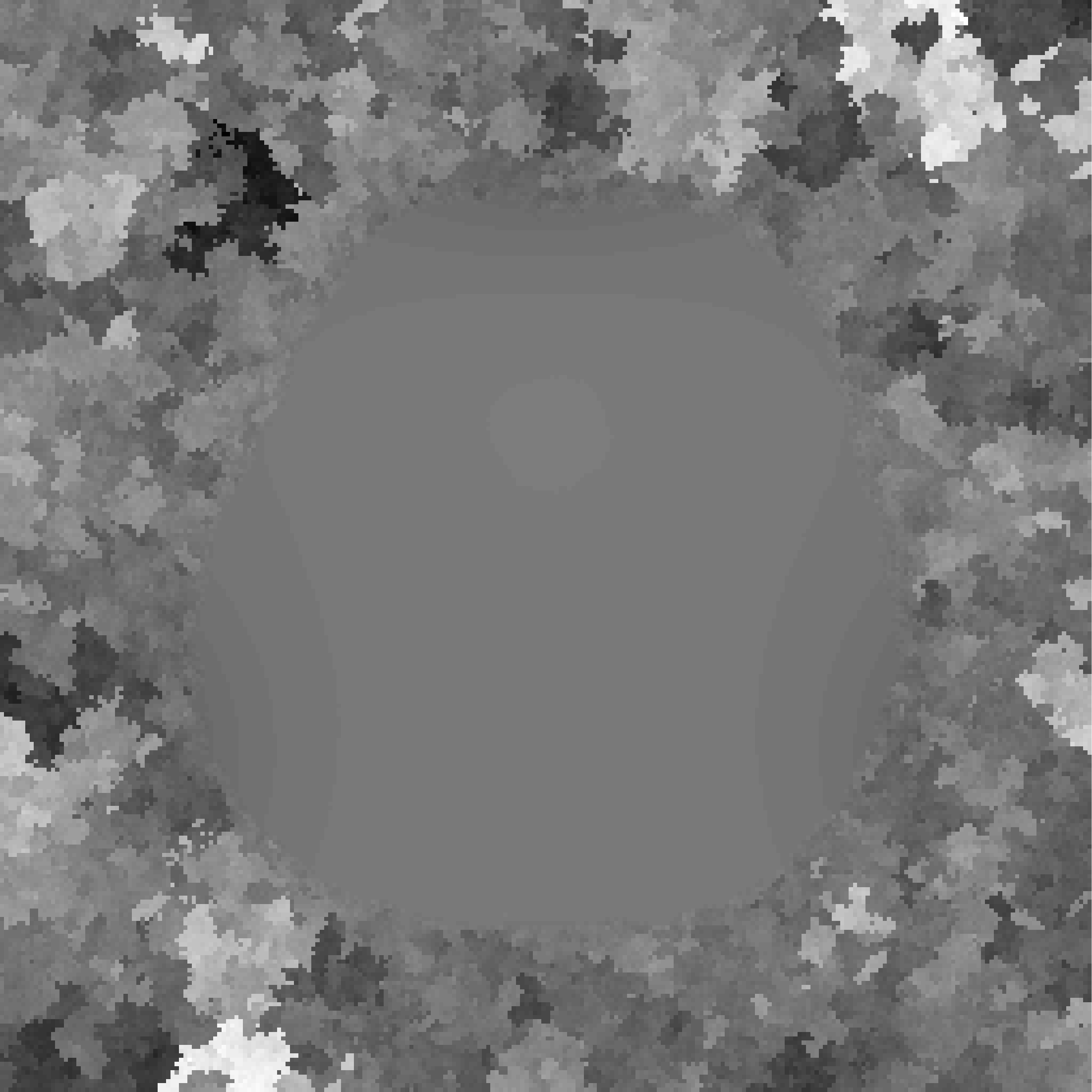}}\hspace{0.005cm}
\subfloat[Local field]{\label{PhantomLocalFieldAx}\includegraphics[width=3.60cm]{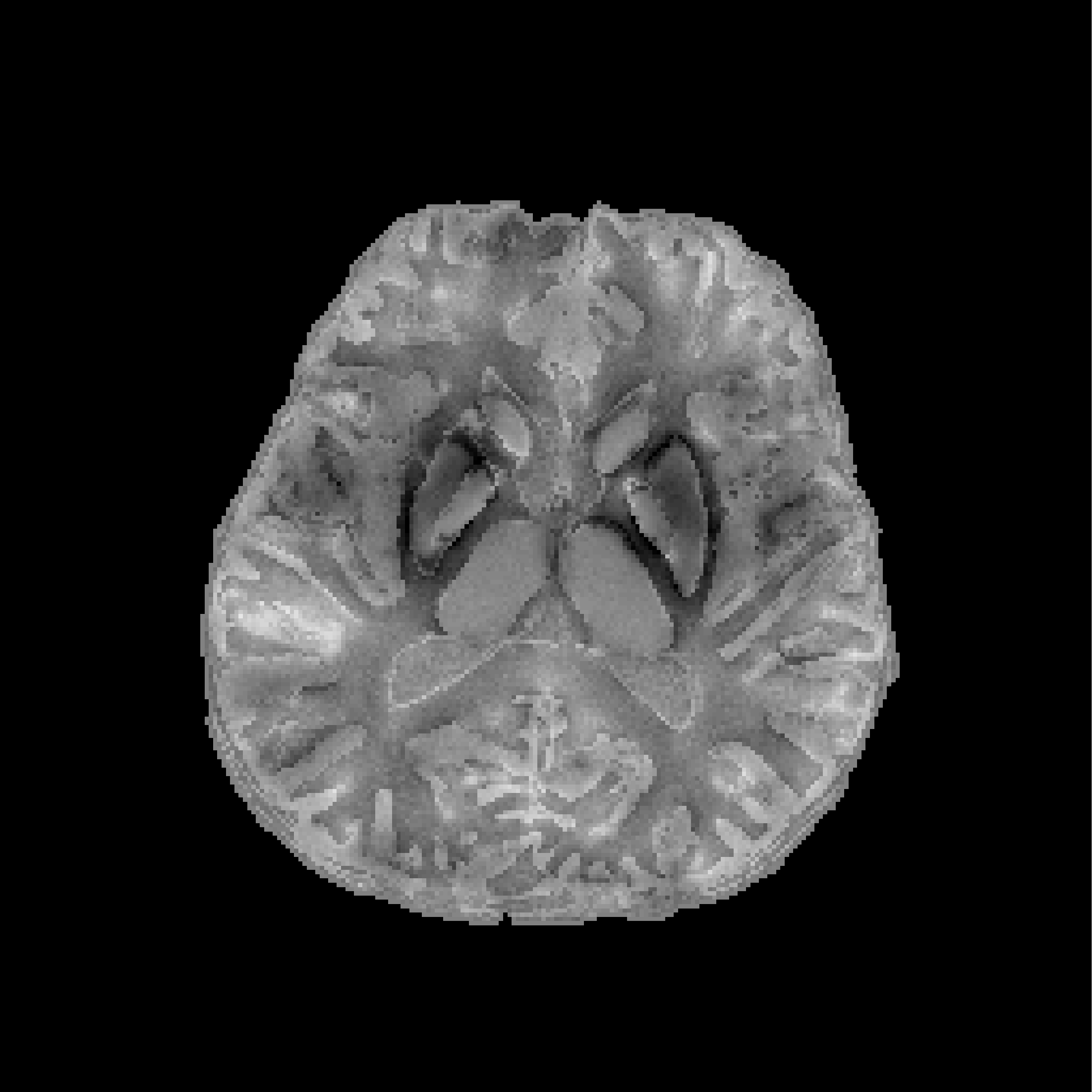}}\vspace{-0.20cm}
\caption{Axial slice images of synthesized data sets for the brain phantom experiments.}\label{BrainPhantomDataSetAxial}
\end{figure}

For the brain phantom experiments, we use $256\times256\times98$ image with spatial resolution $0.9375\times0.9375\times1.5\mathrm{mm}^3$ to simulate the $11$ equispaced multi echo GREs at $3\mathrm{T}$ with $TE$ ranging from $2.6\mathrm{msec}$ to $28.6\mathrm{msec}$. We first simulate the true total field by adding four background susceptibility sources in the true susceptibility image to provide the background field. Then we generate the multi echo complex GRE signal by
\begin{align*}
I(\bk,t)=\wt{m}(\bk)\exp\big\{-i\wt{b}(\bk)\omega_0 B_0TE(t)\big\},~~~~~~\bk\in\OO,~\&~t=1,\cdots,11
\end{align*}
with a given true magnitude image $\wt{m}$ and the true total field $\wt{b}$. Then the white Gaussian noise with standard deviation $0.02$ is added to both real and complex part of each GRE signal. Using the simulated noisy multi echo GRE signal, we estimate the magnitude image and phase data using the method in \cite{L.deRochefort2008}, and the phase is further unwrapped by the method in \cite{D.C.Ghiglia1998} to obtain the noisy and incomplete total field $b$. Finally, we solve the Poisson's equation \cref{LBV} using the method in \cite{D.Zhou2014} to obtain the noisy local field data $b_l$. (See \cref{BrainPhantomDataSet,BrainPhantomDataSetAxial}.)

All regularization based models are initialized with $\chi^0=0$, and both the Frame-HIRE and the TGV-HIRE are also initialized with $v^0=0$. For the parameters, we choose $\hbar=0.125$ for \cref{TKD}, $\eps=0.01$ for \cref{Tikhonov}, $\nu=0.0005$ for the Frame-Int and the Frame-HIRE, $\nu=0.004$ for the Frame-Diff, $\alpha_1=0.00025$ for the TGV-Int and the TGV-HIRE, and $\alpha_1=0.002$ for the TGV-Diff. In addition, we choose $\beta=0.05$ for all split Bregman algorithms to solve the regularization based models including \cref{Algorithm1}.

\begin{figure}[tp!]
\centering
\hspace{-0.1cm}\subfloat[True $\chi$]{\label{PhantomQSMFra}\includegraphics[width=3.60cm]{PhantomQSM.pdf}}\hspace{0.005cm}
\subfloat[TKD]{\label{PhantomTKDFra}\includegraphics[width=3.60cm]{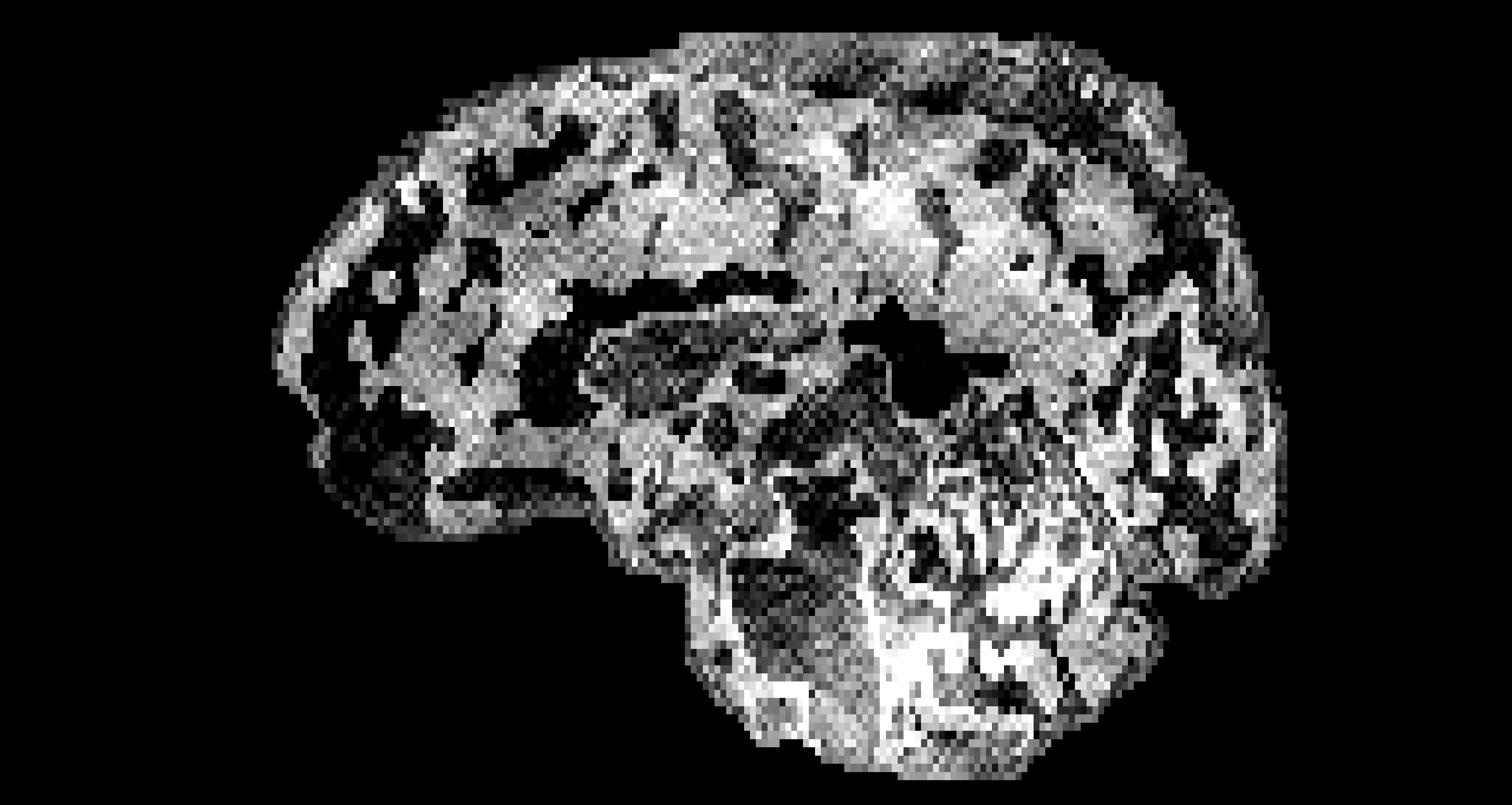}}\hspace{0.005cm}
\subfloat[Tikhonov]{\label{PhantomTikhonovFra}\includegraphics[width=3.60cm]{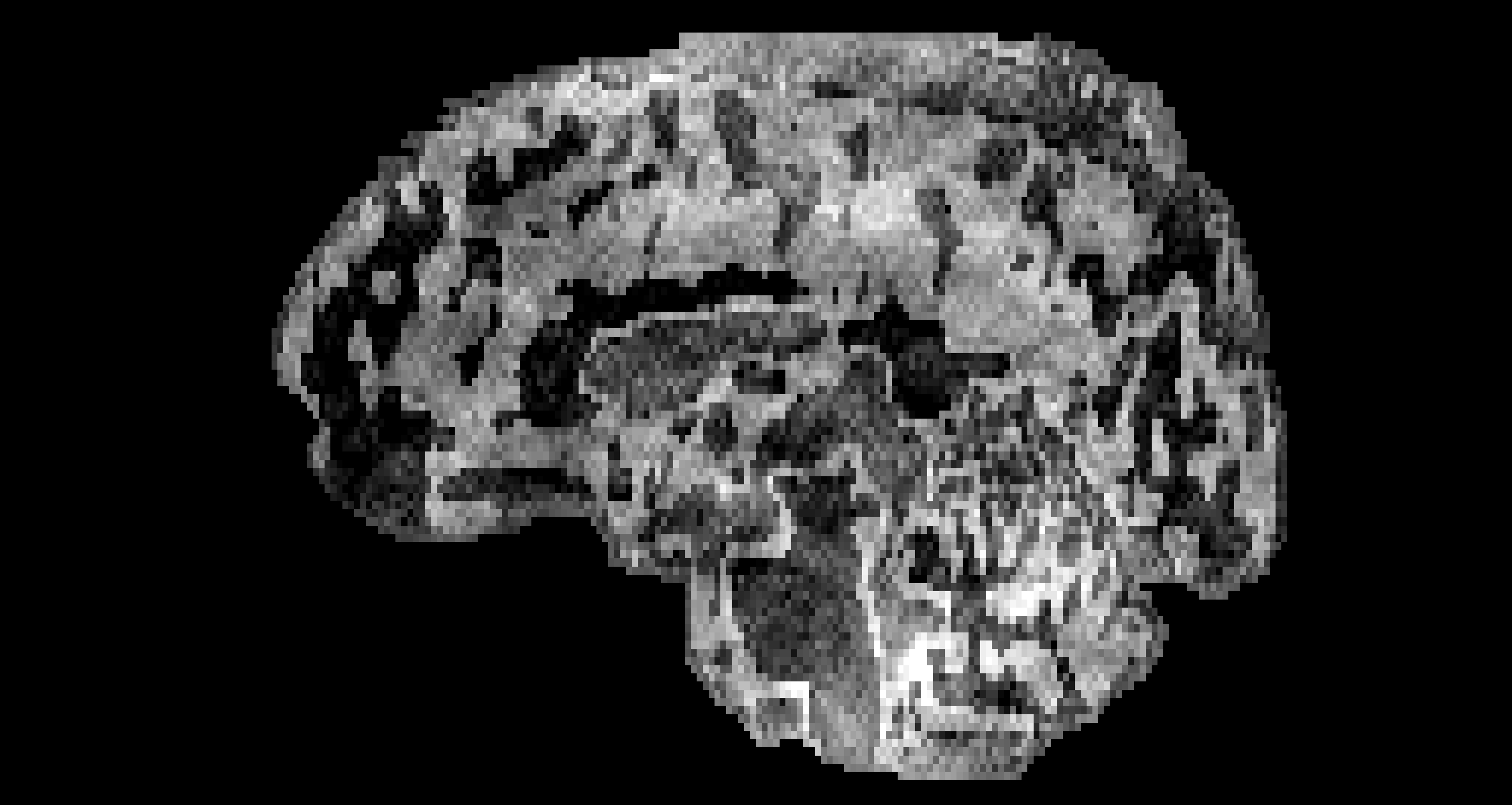}}\vspace{-0.20cm}\\
\subfloat[Frame-Int]{\label{PhantomFraInt}\includegraphics[width=3.60cm]{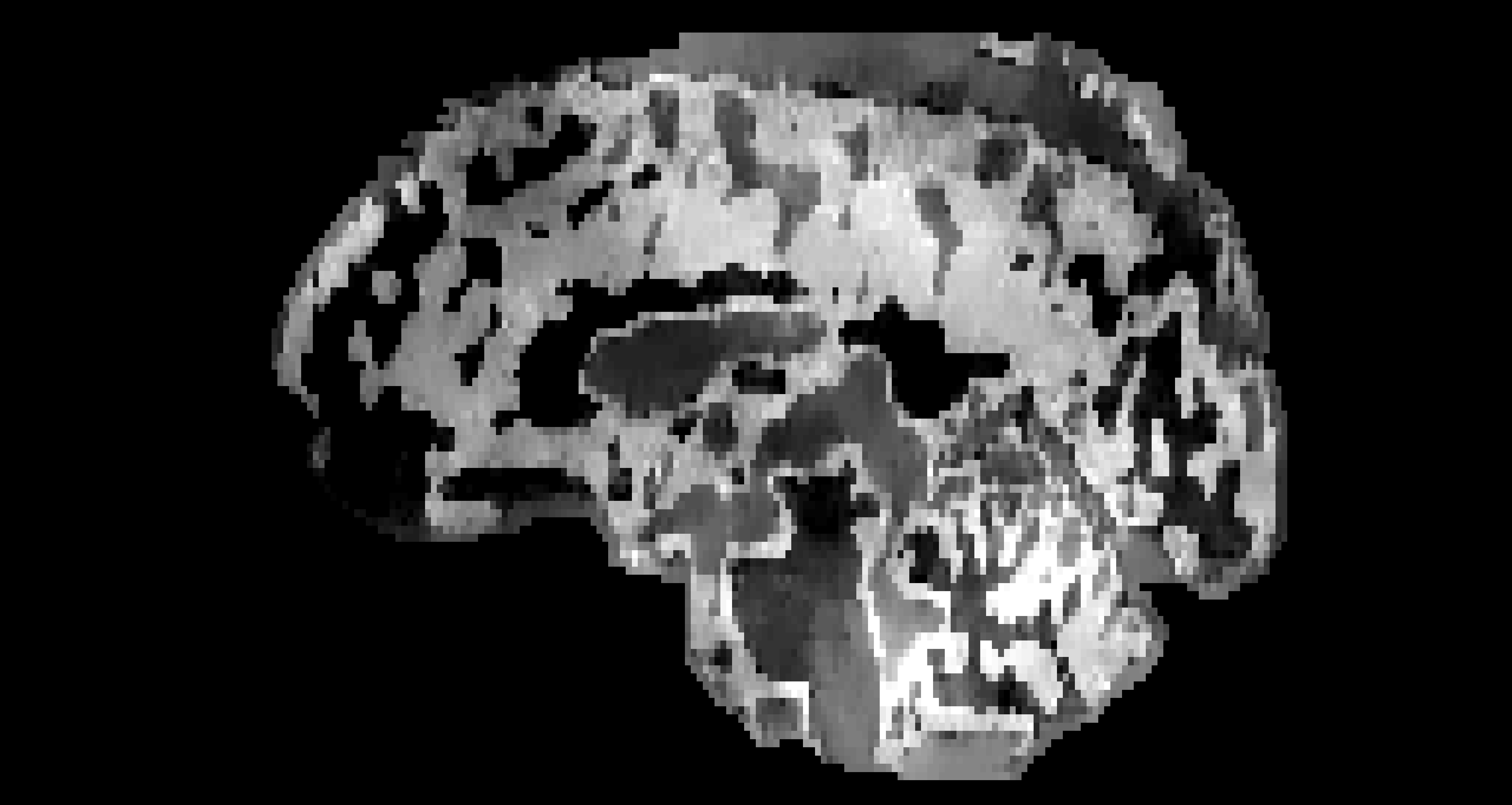}}\hspace{0.005cm}
\subfloat[Frame-Diff]{\label{PhantomFraDiff}\includegraphics[width=3.60cm]{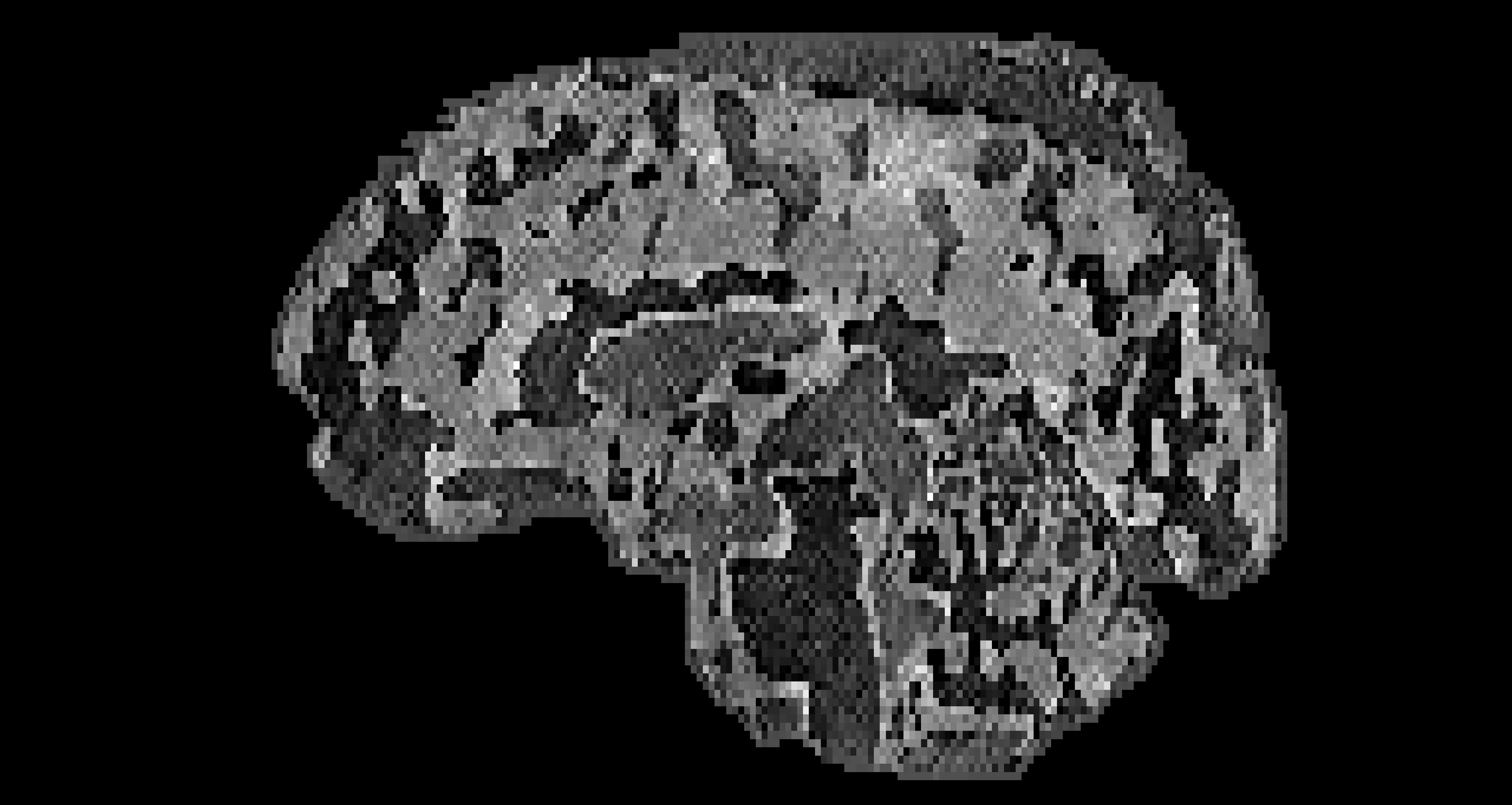}}\hspace{0.005cm}
\subfloat[Frame-HIRE]{\label{PhantomFraHIRE}\includegraphics[width=3.60cm]{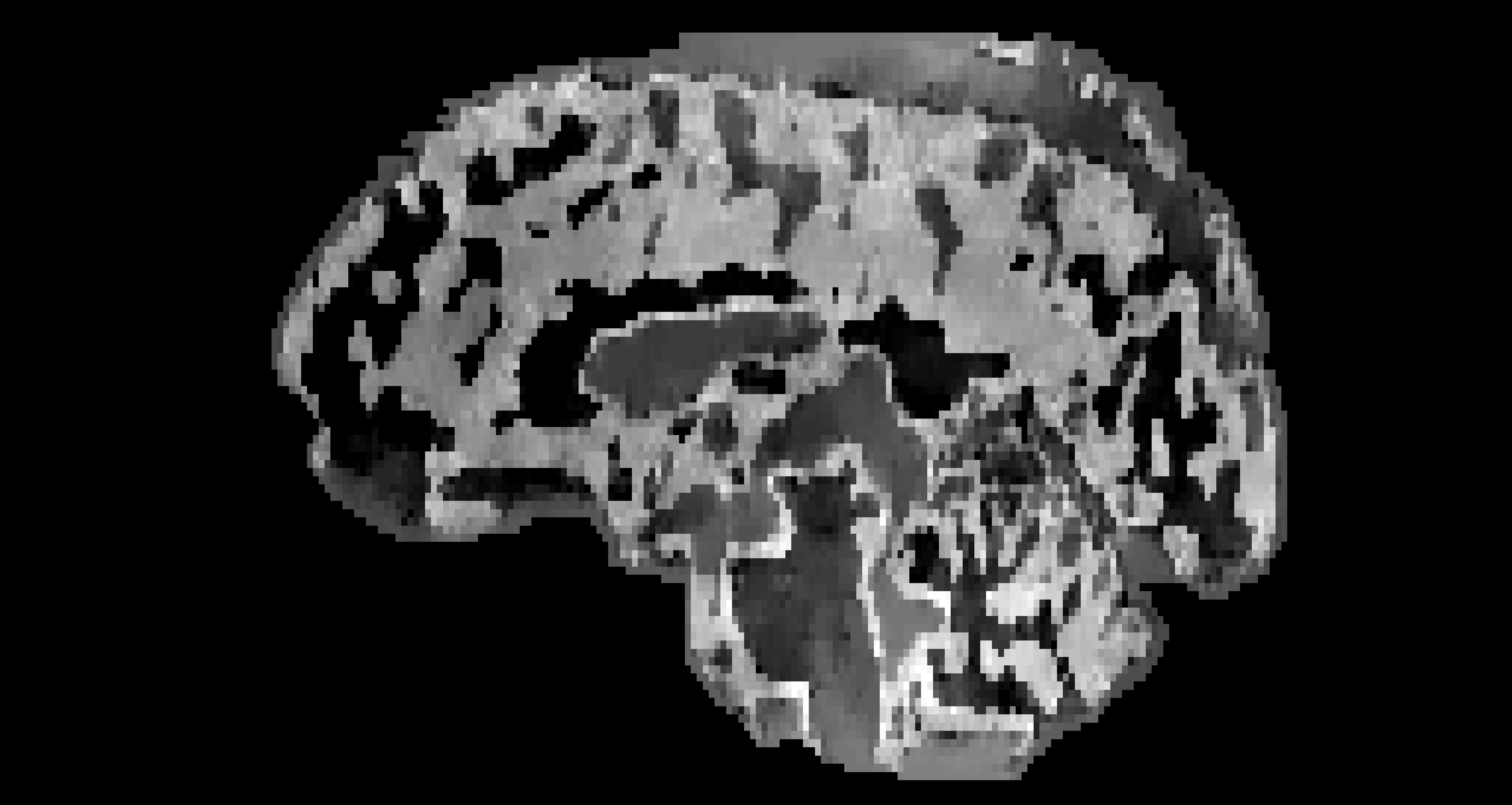}}\vspace{-0.20cm}
\caption{Sagittal slice images comparing QSM reconstruction methods for the brain phantom experiments. All sagittal slice images of brain phantom experimental results are displayed in the window level $[-0.03,0.07]$ for the fair comparison.}\label{BrainPhantomFra}
\end{figure}

\begin{figure}[tp!]
\centering
\hspace{-0.1cm}\subfloat[True $\chi$]{\label{PhantomQSMFraAx}\includegraphics[width=3.60cm]{PhantomQSMAxial.pdf}}\hspace{0.005cm}
\subfloat[TKD]{\label{PhantomTKDFraAx}\includegraphics[width=3.60cm]{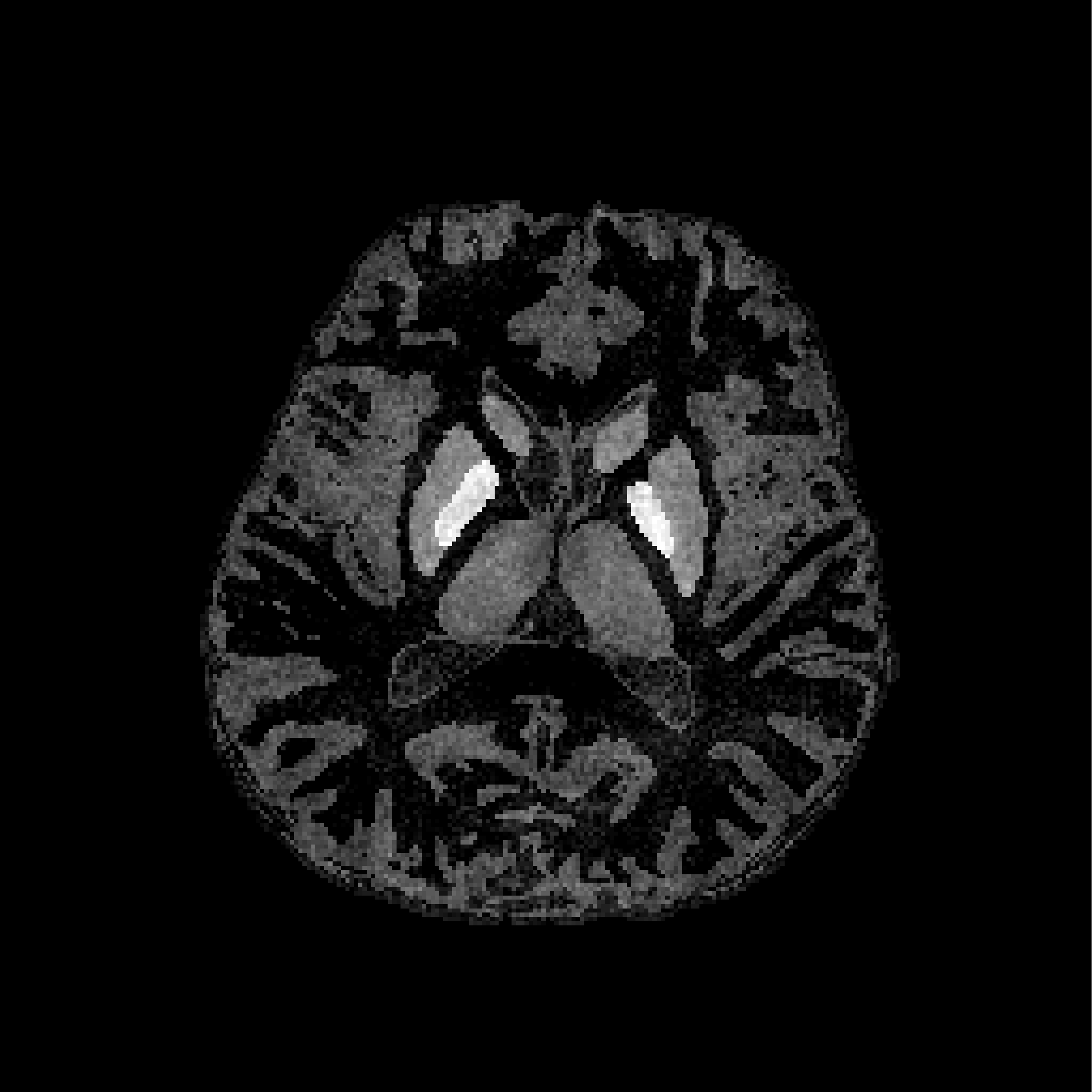}}\hspace{0.005cm}
\subfloat[Tikhonov]{\label{PhantomTikhonovFraAx}\includegraphics[width=3.60cm]{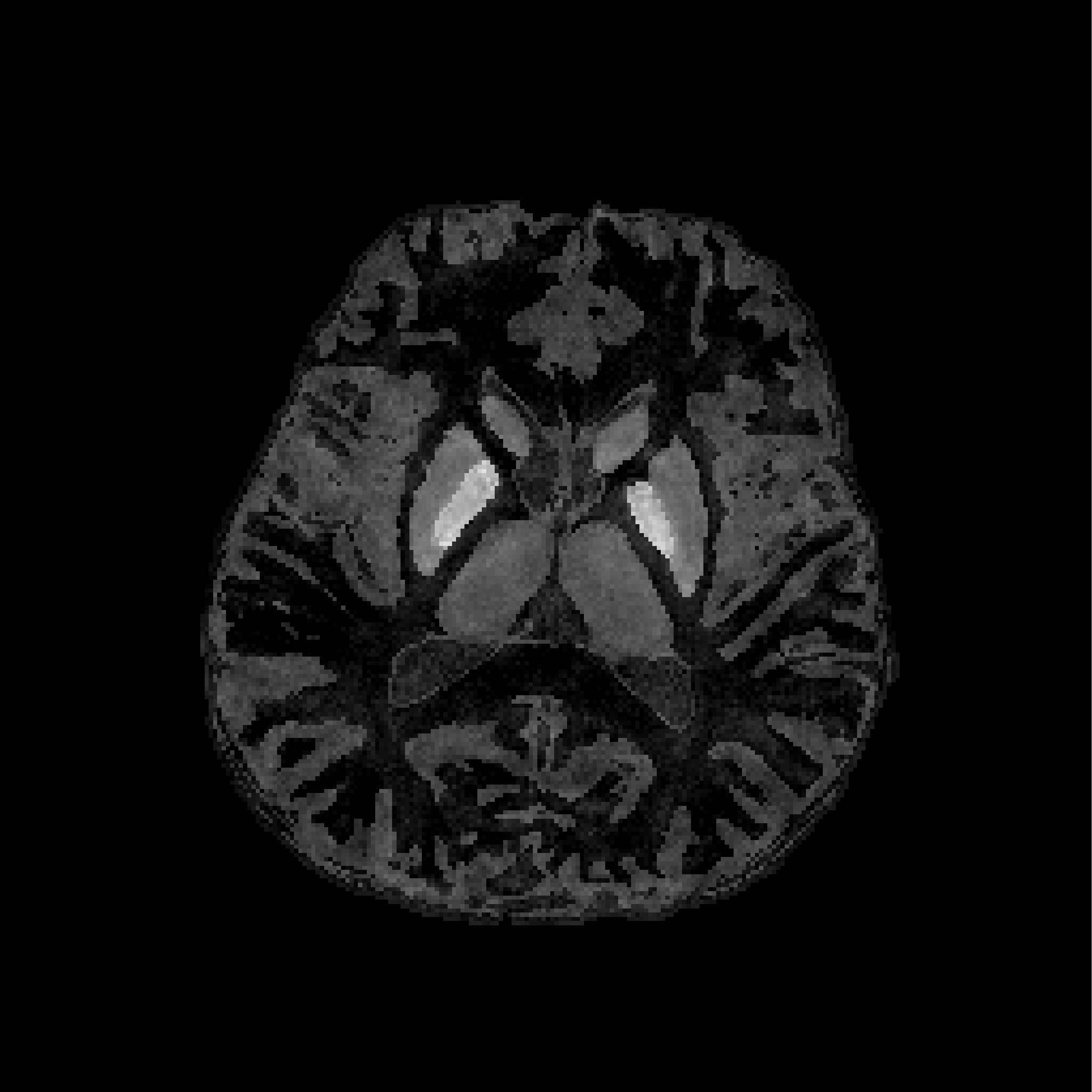}}\vspace{-0.20cm}\\
\subfloat[Frame-Int]{\label{PhantomFraIntAx}\includegraphics[width=3.60cm]{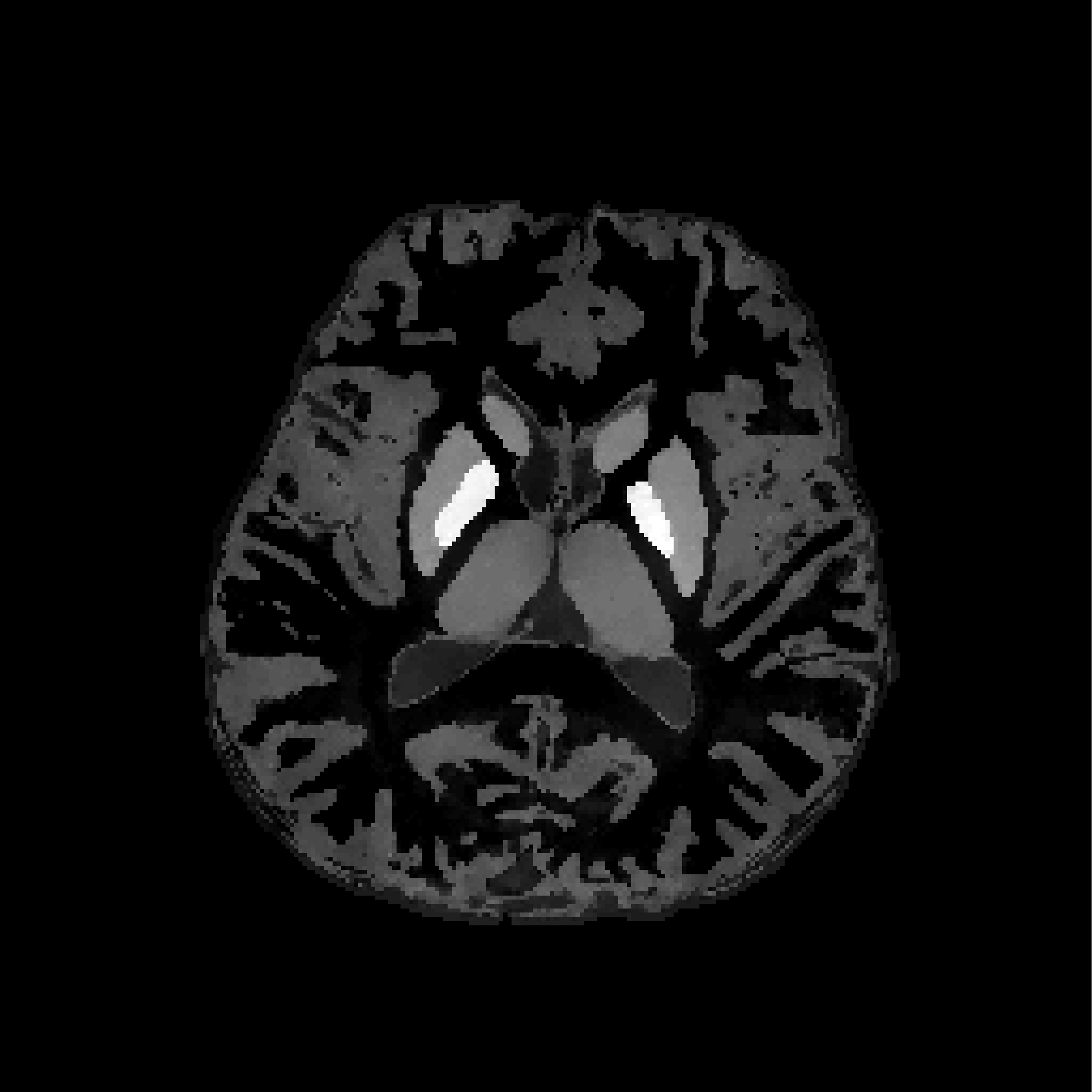}}\hspace{0.005cm}
\subfloat[Frame-Diff]{\label{PhantomFraDiffAx}\includegraphics[width=3.60cm]{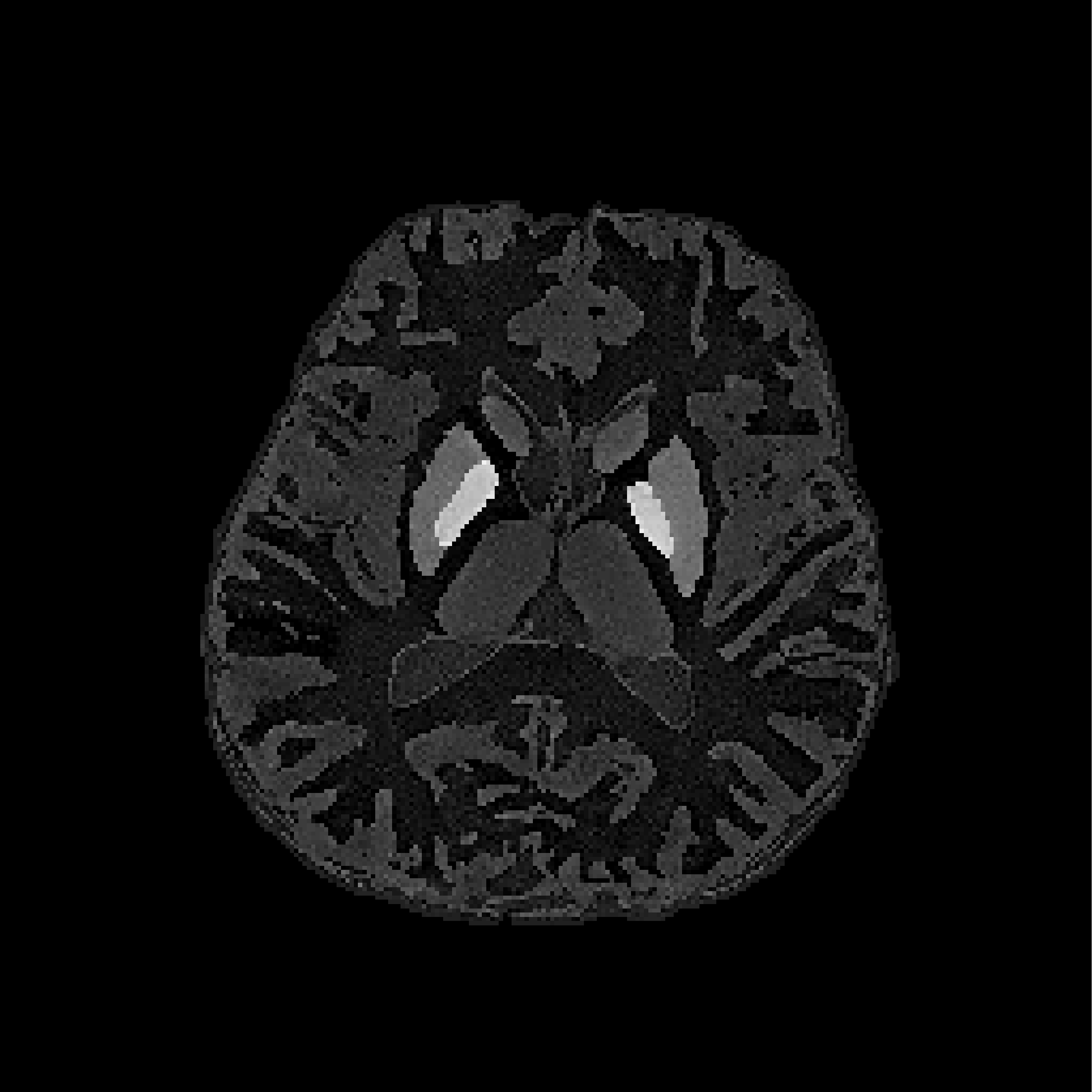}}\hspace{0.005cm}
\subfloat[Frame-HIRE]{\label{PhantomFraHIREAx}\includegraphics[width=3.60cm]{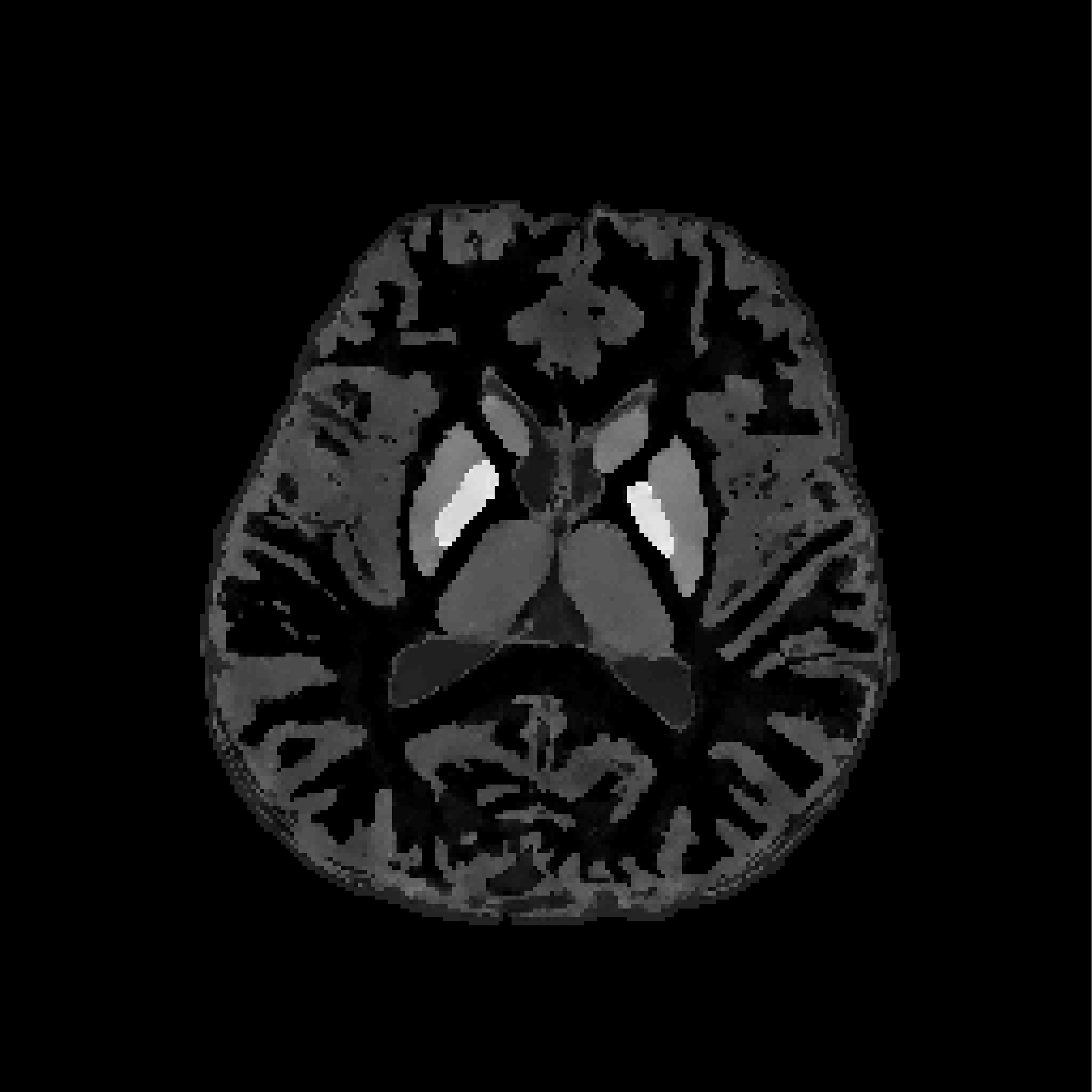}}\vspace{-0.20cm}
\caption{Axial slice images comparing QSM reconstruction methods for the brain phantom experiments with the wavelet frame regularization. All axial slice images of brain phantom experimental results are displayed in the window level $[-0.03,0.19]$ for the fair comparison.}\label{BrainPhantomFraAxial}
\end{figure}

\begin{figure}[tp!]
\centering
\hspace{-0.1cm}\subfloat[True $\chi$]{\label{PhantomQSMTGV}\includegraphics[width=3.60cm]{PhantomQSM.pdf}}\hspace{0.005cm}
\subfloat[TKD]{\label{PhantomTKDTGV}\includegraphics[width=3.60cm]{PhantomTKD.pdf}}\hspace{0.005cm}
\subfloat[Tikhonov]{\label{PhantomTikhonovTGV}\includegraphics[width=3.60cm]{PhantomTikhonov.pdf}}\vspace{-0.20cm}\\
\subfloat[TGV-Int]{\label{PhantomTGVInt}\includegraphics[width=3.60cm]{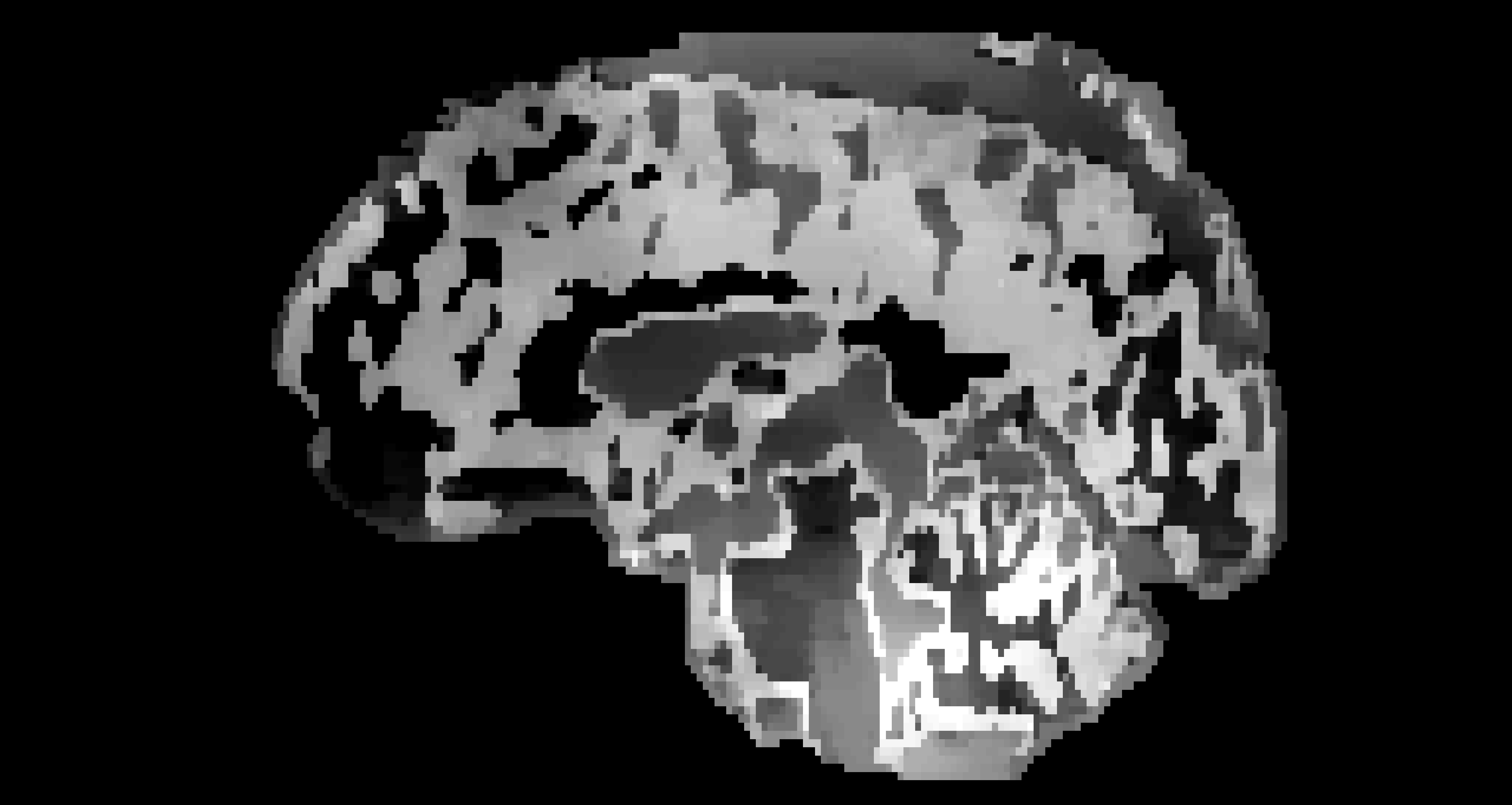}}\hspace{0.005cm}
\subfloat[TGV-Diff]{\label{PhantomTGVDiff}\includegraphics[width=3.60cm]{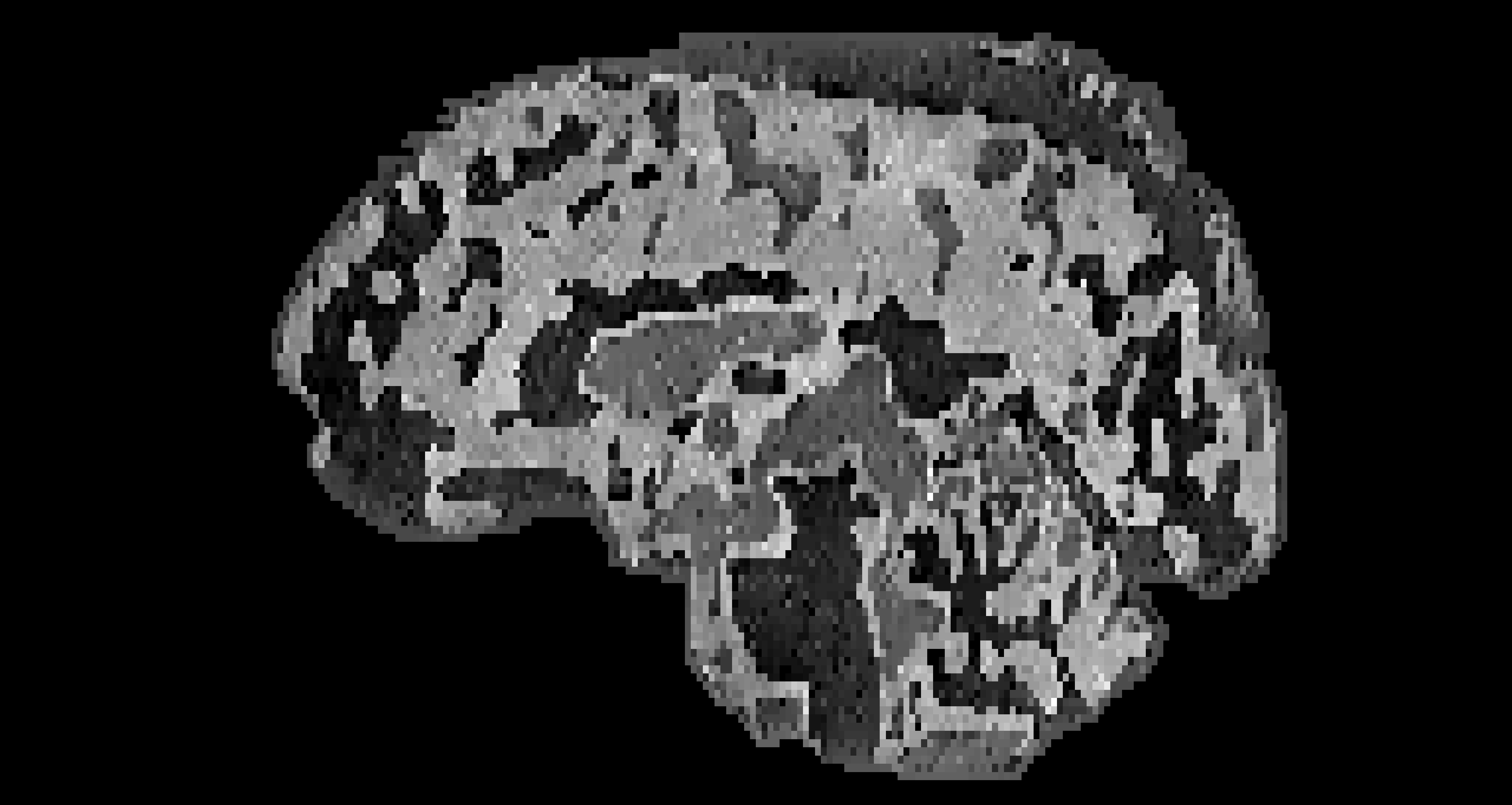}}\hspace{0.005cm}
\subfloat[TGV-HIRE]{\label{PhantomTGVHIRE}\includegraphics[width=3.60cm]{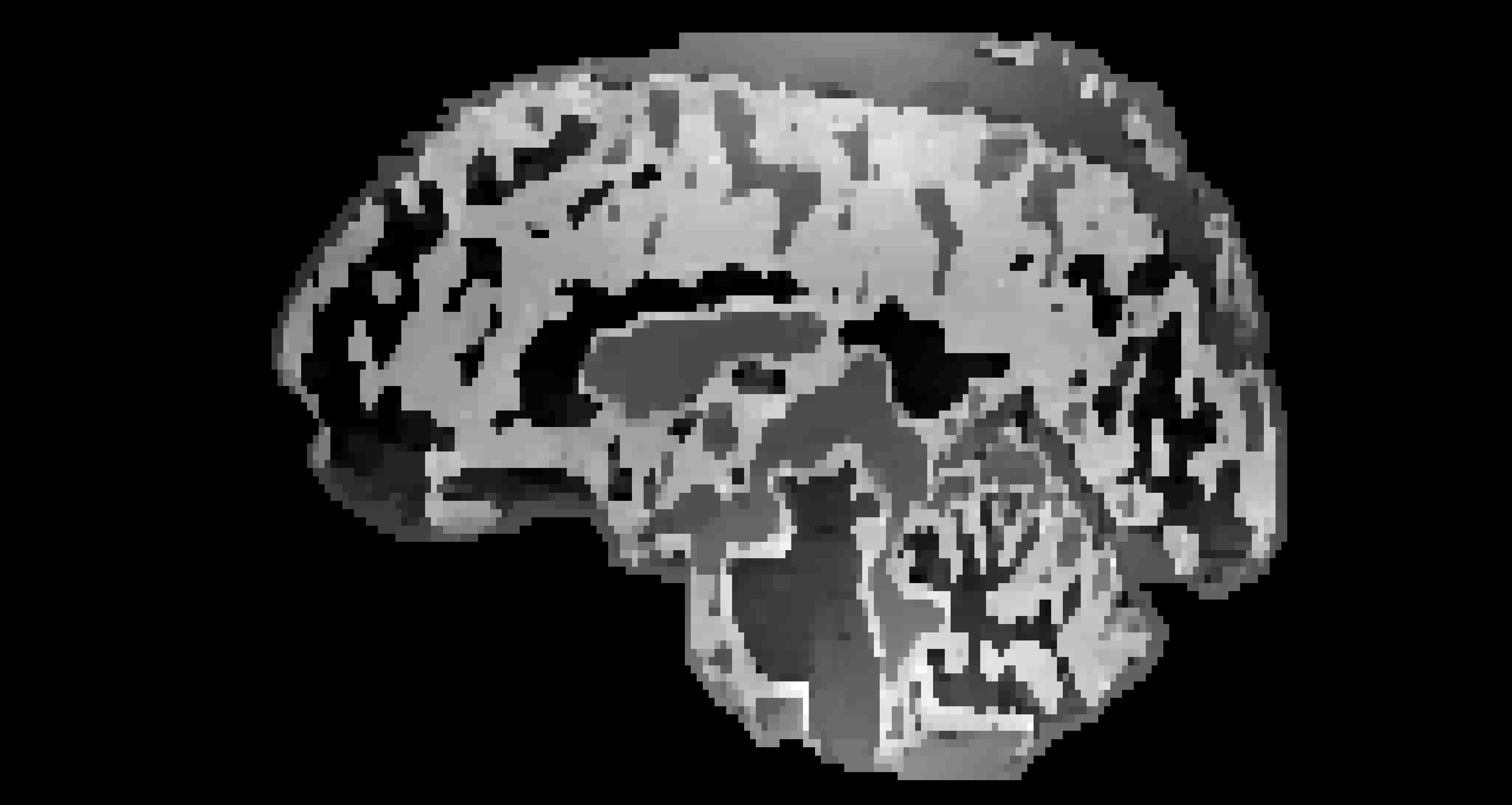}}\vspace{-0.20cm}
\caption{Sagittal slice images comparing QSM reconstruction methods for the brain phantom experiments with the TGV regularization. All sagittal slice images of brain phantom experimental results are displayed in the window level $[-0.03,0.07]$ for the fair comparison.}\label{BrainPhantomTGV}
\end{figure}

\begin{figure}[tp!]
\centering
\hspace{-0.1cm}\subfloat[True $\chi$]{\label{PhantomQSMTGVAx}\includegraphics[width=3.60cm]{PhantomQSMAxial.pdf}}\hspace{0.005cm}
\subfloat[TKD]{\label{PhantomTKDTGVAx}\includegraphics[width=3.60cm]{PhantomTKDAxial.pdf}}\hspace{0.005cm}
\subfloat[Tikhonov]{\label{PhantomTikhonovTGVAx}\includegraphics[width=3.60cm]{PhantomTikhonovAxial.pdf}}\vspace{-0.20cm}\\
\subfloat[TGV-Int]{\label{PhantomTGVIntAx}\includegraphics[width=3.60cm]{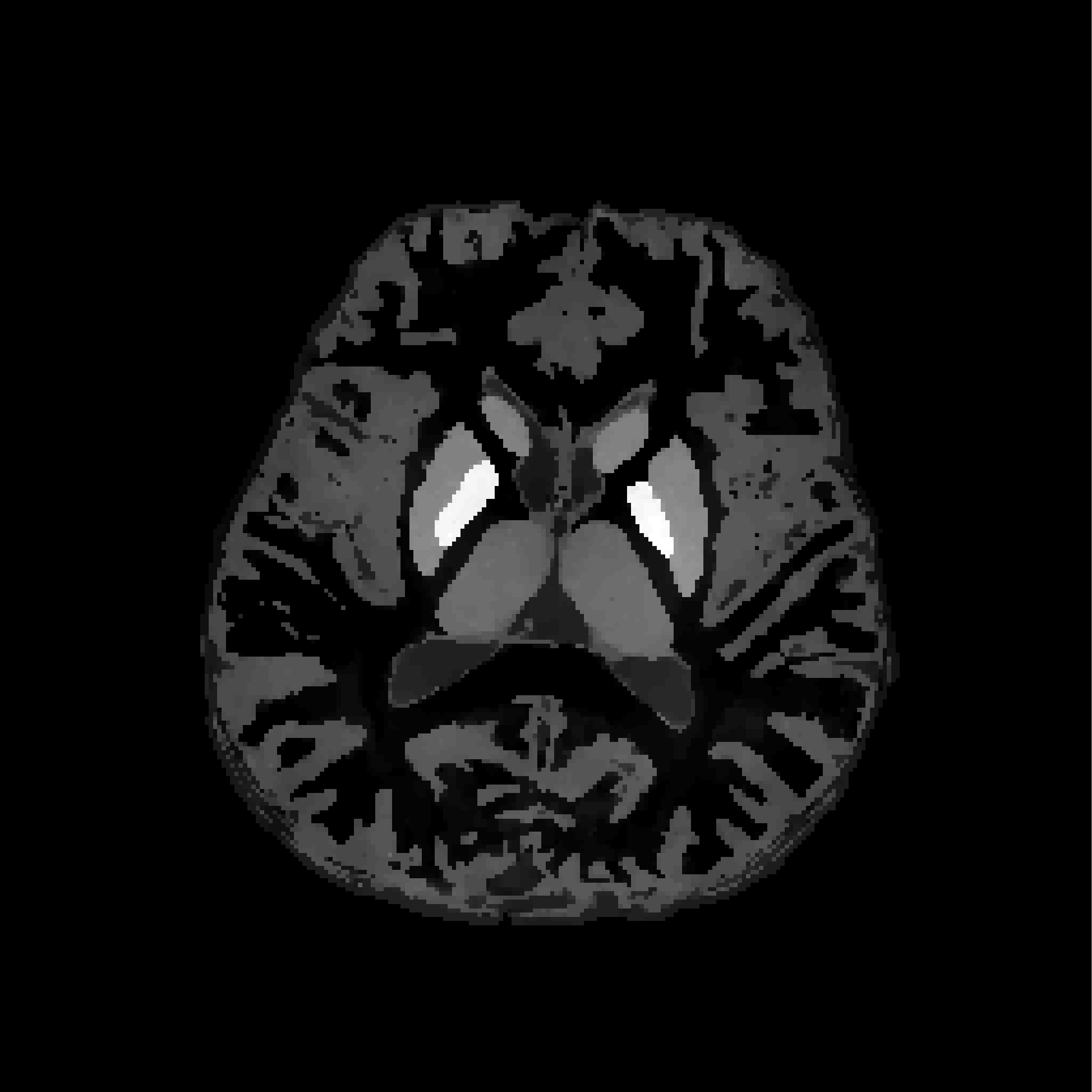}}\hspace{0.005cm}
\subfloat[TGV-Diff]{\label{PhantomTGVDiffAx}\includegraphics[width=3.60cm]{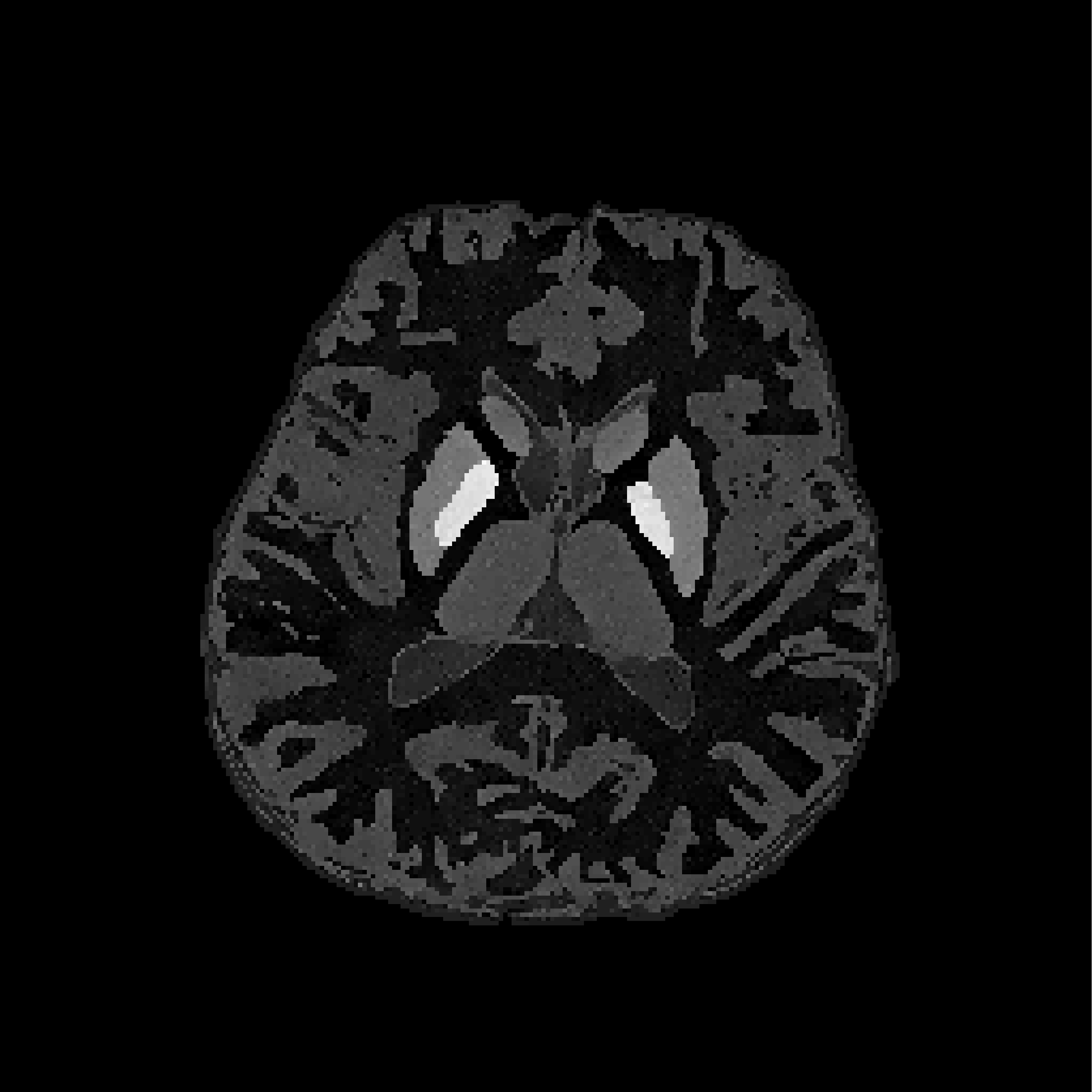}}\hspace{0.005cm}
\subfloat[TGV-HIRE]{\label{PhantomTGVHIREAx}\includegraphics[width=3.60cm]{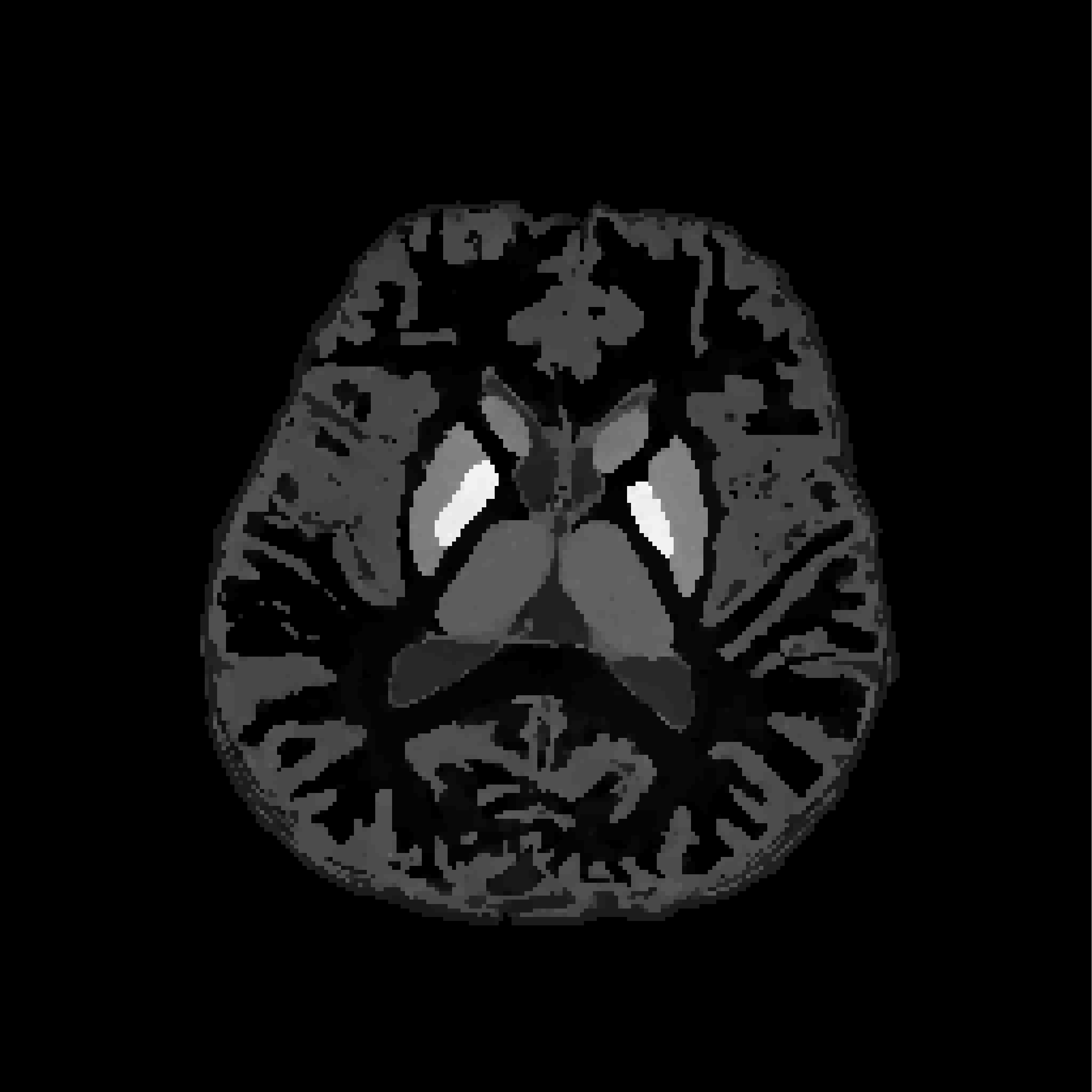}}\vspace{-0.20cm}
\caption{Axial slice images comparing QSM reconstruction methods for the brain phantom experiments with the TGV regularization. All axial slice images of brain phantom experimental results are displayed in the window level $[-0.03,0.19]$ for the fair comparison.}\label{BrainPhantomTGVAxial}
\end{figure}

\cref{ComparisonTableFra} summarizes the relative error and the SSIM of the direct approaches (\cref{TKD,Tikhonov}) and the wavelet frame regularization approaches, and \cref{BrainPhantomFra,BrainPhantomFraAxial} present visual comparisons of the results. In addition, \cref{ComparisonTableTGV} summarizes the aforementioned indices of the direct approaches and the TGV regularization approaches, and \cref{BrainPhantomTGV,BrainPhantomTGVAxial} depict the visual comparisons. We can see that both the Frame-HIRE and the TGV-HIRE consistently outperform the existing direct approaches, the integral approaches, and the differential approaches in both cases. At first glance, this verifies the convention that the regularization based models in general performs better in solving the ill-posed inverse problem of QSM than the direct methods \cite{Y.Kee2017,S.Wang2013}. Most importantly, this result demonstrates that the measured local field data obtained from the phase of a complex GRE MR signal contains the harmonic incompatibility other than the noise, which agrees with our theoretical discovery, and the performance gain mainly comes from taking both the noise in the measured data and the harmonic incompatibility (the incompatibility other than the noise) at the same time. Meanwhile, since this harmonic incompatibility is not taken into account in the integral approaches, the reconstructed susceptibility images contain the shadow artifacts as shown in \cref{PhantomFraInt,PhantomFraIntAx,PhantomTGVInt,PhantomTGVIntAx}. The differential approaches can remove the harmonic incompatibility in the measured data in advance, leading to the shadow artifact removal compared to the integral approach. However, since the noise in $b_l$ was amplified by $\msL$, the final reconstructed images contain the streaking artifacts as shown in \cref{PhantomFraDiff,PhantomFraDiffAx,PhantomTGVDiff,PhantomTGVDiffAx}, leading to the degradation in indices at the same time.

Finally, we mention that even though the TGV-HIRE performs slightly better than the Frame-HIRE from the viewpoint of indices, compared to the Frame-HIRE in \cref{PhantomFraHIRE,PhantomFraHIREAx}, the TGV-HIRE yields an overly smoothed restoration results as shown in \cref{PhantomTGVHIRE,PhantomTGVHIREAx}. In addition, since $\mD(\0)=0$, the $\chi$ subproblem of the TGV-HIRE has a rank deficient system matrix due to the constant offset, unlike the Frame-HIRE whose system matrix has a full column rank due to $W^TW=I$. As a consequence, the CPU time of the TGV-HIRE becomes approximately $3$ times longer than approximately $11.5\mathrm{min}$ of the Frame-HIRE as shown in \cref{ComparisonTimePhantom}, which shows that the TGV regularization approach may not be suitable for the real clinical applications. Therefore, even though it is approximately $1.9$ times slower than the Frame-Int, we can nevertheless conclude that compared to the TGV-HIRE, the Frame-HIRE is able to achieve the efficiency of its split Bregman algorithm as well as the shadow and streaking artifact removal.

\begin{table}[tp!]
\centering
\caption{Comparison of relative error, and structural similarity index map, for the direct approaches and the wavelet frame regularization approach in the brain phantom experiments. The bold-faced numbers indicate the best result.}\label{ComparisonTableFra}
\vspace{-0.2cm}
\begin{tabular}{|c||c|c|c|c|c|}
\hline
\multirow{2}{*}{Indices}&\multicolumn{2}{|c|}{Direct Approach}&\multicolumn{3}{|c|}{Regularization}\\ \cline{2-6}
&TKD&Tikhonov&Integral&Differential&HIRE\\ \hline
RMSE&$0.5579$&$0.5546$&$0.4516$&$0.6143$&$\textbf{0.4183}$\\ \hline
SSIM&$0.6546$&$0.6474$&$0.7485$&$0.6188$&$\textbf{0.7586}$\\ \hline
\end{tabular}
\end{table}

\begin{table}[tp!]
\centering
\caption{Comparison of relative error, and structural similarity index map for the direct approaches and the TGV regularization approach in the brain phantom experiments. The bold-faced numbers indicate the best result.}\label{ComparisonTableTGV}
\vspace{-0.2cm}
\begin{tabular}{|c||c|c|c|c|c|}
\hline
\multirow{2}{*}{Indices}&\multicolumn{2}{|c|}{Direct Approach}&\multicolumn{3}{|c|}{Regularization}\\ \cline{2-6}
&TKD&Tikhonov&Integral&Differential&HIRE\\ \hline
RMSE&$0.5579$&$0.5546$&$0.4129$&$0.4568$&$\textbf{0.3589}$\\ \hline
SSIM&$0.6546$&$0.6474$&$0.7861$&$0.7147$&$\textbf{0.7903}$\\ \hline
\end{tabular}
\end{table}

\begin{table}
\centering
\caption{Comparison of the CPU time for the brain phantom w.r.t. the choice of regularization term.}\label{ComparisonTimePhantom}
\vspace{-0.20cm}
\begin{tabular}{|c||c|c|c|c|c|c|}
\hline
\multirow{2}{*}{Indices}&\multicolumn{3}{|c|}{Wavelet Frame}&\multicolumn{3}{|c|}{TGV}\\ \cline{2-7}
&Integral&Differential&HIRE&Integral&Differential&HIRE\\ \hline
CPU Time&$366.55$&$350.33$&$685.32$&$1327.66$&$365.42$&$2020.29$\\ \hline
\end{tabular}
\end{table}

%
%

\subsection{Experiments on In Vivo MR Data}\label{RealMRData}

\begin{figure}[tp!]
\centering
\hspace{-0.1cm}\subfloat[Magnitude]{\label{RealCornellMag}\includegraphics[width=3.00cm]{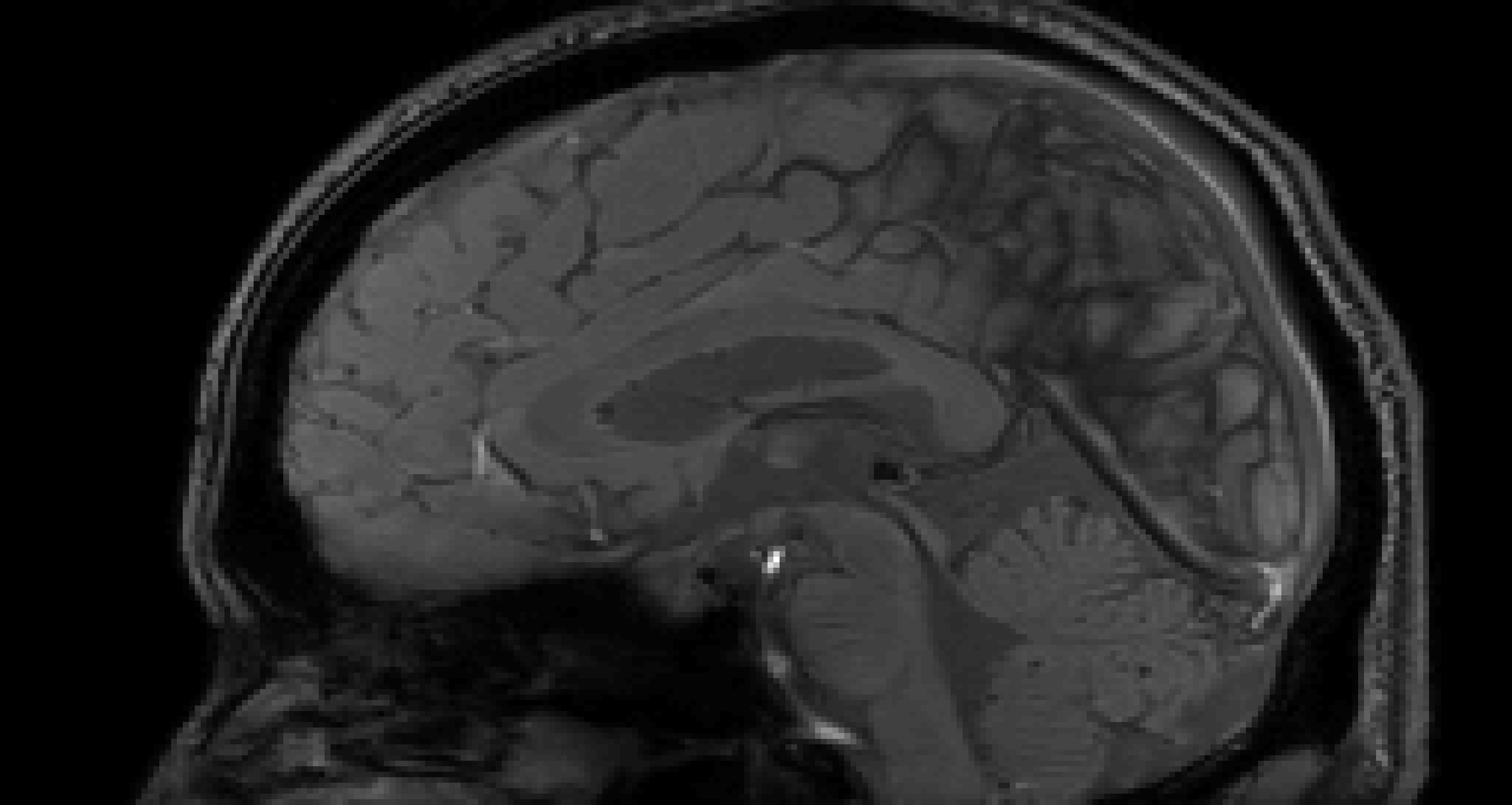}}\hspace{0.005cm}
\subfloat[ROI]{\label{RealCornellMask}\includegraphics[width=3.00cm]{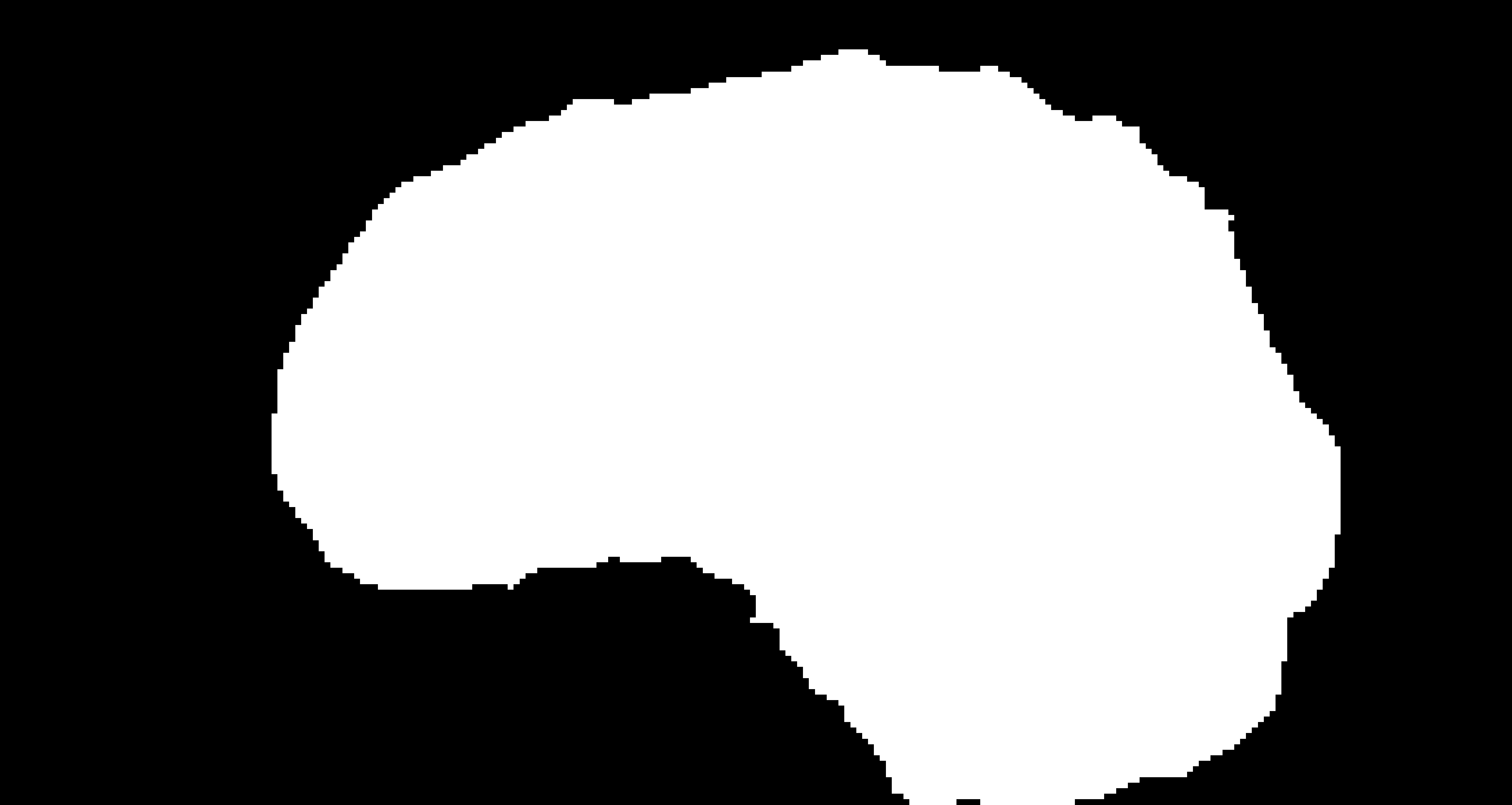}}\hspace{0.005cm}
\subfloat[Phase]{\label{RealCornellPhase}\includegraphics[width=3.00cm]{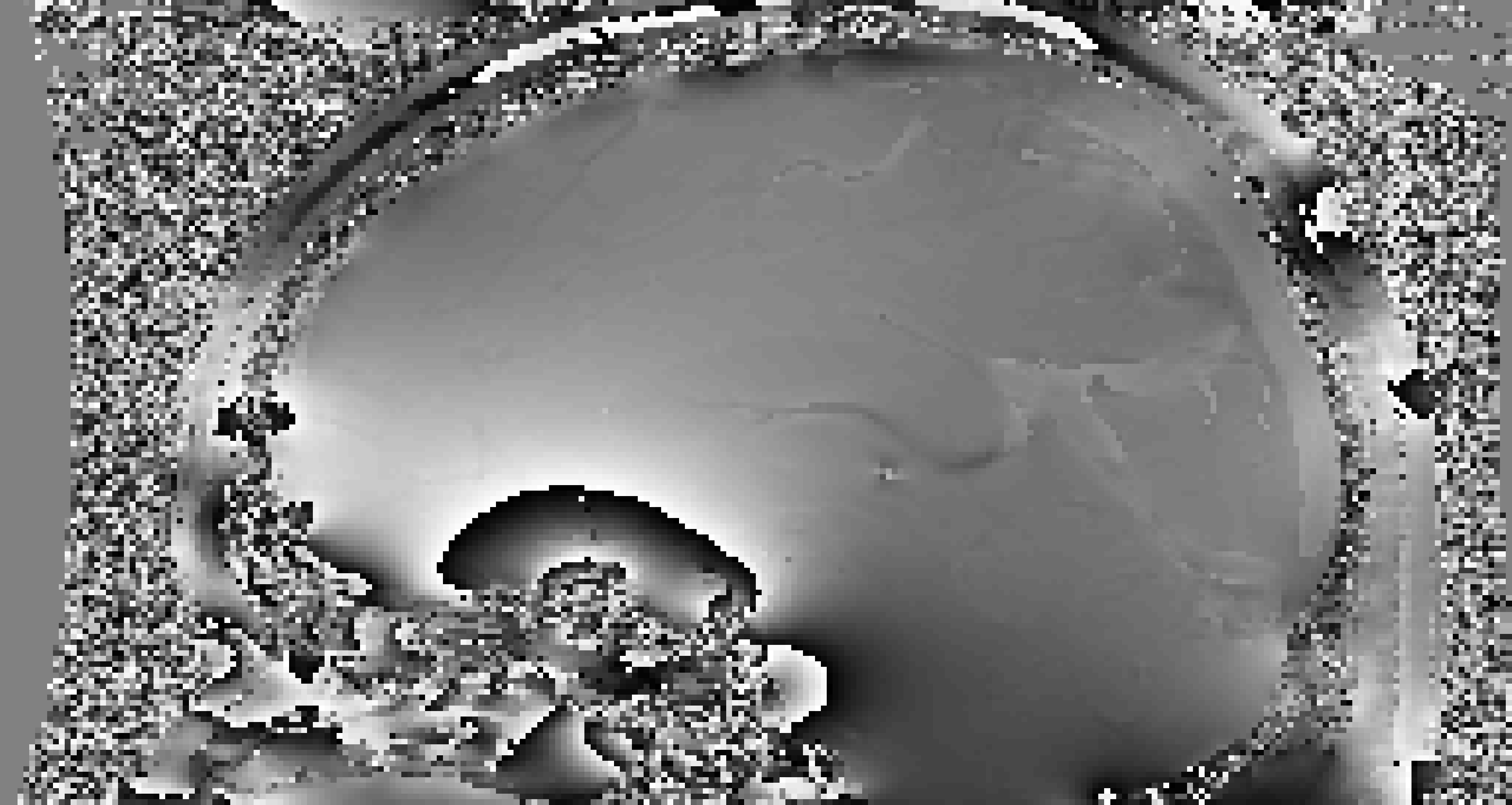}}\hspace{0.005cm}
\subfloat[Total field]{\label{RealCornellTotalField}\includegraphics[width=3.00cm]{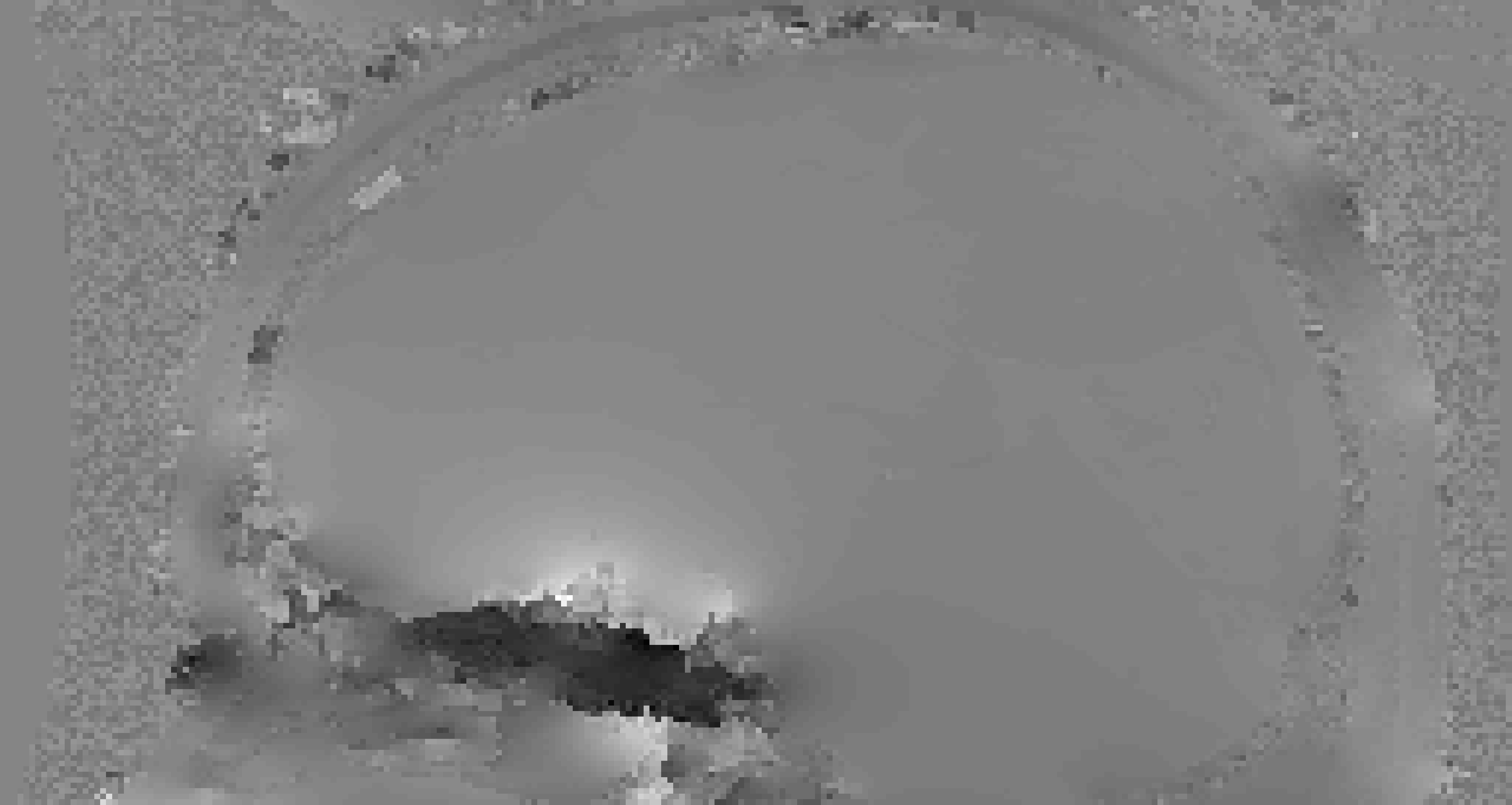}}\hspace{0.005cm}
\subfloat[Local field]{\label{RealCornellLocalField}\includegraphics[width=3.00cm]{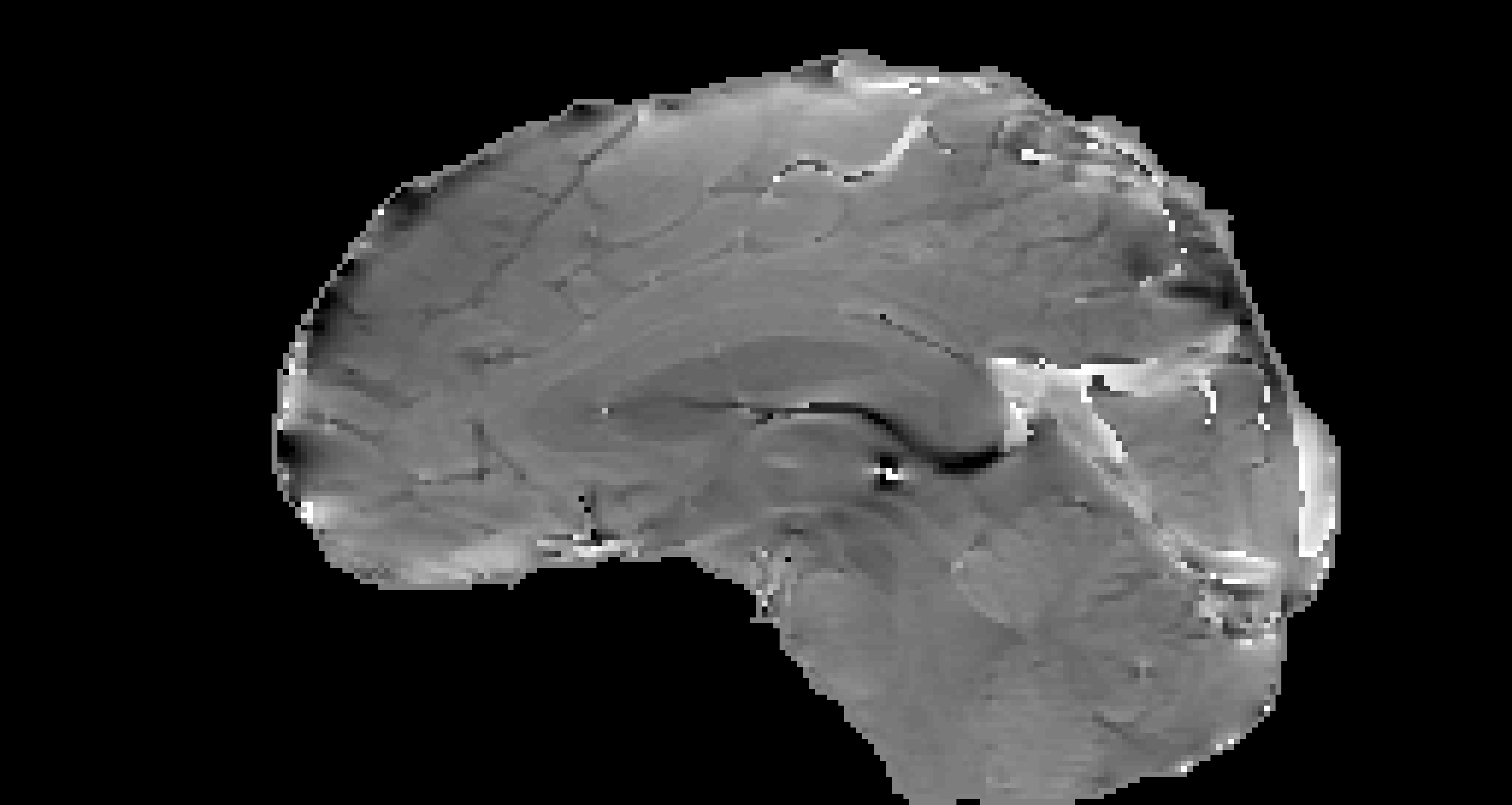}}\vspace{-0.20cm}
\caption{Sagittal slice images of data sets for the in vivo MR data experiments.}\label{InVivoDataSet}
\end{figure}

\begin{figure}[tp!]
\centering
\hspace{-0.1cm}\subfloat[Magnitude]{\label{RealCornellMagAx}\includegraphics[width=3.00cm]{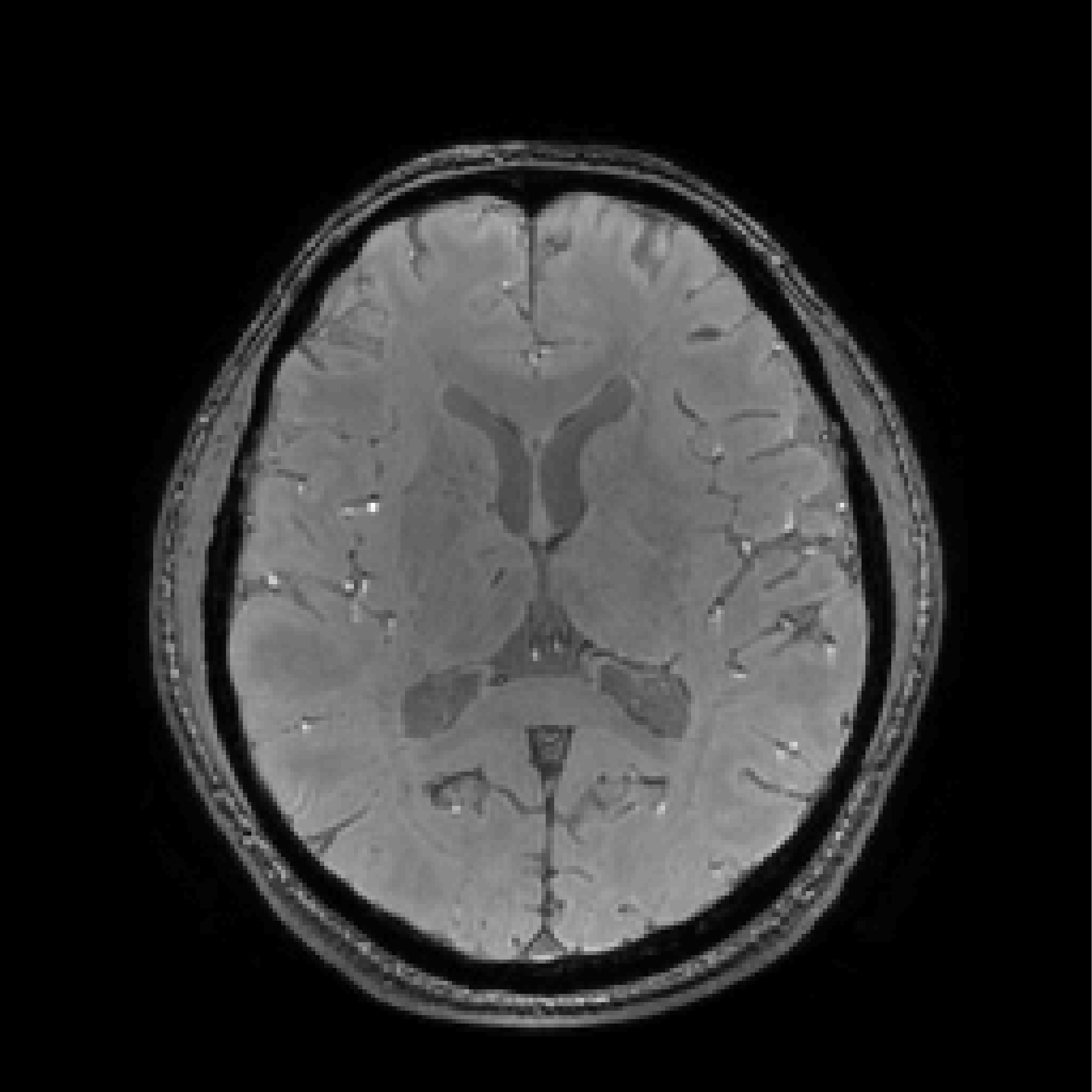}}\hspace{0.005cm}
\subfloat[ROI]{\label{RealCornellMaskAx}\includegraphics[width=3.00cm]{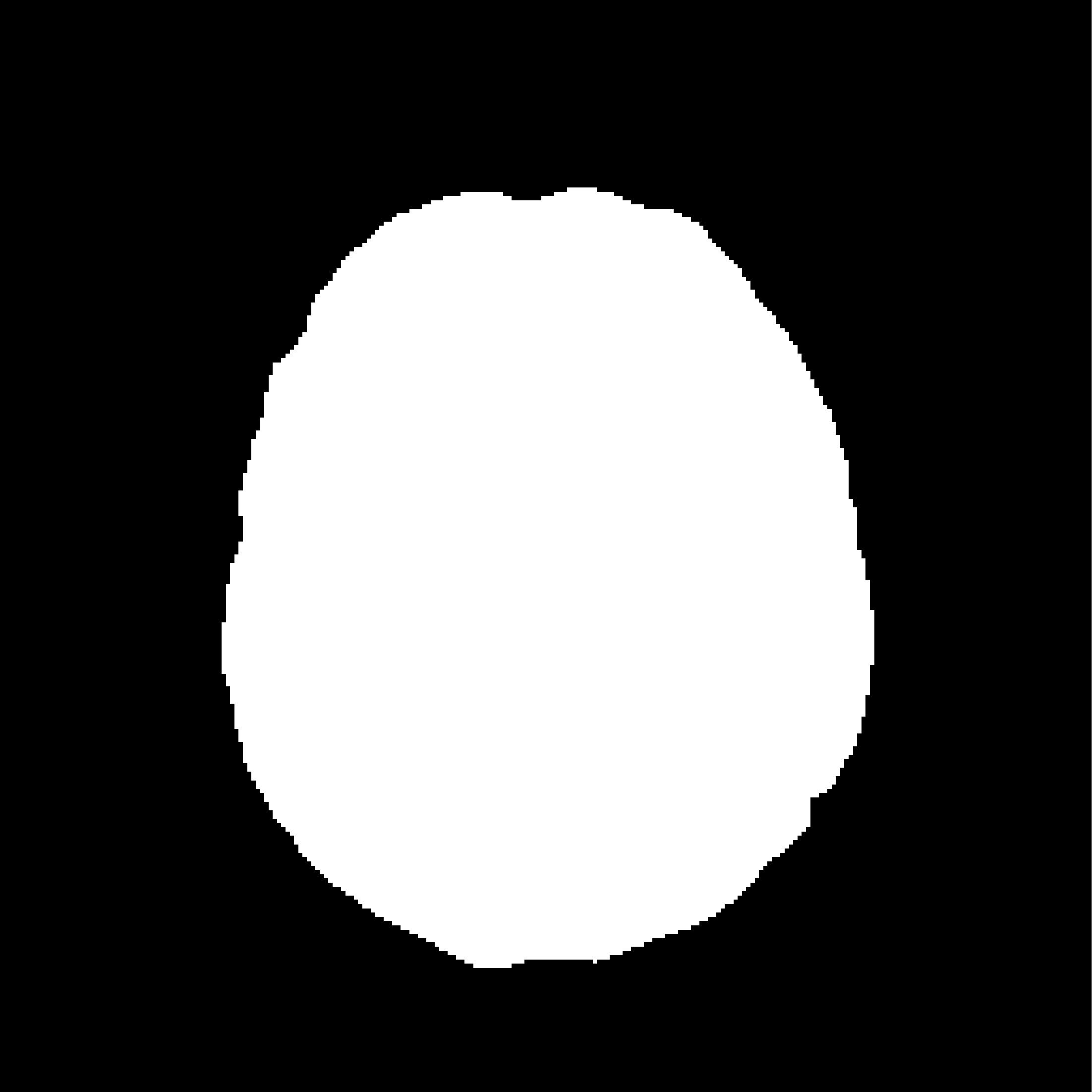}}\hspace{0.005cm}
\subfloat[Phase]{\label{RealCornellPhaseAx}\includegraphics[width=3.00cm]{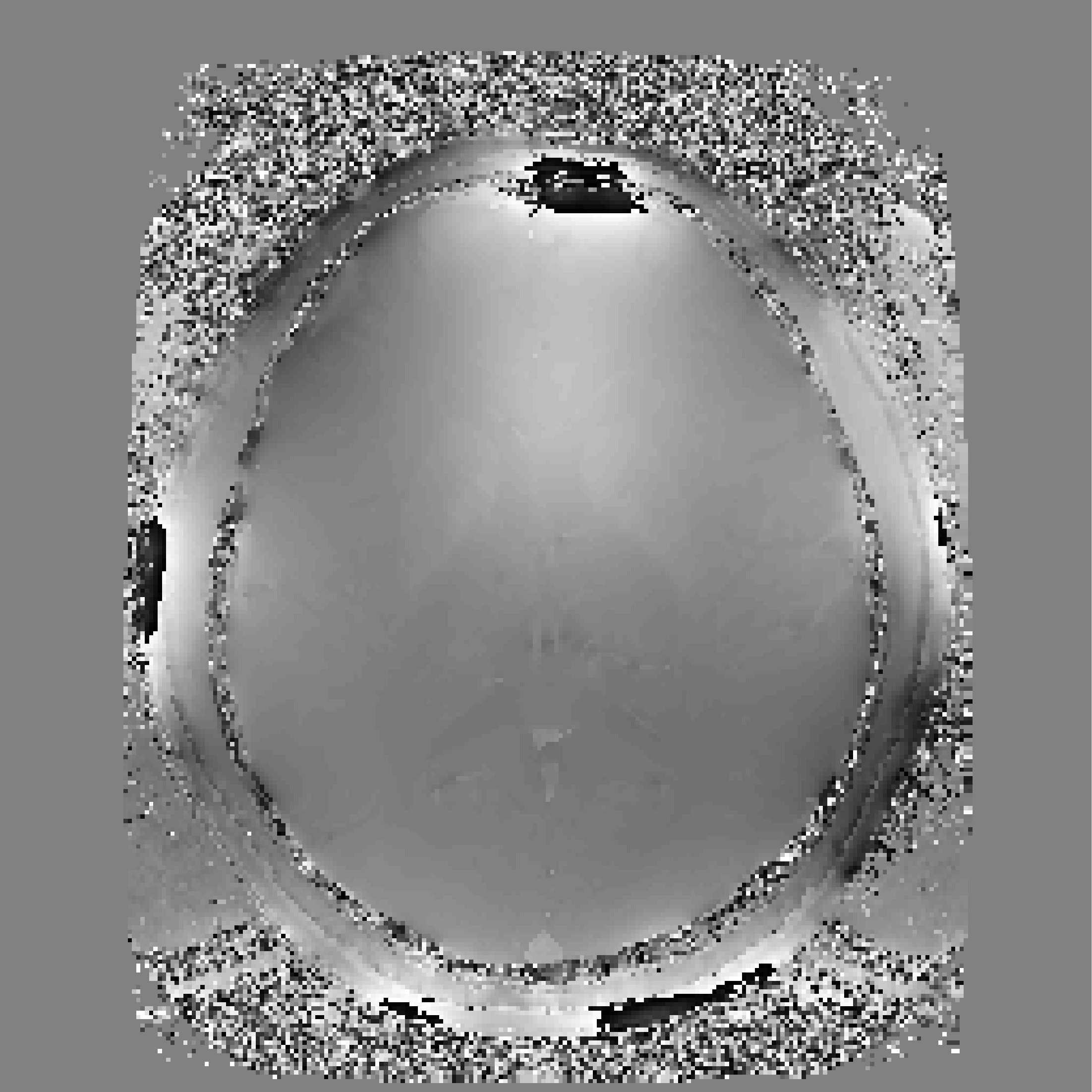}}\hspace{0.005cm}
\subfloat[Total field]{\label{RealCornellTotalFieldAx}\includegraphics[width=3.00cm]{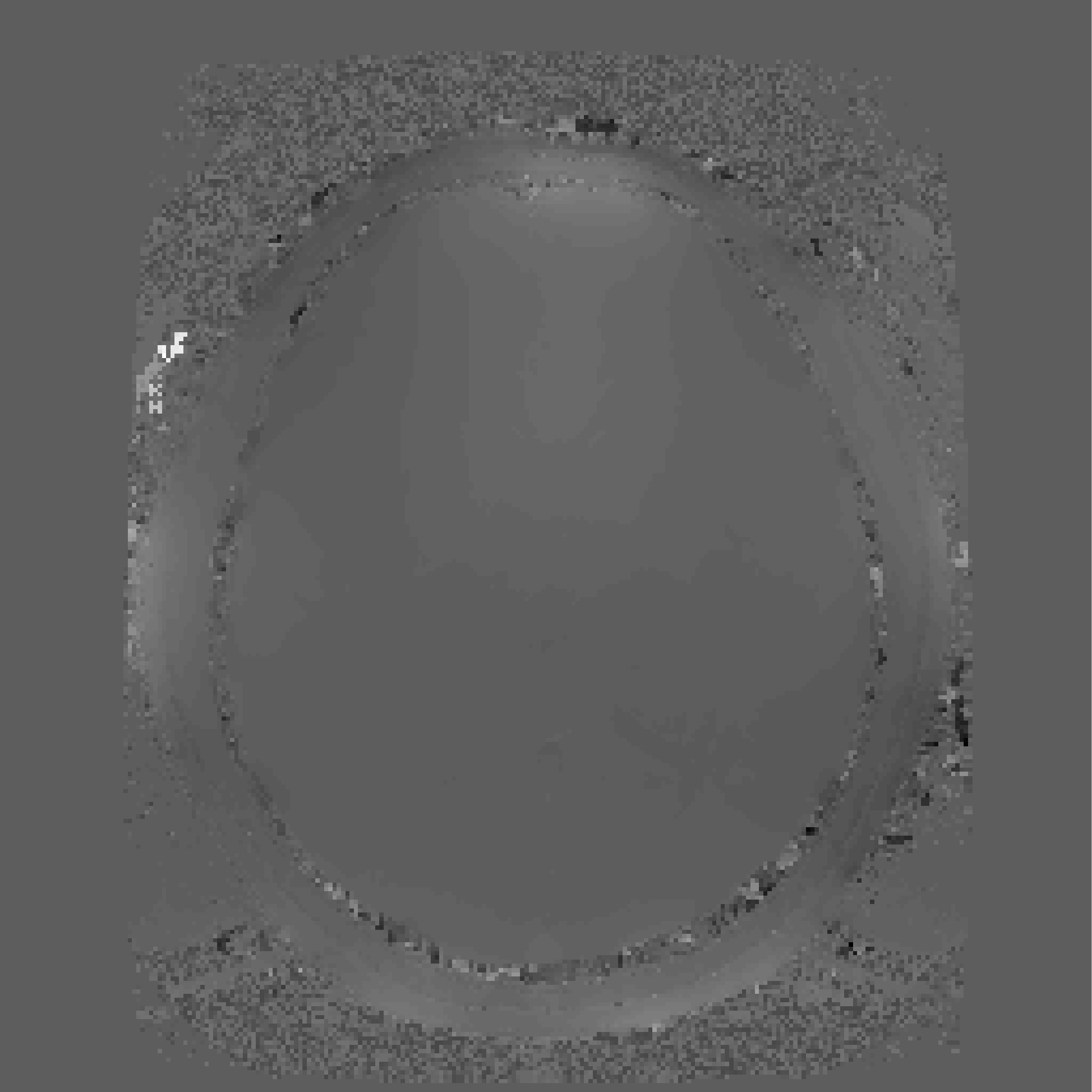}}\hspace{0.005cm}
\subfloat[Local field]{\label{RealCornellLocalFieldAx}\includegraphics[width=3.00cm]{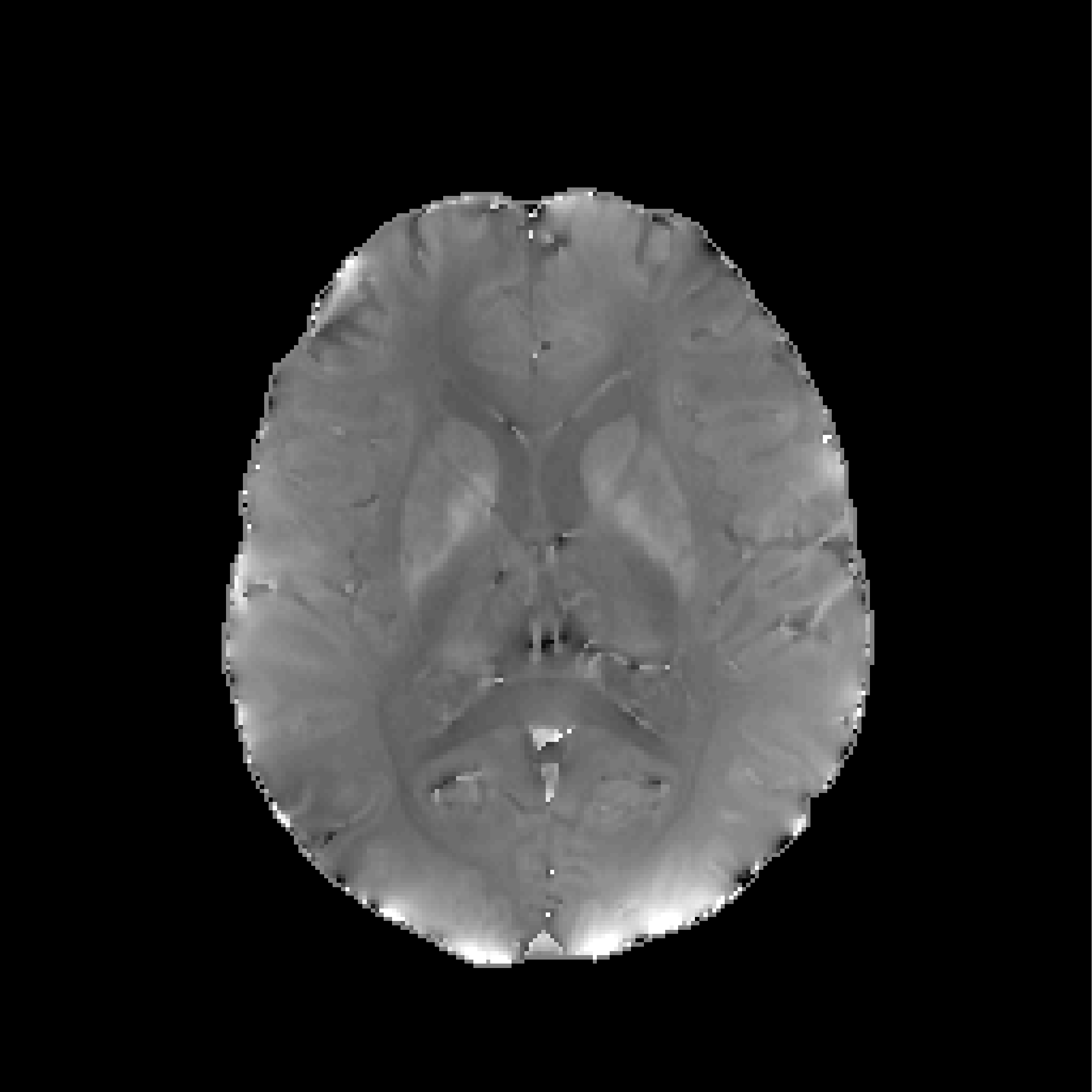}}\vspace{-0.20cm}
\caption{Axial slice images of data sets for the in vivo MR data experiments.}\label{InVivoDataSetAxial}
\end{figure}

The in vivo MR data experiments are conducted using $256\times256\times146$ image with spatial resolution $0.9375\times0.9375\times1\mathrm{mm}^3$ which can be downloaded on Cornell MRI Research Lab webpage. Using the wrapped phase image presented in \cref{RealCornellPhase,RealCornellPhaseAx}, we unwrap the phase using the method in \cite{D.C.Ghiglia1998} to obtain the total field $b$ in \cref{RealCornellTotalField,RealCornellTotalFieldAx} Then the measured local field data $b_l$ in \cref{RealCornellLocalField,RealCornellLocalFieldAx} is obtained by solving the Poisson's equation \cref{LBV} using the method in \cite{D.Zhou2014}.

As in \cref{BrainPhantom}, all regularization based models are initialized with $\chi^0=0$, and both the Frame-HIRE and the TGV-HIRE are also initialized with $v^0=0$. For the parameters, we choose $\hbar=0.1$ for \cref{TKD}, $\eps=0.01$ for \cref{Tikhonov}, $\nu=0.0005$ for the Frame-Int and the Frame-HIRE, $\nu=0.005$ for the Frame-Diff, $\alpha_1=0.00025$ for the TGV-Int and the TGV-HIRE, and $\alpha_1=0.0025$ for the TGV-Diff. In addition, we choose $\beta=0.05$ for all split Bregman algorithms to solve the regularization based models including \cref{Algorithm1}.

\begin{figure}[tp!]
\centering
\hspace{-0.1cm}\subfloat[TKD]{\label{RealCornellTKDFra}\includegraphics[width=3.00cm]{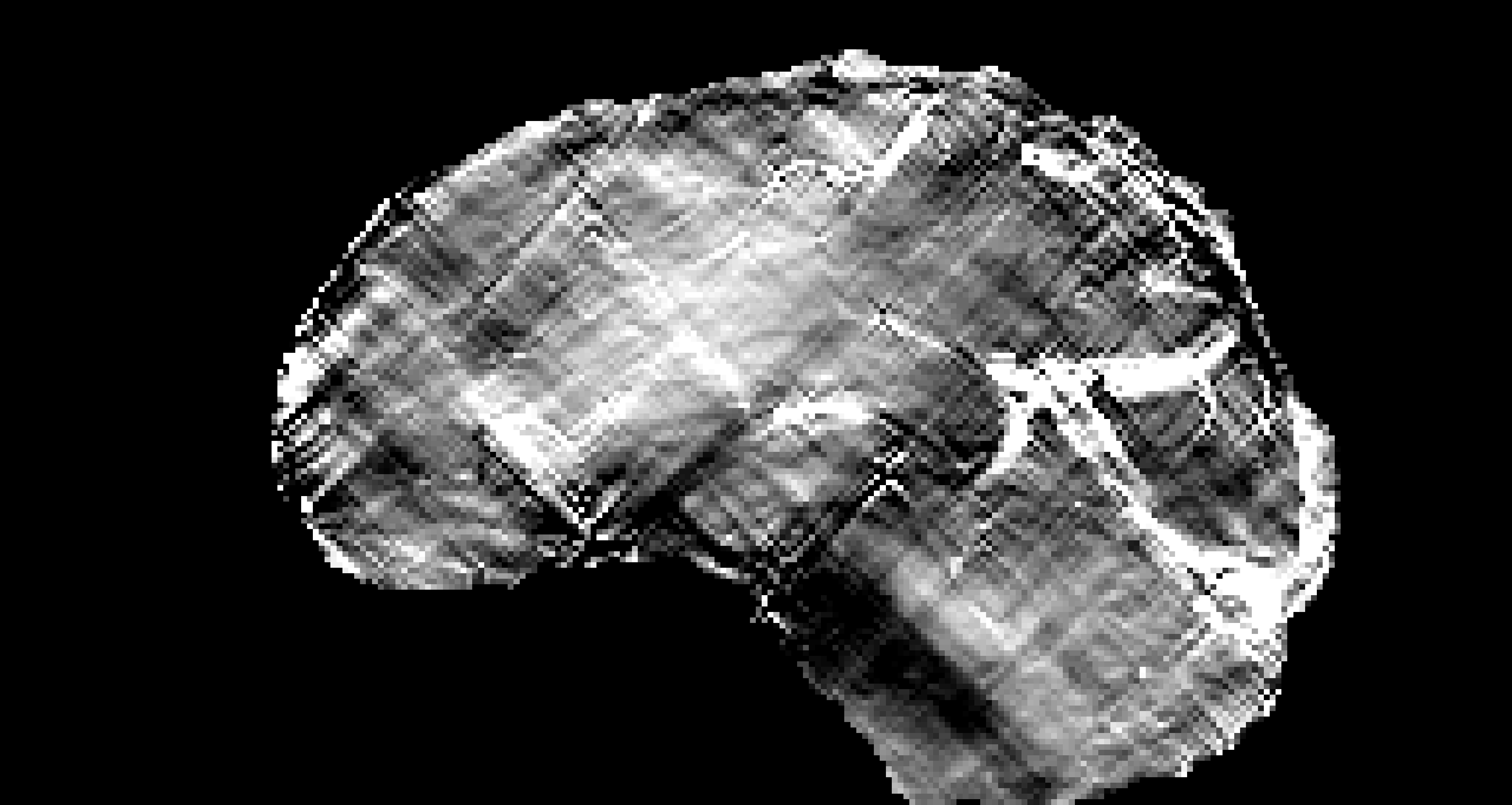}}\hspace{0.005cm}
\subfloat[Tikhonov]{\label{RealCornellTikhonovFra}\includegraphics[width=3.00cm]{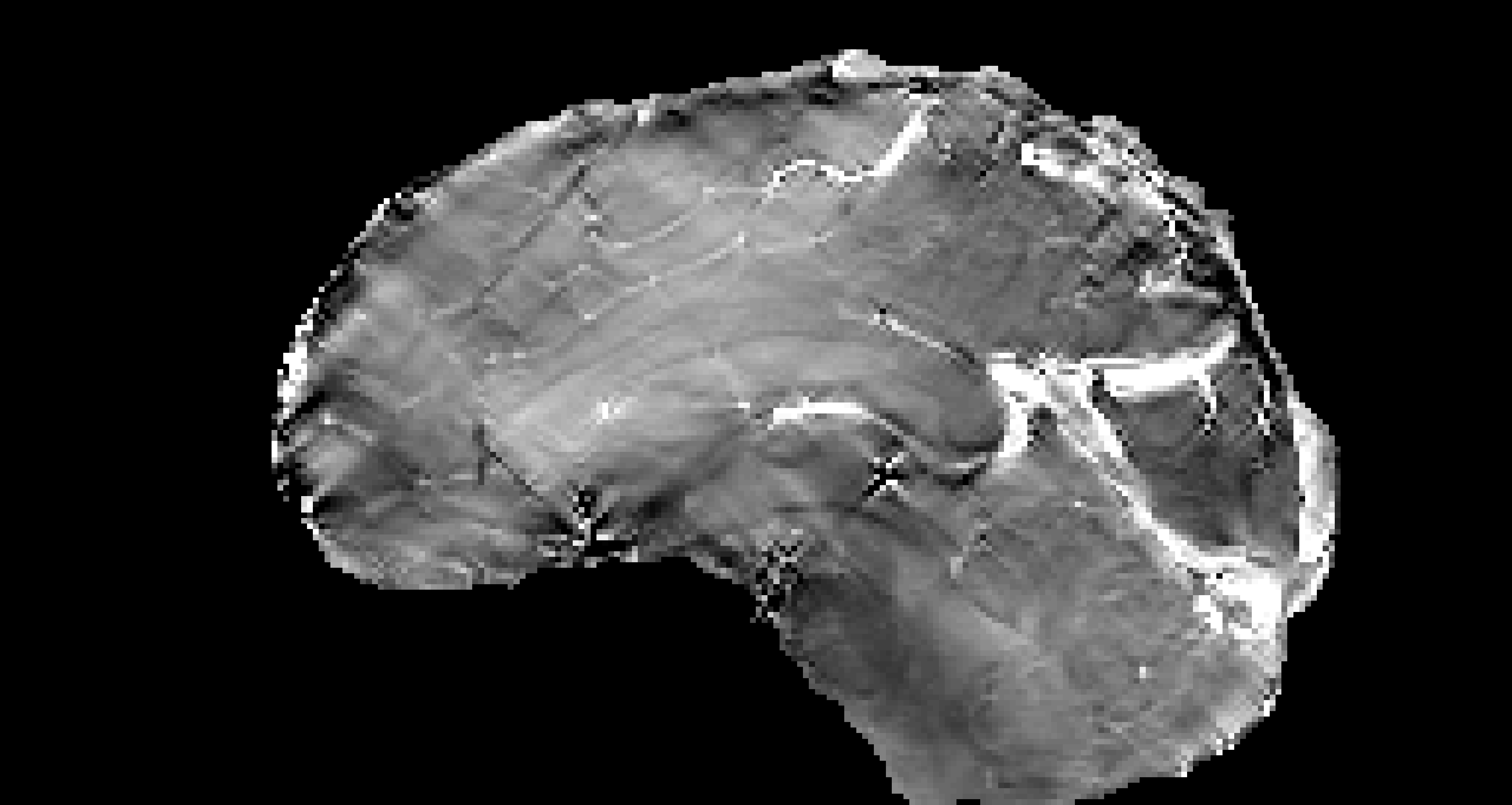}}\hspace{0.005cm}
\subfloat[Frame-Int]{\label{RealCornellFraInt}\includegraphics[width=3.00cm]{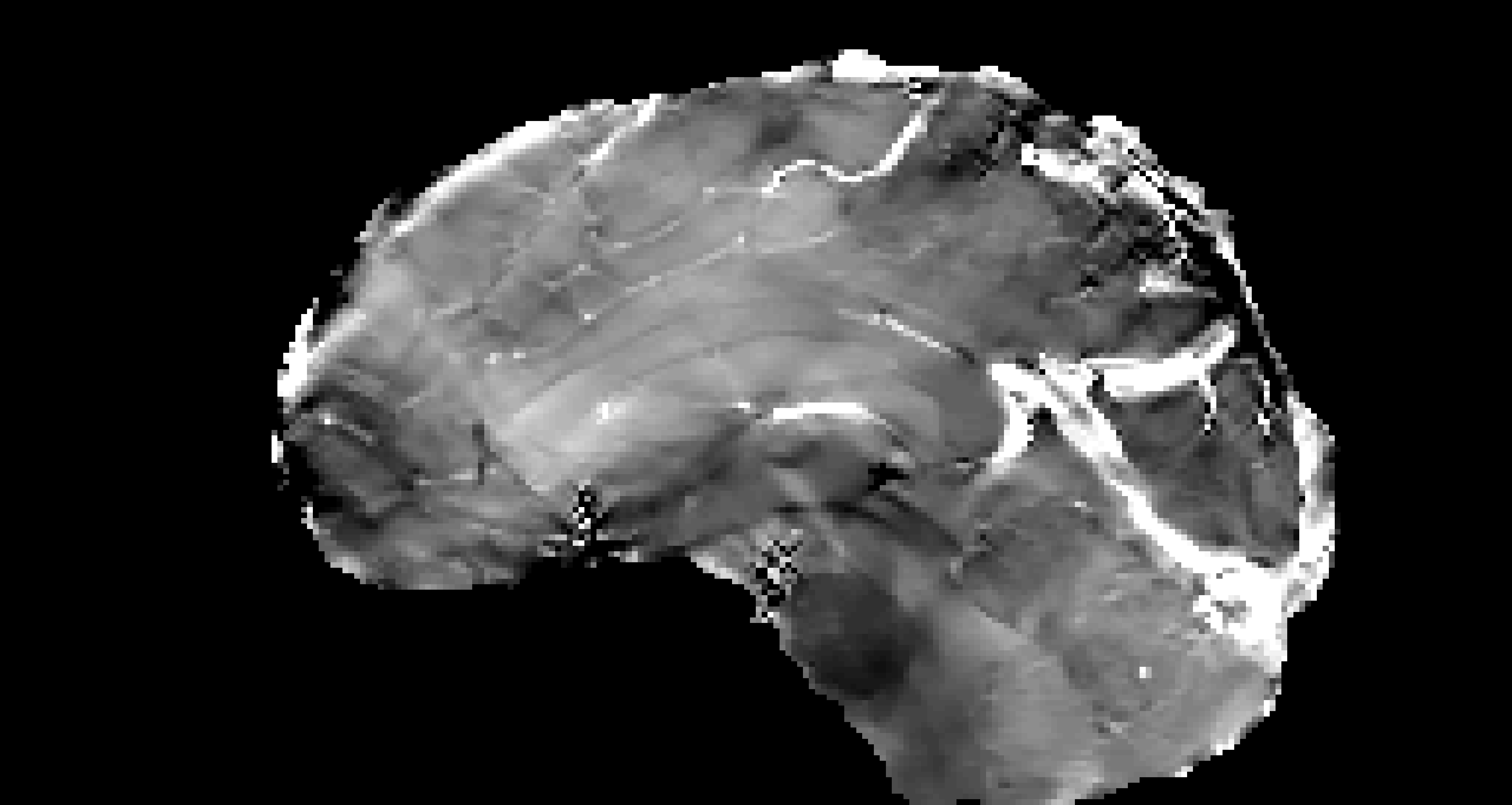}}\hspace{0.005cm}
\subfloat[Frame-Diff]{\label{RealCornellFraDiff}\includegraphics[width=3.00cm]{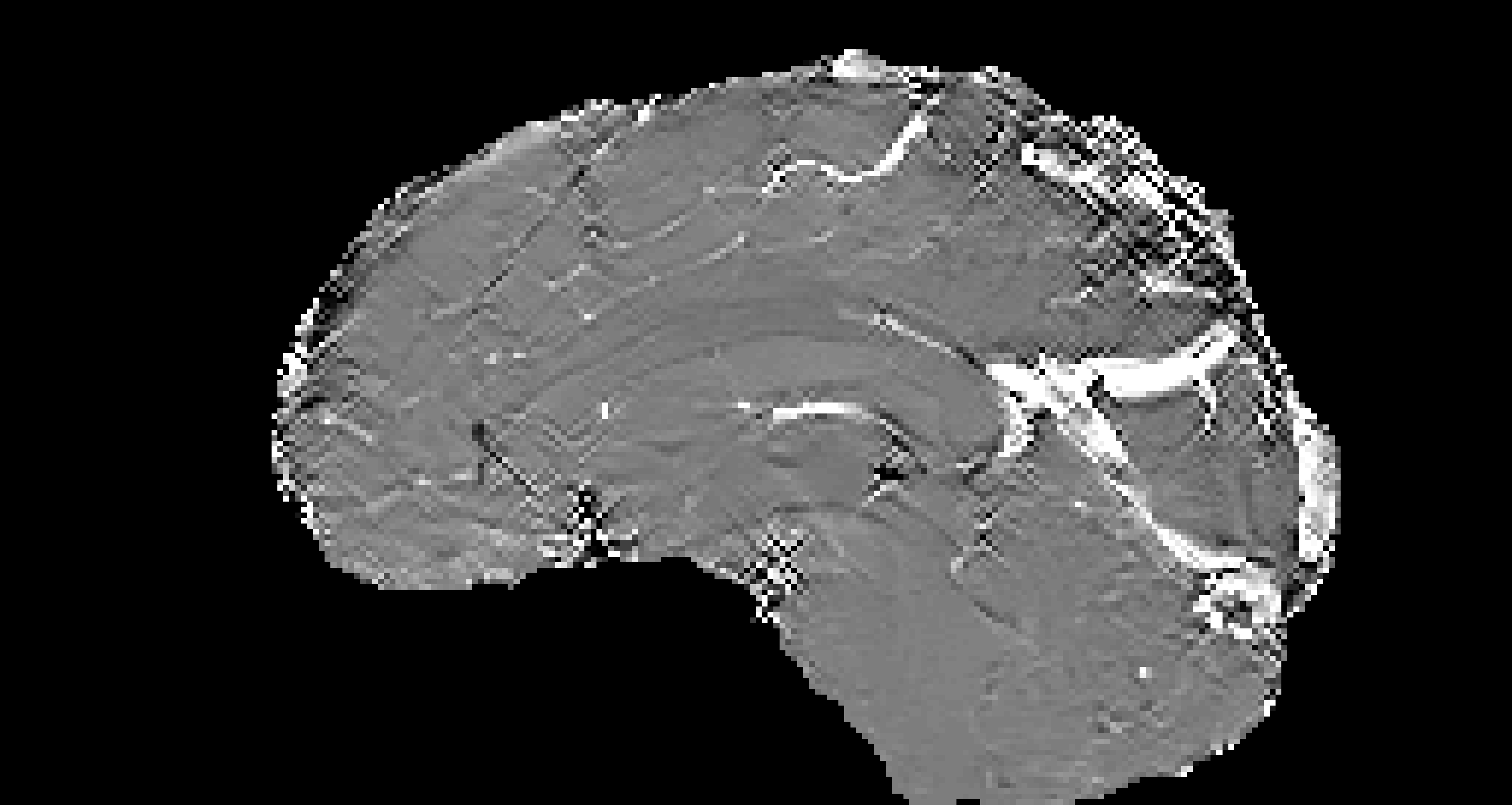}}\hspace{0.005cm}
\subfloat[Frame-HIRE]{\label{RealCornellFraHIRE}\includegraphics[width=3.00cm]{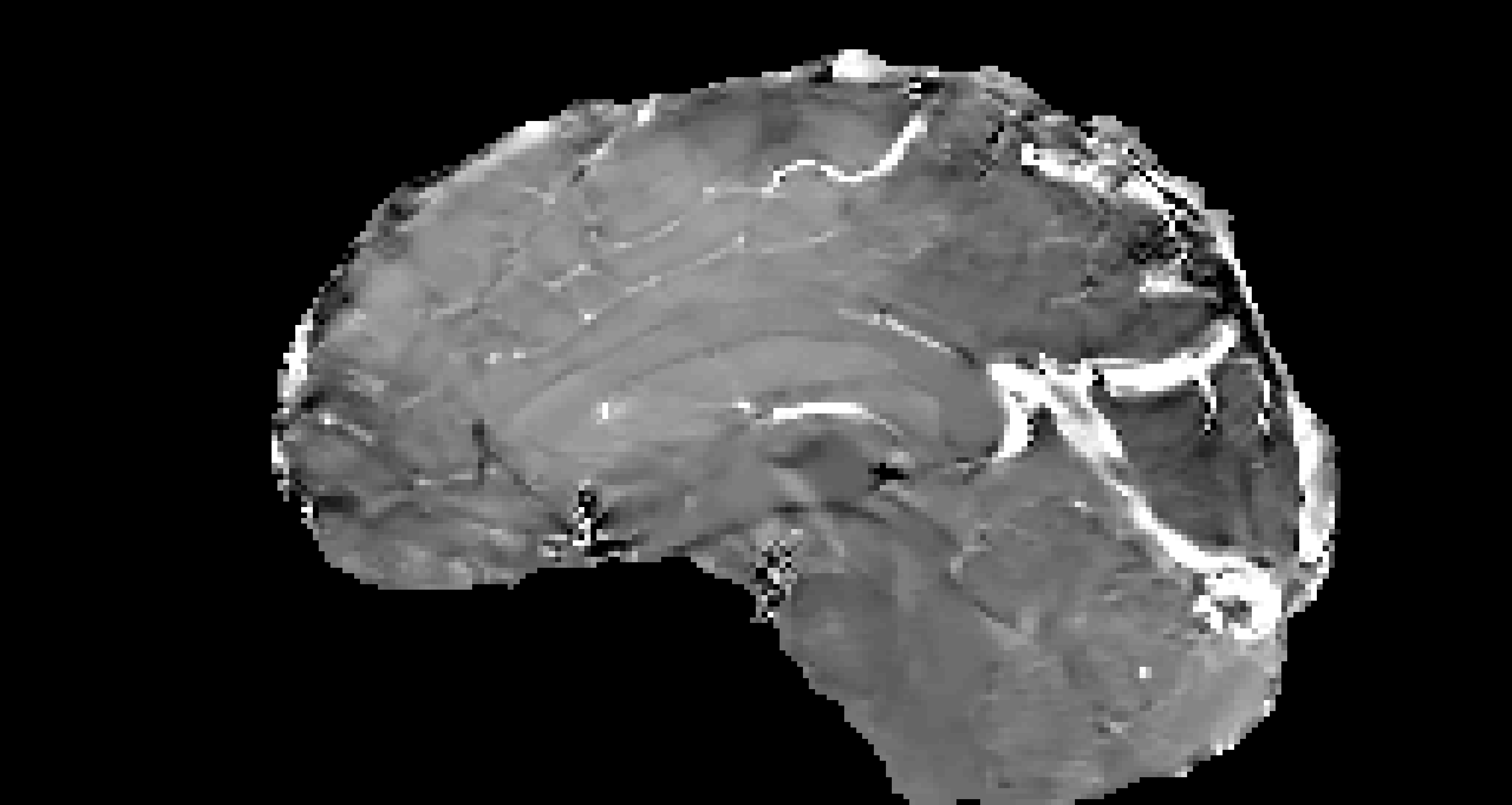}}\vspace{-0.20cm}
\caption{Sagittal slice images comparing QSM reconstruction methods for the in vivo MR data experiments with the wavelet frame regularization. All images of in vivo MR data experimental results are displayed in the window level $[-0.2,0.2]$ for the fair comparison.}\label{RealCornellFra}
\end{figure}

\begin{figure}[tp!]
\centering
\hspace{-0.1cm}\subfloat[TKD]{\label{RealCornellTKDFraAx}\includegraphics[width=3.00cm]{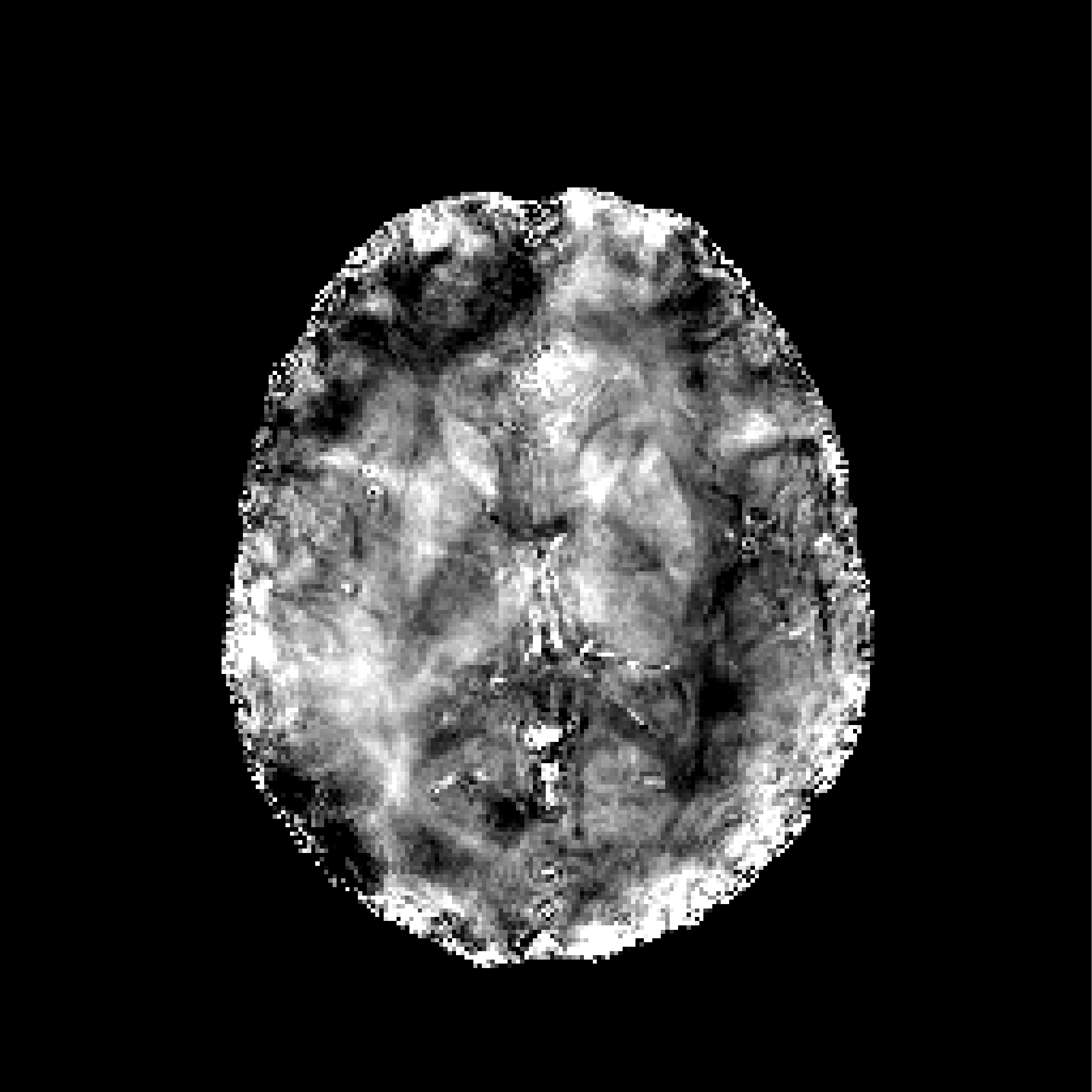}}\hspace{0.005cm}
\subfloat[Tikhonov]{\label{RealCornellTikhonovFraAx}\includegraphics[width=3.00cm]{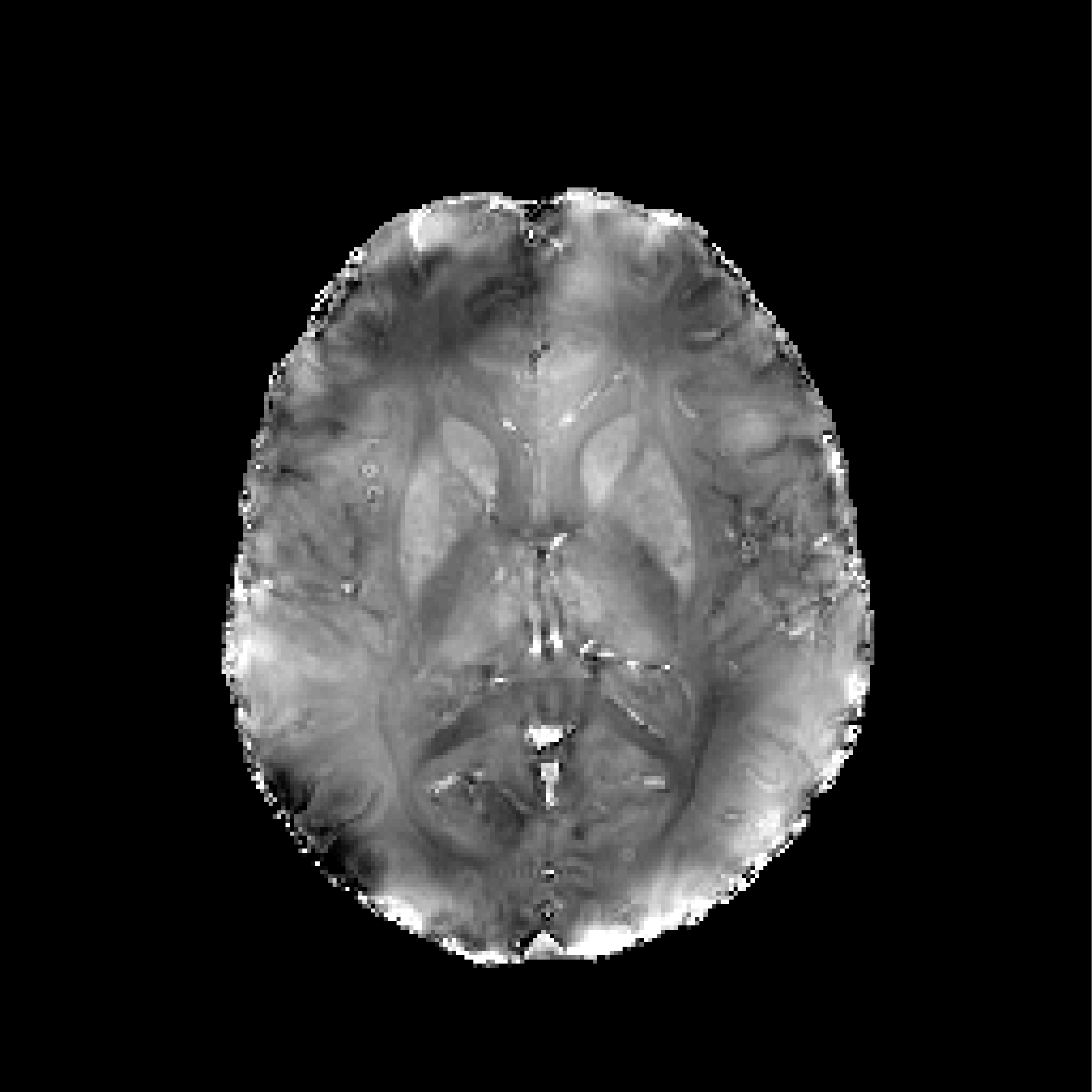}}\hspace{0.005cm}
\subfloat[Frame-Int]{\label{RealCornellFraIntAx}\includegraphics[width=3.00cm]{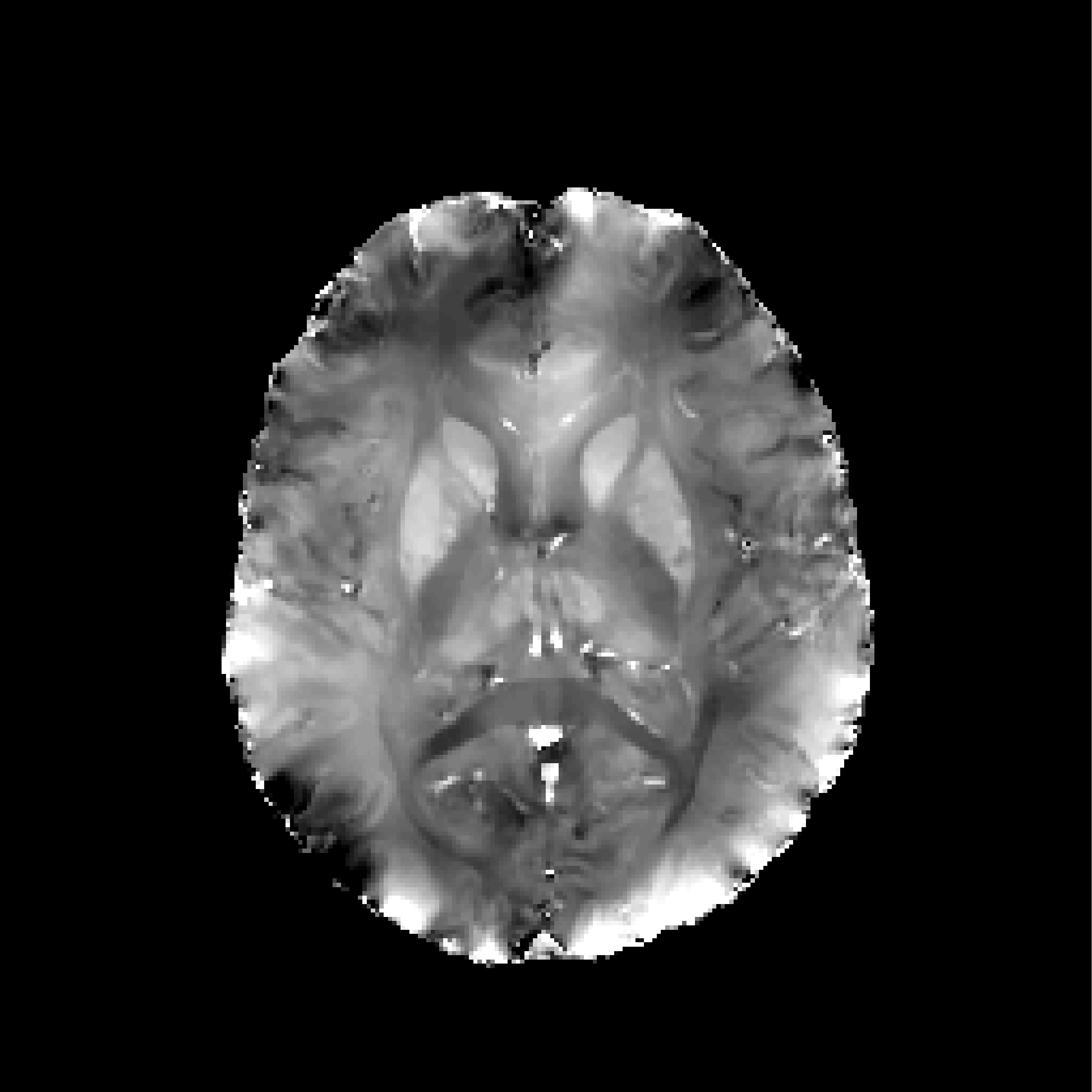}}\hspace{0.005cm}
\subfloat[Frame-Diff]{\label{RealCornellFraDiffAx}\includegraphics[width=3.00cm]{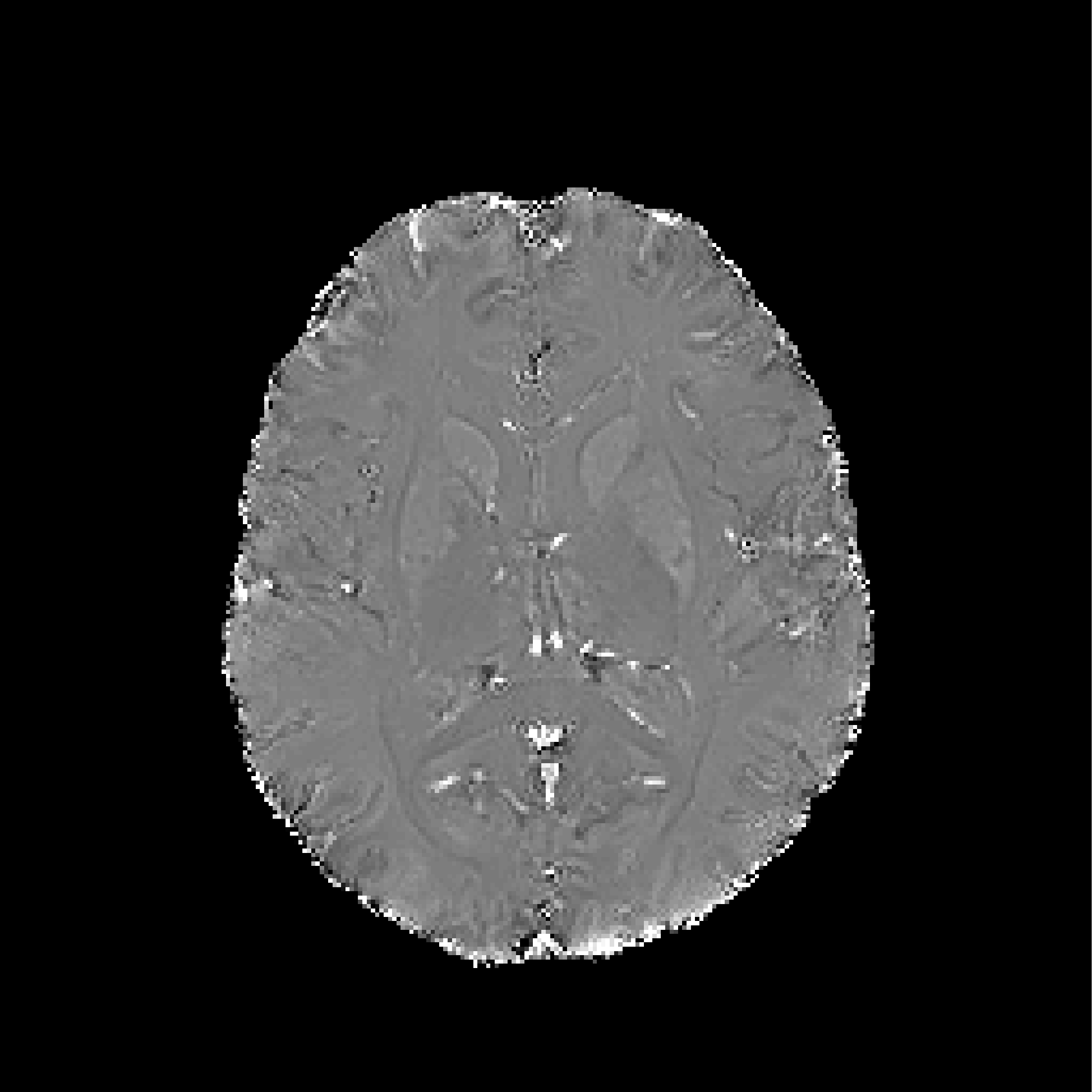}}\hspace{0.005cm}
\subfloat[Frame-HIRE]{\label{RealCornellFraHIREAx}\includegraphics[width=3.00cm]{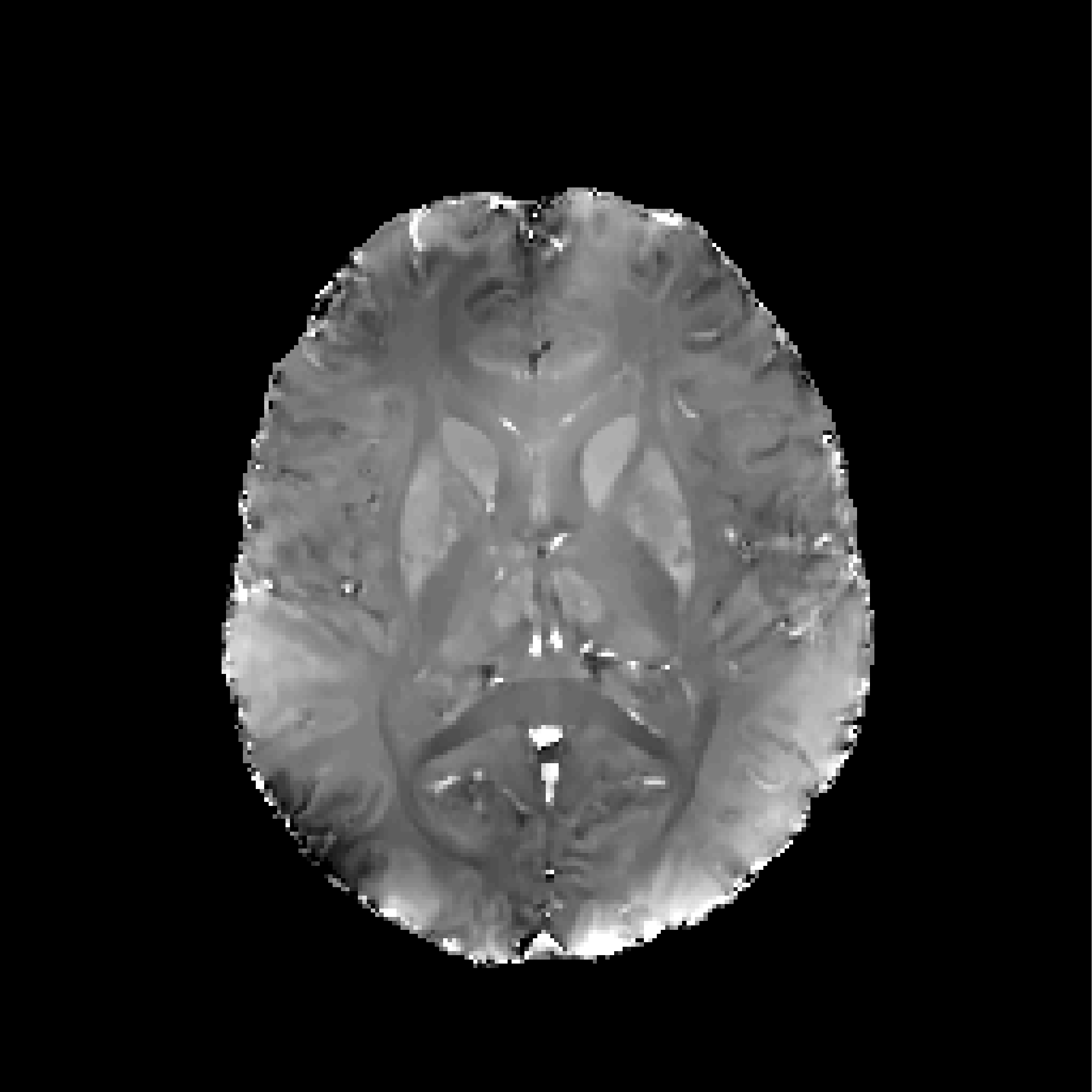}}\vspace{-0.20cm}
\caption{Axial slice images comparing QSM reconstruction methods for the in vivo MR data experiments with the wavelet frame Regularization.}\label{RealCornellFraAxial}
\end{figure}

\begin{figure}[tp!]
\centering
\hspace{-0.1cm}\subfloat[TKD]{\label{RealCornellTKDFraZoom}\includegraphics[width=3.00cm]{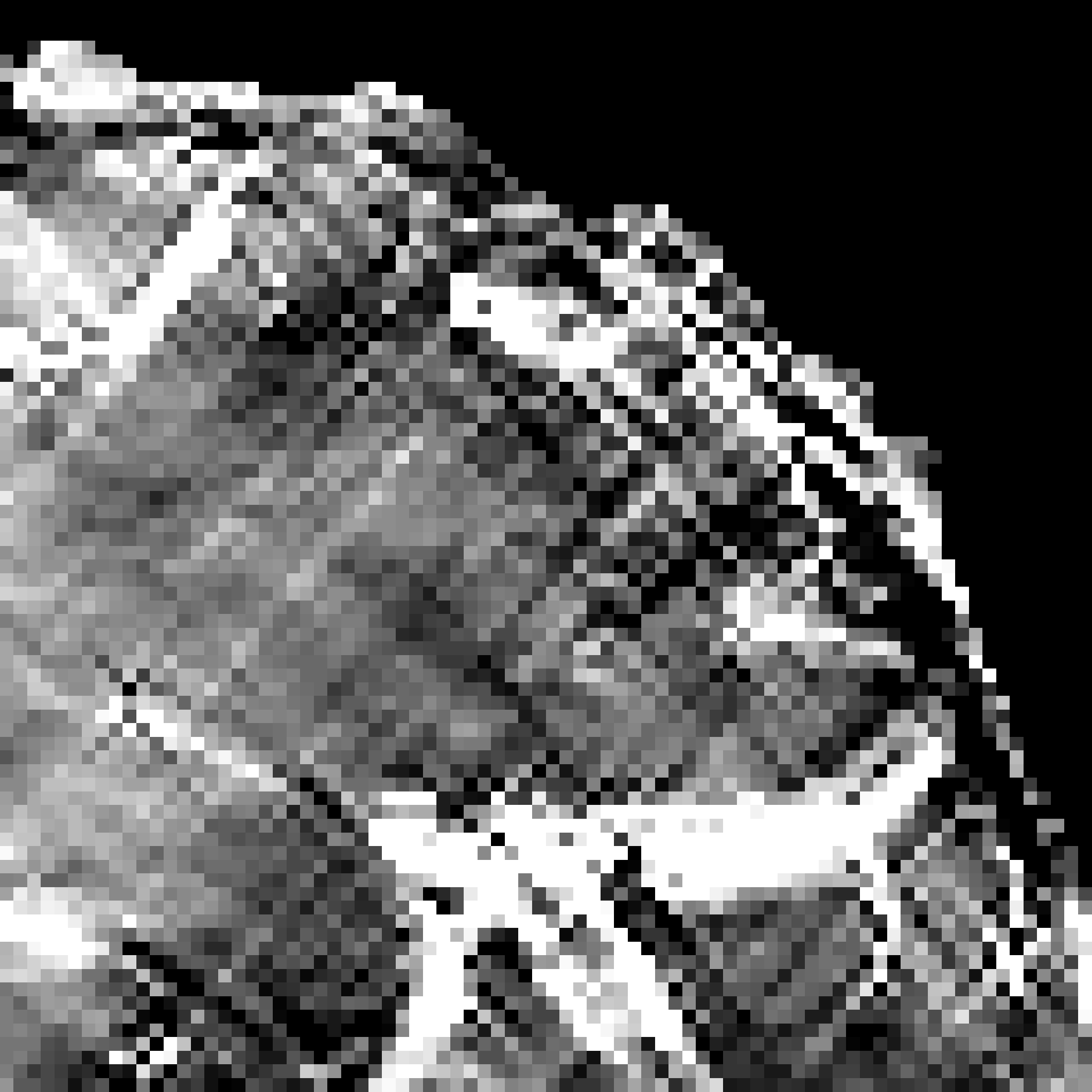}}\hspace{0.005cm}
\subfloat[Tikhonov]{\label{RealCornellTikhonovFraZoom}\includegraphics[width=3.00cm]{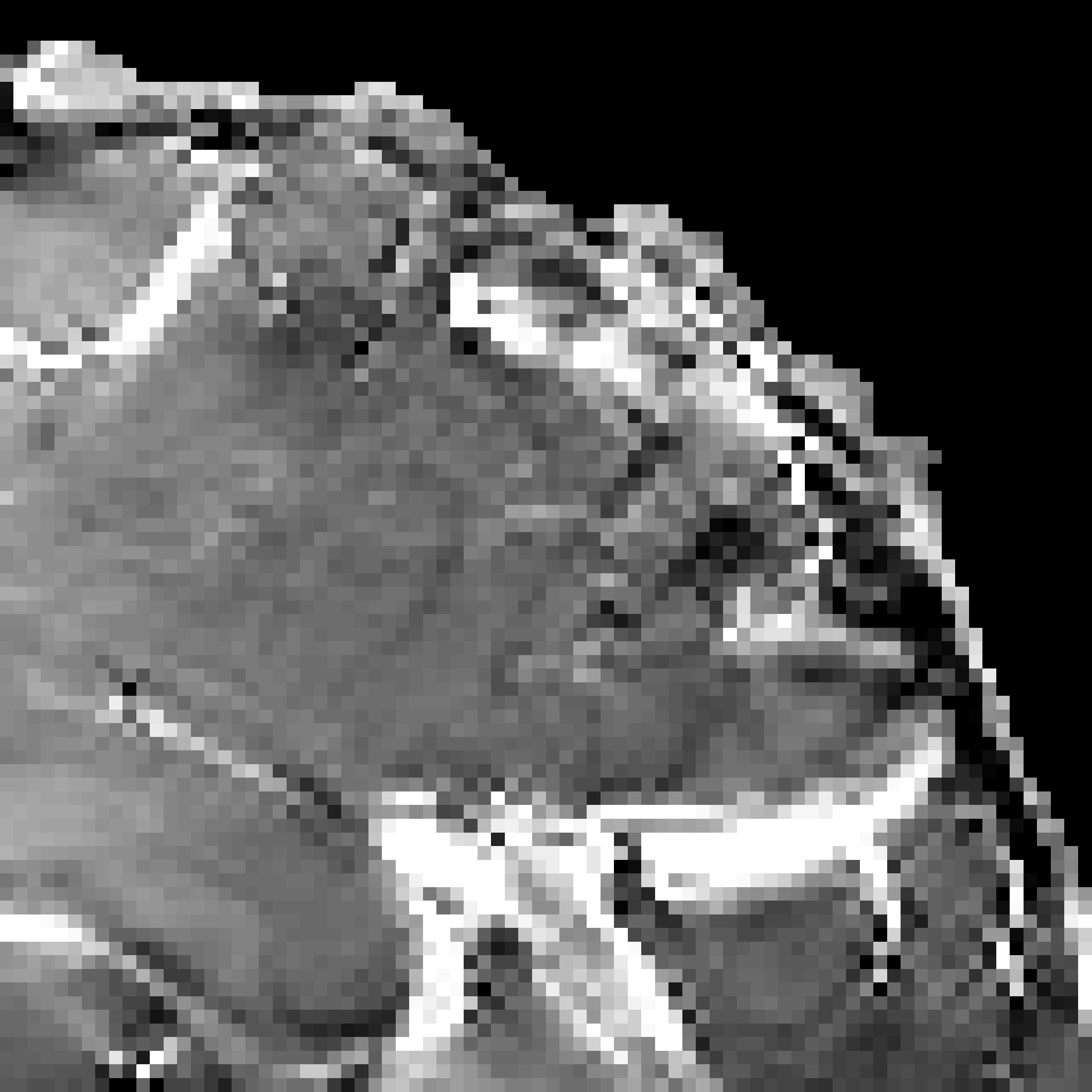}}\hspace{0.005cm}
\subfloat[Frame-Int]{\label{RealCornellFraIntZoom}\includegraphics[width=3.00cm]{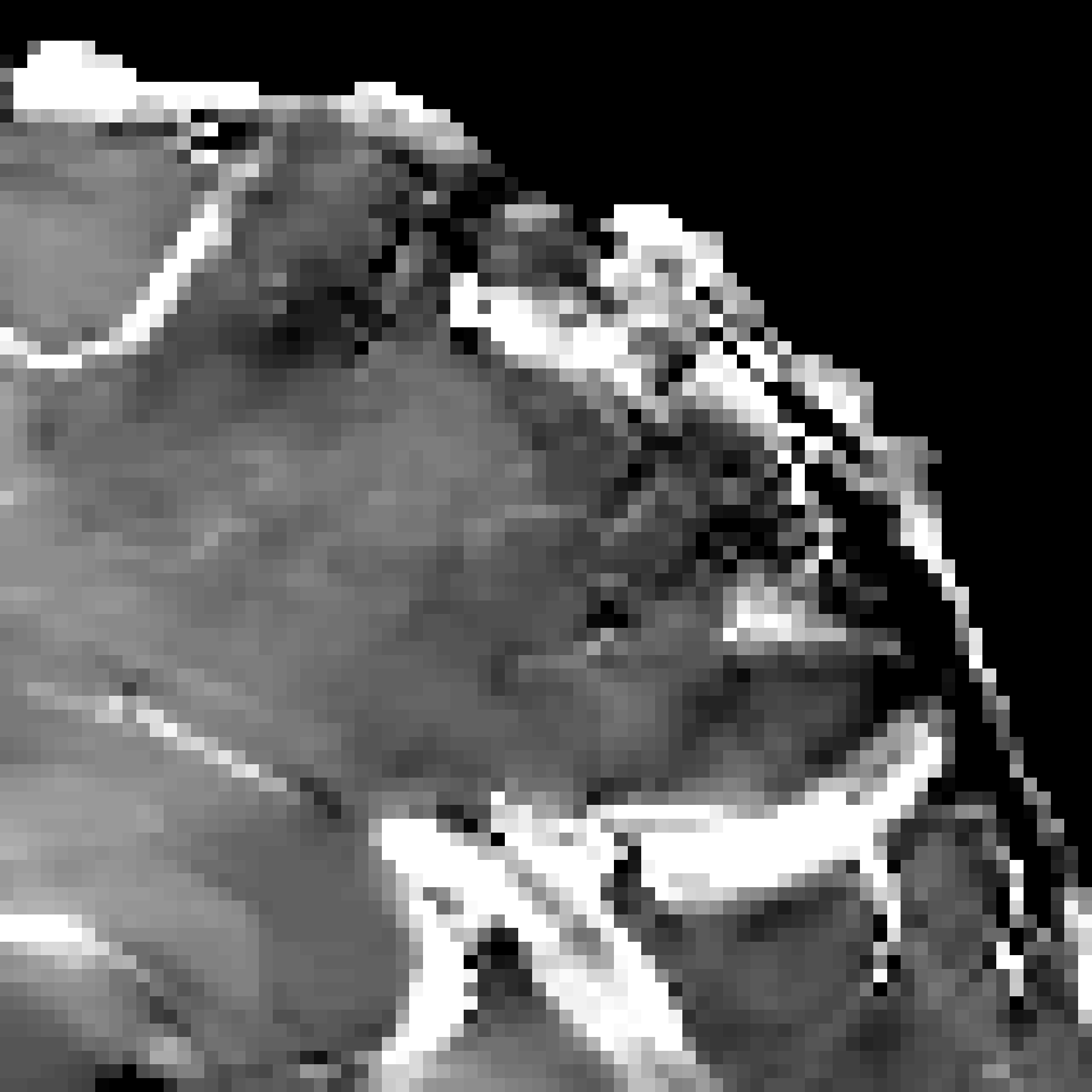}}\hspace{0.005cm}
\subfloat[Frame-Diff]{\label{RealCornellFraDiffZoom}\includegraphics[width=3.00cm]{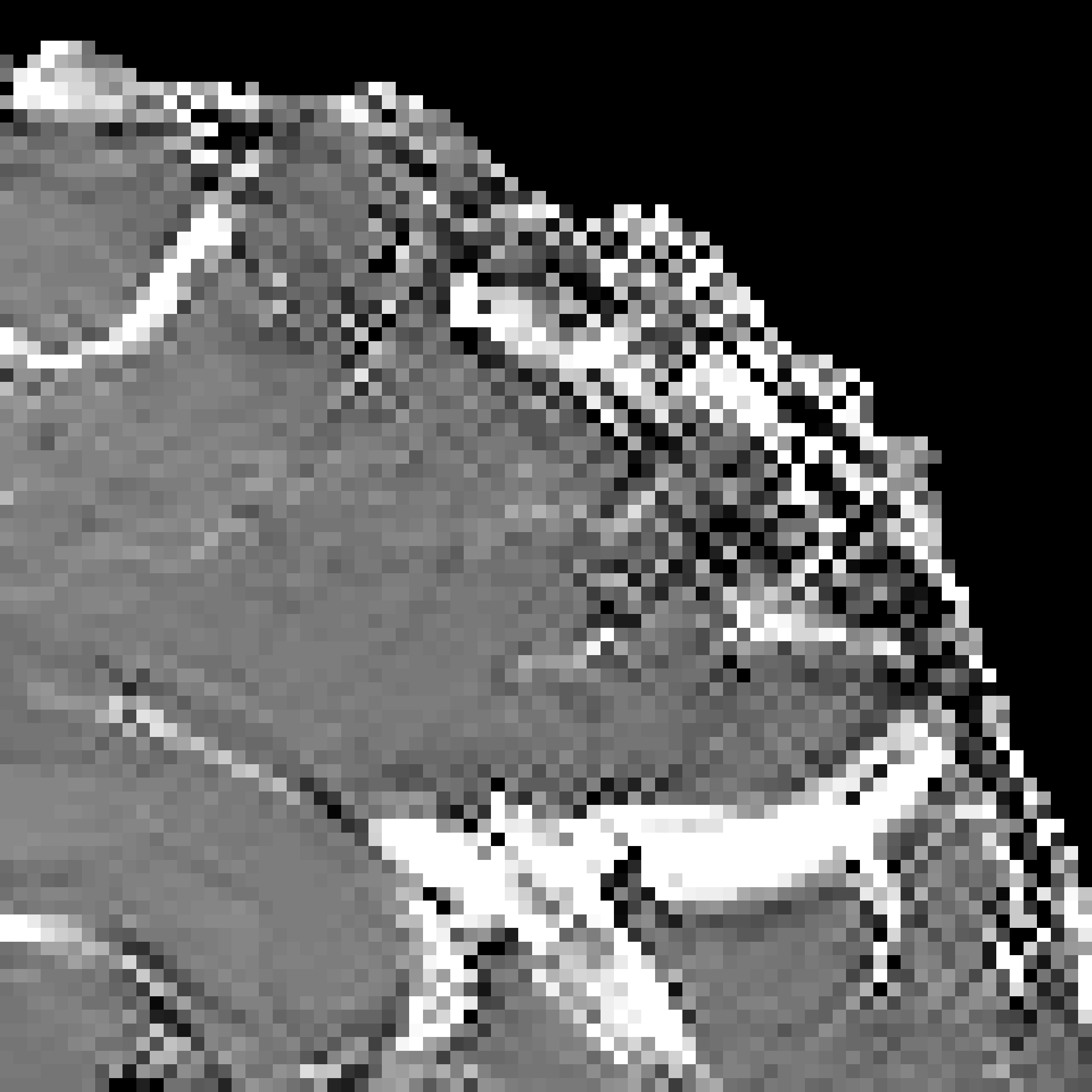}}\hspace{0.005cm}
\subfloat[Frame-HIRE]{\label{RealCornellFraHIREZoom}\includegraphics[width=3.00cm]{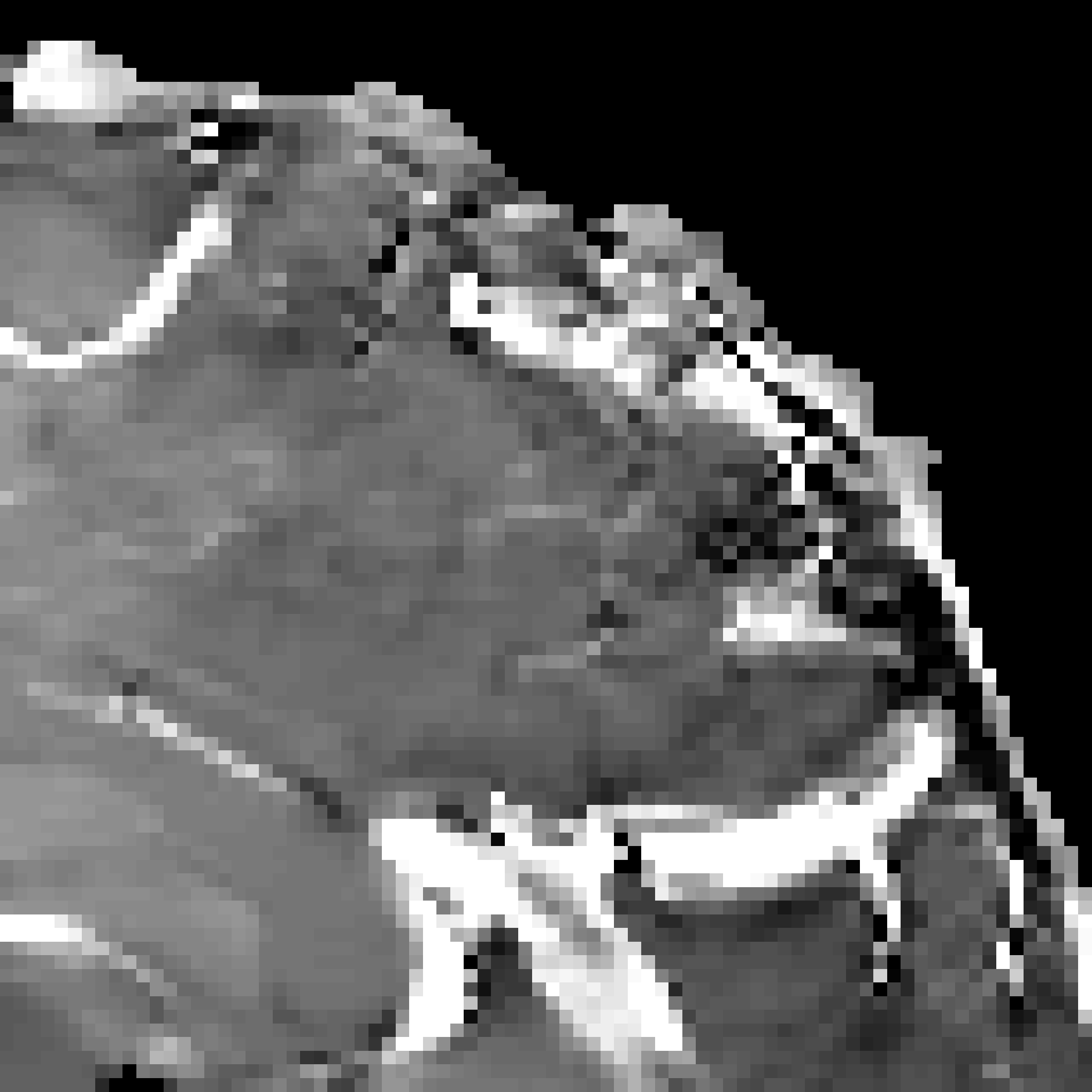}}\vspace{-0.20cm}
\caption{Zoom-in views of \cref{RealCornellFra}.}\label{RealCornellFraZoom}
\end{figure}

\begin{figure}[tp!]
\centering
\hspace{-0.1cm}\subfloat[TKD]{\label{RealCornellTKDTGV}\includegraphics[width=3.00cm]{RealCornellTKD.pdf}}\hspace{0.005cm}
\subfloat[Tikhonov]{\label{RealCornellTikhonovTGV}\includegraphics[width=3.00cm]{RealCornellTikhonov.pdf}}\hspace{0.005cm}
\subfloat[TGV-Int]{\label{RealCornellTGVInt}\includegraphics[width=3.00cm]{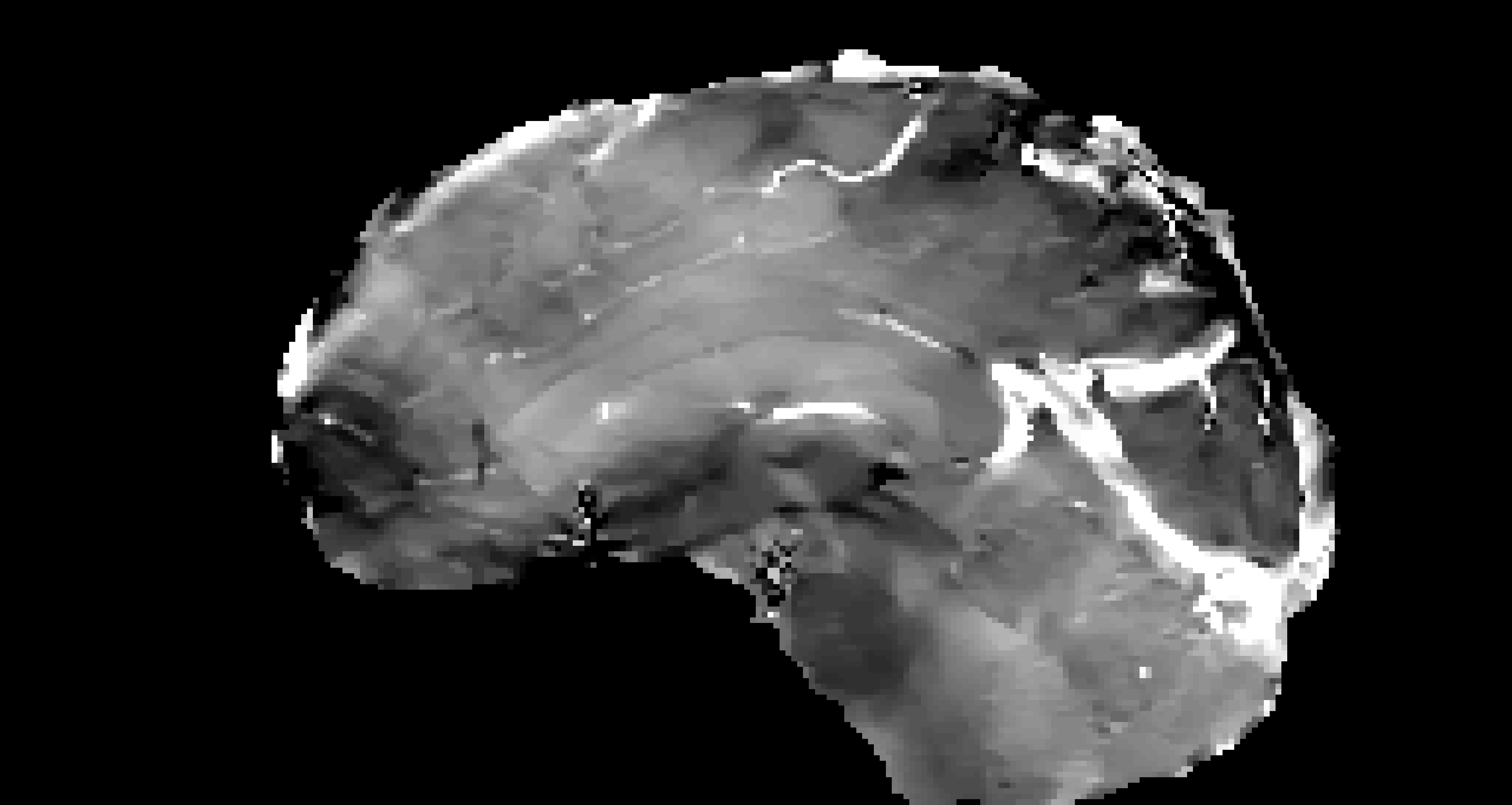}}\hspace{0.005cm}
\subfloat[TGV-Diff]{\label{RealCornellTGVDiff}\includegraphics[width=3.00cm]{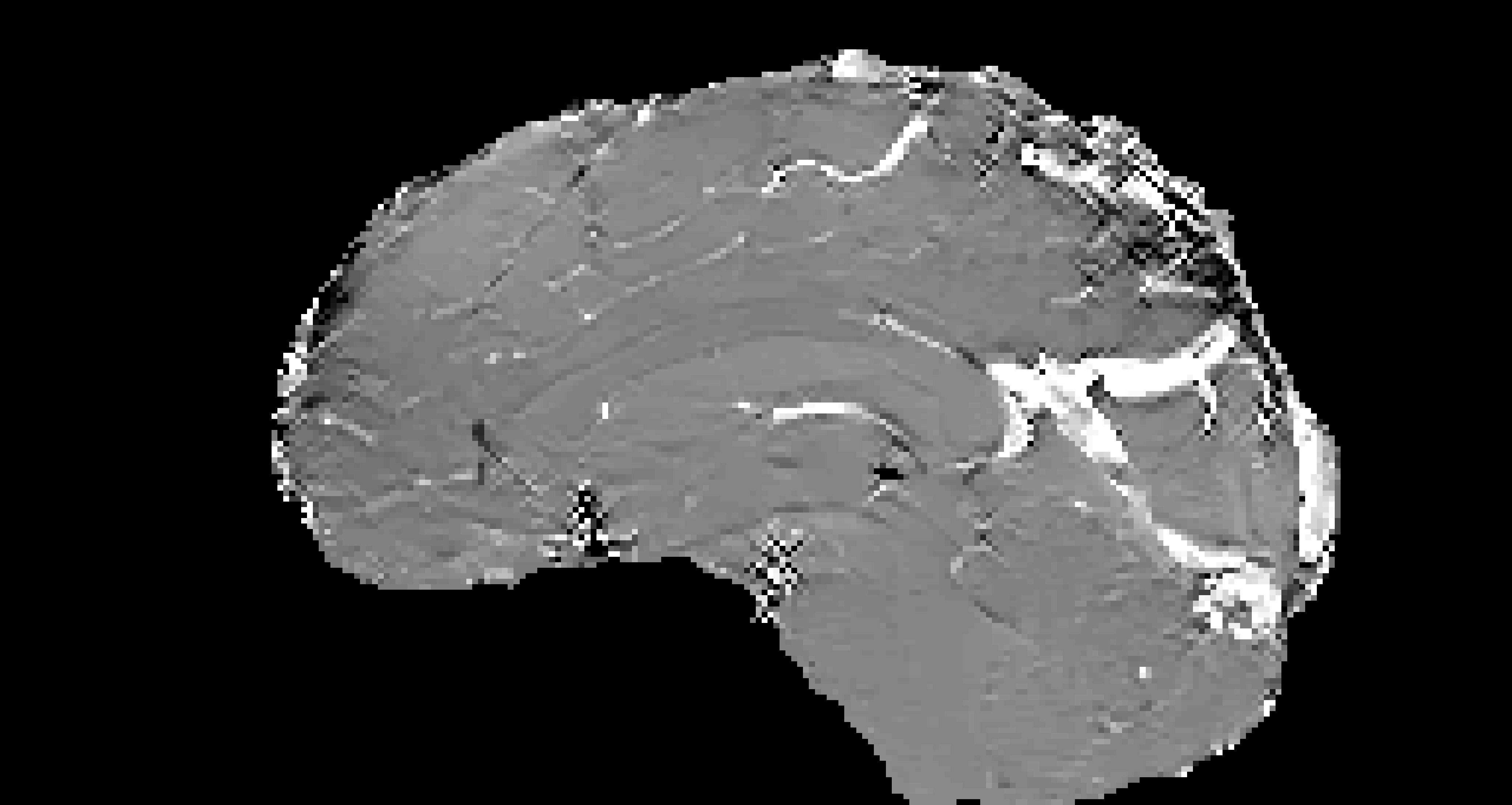}}\hspace{0.005cm}
\subfloat[TGV-HIRE]{\label{RealCornellTGVHIRE}\includegraphics[width=3.00cm]{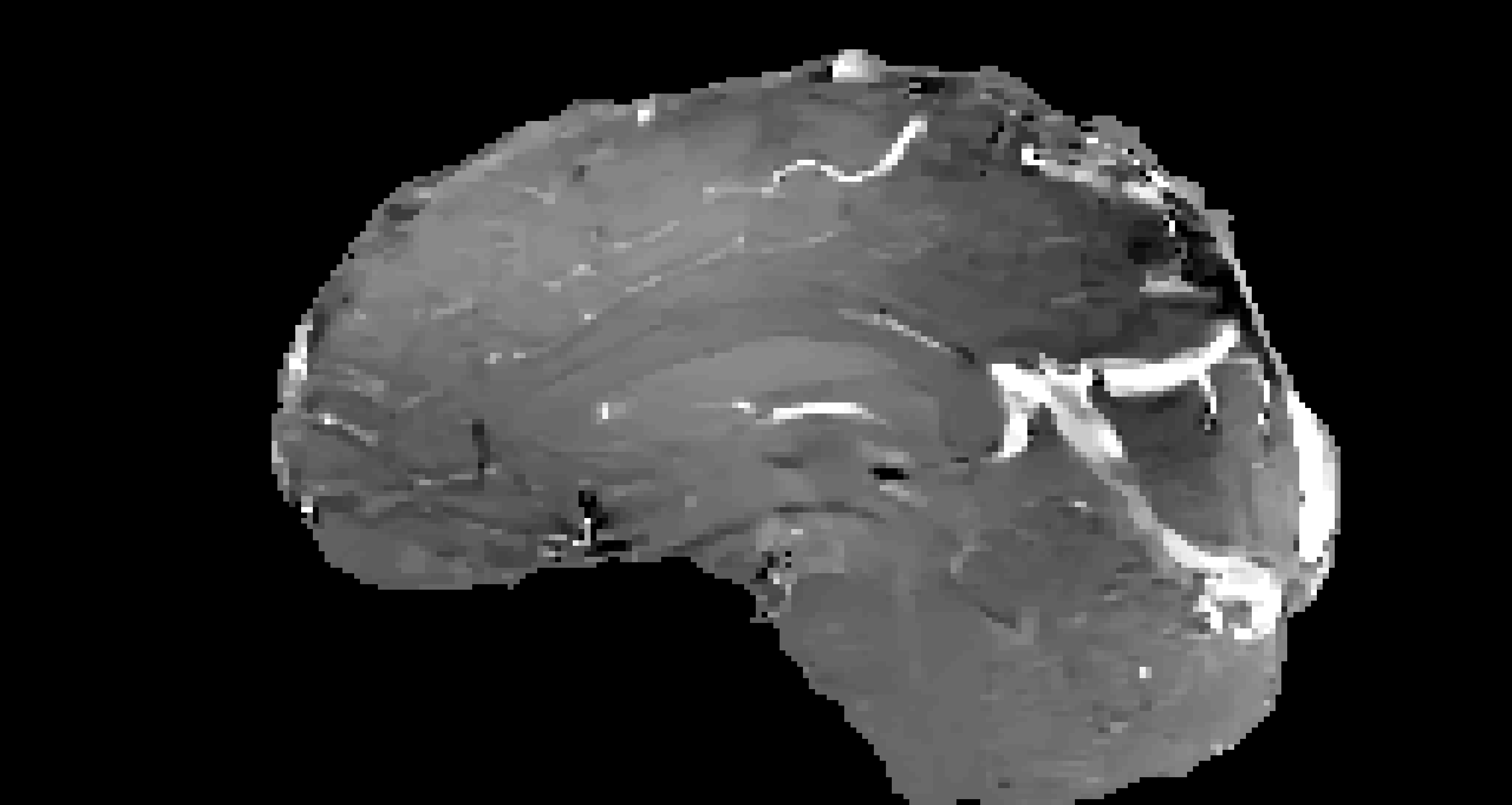}}\vspace{-0.20cm}
\caption{Sagittal slice images comparing QSM reconstruction methods for the in vivo MR data experiments with the TGV regularization. All images of in vivo MR data experimental results are displayed in the window level $[-0.2,0.2]$ for the fair comparison.}\label{RealCornellTGV}
\end{figure}

\begin{figure}[tp!]
\centering
\hspace{-0.1cm}\subfloat[TKD]{\label{RealCornellTKDTGVAx}\includegraphics[width=3.00cm]{RealCornellTKDAxial.pdf}}\hspace{0.005cm}
\subfloat[Tikhonov]{\label{RealCornellTikhonovTGVAx}\includegraphics[width=3.00cm]{RealCornellTikhonovAxial.pdf}}\hspace{0.005cm}
\subfloat[TGV-Int]{\label{RealCornellTGVIntAx}\includegraphics[width=3.00cm]{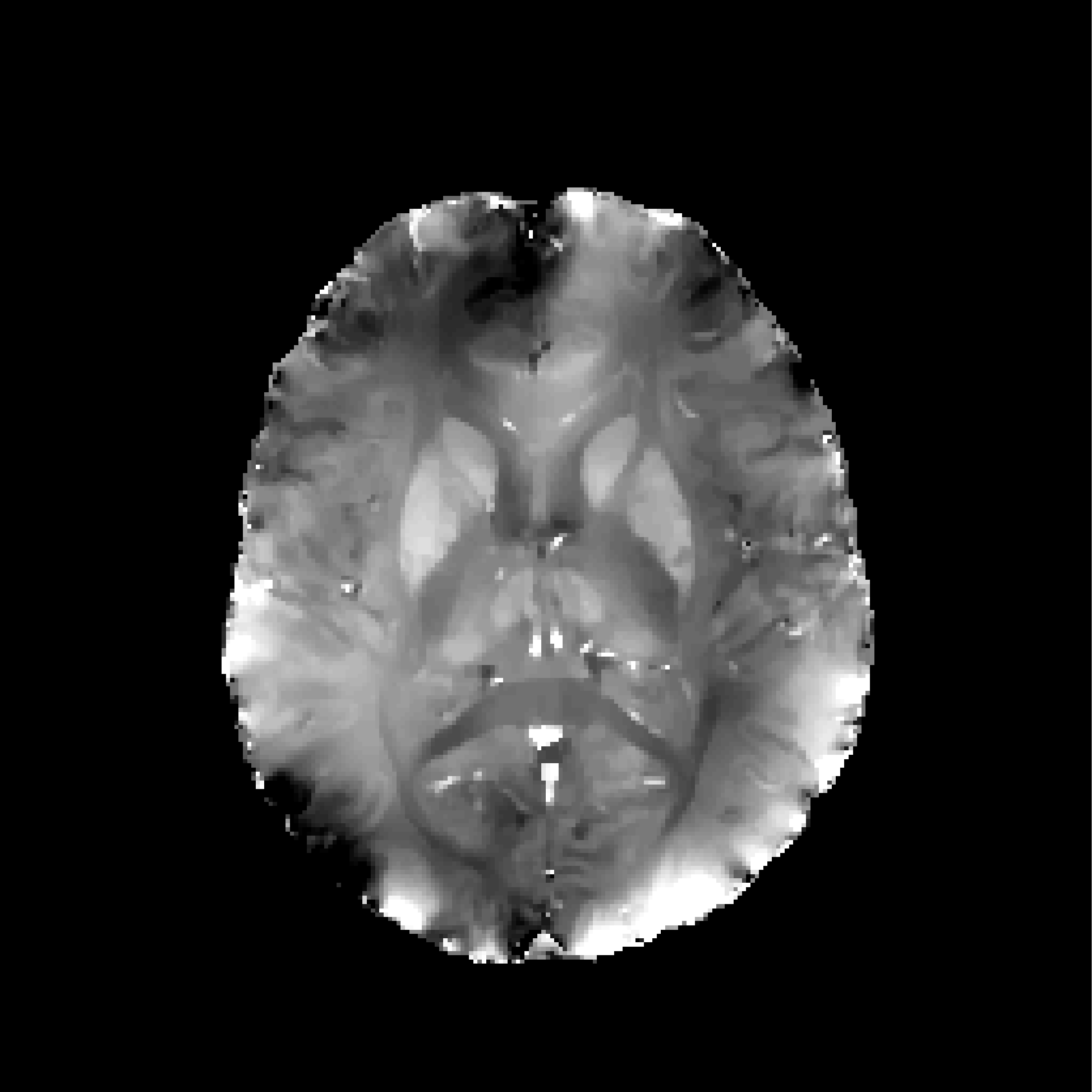}}\hspace{0.005cm}
\subfloat[TGV-Diff]{\label{RealCornellTGVDiffAx}\includegraphics[width=3.00cm]{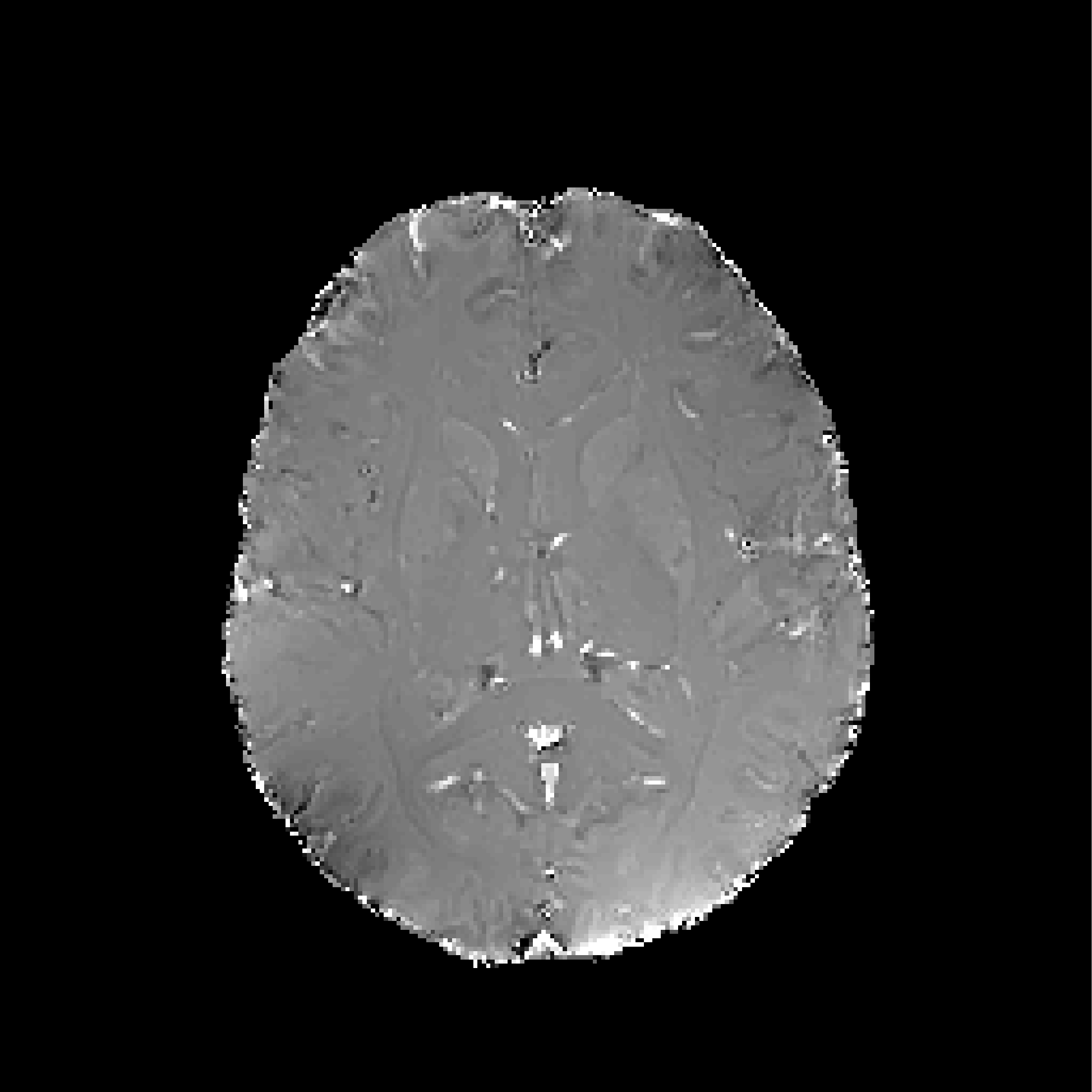}}\hspace{0.005cm}
\subfloat[TGV-HIRE]{\label{RealCornellTGVHIREAx}\includegraphics[width=3.00cm]{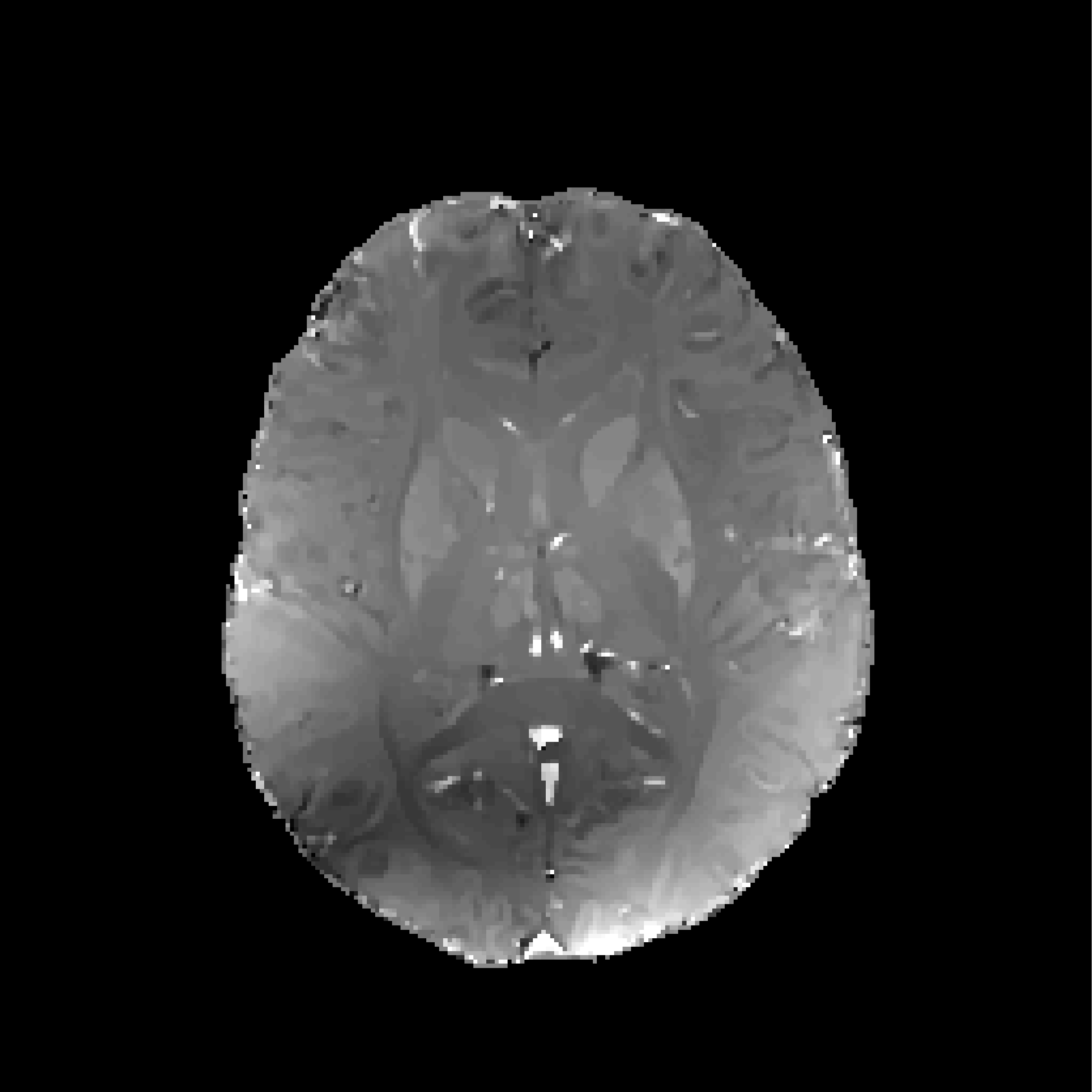}}\vspace{-0.20cm}
\caption{Axial slice images comparing QSM reconstruction methods for the in vivo MR data experiments with the TGV Regularization.}\label{RealCornellTGVAxial}
\end{figure}

\begin{figure}[tp!]
\centering
\hspace{-0.1cm}\subfloat[TKD]{\label{RealCornellTKDTGVZoom}\includegraphics[width=3.00cm]{RealCornellTKDZoom.pdf}}\hspace{0.005cm}
\subfloat[Tikhonov]{\label{RealCornellTikhonovTGVZoom}\includegraphics[width=3.00cm]{RealCornellTikhonovZoom.pdf}}\hspace{0.005cm}
\subfloat[TGV-Int]{\label{RealCornellTGVIntZoom}\includegraphics[width=3.00cm]{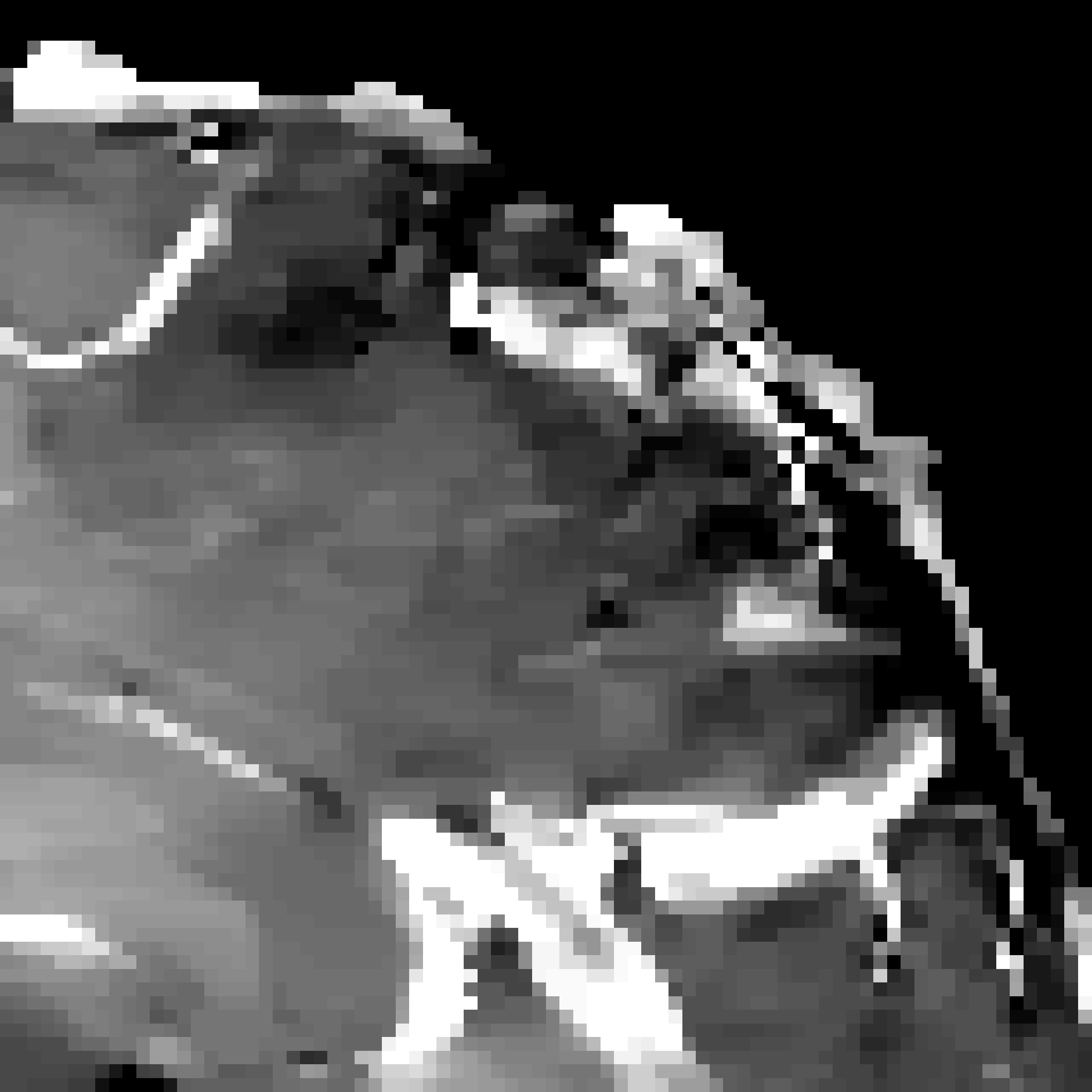}}\hspace{0.005cm}
\subfloat[TGV-Diff]{\label{RealCornellTGVDiffZoom}\includegraphics[width=3.00cm]{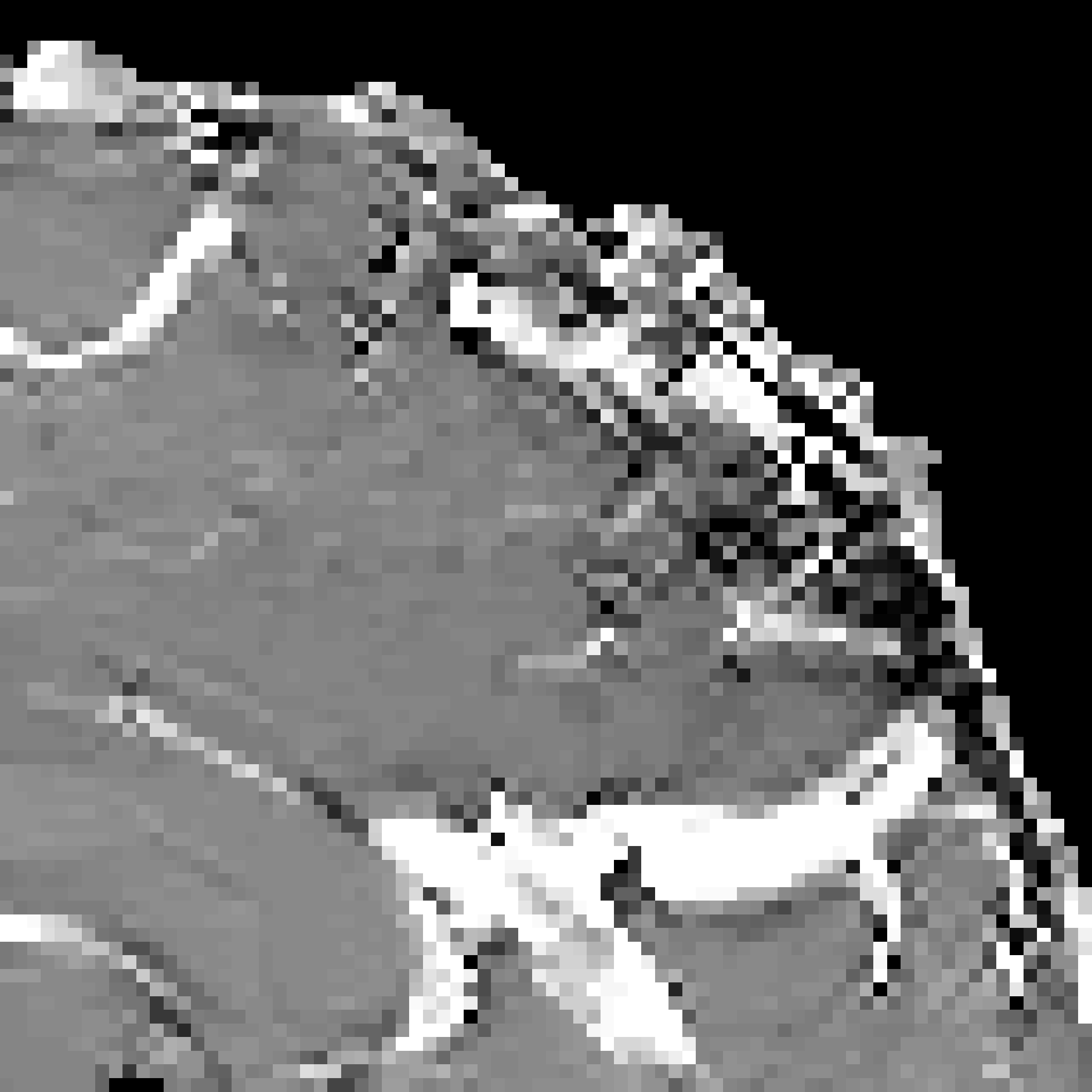}}\hspace{0.005cm}
\subfloat[TGV-HIRE]{\label{RealCornellTGVHIREZoom}\includegraphics[width=3.00cm]{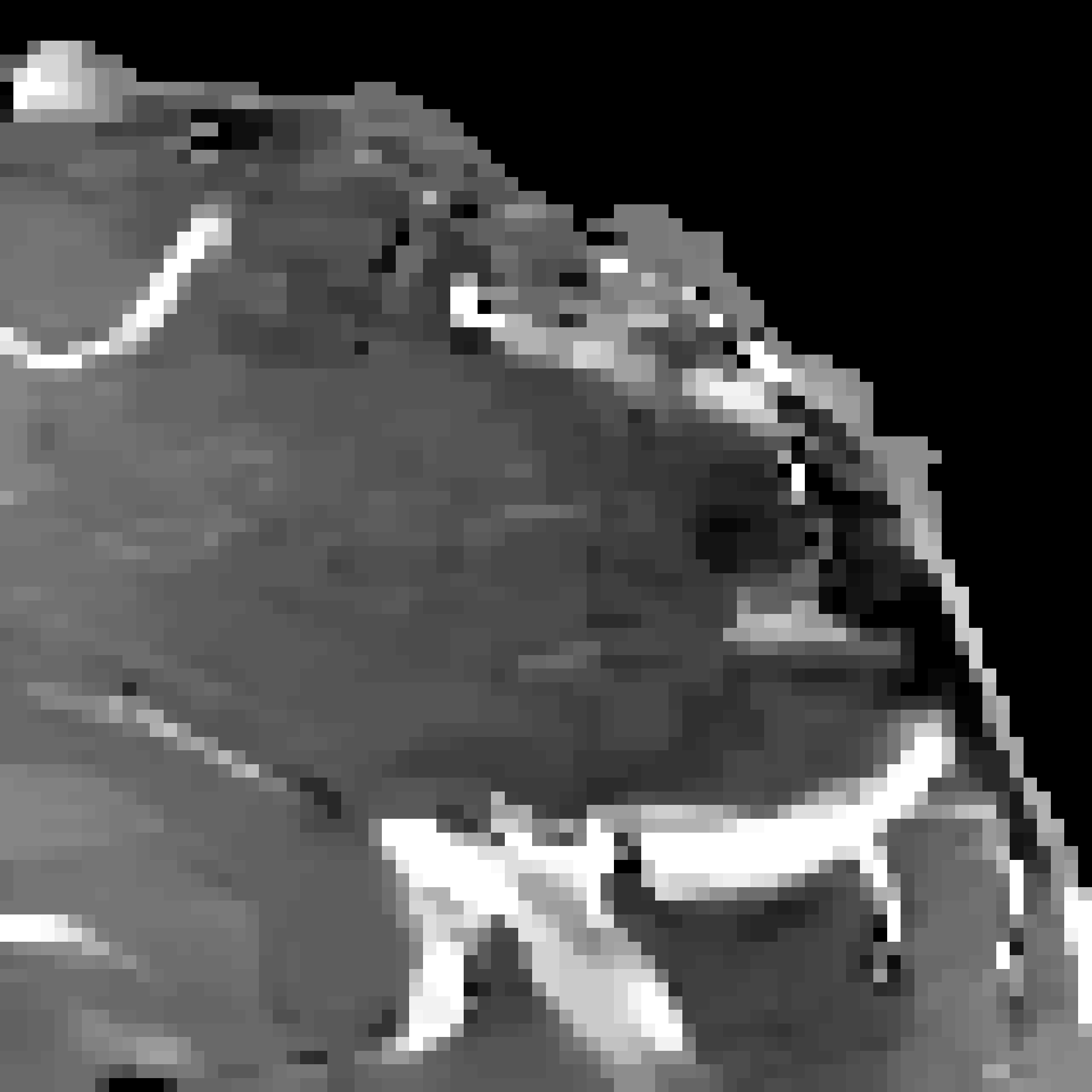}}\vspace{-0.20cm}
\caption{Zoom-in views of \cref{RealCornellTGV}.}\label{RealCornellTGVZoom}
\end{figure}

\cref{RealCornellFra,RealCornellFraAxial} display the visual comparisons of the direct approaches and the wavelet frame regularization approaches, and the zoom-in views of \cref{RealCornellFra} are provided in \cref{RealCornellFraZoom}. We also provide the visual comparisons of the direct approaches and the TGV regularization approaches in \cref{RealCornellTGV,RealCornellTGVAxial,RealCornellTGVZoom}. Since the reference image is not available for in vivo MR data, it is in general more difficult to provide quantitative evaluations than the numerical brain phantom. Nonetheless, we can see from the viewpoint of visual comparison that the pros and cons are almost the same as the numerical brain phantom experiments. It is also worth noting that the HIRE models can reduce the streaking artifacts which propagate from $\p\Om$ into $\Om$ as well as the shadow artifacts. As pointed out in \cite{Y.Wang2015}, the in vivo local field data is prone to the outliers near $\p\Om$ because the GRE signal lacks information outside $\Om$. Hence, we can see that most streaking artifacts propagate from these outliers near $\p\Om$ into the ROI. However, thanks to the sparsity promoting property of $\ell_1$ norm, the term $\lambda\left\|\msL v\right\|_1$ in the HIRE approaches can somehow capture and remove them, leading to the suppression of artifacts propagating from $\p\Om$ into $\Om$ as well as the shadow and streaking artifact removal. Finally, even though we can also note that the Tikhonov regularization can somehow reduce the artifacts, there are some losses of features due to the smoothness prior of the susceptibility image.

Finally, similar to the brain phantom experiments, the TGV-HIRE restores an overly smoothed susceptibility image as shown in \cref{RealCornellTGVHIRE,RealCornellTGVHIREAx,RealCornellTGVHIREZoom} compared to the Frame-HIRE in \cref{RealCornellFraHIRE,RealCornellFraHIREAx,RealCornellFraHIREZoom}. In addition, the split Bregman algorithm of the TGV-HIRE is approximately $7$ times slower than the Frame-HIRE, as shown in \cref{ComparisonTimeRealCornell}, which again shows that the TGV regularization approach may not be suitable for the real clinical applications. Hence, as in the brain phantom experiments, we can conclude that compared to the TGV-HIRE, the Frame-HIRE is able to achieve the efficiency of its split Bregman algorithm as well as the shadow and streaking artifact removal.


\begin{table}[tp!]
\centering
\caption{Comparison of the CPU time for the in vivo data w.r.t. the choice of regularization term.}\label{ComparisonTimeRealCornell}
\vspace{-0.2cm}
\begin{tabular}{|c||c|c|c|c|c|c|}
\hline
\multirow{2}{*}{Indices}&\multicolumn{3}{|c|}{Wavelet Frame}&\multicolumn{3}{|c|}{TGV}\\ \cline{2-7}
&Integral&Differential&HIRE&Integral&Differential&HIRE\\ \hline
CPU Time&$628.25$&$353.32$&$953.10$&$3330.51$&$775.60$&$6776.12$\\ \hline
\end{tabular}
\end{table}

\section{Conclusion}\label{Conclusion}

In this paper, we proposed a new regularization based susceptibility reconstruction model. The proposed HIRE model is based on the identification of the harmonic incompatibility in the measured local field data arising from the underlying PDE \cref{QSM_PDE}. The harmonic property is imposed as a prior of incompatibility via the sparsity under the Laplacian into the integral approach so that we can apply the idea of two system regularization model. By doing so, we can take the incompatibility in the data which is other than the additive noise into account, achieving the susceptibility image reconstruction with less artifacts. Finally, the experimental results show that our proposed approach \cref{WaveletFrameHIRE} outperforms the existing approaches in both brain phantom and in vivo MR data.

\appendix
\section{Preliminaries on Wavelet Frame}\label{WaveletFrame} Provided here is a brief introduction on the tight wavelet frame. Briefly speaking, it is a generalization of the orthogonal wavelet basis (e.g. \cite{Mallat2008}) into the redundant system, and due to the redundancy, it is more robust to errors than the traditional orthonormal basis \cite{B.Dong2017a}. Interested readers may consult \cite{Daubechies1992,I.Daubechies2003,A.Ron1997} for theories of frame and wavelet frame, \cite{Shen2010} for a short survey on the theory and applications of frames, and \cite{B.Dong2013,B.Dong2015} for more detailed surveys.


For a given $\Ps=\left\{\psi_1,\cdots,\psi_r\right\}\subseteq L_2(\R^d)$ with $d\in\N$, a quasi-affine system $\msX(\Ps)$ generated by $\Ps$ is the collection of the dilations and the shifts of the elements in $\Ps$:
\begin{align}\label{QASystem}
\msX(\Ps)=\left\{\psi_{\alpha,n,\bk}:1\leq\alpha\leq r,~n\in\Z,~\bk\in\Z^d\right\},
\end{align}
where $\psi_{\alpha,n,\bk}$ is defined as
\begin{align}\label{QAFramelet}
\psi_{\alpha,n,\bk}=\left\{\begin{array}{cl}
2^{\f{nd}{2}}\psi_{\alpha}(2^n\cdot-\bk)~~&~~n\geq 0;\vspace{0.125em}\\
2^{nd}\psi_{\alpha}(2^n\cdot-2^{n}\bk)~~&~~n<0.
\end{array}\right.
\end{align}
We say that $\msX(\Ps)$ is a tight wavelet frame on $L_2(\R^d)$ if we have
\begin{align}\label{TightFrame}
\left\|f\right\|_{L_2(\R^d)}^2=\sum_{\alpha=1}^r\sum_{n\in\Z}\sum_{\bk\in\Z^d}\left|\left\la f,\psi_{\alpha,n,\bk}\right\ra\right|^2
\end{align}
for every $f\in L_2(\R^d)$. In this case, each $\psi_{\alpha}$ is called a (tight) framelet, and $\left\la f,\psi_{\alpha,n,\bk}\right\ra$ is called the canonical coefficient of $f$.


The constructions of (anti-)symmetric and compactly supported framelets $\Ps$ are usually based on a multiresolution analysis (MRA); we first find some compactly supported refinable function $\phi$ with a refinement mask $q_0$ such that
\begin{align}\label{MRA-RF}
\phi=2^d\sum_{\bk\in\Z^d}q_0[\bk]\phi(2\cdot-\bk).
\end{align}
Then the MRA based construction of $\Ps=\left\{\psi_1,\cdots,\psi_r\right\}\subseteq L_2(\R^d)$ is to find finitely supported masks $q_{\alpha}$ such that
\begin{align}\label{MRA-Fra}
\psi_{\alpha}=2^d\sum_{\bk\in\Z^d}q_{\alpha}[\bk]\phi(2\cdot-\bk),~~~~~\alpha=1,\cdots,r.
\end{align}
The sequences $q_1,\cdots,q_r$ are called wavelet frame mask or the high pass filters of the system, and the refinement mask $q_0$ is also called the low pass filter.

The unitary extension principle (UEP) of \cite{A.Ron1997} provides a general theory of the construction of MRA based tight wavelet frames. Briefly speaking, as long as $\left\{q_0,q_1,\cdots,q_r\right\}$ are compactly supported and their Fourier series $\wh{q}_{\alpha}(\xxi)=\sum_{\bk\in\Z^d}q_{\alpha}[\bk]e^{-i\xxi\cdot\bk}$ satisfy
\begin{align}\label{UEP}
\sum_{\alpha=0}^r\left|\wh{q}_{\alpha}(\xxi)\right|^2=1~~~~\text{and}~~~~\sum_{\alpha=0}^r\wh{q}_{\alpha}(\xxi)\overline{\wh{q}_{\alpha}(\xxi+\nnu)}=0
\end{align}
for all $\nnu\in\left\{0,\pi\right\}^d\setminus\left\{\0\right\}$ and $\xxi\in[-\pi,\pi]^d$, the quasi-affine system $\msX(\Ps)$ with $\Ps=\left\{\psi_1,\cdots,\psi_r\right\}$ defined by \cref{MRA-Fra} forms a tight frame of $L_2(\R^d)$, and the filters $\left\{q_0,q_1,\cdots,q_r\right\}$ form a discrete tight frame on $\ell_2(\Z^d)$ \cite{B.Dong2013}.

\begin{figure}[tp!]
\centering
\subfloat[Refinable function $\phi$]{\label{HaarScale}\includegraphics[width=6.5cm]{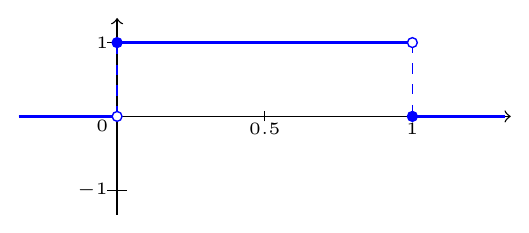}}\hspace{0.005cm}
\subfloat[Framelet function $\psi_1$]{\label{HaarFramelet}\includegraphics[width=6.5cm]{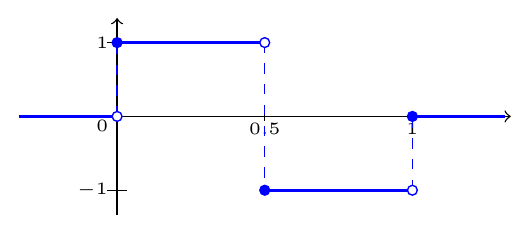}}\vspace{-0.20cm}
\caption{Illustration of $\phi$ and $\psi_1$ in \cref{ExA1}.}\label{IllustrateExA1}
\end{figure}

\begin{example}\label{ExA1} The piecewise constant B-spline (or the Haar framelet) \cite{I.Daubechies2003} for $L_2(\R)$ has one refinable function and one framelet
\begin{align*}
\phi(x)=\left\{\begin{array}{ccl}
1&\text{if}&x\in[0,1)\vspace{0.125em}\\
0&\text{if}&x\notin[0,1)
\end{array}\right.~~~\text{and}~~~\psi_1(x)=\left\{\begin{array}{ccl}
1&\text{if}&x\in[0,1/2)\vspace{0.125em}\\
-1&\text{if}&x\in[1/2,1)\vspace{0.125em}\\
0&\text{if}&x\notin[0,1)
\end{array}\right.
\end{align*}
as shown in \cref{IllustrateExA1}. Here, the associated filters are
\begin{align*}
q_0=\f{1}{2}\big[\begin{array}{cc}
1&1
\end{array}\big]~~~\text{and}~~~q_1=\f{1}{2}\big[\begin{array}{cc}
1&-1
\end{array}\big].
\end{align*}
Since this $\left\{q_0,q_1\right\}$ satisfies \cref{UEP}, $\msX\left(\Ps\right)$ with $\Ps=\left\{\psi_1\right\}$ forms a tight frame on $L_2(\R)$.
\end{example}

The tight frame on $L_2(\R^d)$ with $d\geq 2$ can be constructed by taking tensor products of univariate tight framelets \cite{J.F.Cai2012,J.F.Cai2016,Daubechies1992,B.Dong2013}. Given a set of univariate masks $\big\{q_0,q_1,\cdots,q_r\big\}$, we define multivariate masks $q_{\aal}[\bk]$ with $\aal=(\alpha_1,\cdots,\alpha_d)$ and $\bk=(k_1,\cdots,k_d)$ as
\begin{align*}
q_{\aal}[\bk]=q_{\alpha_1}[k_1]\cdots q_{\alpha_d}[k_d],~~~~~~0\leq\alpha_1,\cdots,\alpha_d\leq r,~~\bk=(k_1,\cdots,k_d)\in\Z^d.
\end{align*}
The corresponding multivariate refinable function and framelets are defined as
\begin{align*}
\psi_{\aal}(\x)=\psi_{\alpha_1}(x_1)\cdots\psi_{\alpha_d}(x_d),~~~~~0\leq\alpha_1,\cdots,\alpha_d\leq r,~~\x=(x_1,\cdots,x_d)\in\R^d
\end{align*}
with $\psi_0=\phi$ for convenience. If the univariate masks $\left\{q_{\alpha}\right\}$ are constructed from UEP, then we can verify that $\left\{q_{\aal}\right\}$ satisfies \cref{UEP} and thus $\msX(\Ps)$ with $\Ps=\left\{\psi_{\aal}:\aal\in\{0,\cdots,r\}^d\setminus\{\0\}\right\}$ forms a tight frame for $L_2(\R^d)$.

In the discrete setting, let $\mI_d\simeq\R^{N_1\times\cdots\times N_d}$ be the space of real valued functions defined on a regular grid $\left\{0,1,\cdots,N_1-1\right\}\times\cdots\times\left\{0,1,\cdots,N_d-1\right\}$. The fast framelet decomposition, or the analysis operator with $L$ levels of decomposition is defined as
\begin{align}
Wu=\left\{W_{l,\aal}u:(l,\aal)\in\left(\left\{0,\cdots,L-1\right\}\times\BB\right)\cup\left\{(L-1,\0)\right\}\right\}
\end{align}
where $\BB=\left\{0,\cdots,r\right\}^d\setminus\left\{\0\right\}$ is the framelet band. Then the frame coefficients $W_{l,\aal}u\in\mI_d$ of $u\in\mI_d$ at level $l$ and band $\aal$ are defined as
\begin{align*}
W_{l,\aal}u=q_{l,\aal}[-\cdot]\circledast u.
\end{align*}
where $\circledast$ denotes the discrete convolution with a certain boundary condition (e.g. the periodic boundary condition), and $q_{l,\aal}$ is defined as
\begin{align}
q_{l,\aal}=\wt{q}_{l,\aal}\circledast\wt{q}_{l-1,\0}\circledast\cdots\circledast\wt{q}_{0,\0}~~\text{with}~~\wt{q}_{l,\aal}[\bk]=\left\{\begin{array}{cl}
q_{\aal}[2^{-l}\bk],&\bk\in 2^l\Z^d\vspace{0.125em}\\
0,&\bk\notin 2^l\Z^d.
\end{array}\right.
\end{align}
We denote by $W^T$, the adjoint of $W$, the fast reconstruction (or the synthesis operator). Then by UEP \cref{UEP}, we have $W^TW=I$.

Finally, we mention that among many different choice of framelets, the ones constructed from the B-spline are the most popular in image processing. This is due to the multiscale structure of the wavelet frame systems, short supports of the (anti-)symmetric framelet functions with varied vanishing moments, and the presence of both low pass and high pass filters in the wavelet frame filter banks, which are desirable in sparsely approximating images \cite{B.Dong2017a}. A tight frame system constructed from the low order B-spline has fewer filters with shorter supports compared to the ones constructed from the high order B-splines. Hence, low order B-spline framelet systems are more computationally efficient while the high order ones are capable of capturing richer image singularities. Moreover, since high order B-spline framelets have larger supports, they may introduce more numerical viscosity, often leading to smoother reconstructions in image restoration tasks. Hence, the choice of framelet systems indeed depends on the the task and the computational cost we can afford \cite{H.Zhang2018}. In this paper, we fix $W$ to be the Haar framelet system for the wavelet frame regularization models as the susceptibility images can be well approximated by piecewise constant functions. Besides, we always fix $L=1$ to avoid the memory storage problem as we solve three dimensional inverse problem. We also note that the choices of $W$ will indeed affect the reconstruction results. For example, the use of data driven tight frames in \cite{J.F.Cai2014} will generate better reconstruction results due to its adaptivity, even though it requires further numerical studies. Nonetheless, we forgo more details on the choice of $W$ in order not to dilute the main focus of this paper.

\section*{Acknowledgments} The authors thank the authors in \cite{Y.Wang2015,C.Wisnieff2013} for making the data sets and the MATLAB toolbox available so that the experiments can be implemented. The authors also thank the anonymous reviewers for their constructive suggestions and comments that helped tremendously with improving the presentations of this paper.

\bibliographystyle{siamplain}

\end{document}